\documentclass[12pt,centertags]{amsart}
\usepackage{amsmath,amstext,amsthm,a4,amssymb,amscd}
\usepackage[mathscr]{eucal}
\usepackage{mathrsfs}
\usepackage{epsfig}
\numberwithin{equation}{section}

\newcommand{\field}[1]{\mathbb{#1}}
\newcommand{\bZ}{\field{Z}}
\newcommand{\bR}{\field{R}}
\newcommand{\bC}{\field{C}}
\newcommand{\bN}{\field{N}}
\newcommand{\bT}{\field{T}}
 \def\cC{\mathscr{C}}
\def\cL{\mathscr{L}}
\def\cO{\mathscr{O}}
 \def\cR{\mathscr{R}}
 \def\cU{\mathscr{U}}
\def\mC{\mathcal{C}}
\def\mA{\mathcal{A}}
\def\mB{\mathcal{B}}
\def\mD{\mathcal{D}}
\def\mF{\mathcal{F}}
\def\mI{\mathcal{I}}
\def\mL{\mathcal{L}}
\def\mO{\mathcal{O}}
\def\mP{\mathcal{P}}
\def\mQ{\mathcal{Q}}
\def\mR{\mathcal{R}}
\def\mT{\mathcal{T}}
\def\mV{\mathcal{V}}
\def\mW{\mathcal{W}}
\newcommand\mS{\mathcal{S}}
\def\kt{\mathfrak{t}}
\def\kg{\mathfrak{g}}
\def\bE{{\bf E}}

\def\bJ{{\bf J}}
\def\bk{{\bf k}}

\def\br{{\bf r}}

\def\Re{{\rm Re}}
\def\Im{{\rm Im}}
\def\Cas{{\rm Cas}}

\newcommand{\til}[1]{\widetilde{#1}}
\newcommand{\tx}{\til{X}}
\newcommand{\tj}{\til{J}}

\newcommand{\tom}{\til{\omega}}
\newcommand{\tl}{\til{L}}
\newcommand{\te}{\til{E}}

\DeclareMathOperator{\End}{End}
\DeclareMathOperator{\Hom}{Hom}
\DeclareMathOperator{\Ker}{Ker}
\DeclareMathOperator{\Coker}{Coker}

\DeclareMathOperator{\Id}{Id}
\DeclareMathOperator{\supp}{supp}
\DeclareMathOperator{\tr}{Tr}
\DeclareMathOperator{\Ind}{Ind}

\DeclareMathOperator{\td}{Td}
\DeclareMathOperator{\vol}{vol}
\DeclareMathOperator{\ch}{ch}
\newcommand{\spin}{$\text{spin}^c$ }
\newcommand{\spec}{\rm Spec}
\newcommand{\norm}[1]{\lVert#1\rVert}

\newcommand{\om}{\omega}

\newtheorem{thm}{Theorem}[section]
\newtheorem{lemma}[thm]{Lemma}
\newtheorem{prop}[thm]{Proposition}
\newtheorem{cor}[thm]{Corollary}
\theoremstyle{definition}
\newtheorem{rem}[thm]{Remark}
\theoremstyle{definition}
\newtheorem{defn}[thm]{Definition}
\newcommand{\be}{\begin{eqnarray}}
\newcommand{\ee}{\end{eqnarray}}
\newcommand{\ov}{\overline}
\newcommand{\wi}{\widetilde}
\newcommand{\var}{\varepsilon}

\newcommand{\comment}[1]{}

\begin{document}
\title{Bergman kernels and symplectic reduction}

\author{Xiaonan Ma}
\address{Centre de Math\'ematiques Laurent Schwartz, UMR 7640 du CNRS,
Ecole Polytechnique, 91128 Palaiseau Cedex,
France (ma@math.polytechnique.fr)}
\
\author{Weiping Zhang}
\address{Chern Institute of Mathematics \& LPMC, Nankai
University, Tianjin 300071, P.R. China.
(weiping@nankai.edu.cn)}

 \begin{abstract} We generalize several recent results concerning
 the asymptotic expansions of Bergman kernels to the framework of
geometric quantization and establish an
asymptotic symplectic identification property. 
More precisely, we study the asymptotic
expansion of the $G$-invariant Bergman kernel of the \spin Dirac
operator associated with high tensor powers of a positive line
bundle on a symplectic manifold.
We also develop a way to compute the coefficients of the
expansion, and compute the first few of them, especially, we
obtain the scalar curvature of the reduction space from the
$G$-invariant Bergman kernel on the total space. These results
generalize the corresponding results in the non-equivariant
setting, which has played a crucial role in the recent work of
Donaldson on stability of projective manifolds, to the geometric
quantization setting. As another kind of application, we
generalize some Toeplitz operator type properties in
semi-classical analysis to the framework of geometric
quantization.  The method we use is inspired by Local Index
Theory, especially by
 the analytic localization techniques developed by Bismut and Lebeau.
\end{abstract}

\maketitle

{
\tableofcontents
}

\setcounter{section}{-1}
\section{Introduction} \label{s0}

The study of the Bergman kernel is a classical subject in the
theory of several complex variables, where usually it concerns the
kernel function of the projection  operator to an infinite
dimensional Hilbert space. The recent interest of the analogue of
this concept in complex geometry mainly started in the paper of
Tian \cite{T90}, which was in turn inspired by a question of Yau.
Here, the projection is onto a finite dimensional space instead.

After \cite{T90}, the Bergman kernel has been studied
extensively in  \cite{Ru98}, \cite{Ze98}, \cite{Catlin99},
\cite{Lu00}, establishing the diagonal asymptotic expansion for
high powers of an ample line bundle. Moreover, the coefficients in
the asymptotic expansion encode geometric information of the
underlying complex projective manifolds. This asymptotic expansion
plays a crucial role in the recent work of Donaldson \cite{Do01},
where the existence of K\"ahler metrics with constant scalar
curvature is shown to be closely related to Chow-Mumford
stability.

In \cite{DLM04a}, \cite{MM04a}, \cite{MM04c},
Dai, Liu, Ma and Marinescu
 studied the {\em full off-diagonal} asymptotic expansion of the
(generalized) Bergman kernel
of the spin$^c$ Dirac operator and the renormalized Bochner--Laplacian
 associated to a positive line bundle
on a compact symplectic manifold. As a by product, they gave a new
proof of the   results mentioned in the previous paragraph. They
find also various applications therein, especially as pointed out
in \cite{MM04c}, the {\em full off-diagonal}
 asymptotic expansion implies Toeplitz operator type properties.
This approach is inspired by the Local Index Theory, especially by
the analytic localization techniques of Bismut-Lebeau \cite[\S 11]{BL91}.
We refer to the above papers
and the recent book \cite{MM05b} for detail
informations of the Bergman kernel on compact symplectic
manifolds.

In this paper, we generalize some of the  results in
\cite{DLM04a}, \cite{MM04a} and \cite{MM04c} to the framework of
geometric quantization, by  studying the asymptotic expansion of
the $G$-invariant Bergman kernel for high powers of an ample line
bundle on symplectic manifolds admitting a Hamiltonian group
action.

To start with, let $(X,\om)$ be a compact symplectic manifold of
real dimension $2n$. Assume that there exists a Hermitian line
bundle $L$ over $X$ endowed with a Hermitian connection $\nabla^L$
with the property that
\be \label{0.2}
\frac{\sqrt{-1}}{2\pi}R^L=\omega,
\ee
 where
$R^L=(\nabla^L)^2$ is the curvature of $(L,\nabla^L)$.

Let $(E,h^E)$ be a Hermitian vector bundle  on $X$ equipped with a
Hermitian connection $\nabla^E$ and $R^E$ denotes the associated
curvature.

Let $g^{TX}$ be a Riemannian metric on $X$.
Let $\bJ:TX\to TX$ be the skew--adjoint linear
map which satisfies the relation
\be \label{0.1}
\om(u,v)=g^{TX}(\bJ u,v)
\ee
 for $u,v \in TX$.

Let $J$ be an almost complex structure such that
\be \label{0.3}
g^{TX}(Ju,Jv)= g^{TX}(u,v), \ \ \ \ \om(Ju,Jv)=\om(u,v)
 \ee
and that $\om(\cdot, J\cdot)$
defines a metric on $TX$. Then $J$ commutes with $\bJ$
and $J=\bJ (-\bJ^2)^{-1/2}$ (cf. \eqref{0a00}).

 Let $\nabla ^{TX}$ be the Levi-Civita connection on $(TX,
g^{TX})$ with curvature $R^{TX}$ and scalar curvature $r^X$. The
connection  $\nabla ^{TX}$ induces a natural connection $\nabla
^{\det}$ on $\det (T^{(1,0)} X)$ with curvature $R^{\det}$, and
the Clifford connection $\nabla^{\rm Cliff}$ on the Clifford
module $\Lambda (T^{*(0,1)}X)$ with curvature $R^{\rm Cliff}$
 (cf. Section \ref{s3.01}).

 The \spin Dirac
operator $D_p$ acts on
$\Omega^{0,{\scriptscriptstyle{\bullet}}}(X,L^p\otimes E)
=\bigoplus_{q=0}^n\Omega^{0,q}(X,L^p\otimes E)$, the direct sum of
spaces of $(0,q)$--forms with values in $L^p\otimes E$.
We denote by  $D_{p}^+$ the restriction of $D_p$ on
$\Omega^{0,{\rm even}}(X,L^p\otimes E)$.
The index of $D_{p}^+$ is defined by
\begin{align}  \label{a0.1}
\Ind (D_{p}^+) = \Ker D_{p}^+- \Coker D_{p}^+.
\end{align}

Let $G$ be a compact connected Lie group with Lie algebra $\kg$
and $\dim G= n_0$.
Suppose that $G$ acts on $X$ and its action on $X$ lifts on $L$ and $E$.
Moreover, we assume the $G$-action preserves the above
 connections and metrics on $TX, L,E$ and $J$.
Then $\Ind (D_{p}^+) $ is a virtual representation of $G$. Denote
by $(\Ker D_p)^G$, $\Ind (D_{p}^+) ^G$ the $G$-trivial components
of $\Ker D_p$, $\Ind (D_{p}^+)$ respectively.

The action of $G$ on $L$ induces naturally a moment map $\mu :
X\to \kg^*$ (cf. (\ref{a6})). We assume that $0\in \kg^*$ is a
regular value of $\mu$.

Set $P=\mu^{-1}(0).$ Then  the Marsden-Weinstein symplectic
reduction $(X_G=P/G, \om_{X_G})$ is a symplectic  orbifold ($X_G$
is smooth if  $G$ acts  freely on $P$).

Moreover, $(L,\nabla ^L)$, $(E,\nabla ^E)$ descend to $(L_G,\nabla
^{L_G})$, $(E_G,\nabla ^{E_G})$ over $X_G$ so that the
corresponding curvature condition $\frac{\sqrt{-1}}{2\pi} R^{L_G}=
\omega_G$ holds (cf. \cite{GuSt82}). The $G$-invariant almost
complex structure $J$ also descends to an almost complex structure
$J_G$ on $TX_G$, and $h^L, h^E, g^{TX}$ descend to $h^{L_G}$,
$h^{E_G}, g^{TX_G}$ respectively.

One can construct the corresponding spin$^c$ Dirac operator
$D_{G,p}$ on $X_G$.

The geometric quantization conjecture of Guillemin-Sternberg
\cite{GuSt82} can be stated as follows: when $E$ is the trivial
bundle $\bC$ on $X$, for any $p>0$,
\begin{align}\label{a0.2}
 \dim\left(\Ind (D_{p}^+) ^G \right)= \dim\left(\Ind (D_{G,p}^+)\right).
\end{align}

When $G$ is abelian, this conjecture was proved by Meinrenken
\cite{Mein95} and Vergne  \cite{Ve96}. The remaining  nonabelian
case was proved
 by  Meinrenken \cite{Mein98} using the symplectic cut
techniques of Lerman, and by Tian and Zhang \cite{TZ98}  using
analytic localization techniques.

More generally, by a result of  Tian and Zhang \cite[Theorem 0.2]{TZ98},
for any general vector bundle $E$ as above,
there exists $p_0>0$ such that for any $p\geq p_0$,
 \eqref{a0.2} still holds.

On the other hand, by \cite[Theorem 2.5]{MM02} (cf.
\eqref{main1}), which is a direct consequence of the Lichnerowicz
formula for $D_{p}$,   for $p$ large enough, $\Coker D_{p}^+$ is
null (cf. also \cite{BorU96}, \cite{Bra99}). Thus there exists
$p_0>0$ such that for any $p\geq p_0$,
\begin{multline}\label{a0.3}
\dim (\Ker D_p)^G = \dim (\Ker D_{G,p})
=\dim\left(\Ind(D_{G,p}^+)\right)
= \int_{X_G} \td (TX_G) \ch (L_G^p\otimes E_G)\\
={\rm rk} (E) \int_{X_G} \frac{(p\, c_1(L_G))^{n-n_0}}{(n-n_0)!}
+ \int_{X_G} \Big(c_1(E_G) + \frac{{\rm rk} (E)}{2} c_1(TX_G)\Big)
\frac{(p \, c_1(L_G))^{n-n_0-1}}{(n-n_0-1)!}\\
 + \cO(p^{n-n_0-2}),
\end{multline}
where $\ch(\cdot), c_1(\cdot), \td(\cdot)$ are the Chern character,
the first Chern class and the Todd class of the corresponding complex
vector bundles
($TX_G$ is a complex vector bundle with complex structure $J_G$).

Set $E_p:=\Lambda (T^{*(0,1)}X)\otimes L^p\otimes E$.
Let $\langle\quad\rangle $ be the $L^2$-scalar product
on $\Omega^{0,\bullet}(X, L^p\otimes E)= \cC^\infty(X, E_p)$
induced by $g^{TX}, h^L, h^E$ as in \eqref{h10}.

Let $P_p^G$ be the orthogonal projection from
$(\Omega^{0,\bullet}(X,L^p\otimes E),\langle\quad\rangle ) $ on $(\Ker D_p)^G$.
The $G$-invariant Bergman kernel is  $P_p^G(x,x')$  $(x,x'\in X)$,
 the smooth kernel of $P_p^G$
with respect to the Riemannian volume form $dv_X(x')$.

 Let ${\rm pr}_1$ and  ${\rm pr}_2$ be the projections from $X\times X$
onto the first and second factor $X$ respectively. Then $P^G_p(x,
x')$ is a smooth section of $ {\rm pr}_1^*( E_p) \otimes {\rm
pr}_2^* (E_p^*) $ on $X\times X$. In particular,  $P_p^G(x,x)\in
\End(E_p)_x =\End(\Lambda (T^{*(0,1)}X)\otimes E)_x$.

The $G$-invariant Bergman kernel $P_p^G(x,x')$ is an analytic
version of $(\Ker D_p)^G$. In view of (\ref{a0.3}),
 it is natural to expect that
  the kernel $P_p^G(x,x')$ should be closely related to  the corresponding
Bergman kernel on the symplectic reduction $X_G$.
The purpose of this paper is to study the asymptotic expansion of
 the $G$-invariant Bergman kernel  $P_p^G(x,x')$ as $p\to \infty$,
and we will relate it to the asymptotic expansion of the Bergman kernel
 on the symplectic reduction $X_G$.

Let $d^X(x,x')$ be the Riemannian distance between  $x,\, x'\in X$.

 In Section \ref{s3.1}, we prove
the following result which allows us to reduce our problem as a
problem near $P=\mu^{-1}(0)$.

\begin{thm} \label{t0.0}     For any open $G$-neighborhood $U$
of $P$ in $X$, $\var_0 >0$,   $l,m \in \bN$,
there exists $C_{l,m}>0$ {\rm (}depend on $U$, $\var_0${\rm )}
such that for $p\geq 1$,  $x,x'\in X, d^X(G x,x')\geq \var_0$
or   $x,x'\in X\setminus U$,
\begin{align}\label{0.5}
|P^G_p(x, x')|_{\cC^m} \leq C_{l,m} p^{-l},
\end{align}
where $\cC^m$ is the $\cC^m$-norm induced by
$\nabla ^L,\nabla ^E$, $\nabla ^{TX}$, $h^L, h^E$ and $g^{TX}$.
\end{thm}

Assume for simplicity that $G$ acts freely on $P$.

Let $U$ be an open $G$-neighborhood of $\mu^{-1}(0)$ such that $G$
acts freely on $U$.

For any $G$-equivariant vector bundle $(F,\nabla^F)$ on $U$, we
denote by $F_{B}$ the  bundle on $U/G=B$ induced naturally by
$G$-invariant sections of $F$ on $U$. The connection $\nabla^F$
induces canonically a connection $\nabla ^{F_B}$ on $F_B$. Let
$R^{F_B}$ be its curvature.
Let \begin{align}\label{a0.5}
\mu^F(K)= \nabla ^F_{K^X}- L_K\in \End(F)
\end{align}
 for $K\in \kg$
and $K^X$ the corresponding vector field on $U$.

Note that $P^G_p\in (\cC^\infty (U\times U, {\rm pr}_1^* E_p
\otimes {\rm pr}_2^* E_p^*))^{G\times G}$, thus we can view
$P^G_p(x,x')$ ($x,\, x'\in U$) as a smooth section of ${\rm
pr}_1^*( E_p)_{B} \otimes {\rm pr}_2^*( E_p^*)_{B}$
 on $B\times B$.

Let $g^{TB}$ be the Riemannian metric on $U/G=B$ induced by $g^{TX}$.
Let $\nabla ^{TB}$ be the Levi-Civita connection on $(TB,g^{TB})$ with
curvature $R^{TB}$.
 Let $N_G$ be the normal bundle to $X_G$ in $B$.
We identify $N_G$ with the orthogonal complement of $TX_G$ in
$(TB|_{X_G}, g ^{TB})$.

Let $g^{TX_G}$, $g^{N_G}$ be the metrics on $TX_G$, $N_G$ induced
by $g^{TB}$ respectively.

Let $P^{TX_G}$, $P^{N_G}$ be the orthogonal projections from
$TB|_{X_G}$ on $TX_G$, $N_{G}$ respectively. Set
\begin{align} \label{a0.6}
&\nabla ^{N_G}= P^{N_G}(\nabla^{TB}|_{X_G})P^{N_G}, \quad
\nabla ^{TX_G}=P^{TX_G}(\nabla^{TB}|_{X_G})P^{TX_G},\\
&{^0\nabla}^{TB} =\nabla ^{TX_G} \oplus\nabla ^{N_G}, \quad
 A= {\nabla}^{TB}|_{X_G}-{^0\nabla}^{TB}. \nonumber
\end{align}
Then $\nabla ^{N_G}, {^0\nabla}^{TB}$ are Euclidean  connections on
$N_G$, $TB|_{X_G}$ on $X_G$,
$\nabla ^{TX_G}$ is the Levi-Civita connection on $(TX_G, g^{TX_G})$, and
 $A$ is the associated second fundamental form.

 Denote by $\vol (Gx)\ (x\in U)$ the
volume of the orbit $Gx$ equipped with the metric induced by $g^{TX}$.
Following \cite[(3.10)]{TZ98}, let $h(x)$ be the function on $U$
defined by
\begin{align}\label{0.6}
h (x)= ( \vol (Gx))^{1/2}.
\end{align}
Then $h$ reduces to a function on $B$.

Denote by $I_{\bC\otimes E}$ the projection from $\Lambda
(T^{*(0,1)}X)\otimes E$ onto $\bC\otimes E$  under the
decomposition $\Lambda (T^{*(0,1)}X)\otimes E= \bC\otimes E \oplus
\Lambda ^{>0} (T^{*(0,1)}X)\otimes E$,
 and  $I_{\bC\otimes E_B}$ the corresponding projection on $B$.

In the whole paper, for any $x_0\in X_G$,  $Z\in T_{x_0}B$, we   write
 $Z=Z^0+Z^\bot$, with $Z^0\in T_{x_0}X_G$, $Z^\bot\in N_{G,x_0}$.

Let $\tau_{Z^0} Z^\bot\in N_{G,\exp^{X_G}_{x_0}(Z^0)}$ be the parallel
 transport of $Z^\bot$ with respect to the connection $\nabla ^{N_G}$
along the geodesic in $X_G$, $[0,1]\ni t \to \exp^{X_G}_{x_0}(t Z^0)$.

 For $\var_0>0$ small enough,
we identify $Z\in T_{x_0}B$, $|Z|<\var_0$ with
$\exp^{B}_{\exp_{x_0}^{X_G}(Z^0)} (\tau_{Z^0} Z^\bot)\in B$.
 Then for $x_0\in X_G$, $Z,Z^{\prime}\in T_{x_0}B$, $|Z|,|Z^{\prime}|<\var_0$,
 the map $\Psi: TB|_{X_G}\times TB|_{X_G} \to B\times B$,
$$\Psi (Z,Z^{\prime})=
(\exp^{B}_{\exp_{x_0}^{X_G}(Z^0)} (\tau_{Z^0} Z^\bot),
\exp^{B}_{\exp_{x_0}^{X_G}(Z^{'0 })} (\tau_{Z^{\prime 0}} Z^{'\bot}))$$
 is well defined.

 We identify $(E_p)_{B, Z}$ to  $(E_p)_{B, x_0}$ by using parallel
transport with respect to $\nabla^{(E_p)_B}$ along $[0,1]\ni u\to uZ$.

Let $\pi_B: TB|_{X_G}\times TB|_{X_G}\to X_G$ be the natural
projection from the fiberwise product of $TB|_{X_G}$ on $X_G$ onto
$X_G$.

From Theorem \ref{t0.0}, we only need to understand
$P^G_p\circ\Psi$, and under our identification, $P^G_p\circ\Psi
(Z,Z^{\prime})$ is a smooth section of $$\pi_B^*(\End (E_{p})_B)=
\pi_B^*(\End ( \Lambda(T^{*(0,1)}X) \otimes E)_B)$$
 on $TB|_{X_G}\times TB|_{X_G} $.

 Let $|\quad|_{\cC ^{m'}(X_G)}$ be the $\cC^{m'}$-norm
on $\cC^\infty(X_G, \End ( \Lambda(T^{*(0,1)}X) \otimes E)_B)$
induced by $\nabla ^{\text{Cliff}_{B}}$, $\nabla ^{E_{B}}$,
$h^E$ and $g^{TX}$. The norm $|\quad|_{\cC ^{m'}(X_G)}$
induces naturally a $\cC^{m'}$-norm along
$X_G$  on $\cC^\infty(TB|_{X_G}\times TB|_{X_G},
\pi_B^*(\End ( \Lambda(T^{*(0,1)}X) \otimes E)_B))$,
we still denote it by $|\quad|_{\cC ^{m'}(X_G)}$.


 Let $dv_B, dv_{X_G}$, $dv_{N_G}$ be the Riemannian volume forms on
$(B, g^{TB})$, $(X_G, g^{TX_G}$), $(N_G,g^{N_G})$ respectively.
Let $\kappa \in \cC^{\infty}(TB|_{X_G},\bR)$,
with  $\kappa =1$ on $X_G$,
be defined by that for $Z \in T_{x_0} B$, $x_0 \in X_G$,
\begin{align}\label{0.9}
dv_B (x_0, Z)= \kappa(x_0, Z)dv_{T_{x_0}B}(Z)
=\kappa(x_0, Z) dv_{X_G}(x_0)dv_{N_{G,x_0}}.
\end{align}

The following result is one of the main results of this paper.

\begin{thm} \label{t0.1}
 Assume  that $G$ acts freely on $\mu^{-1}(0)$ and $\bJ=J$ on $\mu^{-1}(0)$.
  Then there exist  $\mQ_{r}(Z,Z^{\prime})\in
\End (\Lambda (T^{*(0,1)}X) \otimes E)_{B,x_0}$ $(x_0\in X_G, r\in \bN)$,
  polynomials in $Z,Z^{\prime}$ with the same parity as $r$, whose coefficients
are polynomials in $A$, $R^{TB}$, $R^{{\rm Cliff}_B}$, $R^{E_B}$,
$\mu^E$, $\mu ^{\rm Cliff}$ {\rm (}resp. $r^X$, $R^{\det}$, $R^E$;
resp. $h$, $R^L$, $R^{L_B}$; resp. $\mu${\rm )}
and their derivatives at $x_0$ to order $r-1$
 {\rm (}resp. $r-2$; resp. $r$, resp. $r+1${\rm )},
 such that if we denote by
\begin{align}\label{0.7}
&P^{(r)}_{x_0}(Z,Z^{\prime})= \mQ_{r}(Z,Z^{\prime})P(Z,Z^{\prime}),
\quad \mQ_{0}(Z,Z^{\prime})= I_{\bC\otimes E_B},
\end{align}
with
\begin{align}\label{a0.7}
P(Z,Z^{\prime})=& \exp\Big(-\frac{\pi}{2}
|Z^0-Z^{\prime 0}|^2- \pi \sqrt{-1}
\left\langle J_{x_0} Z^0,Z^{\prime 0} \right\rangle \Big)\\
&\times 2 ^{\frac{n_0}{2}}
\exp\Big(- {\pi}\big(|Z^\bot|^2+|Z^{\prime\bot}|^2\big) \Big), \nonumber
\end{align}
then  there exists $C''>0$ such that for any $k,m,m', m''\in \bN$,
there exists $C>0$ such that for
 $x_0\in X_G$, $Z,Z^{\prime}\in T_{x_0}B$, $|Z|, |Z^{\prime}|\leq  \var_0$,
\footnote{In the exponential factor of \cite[(7)]{MZ05a}, 
we missed $m^\prime$ as in the last line of \eqref{0.8} here.}
\begin{multline}\label{0.8}
(1+\sqrt{p}|Z^\bot|+\sqrt{p}|Z^{\prime\bot}|)^{m''}
\sup_{|\alpha|+|\alpha'|\leq m} \left
|\frac{\partial^{|\alpha|+|\alpha'|}}
{\partial Z^{\alpha} {\partial Z^{\prime}}^{\alpha'}} \right.\\
\left. \left (p^{-n+\frac{n_0}{2}}
(h\kappa^{\frac{1}{2}} )(Z)(h\kappa^{\frac{1}{2}} )(Z^{\prime})
P^G_p\circ\Psi(Z,Z^{\prime}) - \sum_{r=0}^k  P^{(r)}_{x_0}
(\sqrt{p} Z,\sqrt{p} Z^{\prime})
p^{-\frac{r}{2}}\right )\right |_{\cC ^{m'}(X_G)}\\
\leq C  p^{-(k+1-m)/2}  (1+\sqrt{p} |Z^0|+\sqrt{p} |Z^{\prime 0}|)
^{2(n+k+m^\prime+2)+m}
\exp (- \sqrt{C''} \sqrt{p} |Z-Z^{\prime}|) +\cO(p^{-\infty}).
\end{multline}
Furthermore, the expansion is uniform in the following sense{\rm:}
for any fixed $k,m,m', m''\in \bN$,
assume that the derivatives of $g^{TX}$, $h^L$, $\nabla ^L$,
$h^E$, $\nabla ^E$,and $J$
with order $\leqslant 2n+k+m+m'+3$
run over a set bounded in the $\cC^{m'}$--\,norm
taken with respect to the parameters and, moreover, $g^{TX}$ runs over a set
bounded below. Then the constant $C$ is independent of $g^{TX}${\rm;}
and the $\cC ^{m'}$-norm in \eqref{0.8} includes also the
derivatives on the parameters.
\end{thm}
In \eqref{0.8},
the term $\cO(p^{-\infty})$ means that for any $l,l_1\in \bN$,
there exists $C_{l,l_1}>0$ such that its $\cC^{l_1}$-norm is dominated
by $C_{l,l_1} p^{-l}$.

It is interesting to see that the kernel $P(Z,Z^{\prime})$ is the
product of two kernels : along $T_{x_0}X_G$, it is the classical
Bergman kernel on $T_{x_0}X_G$ with complex structure $J_{x_0}$,
while along $N_G$, it is the kernel of a harmonic oscillator on
$N_{G,x_0}$.

\begin{rem}\label{t0.3}
 i) Theorem \ref{t0.1} is a special case of Theorem \ref{tue17}
where we do not assume $\bJ=J$ on $P=\mu^{-1}(0)$.
 In Theorem \ref{0t3.6}, we get explicit informations on $P^{(r)}$
when $\bJ$ verifies \eqref{g2}.

ii) If $G$ does not act freely on $P$, then $X_G$ is an orbifold.
In Section \ref{s6.1}, we explain how to modify our arguments to
get the asymptotic expansion, Theorem \ref{t6.1}. Analogous  to
the usual orbifold case \cite[(5.27)]{DLM04a}, $P^G_p(x, x) (x\in
P)$ does not have a uniform asymptotic expansion if the singular
set of $X_G$ is not empty.

iii) Let $\mV$ be an irreducible representation of $G$,
let $P^\mV_p$ be the orthogonal projection from
$\Omega^{0,{\scriptscriptstyle{\bullet}}}(X,L^p\otimes E)$
on $\Hom_G(\mV, \Ker D_{p})\otimes \mV\subset \Ker D_{p}$.
In Section \ref{s6.2},
we get the asymptotic expansion of
the kernel $P^\mV_p(x,x')$ from  Theorems \ref{t0.0}, \ref{t0.1}.

iv) When $G=\{1\}$, Theorem \ref{t0.1} is \cite[Theorem 4.18$^\prime$]{DLM04a}.

v) If we take $Z=Z^{\prime}=0$ in \eqref{0.8}, then we get for $x_0\in X_G$,
\begin{align}\label{aa0.8}
P^{(0)}_{x_0}(0,0)= 2^{\frac{n_0}{2}}I_{\bC\otimes E_B},
\end{align}
and
\begin{align}\label{a0.8}
\Big | p^{-n+\frac{n_0}{2}} h^2(x_0)P^G_p(x_0, x_0)
- \sum_{r=0}^k  P^{(2r)}_{x_0}(0,0)p^{-r}\Big|_{\cC ^{m'}(X_G)}
\leq C  p^{-k-1}.
\end{align}
In  Section \ref{s6.3}, we show that \eqref{aa0.8} and
 \eqref{a0.8} are direct consequences of
the full off-diagonal asymptotic expansion of the Bergman kernel
\cite[Theorem 4.18$^\prime$]{DLM04a}.
In fact, one possible way to get Theorem \ref{t0.1} is to average
the full off-diagonal asymptotic expansion of the Bergman kernel on $X$
 \cite[Theorem 4.18$^\prime$]{DLM04a} with respect to a Haar measure on $G$.
However, we do not know how to get  the full off-diagonal
expansion, especially the fast decay along $N_G$ in \eqref{0.8} in
this way.
\end{rem}

In this paper we will apply the analytic localization techniques
to get Theorem \ref{t0.1}, and this method also gives us an
effective way to compute the coefficients in the asymptotic
expansion (cf. \S \ref{s3.11}). The key observation   is that the $G$-invariant Bergman
kernel is exactly the kernel of the orthogonal projection to the
zero space of a deformation of $D_p^2$ by the Casimir operator
(i.e., to consider $D_p^2-p{\rm Cas}$). This plays an essential
role in proving Theorems \ref{t0.0}, \ref{t0.1}.

 Let $\mathscr{I}_p$ be a section of
$\End (\Lambda (T^{*(0,1)}X)\otimes E)_B$ on $X_G$ defined by
\begin{align}\label{0.10}
 \mathscr{I}_p(x_0)=  \int_{Z\in N_{G}, |Z|\leq \var_0}
h^2 (x_0,Z) P^G_p\circ \Psi((x_0,Z),(x_0,Z)) \kappa (x_0, Z) dv_{N_G}(Z).
\end{align}
By Theorem \ref{t0.0}, modulo $\cO(p^{-\infty})$,
$\mathscr{I}_p(x_0)$ does not depend on $\var_0$, and
\begin{multline}\label{0.11}
\dim (\Ker D_p)^G = \int_X \tr [P^G_p(y,y) ] dv_X(y)
 = \int_U \tr [P^G_p(y,y) ] dv_X(y) +\cO(p^{-\infty})\\
=  \int_B h^2(y)  \tr [P^G_p(y,y) ] dv_B(y) +\cO(p^{-\infty})\\
= \int_{X_G} \tr [\mathscr{I}_p(x_0)] dv_{X_G}(x_0)  +
\cO(p^{-\infty}).
\end{multline}

A direct consequence of Theorem \ref{t0.1} is the following
corollary.
\begin{cor}\label{t0.5}
 Taken $Z=Z^{\prime}\in N_{G,x_0}$, $m=0$ in \eqref{0.8}, we get
\begin{multline}\label{0.13}
\Big|p^{-n+\frac{n_0}{2}} (h^2 \kappa)(Z) P^G_p(Z,Z) -\sum_{r=0}^k  P^{(r)}
_{x_0}(\sqrt{p}Z,\sqrt{p}Z)p^{-r/2}
\Big|_{\cC ^{m'}(X_G)}\\
\leq C  p^{-(k+1)/2} (1+\sqrt{p} |Z|)^{-m''} +\cO(p^{-\infty})   .
\end{multline}
In particular, there exist
$\Phi_r\in \End (\Lambda(T^{*(0,1)}X)\otimes E)_{B,x_0}$ $(r\in \bN)$
 which are polynomials in  $A$, $R^{TB}$, $R^{{\rm Cliff}_B}$, $R^{E_B}$,
$\mu^E$, $\mu ^{\rm Cliff}$, {\rm (}resp. $r^X$, $R^{\det}$,
$R^{E}$; resp.  $h$, $R^{L_B}$, $R^{L}$; resp. $\mu${\rm )}, and
their derivatives at $x_0$ up to order  $2r-1$
 {\rm (}resp. $2r-2$; resp. $2r$; resp. $2r+1${\rm )},
 and $\Phi_0= I_{\bC\otimes E_B}$,
such that for any $k,m'\in \bN$, there exists $C_{k,m'}>0$ such that
for any $x_0\in X_G$, $p\in \bN$,
\begin{align}\label{0.14}
&\Big |p^{-n+n_0}\mathscr{I}_p(x_0)- \sum_{r=0}^{k} \Phi_r(x_0)
p^{-r} \Big |_{\cC^{m'}} \leq C_{k,m'} p^{-k-1}.
\end{align}
\end{cor}

In the rest of Introduction, we will specify our results in the K\"ahler case.

 We suppose now that $(X,\om, J)$ is a compact K\"ahler manifold
and $\bJ=J$ on $X$.
 Assume also that $(L,h^L, \nabla^L)$, $(E,h^E, \nabla^E)$
are holomorphic Hermitian vector bundles with  holomorphic Hermitian
connections, and the action of $G$ on $X,L,E$ is holomorphic.

Let $\ov{\partial}^{L^p\otimes E,*}$
be the formal adjoint of the Dolbeault operator
$\ov{\partial}^{L^p\otimes E}$, then
\begin{align}\label{60.21}
D_p=\sqrt{2}(\ov{\partial}^{L^p\otimes E}+ \ov{\partial}^{L^p\otimes E,*}),
\end{align}
and 
\begin{align}\label{60.20}
D_p^2 = 2 \left(\ov{\partial}^{L^p\otimes
E}\ov{\partial}^{L^p\otimes E,*} + \ov{\partial}^{L^p\otimes
E,*}\ov{\partial}^{L^p\otimes E}\right)
\end{align}
 preserves the
$\bZ$-grading of $\Omega^{0, \bullet}(X,L^p\otimes E)$.

 By the Kodaira vanishing theorem, for $p$ large enough,
\begin{align}\label{60.22}
(\Ker D_p)^G = H^0(X, L^p\otimes E)^G.
\end{align}
Thus for $p$ large enough,
$P^G_p(x, x')\in (L^p\otimes E)_x\otimes (L^p\otimes E)_{x'}^*$ and so
$P^G_p(x, x)\in \End(E_x)$, $\mathscr{I}_p(x_0)\in \End(E_{x_0})$.
In particular, in \eqref{aa0.8},
\begin{align}\label{a0.15}
P^{(0)}_{x_0}(0,0)=2^{\frac{n_0}{2}}{\rm Id}_{E_G}.
\end{align}

\begin{rem}\label{t0.4}
In the special case of $E=\bC$, $P^G_p(x_0, x_0)$ is a
non-negative function on $X_G$, and \eqref{a0.8} has been proved
in \cite[Theorem 1]{Pao03}
 without knowing the informations on $P^{(2r)}_{x_0}(0,0)$,
 while in \cite[Theorem 1]{Pao04}, it was claimed that
$P^{(0)}_{x_0}(0,0)=1$. In \cite[Prop. 1]{Pao03}, Paoletti knew
that for any $l\in \bN$, there is $C>0$ such that for any $p$,
$|P^G_p(x,x)|\leq C p^{-l}$ uniformly on any compact subset of
$X\setminus (\mu^{-1}(0)\cup R)$, with $R$ the subset of unstable
points of the action of $G$.
In \cite{Pao04},  some Toeplitz operator type properties on $X_G$
was also claimed from the analysis of Toeplitz structures
of Boutet de Monvel--Guillemin \cite{BouGu81}, Boutet de Monvel-Sj\"ostrand
\cite{BouSj76} and
Shiffman-Zelditch \cite{ShZe02}.
If we suppose moreover that $G$ is a torus,
Charles \cite{Charles04} has also a different version
on the Toeplitz operator type properties on $X_G$.

 In Section \ref{s6.4}, 
  we will show that  Theorem \ref{t0.1} implies properties of
Toeplitz operators on $X_G$ (which also hold in the symplectic case).
In particular, we recover the results on
Toeplitz operators \cite{Charles04}, \cite{Pao04}.
\end{rem}

\comment{Let $g^{TX_G}$, $g^{N_G}$ be the metric on $TX_G, N_G$ induced by
$g^{TX}$. Let $dv_{X_G}$, $dv_{N_G}$ be the Riemannian volume form on
 $(X_G, g^{TX_G}$), $(N_G,g^{N_G})$. Let $\kappa_0\in \cC^{\infty}(B,\bR)$
be defined by
\begin{align}\label{0.9}
dv_B (x_0, Z^\bot)= \kappa_0(x_0, Z^\bot)
 dv_{X_G}(x_0)dv_{N_G}(Z^\bot).
\end{align}
Then  $\kappa_0=1$ on $X_G$.
}

Let $\wi{h}$ denote the restriction to $X_G$ of the function $h$
defined in (\ref{0.6}).

The second main result of this paper is that we can in fact obtain the
scalar curvature $r^{X_G}$ on the symplectic reduction $X_G$
from $\mathscr{I}_p$. 

\begin{thm} \label{t0.6} If $(X,\om)$ is a compact K\"ahler manifold and
$L,E$ are holomorphic vector bundles with holomorphic Hermitian
connections $\nabla ^L, \nabla ^E$, $\bJ=J$,
and $G$ acts freely on $\mu^{-1}(0)$, then
for $p$ large enough,
$\mathscr{I}_p(x_0)\in \End (E_G)_{x_0}$, and
in \eqref{0.14}, $\Phi_r(x_0)\in  \End (E_G)_{x_0}$
are polynomials in $A$, $R^{TB}$, $R^{E_B}$,
$\mu^E$, $R^E$ {\rm (}resp. $h$, $R^{L_B}$; resp. $\mu${\rm )}
and their derivatives at $x_0$ to order $2r-1$
 {\rm (}resp. $2r$, resp. $2r+1${\rm )},
 and $\Phi_0= {\rm Id}_{E_G}$. Moreover
\begin{align}\label{0.15}
\Phi_1(x_0)= \frac{1}{8\pi} r^{X_G}_{x_0}
+ \frac{3}{4\pi} \Delta_{X_G}\log \wi{h}
+\frac{1}{2\pi}R^{E_G}_{x_0}(w^0_j, \ov{w}^0_j).
\end{align}
Here  $r^{X_G}$ is the Riemannian scalar curvature of $(TX_G,g^{TX_G})$,
$\Delta_{X_G}$ is the Bochner-Laplacian on $X_G$ {\rm (}cf. \eqref{r6}{\rm )},
and $\{w^0_j\}$ is an orthonormal basis of $T^{(1,0)}X_G$.
\end{thm}

Since the non-equivariant version of this result has
already played a crucial role in the work of Donaldson mentioned
before, we have reason to believe that Theorem \ref{t0.6} might 
also play a role in the study of stability of projective
manifolds. Indeed, as Donaldson usually interprets his results in
the framework of geometric quantization, this seems likely to be so.

We recover \eqref{a0.3} from \eqref{0.15} after taking  the trace,
and then the integration  on $X_G$. Thus  \eqref{0.15} is a local
version of
 \eqref{a0.3} in the spirit of the Local Index Theory. The appearance
of the term $\frac{3}{4\pi}\Delta_{X_G}\log \wi{h}$
is unexpected.\\

Let $T$ be the torsion of the connection ${^0\nabla}^{TX}$ in
\eqref{h01} on $U$. The curvature $\Theta$ of the principal bundle
$U\to  B$ relates to the torsion $T$ by \eqref{h3}.

 Following \eqref{0g4} and \eqref{a33.16}, we choose $\{e^\bot_j\}$
to be an orthonormal basis of $N_{G,x_0}$ and
$\{\tfrac{\partial}{\partial z^0_j}\}\in T_{x_0}^{(1,0)}X_G$ 
to be the holomorphic basis of the normal
coordinate on $X_G$, and define $\mT_{klm}, \wi{\mT}_{jkl}$ as in
\eqref{34.21}. In particular, by Remark \ref{r4.4},
 $\wi{\mT}_{jkl}=0$ if $G$ is abelian.

 The $G$-invariant section $\wi{\mu}^E$ of $TY\otimes \End(E)$ on $U$
is defined by \eqref{h6} and \eqref{ah7}.

 If there is no other
specific notification in the next formula \eqref{4.74},  when we
meet the operation $|\quad|^2$,
 we will first do this operation, then take the sum of the indices.

\begin{thm}\label{t4.15} Under the assumption of Theorem
\ref{t0.6},  for $p>0$ large enough, $P^G_p(x, x)\in \End(E_x)$
and $P^{(r)}_{x_0}(0,0)\in \End(E_{x_0})$. Moreover,
\begin{multline}\label{4.74}
P^{(2)}_{x_0}(0,0)=2 ^{\frac{n_0}{2}}\left\{\frac{1}{8\pi} r^{X_G}_{x_0}
 +\frac{1}{\pi} R^{E_G}(\tfrac{\partial}{\partial z^0_j},
\tfrac{\partial}{\partial \ov{z}^0_j}) + \frac{1}{\pi}\Delta_{X_G}
\log \widetilde{h}
\right. \\
-\frac{3}{8\pi} \nabla_{e^\bot_k}\nabla_{e^\bot_k}\log h
- \frac{2}{\pi}\sqrt{-1}  \nabla_{JT(\tfrac{\partial}{\partial z^0_j},
\tfrac{\partial}{\partial \ov{z}^0_l})} \log h
- \frac{3}{\pi}\Big|\nabla_{\tfrac{\partial}{\partial \ov{z}^0_j}}
 \log h\Big|^2\\
- \frac{5}{4\pi}\Big|\nabla_{e^\bot_j} \log h\Big|^2
+ \frac{1}{2 \pi}
|T(e^\bot_k,\tfrac{\partial}{\partial \ov{z}^0_j})|^2
-\frac{1}{2\pi} \Big| \sum_j T(\tfrac{\partial}{\partial z^0_j},
\tfrac{\partial}{\partial \ov{z}^0_j})\Big|^2 \\
+\frac{1}{2\pi} \Big| T(\tfrac{\partial}{\partial z^0_i},
\tfrac{\partial}{\partial \ov{z}^0_j})\Big|^2
+ \frac{1}{24 \pi}\mT_{klm}^2
+ \frac{1}{2^6 \pi} \wi{\mT}_{ijk} (3 \wi{\mT}_{kji} - \wi{\mT}_{ijk})
+ \frac{1}{2\pi}\left\langle \wi{\mu}^E_{x_0},
\wi{\mu}^E_{x_0}\right\rangle_{g^{TY}} \\
+ \frac{1}{\pi} \left\langle \wi{\mu}^E,
T(\tfrac{\partial}{\partial z^0_l},
\tfrac{\partial}{\partial \ov{z}^0_l})\right\rangle
+\frac{3\sqrt{-1}}{2\pi}\left\langle \wi{\mu}^E,  J e^\bot_k \right\rangle
\nabla_{e^\bot_k}\log h
\left.+ \frac{\sqrt{-1}}{4 \pi} \left\langle J e^\bot_k,
\nabla^{TY}_{e^\bot_k} \wi{\mu}^E\right\rangle
 \right\}.
\end{multline}
\end{thm}

\begin{rem} \label{at0.6}
Certainly, if we only assume that $\bJ=J$ on $U$, a neighborhood
of $P=\mu^{-1}(0)$, then we still have  $\Phi_r(x_0)\in  \End
(E_G)_{x_0}$, as we work on the kernel of the Dirac operator
$D_p$. Set $\mathscr{I}_{p,0}= I_{\bC \otimes E_G}\mathscr{I}_p
I_{\bC \otimes E_G}$, the component of $\mathscr{I}_p$ on $\bC
\otimes E_G$. As the computation is local, we still have Theorem
\ref{t0.6} with $\mathscr{I}_p$ replaced by $\mathscr{I}_{p,0}$
and $\mathscr{I}_p-\mathscr{I}_{p,0}=\cO(p^{-\infty})$ (cf.
\eqref{3.19}). If we only work on the $\ov{\partial}$-operator,
 i.e. the holomorphic sections, in Section \ref{s8.1},
we explain how to reduce the case of  general $\bJ$ to the case
$\bJ=J$. Same remark holds for $P^G_p(x_0,x_0)$.
\end{rem}


Let $i:P\hookrightarrow X$ be the natural injection.

Let $\pi_G:  \cC^{\infty}(P, L^p\otimes E)^G
 \to \cC^{\infty} (X_G,L^p_G\otimes E_G)$ be the natural identification.

By a result of Zhang \cite[Theorem 1.1 and Proposition 1.2]{Z99},
one sees that for $p$ large enough, the map
$$ \pi_G\circ i^*: \cC^{\infty}(X, L^p \otimes E)^G
 \to \cC^{\infty} (X_G,L^p_G\otimes E_G) $$
induces a natural isomorphism
\begin{align}\label{60.23}
\sigma_p=  \pi_G \circ i^* :  H^0(X, L^p\otimes E)^G
\to H^0(X_G, L^p_G\otimes E_G).
\end{align}
(When $E=\bC$, this result was first proved in \cite[Theorem
3.8]{GuSt82}.)

The following result is a symplectic version of the above isomorphism
which is proved in Corollary \ref{at6.11}, 
 as a simple application of the Toeplitz
operator type properties proved in that subsection.
It might  be regarded as an ``asymptotic symplectic
quantization identification'', generalizing the corresponding
holomorphic identification \eqref{60.23}.

\begin{thm}\label{at100.6} If $X$ is a compact symplectic manifold
and $\bJ=J$, then the natural map 
$\sigma_p: (\Ker D_p)^G \to \Ker D_{G,p}$ defined in \eqref{6.45}
is an isomorphism for $p$ large enough.
\end{thm}



Let $\langle \, , \, \rangle_{L^p_G\otimes E_G}$ be the metric on
$L^p_G\otimes E_G$ induced by $h^{L_G}$ and $h^{E_G}$.

In view of \cite[(3.54)]{TZ98},
the natural Hermitian product on $\cC^{\infty} (X_G,L^p_G\otimes E_G)$
is the following weighted Hermitian product
$\langle \, , \, \rangle_{\widetilde{h}}$:
\begin{align}\label{0.26}
\langle s_1,s_2\rangle_{\widetilde{h}}
=\int_{X_G} \langle s_1, s_2\rangle_{L^p_G\otimes E_G} (x_0)
\widetilde{h}^2 (x_0)\,
dv_{X_G}(x_0).
\end{align}
In fact,  $\pi_G:  (\cC^{\infty}(P, L^p\otimes E)^G,\langle \, , \, \rangle)
 \to (\cC^{\infty} (X_G,L^p_G\otimes E_G),
\langle \, , \, \rangle_{\widetilde{h}})$ is an isometry.

We still denote by $\langle \quad \rangle $ the scalar product on
$ H^0(X, L^p\otimes E)^G$ induced by \eqref{60.22}.

\begin{thm}\label{t0.9} The isomorphism
$(2 p) ^{-\frac{n_0}{4}}\sigma_p$ is an asymptotic isometry from
$(H^0(X, L^p\otimes E)^G, \langle \, , \,\rangle)$ onto $(H^0(X_G,
L^p_G\otimes E_G), \langle \, , \, \rangle_{\widetilde{h}})$,
 i.e.,
if  $\{s_i^p\}_{i=1}^{d_p}$ is an orthonormal basis of
$(H^0(X, L^p\otimes E)^G, \langle \, ,\,\rangle)$, then
\begin{align}\label{0.27}
(2p)^{-\frac{n_0}{2}} \langle\sigma_p s_i^p,\sigma_p s_j^p\rangle
_{\widetilde{h}} = \delta_{ij} +\cO\left(\frac{1}{p}\right).
\end{align}
\end{thm}

From the explicit formula \eqref{4.74}, one can also get the
coefficient of $p^{-1}$ in the expansion \eqref{0.27}. We leave it
to the interested readers.


\comment{Let $\overline{\partial}^{L^p_G\otimes E_G,*}$
be the formal adjoint of $\overline{\partial}^{L^p_G\otimes E_G}$
associated to $g^{TX_G}$, $h^{L_G}$, $h^{E_G}$. Set
\begin{align}\label{0.16}
\widetilde{D}_{X_G,p}=\sqrt{2}\left(\widetilde{h}\overline{\partial}^{L^p_G\otimes
E_G}\widetilde{h}^{-1}+\widetilde{h}^{-1}\overline{\partial}^{L^p_G\otimes
E_G,*}\widetilde{h}\right)
\end{align}
be the deformed \spin Dirac operator on $X_G$ constructed in
\cite[(3.54)]{TZ98}, which appears there naturally through the
consideration of geometric quantization.
Let $\wi{P}^{X_G}_{p}$ denote the orthogonal projection from
$\Omega^{0,*}(X_G,L_G^p\otimes E_G)$ to $\Ker
\widetilde{D}_{X_G,p}$. Let $\wi{P}^{X_G}_{p}(x_0,x_0')$
 $(x_0,\, x_0'\in X_G)$
denote the corresponding Bergman kernel with respect to $dv_{X_G}(x_0')$.

The following result is a variant of \cite[Theorem 1.3]{DLM04a}.
}

Let $\wi{P}^{X_G}_{p}$ denote the orthogonal projection from
$(\cC^\infty (X_G,L_G^p\otimes E_G),
\langle \, , \, \rangle_{\widetilde{h}})$ onto $H^0(X,L_G^p\otimes E_G)$.
Let $\wi{P}^{X_G}_{p}(x_0,x_0')$
 $(x_0,\, x_0'\in X_G)$ be the smooth kernel of the operator
$\wi{P}^{X_G}_{p}$ with respect to $\widetilde{h}^2(x_0')dv_{X_G}(x_0')$.

The following result is an easy consequence of  \cite[Theorem
1.3]{DLM04a}.
\begin{thm} \label{t0.7} Under the assumption of Theorem
\ref{t0.6}, there
   exist smooth coefficients $\widetilde{\Phi}_r(x_0)\in \End (E_G)_{x_0}$
 which
are polynomials in $R^{TX_G}$, $R^{E_G}$ {\rm (}resp. $\wi{h}${\rm )},
  and their derivatives at $x_0$ to order $2r-1$ {\rm (}resp. $2r${\rm )},
 and $\widetilde{\Phi}_0= {\rm Id}_{E_G}$,
such that for any $k,l\in \bN$, there exists $C_{k,l}>0$ such that
for any $x_0\in X_G$, $p\in \bN$,
\begin{align}\label{0.22}
&\Big |p^{-n+n_0}\widetilde{h}^2(x_0)\wi{P}^{X_G}_{p}(x_0,x_0)- \sum_{r=0}^{k}
\widetilde{\Phi}_r(x_0) p^{-r} \Big |_{\cC^l} \leq C_{k,l} p^{-k-1}.
\end{align}
Moreover, the following identity holds,
\begin{align}\label{0.25}
\widetilde{\Phi}_1(x_0)= \frac{1}{8\pi} r^{X_G}_{x_0} + \frac{1}{2\pi}
\Delta_{X_G}\log \widetilde{h}
+\frac{1}{2\pi}R^{E_G}_{x_0}(w^0_j, \ov{w}^0_j) .
\end{align}
\end{thm}

\begin{rem} From (\ref{0.15}) and (\ref{0.25}), one sees that in general
$\Phi_1\neq \widetilde{\Phi}_1$, if $\wi{h}$ is not constant on $X_G$.
This reflects a subtle difference between the Bergman kernel
and the geometric quantization.
\end{rem}

From the works \cite{DLM04a}, \cite{MM04a} and the
present paper, we see clearly that the asymptotic expansion of
Bergman kernel is parallel to the small time asymptotic expansion
of the heat kernel. To localize the problem, the spectral gap
property \eqref{main1}
 and the finite propagation speed of solutions of hyperbolic equations play
essential roles.

Let $U$ be a $G$-neighborhood of $\mu^{-1}(0)$ as in Theorem
\ref{t0.1}, in this paper, we will then work on $U/G$.

Indeed, after doing suitable  rescaling on the coordinate, we get
the limit operator $\cL^0_2$ (cf. \eqref{g9}) which is the sum of
two parts, along $T_{x_0}X_G$, its kernel is infinite dimensional
and gives us the classical Bergman kernel as in $\bC^{n-n_0}$,
while along $N_G$, it is a harmonic oscillator and its kernel is
one dimensional. This explains well why we can expect  to get the
fast decay estimate along $N_G$ in \eqref{0.8}.

This paper is organized as follows. In  Section \ref{s1},
we study connections and Laplacians associated to a principal bundle.
  In  Section \ref{s3}, we localize the problem by using
the spectral gap property and finite propagation speed,
then we use the rescaling technique in local index theory to prove
Theorem \ref{tue17} which is a version of Theorem \ref{t0.1}
without assumption on $\bJ$. We assume $G$ acts freely on
$P=\mu^{-1}(0)$ in Sections \ref{0s3.2}-\ref{s3.5}, and in Section
\ref{s6.1} we explain Theorem \ref{t6.1}, the version of Theorem
\ref{t0.1} where we only assume that $\mu$ is regular at $0$.
 In  Section \ref{s30}, we get explicit informations on the coefficients
$P^{(r)}$ when $\bJ$ verifies \eqref{g2},
thus we get an effective way to compute
its first coefficients of the asymptotic expansion \eqref{0.8}.
Especially, we establish \eqref{0.7} and \eqref{a0.7}.
 In  Section \ref{s6}, we explain various applications of
our Theorem \ref{t0.1}, including Toeplitz properties, etc.
 In  Section \ref{s4}, we  compute the coefficients
$\Phi_1$ in Theorem \ref{t0.6} and
$P^{(2)}_{x_0}(0,0)$ in Theorem \ref{t4.15}
and in the general case: $\bJ\neq J$.
In  Section \ref{s8}, we prove Theorems  \ref{t0.9}, \ref{t0.7}.

Some results of this paper have been announced
in \cite{MZ05a, MZ05e}.

Notation : In the whole paper, if there is no other specific notification,
when in a formula a subscript index appears two times,  we sum up
with this index.

\newpage

\section{Connections and Laplacians associated to a principal bundle}\label{s1}

In this Section, for $\pi : X\to B=X/G$ a $G$-principal bundle, we
will study the associated connections and Bochner-Laplacians. The
results in this Section extend the corresponding ones in \cite[\S
1d)]{B86} and \cite[\S 5.1, 5.2]{BeGeVe} where the  metric along
the fiber is parallel along the horizontal direction.
These results will be used in Proposition \ref{0t3.2}
and in Sections \ref{s3.111}, \ref{s4}.

If $G$ acts only infinitesimal freely on $X$, then $B=X/G$ is an
orbifold. The results in this Section can be extended easily to
this situation, as will be  explained in Section \ref{s6.1}.

This Section is organized as follows. In  Section \ref{s4.2}, we
study the Levi-Civita connection for a principal bundle which
extends the results of  \cite[\S 1d)]{B86}. In  Section
\ref{as4.2}, we study the relation of the  Laplacians on the total
and base manifolds.

\subsection{Connections associated to a principal bundle}\label{s4.2}
Let a compact connected Lie group $G$ acts smoothly on the left on
a smooth manifold $X$ and $\dim G=n_0$. We suppose temporary that
$G$ acts freely on $X$. Then  $$\pi : X\to B=X/G$$ is a
$G$-principal bundle.
We denote by $TY$ the relative tangent bundle for the fibration
$\pi : X\to B$.

Let $g^{TX}$ be a $G$-invariant metric on $TX$. Let $\nabla ^{TX}$
be the Levi-Civita connection on $TX$.
By the explicit equation for $\left\langle\nabla ^{TX}_\cdot
\cdot,\cdot \right\rangle$ in \cite[(1.18)]{BeGeVe}, for
$W,Z,Z^{\prime}$ vector fields on $X$,
\begin{multline} \label{ah1}
2\left\langle{\nabla}^{TX}_{W} Z,Z^{\prime} \right\rangle=
W \left\langle Z,Z^{\prime} \right\rangle +Z \left\langle W, Z^{\prime} \right\rangle
-Z^{\prime} \left\langle W, Z\right\rangle \\
-\left\langle W, [Z,Z^{\prime}]\right\rangle
-\left\langle Z, [W, Z^{\prime}] \right\rangle + \left\langle Z^{\prime},
[W,Z]\right\rangle.
\end{multline}

Let $T^H X$ be the orthogonal complement of $TY$ in $TX$.

For $U\in TB$, let $U^H\in T^HX$ be the lift of $U$.

Let $g^{TY}, g^{T^HX}$ be $G$-invariant metrics on
$TY, T^HX$ induced by $g^{TX}$.
Let $P^{TY}$, $P^{T^HX}$ be the orthogonal projections from $TX$ onto
$TY$, $T^HX$.

Let $g^{TB}$ be the metric on $TB$ induced by $g^{T^HX}$.
Let $\nabla ^{TB}$ be the Levi-Civita connection on $(TB, g^{TB})$
with curvature $R^{TB}$.
Set
 \begin{align}\label{h01}
\nabla ^{T^HX}=\pi ^*\nabla ^{TB},\quad \nabla ^{TY}= P^{TY}\nabla ^{TX}P^{TY},
\quad {^0\nabla}^{TX}= \nabla ^{TY}\oplus \nabla ^{T^HX}.
\end{align}
 Then $\nabla ^{T^HX}$, ${^0\nabla}^{TX}$ define Euclidean connections
on $T^HX$, $TX$, and $\nabla ^{TY}$ is the connection on $TY$
induced by $\nabla ^{TX}$ (cf. \cite[Def. 1.6]{B86}).
\comment{
defines a connection on $T^HX$.
}

Let $T$ be the torsion of  ${^0\nabla}^{TX}$, and let
$S\in T^*X \otimes \End(TX)$, $\dot{g}^{TY}_\cdot\in T^*B\otimes \End(TY)$
be defined by
 \begin{align}\label{h0}
&S= \nabla ^{TX}-{^0\nabla}^{TX},
&\dot{g}^{TY}_{U}= (g^{TY})^{-1} (L_{U^H}g^{TY}) \quad {\rm for}
\,\, U\in TB.
\end{align}
Then $S$ is a $1$-form on $X$ taking values in the skew-adjoint
endomorphisms of $TX$.

By \cite[Theorem 1.2]{BGS88b} (cf. \cite[Theorems 1.1 and 1.2]{B97})
the proof of which can also be found in \cite[Prop. 10.2]{BeGeVe}
where one applies directly \eqref{ah1}, we know that
$\nabla ^{TY}$ is the Levi-Civita connection on $TY$ along the fiber $Y$, and
for $U\in TB$,
\begin{align}\label{h1}
\nabla ^{TY}_{U^H}= L_{U^H} + \frac{1}{2} (g^{TY})^{-1}(L_{U^H}g^{TY})
= L_{U^H} + \frac{1}{2}\dot{g}^{TY}_{U} .
\end{align}

Let $\kg$ be the Lie algebra of $G$.
For $K\in \kg$, we denote by $K^X_x= \frac{\partial}{\partial t}e
^{-tK}x|_{t=0}$ the corresponding vector field on $X$, then $g
K^X_x= ({\rm Ad}_g (K))^X_{gx}$. Thus we can identify the trivial
bundle $X\times \kg$ with ${\rm Ad}$-action of $G$ on $\kg$ to the
$G$-equivariant bundle $TY$ by the map $K\to K^X$.

 Let $\theta:
TX\to \kg$ be the connection form of the principal bundle $\pi:
X\to B$ such that $T^HX=\Ker \theta$, and $\Theta$ its curvature.

\comment{
$Y_1,Y_2$ two $G$-invariant sections of $TY$,
then $[Y_i,K_1^X ]=0$, $L_{Y_i}g^{TY}=0$ for $i=1,2$.
As $g^{TX}$ is $G$-invariant,
by \eqref{ah1}, we get
\begin{align}\label{ah2}
2 \left\langle \nabla ^{TY}_{Y_1}Y_2, K^X_1\right\rangle
= Y_1\left\langle Y_2, K^X_1\right\rangle
+Y_2\left\langle Y_1, K^X_1\right\rangle
+ \left\langle [Y_1, Y_2], K^X_1\right\rangle
= \left\langle [Y_1, Y_2], K^X_1\right\rangle .
\end{align}
By \eqref{ah2}, we get
&\nabla ^{TY}_{Y_1}Y_2 =\frac{1}{2}[Y_1,Y_2];
}

For $K_1,K_2\in \kg$, $U,V\in TB$,  as $U^H$ is $G$-invariant, we have
\begin{align}\label{h2b}
&L_{U^H}K_1^X= -[K_1^X, U^H]=0.
\end{align}
By (\ref{h1}), (\ref{h2b}), we get $T\in \Lambda ^2(T^*X)\otimes TY$ and
\begin{align}\label{h3}
&T(U^H, V^H)=\Theta(U^H, V^H)= -P^{TY}[U^H, V^H],
\quad   T(K_1^X,K_2^X ) =0, \\
&T(U^H, K_1^X) =  \frac{1}{2} (g^{TY})^{-1}(L_{U^H}g^{TY}) K_1^X
=\frac{1}{2}\dot{g}^{TY}_{U}K_1^X.\nonumber
\end{align}
And by \eqref{ah1}, (\ref{h1}), (\ref{h2b}) and (\ref{h3}), for $W\in TX$, 
we have (cf. also \cite[(1.28)]{B86}, \cite[Prop. 10.6]{BeGeVe}),
\begin{align}\label{h4}
S(W)(TY)\subset T^HX,& \quad S(U^H)V^H\in TY,\nonumber\\
2 \left\langle S(U^H)K_1^X, V^H  \right\rangle
&= 2 \left\langle S(K_1^X)U^H, V^H\right\rangle
= \left\langle T(U^H, V^H), K_1^X\right\rangle,\\
\left\langle S(K_2^X)U^H, K_1^X \right\rangle
&=-\left\langle S(K_2^X)K_1^X, U^H  \right\rangle \nonumber\\
&=  \frac{1}{2}U^H \left\langle  K_2^X,K_1^X \right\rangle
=\left\langle  T(U^H, K_1^X),K_2^X \right\rangle .\nonumber
\end{align}

Let $\{e_i\}$ be an orthonormal basis of $TB$. By (\ref{h0}) and (\ref{h4}),
for $Y$ a section of $TY$,
 \begin{align}\label{ah4}
\nabla ^{TX}_{U^H}Y= \nabla ^{TY}_{U^H}Y
+\frac{1}{2} \left\langle T(U^H, e^H_i), Y\right\rangle e^H_i.
\end{align}
\begin{prop}\label{p4.01} Let $\{f_l\}_{l=1}^{n_0}$
 be a $G$-invariant orthonormal frame of $TY$, then
\begin{align}\label{bh4}
\sum_{l=1}^{n_0} \nabla ^{TY}_{f_l} f_l=0.
\end{align}
\end{prop}
\begin{proof} \eqref{bh4} is analogue to the fact that
any left invariant volume form on $G$ is also right invariant.
We only need to work on a fiber $Y_b$, $b\in B$.

Let $dv_Y$ be the Riemannian volume form on $Y_b$.

By using $L_{f_k}f_l= \nabla ^{TY}_{f_k}f_l- \nabla ^{TY}_{f_l}f_k$
and $dv_Y$ is preserved by $\nabla^{TY}$ on $Y_b$, we get
\begin{align}\label{bh5}
L_{f_k} dv_Y = \sum_{l=1}^{n_0} \left\langle \nabla ^{TY}_{f_l} f_k,
 f_l\right\rangle dv_Y.
\end{align}
Now from $L_{f_k}= i_{f_k} d^Y +  d^Y i_{f_k}$ and
$ \left\langle \nabla ^{TY}_{f_l} f_k, f_l\right\rangle$
is $G$-invariant and \eqref{bh5}, we get
\begin{align}\label{bh6}
0= \int_{Y_b} L_{f_k} dv_Y
= \sum_{l=1}^{n_0} \left\langle \nabla ^{TY}_{f_l} f_k,
 f_l\right\rangle \int_{Y_b} dv_Y.
\end{align}
From \eqref{bh6}, we get  \eqref{bh4}.
\end{proof}

\begin{rem} \label{p4.02} If $g^{TY}$ is induced by a family of
$Ad$-invariant metric on $\kg$ under the isomorphism from $X\times \kg$
to $TY$ defined by $K\to K^X$, then \eqref{bh4} is trivial.
In this case, as in \cite[Theorem 11.3]{Fegan91},
 for $Y_1,Y_2$ two $G$-invariant sections of $TY$, by \eqref{ah1}, we have
\begin{align}\label{bh7}
\nabla ^{TY}_{Y_1} Y_2 = \frac{1}{2} [Y_1,Y_2].
\end{align}
\end{rem}

\subsection{Curvatures and Laplacians
 associated to a principal bundle}\label{as4.2}

Let $(F,h^F)$ be a $G$-equivariant Hermitian  vector bundle on $X$
with a $G$-invariant Hermitian connection $\nabla ^F$ on $X$. For
any $K\in \kg$, denote by $L_K$ the infinitesimal action induced
by $K$ on the corresponding vector bundles.

Let $\mu ^F$ be the section of $\kg ^* \otimes \End(F)$ on $X$
 defined by,
\begin{align}\label{h6}
\mu ^F(K)= \nabla ^F_{K^X}  - L_K \quad \, {\rm for} \, \, K\in \kg.
\end{align}
By using the identification $X\times \kg \to TY$, $\mu^F$ defines a
$G$-invariant section $\wi{\mu}^F$ of $TY \otimes \End(F)$ on $X$
such that
 \begin{align}\label{ah7}
\left\langle \wi{\mu} ^F, K^X\right\rangle = \mu ^F (K).
\end{align}

The curvature $R^F_\mu$ of the Hermitian connection $\nabla ^F-\mu
^F(\theta)$ on $F$ is $G$-invariant. Moreover as $\nabla ^F$ is
$G$-invariant, by \eqref{h6},
 \begin{align}\label{ah8}
 R^F_\mu(K^X,v)=[L_K,\nabla ^F-\mu ^F(\theta)] (v)=0
\end{align}
for $K\in \kg$, $v\in TX$, and
\begin{align}\label{h7}
R^F_\mu = R^F-\nabla ^F (\mu ^F(\theta))
+\mu^F(\theta)\wedge \mu^F(\theta).
\end{align}

The Hermitian vector bundle $(F,h^{F})$ induces a Hermitian vector
bundle $(F_{B},h^{F_{B}})$ on $B$ by identifying $G$-invariant
sections of $F$ on $X$.

For $s\in \cC^\infty(B, F_{B})\simeq \cC^\infty(X, F) ^G$, we
define
\begin{align}\label{h5}
\nabla ^{F_{B}}_U s= \nabla ^F_{U^H}s.
\end{align}
Then  $\nabla ^{F_{B}}$ is a Hermitian connection on $F_{B}$ with
curvature $R^{F_{B}}$.

Observe that $\nabla ^{F_{B}}$ is the restriction of the
connection
 $\nabla ^F-\mu ^F(\theta)$ to $\cC^\infty(X, F) ^G$,
and $R^{F_{B}}$ is the section induced by
$R^F_\mu$.
From (\ref{h7}), for $U_1,U_2\in TB$, we get
\begin{align}\label{h8}
R^{F_{B}}(U_1,U_2) = R^F(U^H_1,U^H_2) -\mu ^F (\Theta)(U_1,U_2) .
\end{align}

Let $dv_X$ be the  Riemannian volume form on $(X, g^{TX})$.
We define a scalar product on $\cC^{\infty}  (X, F)$ by
\begin{align}\label{h10}
\langle s_1,s_2 \rangle =\int_{X}\langle s_1,
s_2\rangle_{F}(x)\,dv_{X}(x)\,.
\end{align}
As in \eqref{h10}, $h^{F_B}$, $g^{TB}$ induce a natural scalar product
 $\langle \quad\rangle$ on $\cC^{\infty}  (B, F_B)$.

Denote by $\vol (Gx)\ (x\in X)$ the volume of the orbit $Gx$
equipped with the metric induced by $g^{TX}$.
The function $$h (x)=\sqrt{\vol (Gx)},\ \ \ \ \ \  x\in X,$$ as in
\eqref{0.6}
 is $G$-invariant and defines a function on $B$.

Denote by $\pi_G: \cC^{\infty} (X,F)^G \to \cC^{\infty} (B,F_{B})$
the natural identification. Then the map
\begin{align}\label{h12}
\Phi = h \pi_G:
(\cC^{\infty}  (X, F)^G,\langle \, ,\,  \rangle)
 \to (\cC^{\infty}  (B, F_B), \langle \, ,\,  \rangle)
\end{align}
is an isometry.

Let $\{e_a\}_{a=1}^m$ be an orthonormal frame of $TX$.

Let $(E,h^E)$ be a Hermitian vector bundle on $X$
and let $\nabla  ^E$ be a Hermitian connection on $E$.
The usual Bochner Laplacians $\Delta  ^E, \Delta_X$ are defined by
\be\label{r6} \Delta  ^E:=-\sum_{a=1}^m \left ((\nabla ^E_{e_a})^2
-\nabla  ^E_{ \nabla ^{TX}_{e_a}e_a}\right ),\quad \Delta_X=
\Delta  ^\bC. \ee

Let $\{f_l\}_{l=1}^{n_0}$ be a $G$-invariant orthonormal frame of $TY$,
and $\{f^l\}$ its dual basis,
and let $\{e_i\}$ be an orthonormal frame of $TB$,
then $\{e^H_i, f_l\}$ is an orthonormal frame of $TX$.

To simplify the notation,
 for $\sigma_1,\sigma_2\in TY\otimes \End(F)$, we denote by
$\langle \sigma_1,\sigma_2\rangle_{g^{TY}}\in \End(F)$ the contraction of
$\sigma_1\otimes \sigma_2$ on the part of $TY$ by $g^{TY}$.
In particular,
\begin{align}\label{ah12}
\langle\wi{\mu}^{F}, \wi{\mu}^{F}\rangle_{g^{TY}}
=  \sum_{l=1}^{n_0}\langle\wi{\mu}^{F}, f_l\rangle^2 \in \End(F).
\end{align}

The following result extends \cite[Prop. 5.6, 5.10]{BeGeVe}
where $F=X\times_G V$ for a $G$-representation $V$,
and where $g^{TY}$ is induced by a fixed $Ad$-invariant metric on $\kg$
under the isomorphism from $X\times \kg$ to $TY$ defined by $K\to K^X$
(Thus $h$ is constant on $B$).

\begin{thm}\label{t4.1} As an operator on $\cC^{\infty}  (B, F_B)$, we have
\begin{align}\label{h13}
\Phi \Delta^F \Phi^{-1} =  \Delta^{F_B}
- \langle\wi{\mu}^{F}, \wi{\mu}^{F}\rangle_{g^{TY}} - \frac{1}{h}\Delta_B h.
\end{align}
\end{thm}
\begin{proof} At first by \eqref{h3} and \eqref{h4},
\begin{multline}\label{h14}
\frac{1}{h}(e_i h) = \frac{1}{2}({L_{e^H_i}dv_Y})/{dv_Y}
= \frac{1}{2}\left\langle L_{e^H_i}f^l, f^l \right\rangle
= -\frac{1}{2}\left\langle L_{e^H_i}f_l, f_l \right\rangle\\
=\frac{1}{4}( L_{e_i^H}g^{TY})(f_l, f_l)
=\frac{1}{2}\left\langle T(e_i^H, f_l), f_l \right\rangle
= -\frac{1}{2}\left\langle S(f_l)f_l,e_i^H \right\rangle.
\end{multline}

 As $\wi{\mu}^{F}$ is $G$-invariant, then $\langle\wi{\mu}^{F}, f_l\rangle$
is also a $G$-invariant section of $\End(F)$.

By \eqref{h6}, $\nabla ^{F}_{f_l}= \langle\wi{\mu}^{F}, f_l\rangle$ on
$\cC^{\infty} (X,F)^G$, and by \eqref{h0},
$\nabla ^{TX}_{f_l}f_l =\nabla ^{TY}_{f_l}f_l + S(f_l)f_l$,
 thus by  \eqref{h12}, we get for $1\leq l\leq n_0$,
\begin{align}\label{h15}
\Phi  [(\nabla ^{F}_{f_l})^2
- \nabla ^{F}_{\nabla ^{TX}_{f_l}f_l} ] \Phi^{-1}
= \langle\wi{\mu}^{F}, f_l\rangle^2
-  \langle\wi{\mu}^{F}, \nabla ^{TY}_{f_l}f_l\rangle
- h \nabla ^{F_B}_{S(f_l)f_l} h^{-1} .
\end{align}
From \eqref{h4}, \eqref{bh4}, \eqref{r6}, \eqref{ah12}, \eqref{h14}
and \eqref{h15}, we have
\begin{multline}\label{h16}
\Phi \Delta ^{F}\Phi^{-1} =
- \sum_{i=1}^{2n-n_0}\Phi \Big [(\nabla ^{F}_{e_i^H})^2
- \nabla ^{F}_{\nabla ^{TX}_{e_i^H}e_i^H}\Big ] \Phi^{-1}
- \sum_{l=1}^{n_0}\Phi \Big [(\nabla ^{F}_{f_l})^2
- \nabla ^{F}_{\nabla ^{TX}_{f_l}f_l}\Big ] \Phi^{-1}\\
 = h \Delta^{F_B} h^{-1}- \sum_{l=1}^{n_0} \langle\wi{\mu}^{F}, f_l\rangle^2
- 2(e_i h) \nabla ^{F_B}_{e_i} h^{-1}
= \Delta^{F_B}-\langle\wi{\mu}^{F}, \wi{\mu}^{F}\rangle_{g^{TY}}
- \frac{1}{h}\Delta_B h.
\end{multline}
\end{proof}


\section{$G$-invariant Bergman kernels}\label{s3}

In this Section, we study the uniform estimate with
its derivatives on $t=\frac{1}{\sqrt{p}}$ of  the $G$-invariant
 Bergman kernel $P^G_p(x,x')$  of $D^2_p$ as $p\to \infty$.

The first main difficulty is to  localize the problem to arbitrary
small neighborhoods of $P=\mu ^{-1}(0)$, so that one can study the
$G$-invariant Bergman kernel in the spirit of \cite{DLM04a}. Our
observation here is that the $G$-invariant Bergman kernel is
exactly the kernel of the orthogonal projection on the zero space
of an operator $\mL_{p}$, which is a deformation of $D_p^2$ by the
Casimir operator. Moreover, $\mL_{p}$ has a spectral gap property
(cf. \eqref{a11}, \eqref{a12}). In the spirit of \cite[\S
3]{DLM04a}, this allows us to localize the problem to a problem
near a $G$-neighborhood of $Gx$. By combining with the
Lichnerowicz formula, we get Theorem \ref{t0.0} in Section
\ref{s3.1}.

After localizing  the problem to a problem near $P$, we first
replace $X$ by $G\times \bR ^{2n-n_0}$, then we reduce it to a
problem on $\bR ^{2n-n_0}$. On $\bR ^{2n-n_0}$, the problem in
Section \ref{s3.3} is similar to a problem on $\bR^{2n}$
considered in \cite[\S 3.3]{DLM04a}.

 Comparing with  the operator in \cite[\S 3.3]{DLM04a},
we have an extra quadratic term along the normal direction of
$X_G$. This allows us to improve the estimate in the normal
direction. After suitable rescaling, we will introduce a family of
Sobolev norms defined by the rescaled connection on $L^p$ and the
rescaled moment map in this situation, then we can extend  the
functional analysis techniques developed in
 \cite[\S 3.3]{DLM04a} and \cite[\S 11]{BL91}.

This section is organized as follows. In Section \ref{s3.00}, we
recall a basic property on the Casimir operator of a compact
connected Lie group. In Section \ref{s3.01}, we recall the
definition of  spin$^c$ Dirac operators for an almost complex
manifold. In Section \ref{s3.0}, we introduce the operator
$\mL_{p}$ to study the $G$-invariant Bergman kernel $P^G_p$ of
$D_p^2$. In Section \ref{s3.1}, we explain that the asymptotic
expansion of $P^G_p(x,x')$ is localized on a $G$-neighborhood of
$Gx$, and we establish Theorem \ref{t0.0}. In Section \ref{0s3.2},
we show that our problem near $P$ is equivalent to a problem on
$U/G$ for any open $G$-neighborhood $U$ of $P$. In Section
\ref{s3.2}, we derive an asymptotic expansion of $\Phi \mL_p\Phi
^{-1}$ in coordinates of $U/G$. In Section \ref{s3.3}, we study
the uniform estimate with its derivatives on $t$ of the Bergman
kernel associated to the rescaled operator $\cL^t_2$ from $\Phi
\mL_p\Phi ^{-1}$ using heat kernel.
 In Theorem \ref{tue14}, we estimate
uniformly the remainder term of the Taylor expansion of  $e
^{-u\cL^t_2}$
 for $u\geq u_0>0,\ 0<t\leq t_0\leq1$. In Section \ref{s3.5}, we identify
 $J_{r,u}$,
the coefficient of the Taylor expansion of  $e ^{-u\cL^t_2}$, with
the Volterra expansion of the heat kernel, thus giving  a way to
compute the coefficient $P^{(r)}_{x_0}$ in  Theorem \ref{t0.1}. In
Section \ref{s3.5},  we prove Theorem \ref{t0.1} except
\eqref{0.7} and \eqref{a0.7}.

We use the notation  in Section \ref{s1}.
In Sections \ref{0s3.2}-\ref{s3.5}, we assume $G$ acts freely
on $P=\mu^{-1}(0)$.


\subsection{Casimir operator}\label{s3.00}

Let $G$ be a  compact connected Lie group with Lie algebra $\mathfrak{g}$ and $\dim G=n_{0}$.
We choose an Ad-invariant metric on $\kg$ such that it is the
minus Killing form on the semi-simple part of $\kg$.

Let $\{K_j\}_{j=1}^{\dim G}$ be an orthogonal basis of $\kg$ and $\{K^j\}$
be its dual basis of $\kg ^*$.

The Casimir operator $\Cas$ of $\kg$ is defined as the following
element of the universal enveloping algebra $U(\kg)$ of $\kg$,
\begin{align}\label{a1}
\Cas := \sum_{j=1}^{\dim G} K_jK_j.
\end{align}
Then $\Cas$ is independent of the choice of $\{K_j\}$
 and belongs to the center of $U(\kg)$.

Let $\kt$ be the Lie algebra of a maximum torus $T$ of $G$,
and $\kt^*$ its dual.
Let $|\quad|$ denote the norm on $\kt^*$ induced by
the Ad-invariant metric on $\kg$.

Let $\mW\subset \kt^*$ be the fundamental Weyl chamber associated to the set of
positive roots $\Delta ^ +$ of $G$, and
its closure $\overline{\mW}\subset \kt^*$.

Let $I=\{K\in \kt; \exp(2\pi K)=1\in T\}$ the integer lattice
such that $T=\kt/2 \pi I$,
 and $P= \{\alpha\in \kt^* ; \alpha(I)\subset  \bZ\}$
the lattice of integral forms.

Let $\varrho_G$ be the half sum of the positive roots of $G$.

By the Weyl character formula \cite[Theorem 8.21]{Fegan91},
the irreducible representations correspond one to one to
$\vartheta \in \ov{\mW}\cap P$, the highest weight of the representation.

Moreover, for any irreducible representation $\rho:G\to \End(E)$
with highest weight $\vartheta\in  \ov{\mW}\cap P$,
classically, the action of $\Cas$
on $E$ is given by (cf. \cite[Theorem 10.6]{Fegan91}),
\begin{align}\label{a2}
\rho(\Cas) = - (|\vartheta+\varrho_G|^2 - |\varrho_G|^2) \Id_ E.
\end{align}

Set
\begin{align}\label{a3}
\nu_1 :=\inf_{0\neq \vartheta \in \ov{\mW}\cap P}( |\vartheta+\varrho_G|^2
 - |\varrho_G|^2)>0.
\end{align}

By (\ref{a2}), for any representation $\rho:G\to \End(E)$,
if the $G$-invariant subspace $E^G$ of $E$ is zero, then
\begin{align}\label{a4}
-\rho(\Cas) \geq \nu_1 \Id_E.
\end{align}

\subsection{Spin$^c$ Dirac operator}\label{s3.01}

Let $(X,\om)$ be a compact symplectic manifold of real dimension
$2n$. Assume that there exists a Hermitian line bundle $L$ over
$X$ endowed with a Hermitian connection $\nabla^L$ with the
property that $$\frac{\sqrt{-1}}{2\pi}R^L=\omega,$$ where
$R^L=(\nabla^L)^2$ is the curvature of $(L,\nabla^L)$.

Let $(E,h^E)$ be a Hermitian vector bundle  on $X$ with Hermitian
connection $\nabla^E$ and its curvature $R^E$.

Let $g^{TX}$ be a Riemannian metric on $X$.

Let $\bJ :TX\longrightarrow TX$ be the skew--adjoint linear
map which satisfies the relation
\begin{align}  \label{0a0}
\om(u,v)=g^{TX}({\bJ}u,v)
\end{align}
 for $u,v \in TX$.

Let $J$ be an almost complex structure such that
\begin{align}  \label{0a01}
g^{TX}(Ju,Jv)= g^{TX}(u,v),\ \ \ \ \ \om(Ju,Jv)=\om(u,v),
\end{align}
and that $\om(\cdot,J\cdot)$ defines a metric on $TX$.
 Then $J$ commutes with $\bJ$ and
 $$-\left \langle J\bJ \cdot, \cdot \right \rangle= \om(\cdot,J \cdot
 )$$
is positive by our assumption. Thus $-J\bJ \in \End(TX)$ is
symmetric and positive, and one verifies easily that
\begin{align}  \label{0a00}
-J\bJ = (-\bJ ^2)^{1/2}, \quad J=\bJ (-\bJ ^2)^{-1/2}.
\end{align}

The almost complex structure $J$ induces a splitting
$$TX\otimes_\bR \bC=T^{(1,0)}X\oplus T^{(0,1)}X,$$ where
$T^{(1,0)}X$ and $T^{(0,1)}X$ are the eigenbundles of $J$
corresponding to the eigenvalues $\sqrt{-1}$ and $-\sqrt{-1}$
respectively.
Let $T^{*(1,0)}X$  and $T^{*(0,1)}X$ be the corresponding dual
bundles.

For any $v\in TX\otimes_\bR \bC$ with decomposition
$v=v_{1,0}+v_{0,1} \in T^{(1,0)}X\oplus T^{(0,1)}X$, let
${\overline v^\ast_{1,0}}\in T^{*(0,1)}X$ be the metric dual of
$v_{1,0}$. Then
\begin{align}  \label{0a2}
c(v):=\sqrt{2}({\overline v^\ast_{1,0}}\wedge- i_{v_{0,1}})
 \end{align}
defines the Clifford action of $v$ on $\Lambda (T^{*(0,1)}X)$,
where $\wedge$ and $i$ denote the exterior and interior
multiplications  respectively.

Set
\begin{equation}\label{0a1}
\nu_0:=\displaystyle\inf_{{u\in T_x^{(1,0)}X,\,x\in X}}
R^L_x(u,\overline{u})/|u|^2_{g^{TX}} >0.
\end{equation}

Let $\nabla^{TX}$ be the Levi-Civita connection
of the metric $g^{TX}$ with curvature $R^{TX}$.
We denote by  $P^{T^{(1,0)}X}$
  the projection from $T X\otimes_\bR \bC$ to $T^{(1,0)}X$.

Let $\nabla^{T^{(1,0)}X}=P^{T^{(1,0)}X}\,\nabla^{TX}P^{T^{(1,0)}X}$
 be the Hermitian connection on $T^{(1,0)}X$ induced by  $\nabla^{TX}$
with curvature  $R^{T^{(1,0)}X}$.

By \cite[pp.397--398]{LaMi89}, $\nabla^{TX}$ induces canonically
a Clifford connection $\nabla^{\text{Cliff}}$ on $\Lambda (T^{*(0,1)}X)$
and its curvature $R^{\text{Cliff}}$ (cf. also \cite[\S 2]{MM02}).

Let $\{e_a\}_a$ be an orthonormal basis of $TX$. Then
\begin{align}  \label{0a3}
R^{\text{Cliff}}=
\frac{1}{4}\sum_{ab}\langle R^{TX}e_{a} ,e_{b}\rangle c(e_{a})c(e_{b})
 +\frac{1}{2}\tr\left[R^{T^{(1,0)}X} \right].
\end{align}

Let $\nabla^{E_p}$ be the connection on
\begin{align}  \label{0a4}
E_p : =\Lambda (T^{*(0,1)}X)\otimes L^p\otimes E
\end{align}
induced by $\nabla^{\text{Cliff}}$, $\nabla^L$ and  $\nabla^E$.

Let $\langle \quad \rangle_{E_p}$ be the metric on $E_p$ induced by
 $g^{TX}$, $h^L$ and $h^E$.


The $L^2$--scalar product $\langle \quad\rangle$  on
$\Omega^{0,\scriptscriptstyle{\bullet}}(X,L^p\otimes E)$,
the space of smooth sections of $E_p$, is given by \eqref{h10}.
We denote the corresponding norm by $\norm{\cdot}_{L^ 2}$.

\begin{defn}\label{Dirac}
The \spin Dirac operator $D_p$ is defined by
\begin{equation}\label{defDirac}
D_p:=\sum_{a=1}^{2n}c(e_a)\nabla^{E_p}_{e_a}:
\Omega^{0,\scriptscriptstyle{\bullet}}(X,L^p\otimes E)\longrightarrow
\Omega^{0,\scriptscriptstyle{\bullet}}(X,L^p\otimes E)\,.
\end{equation}
\end{defn}
\noindent Clearly, $D_p$ is a formally self--adjoint, first order
elliptic differential operator on
$\Omega^{0,\scriptscriptstyle{\bullet}}(X,L^p\otimes{E})$, which
interchanges $\Omega^{0,\text{even}}(X,L^p\otimes E)$ and
$\Omega^{0,\text{odd}}(X,L^p\otimes E)$.

If $A$ is any operator, we denote by $\spec (A)$ the spectrum of $A$.

The following  result was proved in \cite[Theorems 1.1, 2.5]{MM02}
by applying directly the Lichnerowicz formula (cf. also
\cite[Theorem 1]{BVa89} in the holomorphic case).

\begin{thm}\label{t2.1}
There exists $C_L>0$ such that for any $p\in \bN$ and
any $s\in\Omega^{>0}(X,L^p\otimes E)
=\bigoplus_{q\geqslant 1}\Omega^{0,q}(X,L^p\otimes E)$,
\begin{equation}\label{main1}
\norm{D_{p}s}^2_{L^ 2}\geqslant(2p\nu_0-C_L)\norm{s}^2_{L^ 2}\,.
\end{equation}
Moreover  ${\spec }(D^2_p )\subset \{0\}\cup [2p\nu_0
-C_L,+\infty[$.
\end{thm}

\subsection{$G$-invariant Bergman kernel}\label{s3.0}

Suppose that the compact connected Lie group $G$
acts on the left of $X$, and the action of $G$ lifts on $L,E$ and
preserves the metrics and connections, $\om$
and the almost complex structure $J$.

Let $\mu : X \to \kg^*$ be defined by
\begin{align}\label{a6}
2 \pi \sqrt{-1}\mu(K) : =\mu ^L(K)
=\nabla ^L_{K^X}- L_K, \, \, K\in \kg.
\end{align}
Then $\mu$ is  the corresponding moment map
(cf. \cite[Example. 7.9]{BeGeVe}), i.e.  for any $K\in \kg$,
\begin{align}\label{a5}
d\mu(K)= i_{K^X}\om.
\end{align}

\comment{
Suppose that $(X,\omega)$ admits a Hamiltonian action of $G$.
If $K\in \kg$, let $K^X$ be the corresponding vector field on $X$.
Let $\mu: X\rightarrow {\kg}^*$ be the corresponding moment map
(cf. \cite[Example. 7.9]{BeGeVe}),
i.e.  for any $K\in \kg$,
\begin{align}\label{a5}
d\mu(K)= i_{K^X}\om.
\end{align}
Then a formula due to Kostant \cite{Kostant70}
 induces a natural $\kg $ action on $L$
\begin{align}\label{a6}
L_K s =\nabla ^L_{K^X} s -2 \pi \sqrt{-1}  \mu(K) s, \quad
s\in \cC^{\infty}(X,L), \, K\in \kg.
\end{align}
We make the assumption that this $\kg $
action can be lifted to a $G$ action on $L$, and the $G$-action on $X$ lifts on $E$.
Then this $G$ action preserves $\nabla^L$. One can also assume,
after an integration over $G$ if necessary, that $G$ preserves the Hermitian
metrics $h^L, h^E$, the Hermitian connection $\nabla ^E$,
 the almost complex structure $J$ and thus also the
Riemannian metric $g^{TX}$.
}

For $V$ a subspace of $\Omega ^{0,\bullet}(X,L^p\otimes E)$,
we denote by $V^\bot$ the orthogonal complement of $V$
in $(\Omega ^{0,\bullet}(X,L^p\otimes E),\langle \quad \rangle)$.

Let $\Omega ^{0,\bullet}(X,L^p\otimes E)^G$, $(\Ker D_p)^G$
be the $G$-invariant subspaces of
 $\Omega ^{0,\bullet}(X,L^p\otimes E)$, $\Ker D_p$.
Let $P_p^G$ be the orthogonal projection from
$\Omega^{0,\bullet}(X,L^p\otimes E)$
on $(\Ker D_p)^G$.

\begin{defn}\label{ta0}
The $G$-invariant Bergman kernel $P_p^G(x,x')$ $(x,x'\in X)$ of $D_p$ is
 the smooth kernel of $P_p^G$
with respect to the Riemannian volume form $dv_X(x')$.
\end{defn}
Let $\{S^p_i\}_{i=1}^{d_p}$
$(d_p = \dim (\Ker D_p)^G)$ be any orthonormal basis of
$(\Ker D_p)^G$ with respect to the norm $\|\, \, \|_{L^2}$, then
\begin{align} \label{a7}
&P_p^G(x,x') = \sum_{i=1}^{d_p} S^p_i (x) \otimes (S^p_i(x'))^*
\in (E_p)_x \otimes (E_p^*)_{x'}.
\end{align}
Especially, $P_p^G(x,x)\in \End(E_p)_x\simeq
\End(\Lambda (T^{*(0,1)}X)\otimes E)_x$.

\comment{
For any $K\in \kg$, denote by $L_K$ the infinitesimal action induced
by $K$ on the corresponding vector bundles.
Let $K^E$ be the vector field induced by the action of $K$ on $E$.
Let $\mu^E(K)\in \End (E)$ be the vertical part of $K^E$ with respect
to $G$-invariant connection $\nabla ^E$ on $E$, then the action $L_K$
of $K$ on smooth sections of $E$ is given by
\begin{align}\label{a8}
L_K= \nabla ^E_{K^X}- \mu^E(K).
\end{align}
}

We use the notation $\mu^F$ in \eqref{h6} now.

Recall that the Lie derivative $L_K$ on $TX$ is given by
\begin{align}\label{a9}
L_K V= \nabla ^{TX}_{K^X} V- \nabla ^{TX}_V {K^X}.
\end{align}
Thus
\begin{align}\label{0a8}
\mu^{TX}(K)= \nabla ^{TX}_\cdot {K^X} \in \End(TX),
\end{align}
and the action on $\Lambda (T^{*(0,1)}X)$ induced by $\mu^{TX}(K)$
is given by
\begin{align}\label{0a10}
\mu^{\rm Cliff}(K) =\frac{1}{4} \sum_{a=1}^{2n} c(e_a)c(\nabla
^{TX}_{e_a} K^X) +\frac{1}{2} \tr [ P^{T^{(1,0)}X}\nabla
^{TX}_\cdot K^X].
\end{align}
Thus the action $L_K$ of $K$ on smooth sections of
$\Lambda (T^{*(0,1)}X)$ is given by  (cf.  \cite[(1.24)]{TZ98})
\begin{align}\label{a10}
L_K = \nabla ^{\text{Cliff}}_{K^X} - \mu^{\rm Cliff}(K).
\end{align}

By  \eqref{a6} and \eqref{a10},
the action $L_K$ of $K$ on $\Omega^{0,\bullet}(X,L^p\otimes E)$
is $\nabla ^{E_p}_{K^X} -\mu^{E_p}(K)$ with
\begin{align}\label{0a9}
\mu^{E_p}(K) = 2\pi \sqrt{-1}p \mu(K) + \mu^E(K)+\mu^{\rm Cliff}(K).
\end{align}

\begin{defn} \label{ta1} The (formally) self-adjoint operator
$\mL_{p}$ acting on
$(\Omega ^{0,\bullet}(X,L^p\otimes E),\langle \, , \, \rangle)$
is defined by,
\begin{align} \label{a11}
\mL_{p} = D_p^2 - p \sum_{i=1}^{\dim G} L_{K_i}L_{K_i}.
\end{align}
\end{defn}

The following  result will play a crucial role
in the whole paper.

\begin{thm}  \label{ta2}
The projection $P_{p}^G$ is the orthogonal projection from $\Omega
^{0,\bullet}(X,L^p\otimes E)$ onto $\Ker( \mL_p)$. Moreover, there
exist $\nu$, $C_L>0$ such that for any $p\in \bN$,
\begin{align} \label{a12}
\begin{split}
&\Ker( \mL_p) = (\Ker D_p)^G,\\
&{\spec (\mL_p )}{\subset} \{0\}{\cup} [2p \nu
-C_L,+\infty[.
\end{split}
\end{align}
\end{thm}
\begin{proof} By (\ref{a11}),
 for any $s\in\Omega ^{0,\bullet}(X,L^p\otimes E)$,
\begin{align}\label{a13}
 \langle \mL_p s, s \rangle = \|D_ps\|^2 _{L^2}
+ p \sum_{i=1}^{\dim G} \|L_{K_i}s\|^2_{L^2}.
\end{align}
Thus $\mL_p s=0$ is equivalent to
\begin{align}\label{a14}
D_ps= L_{K_i}s=0.
\end{align}
This means $s$ is fixed by the $G$-action.
Thus we get the first equation of (\ref{a12}).

For $s\in  (\Ker \mL_p)^\bot$, there exist $s_1\in \Omega
^{0,\bullet}(X,L^p\otimes E)^G \cap (\Ker D_p)^\bot$, $s_2\in
(\Omega ^{0,\bullet}(X,L^p\otimes E)^G)^\bot$, such that
$s=s_1+s_2$.
Clearly, $$D_p s_1\in \Omega ^{0,\bullet}(X,L^p\otimes E)^G, \ \ \
D_ps_2\in (\Omega ^{0,\bullet}(X,L^p\otimes E)^G)^\bot.$$

By Theorem \ref{t2.1} and (\ref{a4}),
\begin{align}\label{a15}
\langle \mL_p s, s \rangle &= - p \langle \rho(\Cas) s_2, s_2 \rangle
 + \|D_ps_2\|^2_{L^2}+ \|D_ps_1\|^2_{L^2}\\
&\geq p\nu_1 \|s_2\|^2_{L^2} + (2p\nu_0-C_L)
\|s_1\|^2_{L^2},\nonumber
\end{align}
from which we get \eqref{a12}.
\end{proof}

We assume that $0\in \kg^*$ is a regular value of $\mu$.
Then $X_G=\mu^{-1}(0)/G$ is an orbifold
($X_G$ is smooth if $G$ acts freely on $P=\mu^{-1}(0)$).
Furthermore, $\omega$ descends to a symplectic form $\omega_G$ on
$X_G$. Thus one gets the Marsden-Weinstein symplectic reduction
$(X_G,\omega_G)$.

Moreover, $(L,\nabla ^L)$, $(E,\nabla ^E)$ descend to $(L_G,\nabla
^{L_G})$, $(E_G,\nabla ^{E_G})$ over $X_G$ so that the
corresponding curvature condition holds \cite{GuSt82} :
\begin{align}\label{a16}
\frac{\sqrt{-1}}{2\pi}
R^{L_G}= \omega_G.
 \end{align}

The $G$-invariant almost complex structure $J$ also descends to an
almost complex structure $J_G$ on $TX_G$, and $h^L, h^E, g^{TX}$
descend to $h^{L_G}$, $h^{E_G}, g^{TX_G}$.

We can construct the corresponding spin$^c$ Dirac operator
$D_{G,p}$ on $X_G$.

Let $P_{G,p}$ be the orthogonal projection from
$\Omega^{0,\bullet}(X_G ,L^p_G\otimes E_G)$ on $\Ker D_{G,p}$, and
let $P_{G,p}(x,x')$ be the smooth kernel of $P_{G,p}$ with respect
to the Riemannian volume form $dv_{X_G}(x')$.

The purpose of this paper is to study the asymptotic expansion of
$P_p^G(x,x')$ when $p\to \infty$, and we will relate it to the
asymptotic expansion of the Bergman kernel $P_{G,p}$ on $X_G$.

\subsection{Localization of the problem}\label{s3.1}

Let $a^X$ be the injectivity radius of $(X, g^{TX})$, and $\var\in
(0,a^X/4)$. If $x\in X$, $Z\in T_x X$,
let $\bR \ni u\to x_u=\exp^X_x(uZ)\in X$ be the geodesic in  $(X,g^{TX})$,
  such that $x_0=x, \frac{d x_u}{d u}|_{u=0}= Z$.

For $x\in X$, we denote by $B^{X}(x,\varepsilon)$ and
$B^{T_xX}(0,\varepsilon)$ the open balls in $X$ and $T_x X$ with
center $x$ and radius $\varepsilon$, respectively.
The map $T_x X\ni Z \to \exp^X_x(Z)\in X$ is a diffeomorphism from
$B^{T_xX}(0,\varepsilon)$  on $B^{X}(x,\varepsilon)$ for
$\varepsilon \leq a^X$.

From now on, we identify $B^{T_xX}(0,\varepsilon)$ with
$B^{X}(x,\varepsilon)$ for $\varepsilon \leq a^X/4$.

Let $f : \bR \to [0,1]$ be a smooth even function such that
 \begin{align} \label{c2}
f(v) = \left \{ \begin{array}{ll}  1 \quad {\rm for}
\quad |v| \leq  \var/2, \\
  0 \quad {\rm for} \quad |v| \geq  \var.
\end{array}\right.
\end{align}
Set
 \begin{align} \label{1c2}
F(a)= \Big(\int_{-\infty}^{+\infty}f(v) dv\Big)^{-1}
 \int_{-\infty}^{+\infty} e ^{i v a} f(v) dv.
\end{align}
Then $F(a)$ is an even function and lies in Schwartz space $\mathcal{S} (\bR)$
 and $F(0)=1$.

 Let $\wi{F}$ be the holomorphic function on
 $\bC$ such that $\wi{F}(a ^2) =F(a)$.
The restriction of $\wi{F}$
 to $\bR$ lies  in the Schwartz space $\mS (\bR)$.

Let $\wi{F}(\mL_p)(x,x')$ be the smooth kernel of $\wi{F}(\mL_p)$
with respect to the volume form $dv_X(x')$.

\begin{prop}\label{0t3.0}
 For any $l,m\in \bN$, there exists $C_{l,m}>0$ such that  for $p\geq C_L/\nu$,
\begin{align}\label{1c3}
|\wi{F}(\mL_p)(x,x') - P_p^G(x,x')|_{\cC^m(X\times X)} \leq C_{l,m} p^{-l}.
\end{align}
Here the $\cC^m$ norm is induced by $\nabla^L, \nabla^E$,
 $\nabla^{\text{Cliff}}$, $h^L, h^E$ and $g^{TX}$.
\end{prop}

\begin{proof} For $a\in \bR$, set
\begin{eqnarray}\label{0c2}
\phi_p(a) = 1_{[p\nu, +\infty[} (a) \wi{F}(a).
\end{eqnarray}
Then by Theorem \ref{ta2}, for $p> C_L/\nu$,
\begin{eqnarray}\label{0c3}
\wi{F}(\mL_p)-P_p^G = \phi_p(\mL_p).
\end{eqnarray}
By (\ref{1c2}), for any $m\in \bN$ there exists $C_{m}>0$ such that
\be\label{1c9}
\sup_{a\in \bR} |a|^m |\wi{F}(a) | \leq C_{m}.
\ee

As $X$ is compact, there exist $\{x_i\}_{i=1}^r\subset X$ such
that $\{U_i = B^X(x_i,\var)\}_{i=1}^r$ is a covering of $X$.
We identify $B^{T_{x_i}X}(0,\var)$ with $B^{X} (x_i,\var)$ by
geodesics  as above.

We identify $(E_p)_Z$ for $Z\in B^{T_{x_i}X}(0,\var)$ to
$(E_p)_{x_i}$  by parallel transport with respect to the
connection $\nabla^{E_p}$ along the curve $\gamma_Z: [0,1]\ni u
\to \exp^X_{x_i} (uZ)$.

Let $\{e_j\}_{j=1}^{2n}$ be an orthonormal basis of $T_{x_i}X$.
Let $\wi{e}_j (Z)$ be the parallel transport of ${e}_j$ with
respect to $\nabla^{TX}$ along the above curve.

Let $\Gamma ^E, \Gamma ^L, \Gamma ^{\text{Cliff}}$ be the
corresponding connection forms of $\nabla^E$, $\nabla^L $ and
$\nabla^{\text{Cliff}}$ with respect to any fixed frame for $E,L$,
$\Lambda (T^{*(0,1)}X)$ which is parallel along the curve
$\gamma_Z$ under the trivialization on $U_i$. Then $\Gamma ^L$ is
a usual $1$-form.

Denote by  $\nabla_U$  the ordinary differentiation
 operator on $T_{x_i}X$ in the direction $U$. Then
\begin{align}\label{c10}
&\nabla ^{E_p} = \nabla +p\, \Gamma ^L
+ \Gamma ^{\text{Cliff}}+ \Gamma ^E, \quad
D_p = c(\wi{e}_j) \nabla ^{E_p} _{\wi{e}_j}.
\end{align}

Let $\{ \varphi_i \}$ be a partition of unity subordinate to
$\{U_i\}$.

For $l\in \bN$, we define a Sobolev norm on the $l$-th Sobolev
space $H^l(X,E_p)$ by \be\label{c11} \| s\| _{H^l_p}^2 = \sum_i
\sum_{k=0}^l \sum_{i_1, \cdots ,i_k=1} ^{2n}
\|\nabla_{e_{i_1}}\cdots  \nabla_{e_{i_k}}(\varphi _i s)\|_{L^2}^2
\ee Then by (\ref{c10}), there exist $C,C', C''>0$ such that for
$p\geq 1$, $s\in H^2(X, E_p)$, \be\label{c12}
  C'\|D_p^2 s\|_{L^2} -C'' p^2\|s\|_{L^2} \leq
\|s\|_{H^2_p}   \leq C(\|D_p^2 s\|_{L^2} + p^2\|s\|_{L^2}).
\ee

Observe that $D_p$ commutes with the $G$-action, thus
\begin{align}\label{0c15}
[D_p,L_{K_j}]=0.
\end{align}

By (\ref{a11}), \eqref{0c15}, and the facts that $D_p$ is self-adjoint and
$L_{K_j}$ is skew-adjoint, we know
\begin{multline}\label{0c11}
\|\mL_p s\|_{L^2}^2 = \|D_p ^2 s\|_{L^2}^2
+p^2 \| \sum_jL_{K_j}L_{K_j}s\|_{L^2}^2
 -2p\, \Re \sum_j\langle D_p^2 s,L_{K_j}L_{K_j}s \rangle\\
= \|D_p^2 s\|_{L^2}^2 +p^2 \|\sum_j L_{K_j}L_{K_j}s\|_{L^2}^2
+2p \sum_j\| L_{K_j}D_ps\|_{L^2}^2.
\end{multline}

From (\ref{c12}), and (\ref{0c11}), there exists $C>0$ such that
\begin{align}\label{0c13}
\| s\| _{H^2_p} \leq C(\|\mL_p s\|_{L^2}+ p^2 \|s\|_{L^2}).
\end{align}

Let $Q$ be a differential operator of order $m\in \bN$ with scalar
principal symbol and with  compact support in $U_i$, then \be
\label{c14} \qquad [\mL_p,Q] =  [D^2_p,  Q ]-p \sum_j
[L_{K_j}L_{K_j}, Q] \ee is a differential operator  of order $m+1$.
Moreover, by \eqref{0a9}, \eqref{c10}, the leading term of order
$m-1$  differential operator in $[L_{K_j}L_{K_j}, Q]$ is $p^2
[((\Gamma^L -2\pi \sqrt{-1}\mu)(K_j))^2,Q]$. Thus by (\ref{0c13})
and (\ref{c14}),
\begin{align}\label{c13}
\|Qs\|_{H^2_p}   &\leq C(\|\mL_pQ s\|_{L^2} + p^2 \|Qs\|_{L^2})\\
 &\leq C(\|Q \mL_p s\|_{L^2} + p\|s\|_{H^{m+1}_p}
+p^2 \|s\|_{H^{m}_p}+ p^3  \|s\|_{H^{m-1}_p}).\nonumber
\end{align}
This means \be\label{c17} \|s\|_{H^{2m+2}_p} \leq C_m p^{2m+2}
\sum_{j=0}^{m+1} \|\mL_p ^js\|_{L^2}. \ee Moreover, from $$\langle
\mL_p^{m'} \phi_p(\mL_p)Q s,s'\rangle =\langle s,Q^ *
\phi_p(\mL_p) \mL_p^{m'} s'\rangle,$$ (\ref{0c2}) and (\ref{1c9}),
we know that for any $l,m'\in \bN$, there exists $C_{l,m'}>0$ such
that for $p\geq 1$, \be\label{c18} \|\mL_p^{m'}
\phi_p(\mL_p)Qs\|_{L^2}  \leq C_{l,m'}
 p^{-l+m}  \|s\|_{L^2}.
\ee

We deduce from (\ref{c17}) and (\ref{c18}) that if $Q_1,Q_2$ are
differential operators of order $m,m'$  with compact support in
$U_i$, $U_j$ respectively, then for any $l>0$, there exists
$C_l>0$ such that for $p\geq 1$,
\be\label{c19} \|Q_1 \phi_p( \mL_p) Q_2
s\|_{L^2}\leq C_l p^{-l} \|s\|_{L^2}. \ee

On $U_i\times U_j$, by
using Sobolev inequality and  (\ref{0c3}),
 we get  Proposition \ref{0t3.0}.
\end{proof}

Observe that $K_j^X$ are vector fields along the orbits of the
$G$-action, thus the contribution of $pL_{K_j}L_{K_j}$ in the wave
operator $\exp(\sqrt{-1}t \sqrt{\mL_p})$ will propagate along the
$G$-orbits, and the principal symbol of $\mL_p$ is given by
$$\sigma (\mL_p)(\xi) =|\xi|^2 +p \sum_j \langle  K_j^X, \xi
\rangle ^2\quad \mbox{for}\, \, \xi\in T^*X.$$

 By the finite propagation speed for solutions of hyperbolic equations
\cite[\S 7.8]{ChPi81}, \cite[\S 4.4]{Taylor81},
\cite[I. \S 2.6, \S 2.8]{Taylor96},
$\wi{F}(\mL_p)(x, x')$ only depends on the restriction of $\mL_p$ to
$G\cdot B^X(x,\var)$ and
\begin{align}\label{c20}
\wi{F}(\mL_p)(x, x')=0,     \quad {\rm if} \, d^X(Gx, x') \geq \var .
\end{align}
(When we apply the proof of \cite[\S 2.6, \S 2.8]{Taylor96}, we
need to suppose that $\Sigma_1,\Sigma_2$ therein are
$G$-space-like surfaces for the operator $\frac{\partial
^2}{\partial t^2}- D_p^2$).

Combining with Proposition \ref{0t3.0}, we know that the
asymptotic of $P_p^G(x,x')$ as $p\to \infty$ is localized on a
neighborhood of $Gx$.

\comment{
\begin{thm}\label{0t3.1} For  $U$ a $G$-open neighborhood
of $P=\mu^{-1}(0)$, then for any $l,m \in \bN$,
there exists $C_{l,m}>0$ (depend on $U$)
such that for   $x,x'\in X\setminus U$,
\begin{align}\label{c21}
|P^G_p(x, x')|_{\cC^m} \leq C_{l,m} p^{-l}.
\end{align}
\end{thm}
}

\begin{proof}[Proof of Theorem \ref{t0.0}]
From Proposition  \ref{0t3.0} and \eqref{c20}, we get \eqref{0.5}
for any $x,x'\in X$, $d^X(Gx,x')\geq \var_0$. Now we will establish
\eqref{0.5} for $x,x'\in X\setminus U$.

Recall that $U$ is a $G$-open neighborhood of $P=\mu^{-1}(0)$.

As $0$ is a regular value of $\mu$, there exists $\epsilon_0>0$
such that $\mu : X_{ 2\epsilon_0}= \mu ^{-1}(B^{\kg^*}(0,
2\epsilon_0)) \to B^{\kg^*}(0, 2 \epsilon_0)$ is a submersion,
$X_{2\epsilon_0}$ is a $G$-open subset of $X$.

Fix $\var, \epsilon_0 >0$ small enough such that
$X_{2 \epsilon_0}\subset U$, and  $d^X(x,y)>4\var$ for any
$x\in  X_{\epsilon_0}$, $y\in X\setminus U$.
Then  $V_{\epsilon_0}= X\setminus X_{\epsilon_0}$
is a smooth $G$-manifold with boundary $\partial V_{\epsilon_0}$.

Consider the operator $\mL_p$ on $V_{\epsilon_0}$ with the
Dirichlet boundary condition. We denote it by $\mL_{p,D}$.
Note that $\mL_{p,D}$ is self-adjoint.

Still from \cite[\S 2.6, \S 2.8]{Taylor96}, the wave operator
 $\exp(\sqrt{-1}t \sqrt{\mL_{p,D}})$ is well defined and
$\exp(\sqrt{-1}t \sqrt{\mL_{p,D}})(x,x')$ only depends on the restriction
of $\mL_p$  to $G\cdot B^X (x,t)\cap V_{\epsilon_0}$, and is zero if
$d^X(Gx,x') \geq t$.
Thus, by (\ref{1c2}),
\begin{align}\label{1c22}
\wi{F}(\mL_p) (x,x')=\wi{F}(\mL_{p,D}) (x,x'),
\quad {\rm if} \, x,x'\in X\setminus U.
\end{align}

Now for $s\in \cC^\infty_0 (V_{\epsilon_0}, E_p)$,
   after taking an
integration over $G$, we can get the decomposition $s=s_1+s_2$
with $s_1\in  \Omega ^{0,\bullet}(X,L^p\otimes E)^G$, $s_2\in
(\Omega ^{0,\bullet}(X,L^p\otimes E)^G)^\bot$ and $\supp s_i
\subset V_{\epsilon_0}\setminus \partial V_{\epsilon_0}$.

Since $\sum_{i=1}^{\dim G} L_{K_i}L_{K_i}$ commutes with
 the $G$-action,
 $\mL_p s_1\in  \Omega ^{0,\bullet}(X,L^p\otimes E)^G$,
$\mL_p s_2\in  (\Omega ^{0,\bullet}(X,L^p\otimes E)^G)^\bot$ and,
by (\ref{a11}), (\ref{a15}),
\begin{multline}\label{c23}
\langle \mL_p s, s \rangle = \langle \mL_p s_1, s_1 \rangle
+\langle \mL_p s_2, s_2 \rangle \\
=\|D_ps_2\|^2_{L^2}
-p \langle \rho(\Cas)s_2,s_2\rangle +  \langle D_p^2 s_1,s_1\rangle\\
\geq p\nu_1 \|s_2\|^2_{L^2} + \langle D_p^2 s_1,s_1\rangle.
\end{multline}

To estimate the term $\langle D_p^2 s_1,s_1\rangle$, we will use
the Lichnerowicz formula.

Recall that the Bochner-Laplacian  $\Delta ^{E_p}$ on $E_p$
is defined by \eqref{r6}.

Let $r^X$ be the Riemannian scalar curvature of $(TX,g^{TX})$.

Let $\{w_a\}$ be an orthonormal frame of $(T^{(1,0)}X, g^{TX})$.
Set
\begin{align}  \label{c24}
\begin{split}
&\om_d=-\sum_{a,b} R^L (w_a,\overline{w}_b)\,\overline{w}^b\wedge
\,i_{\overline{w}_a}\,,\\
& \tau(x)=\sum_a R^L (w_a,\overline{w}_a)\, ,
\qquad  R^E_\tau = \sum_a R^E(w_a,\overline{w}_a)\, ,    \\
&\mathbf{c}(R)=\sum_{a<b}\left(R^E+\tfrac{1}{2}
\tr[R^{T^{(1,0)}X}]\right)(e_a,e_b)\,c(e_a)\,c(e_b)\,.
\end{split}
\end{align}


The Lichnerowicz formula \cite[Theorem 3.52]{BeGeVe} (cf.
\cite[Theorem 2.2]{MM02}) for $D_p^2$ is
\begin{equation}\label{Lich}
D^2_p=\Delta ^{E_p}-2p\om_d-
p\tau+\tfrac{1}{4}r^X+ \mathbf{c}(R).
\end{equation}
Especially, as $\supp s_i \subset V_{\epsilon_0}
\setminus \partial V_{\epsilon_0}$, from \eqref{Lich}, we get
\begin{align}\label{c25}
 \langle D_p^2 s_1,s_1\rangle= \|\nabla ^{E_p} s_1\|^2_{L^2}
-p  \langle (2 \om_d+ \tau)s_1,s_1\rangle
+   \langle (\tfrac{1}{4}r^X+ \mathbf{c}(R))s_1,s_1\rangle .
\end{align}

Since $s_1\in \Omega ^{0,\bullet}(X,L^p\otimes E)^G$,
   from (\ref{h6}), for any $K\in\kg$,
\begin{align}\label{c26}
 \nabla ^{E_p}_{K^X} s_1 =(L_{K}+\mu^{E_p}(K)) s_1
=\mu^{E_p}(K)s_1.
\end{align}

From \eqref{0a9} and \eqref{c26},  there exist $C,C'>0$ such that
\begin{align}\label{0c20}
\|\nabla ^{E_p} s_1\|^2_{L^2}
&\geq C\sum_j \|\nabla ^{E_p} _{K^X_j} s_1\|^2_{L^2}
= C\sum_j\| \mu^{E_p}(K_j)s_1\|^2_{L^2}\\
&\geq Cp^2 \| |\mu|s_1\|^2_{L^2} - C'\|s_1\|^2_{L^2}
\geq C\epsilon_0^2 p^2  \|s_1\|^2_{L^2}- C'\|s_1\|^2_{L^2}.\nonumber
\end{align}

From \eqref{c23}-\eqref{0c20}, for $p$ large enough,
\begin{align}\label{0c21}
\langle \mL_p s, s \rangle \geq p\nu_1  \|s_2\|^2_{L^2}+ Cp^2  \|s_1\|^2_{L^2}.
\end{align}
Thus there are $C,C'>0$ such that for $p\geq 1$,
\begin{align}\label{0c22}
{\spec} ({\mL_{p,D}} )\subset [Cp-C',\infty[.
\end{align}

At first as $K^X_j\in T\partial V_{\epsilon_0}$ for any $j$,
thus $L_{K_j}$ preserves the Dirichlet boundary condition.
We get for $l\in \bN$,
\begin{align}\label{2c20}
L_{K_j} \phi_p(\mL_{p,D})= \phi_p(\mL_{p,D})L_{K_j},\quad
(\mL_{p,D})^l  \phi_p(\mL_{p,D})= \phi_p(\mL_{p,D}) (\mL_{p,D})^l.
\end{align}
Thus from \eqref{a11} and \eqref{2c20},
\begin{align}\label{2c23}
D^2_{p,D} \leq \mL_{p,D},
\end{align}
and for $l\in \bN$, $(D^2_{p,D})^l$ commutes with
the operator $\phi_p(\mL_{p,D})$.

Let $\phi_p(\mL_{p,D})(x,x')$ be the smooth kernel of $\phi_p(\mL_{p,D})$
with respect to $dv_X(x')$.

Then from the above argument
 we get that for any $l,k\in \bN$, $(D^2_{p,x})^l (D^2_{p,x'})^k
\phi_p(\mL_{p,D})(x,x')$
verifies the Dirichlet boundary condition for $x,x'$ respectively.

By \eqref{c10} and the elliptic estimate for  Laplacian with
Dirichlet boundary condition \cite[Theorem 5.1.3]{Taylor96}, there exists
$C>0$ such that for $s\in H^{2m+2}(X,E_p)\cap H^1_0(X,E_p)$, $p\in \bN$,
we have
\begin{align}\label{2c21}
\|s\|_{H^{2m+2}}\leq C(\|D^2_p s\|_{H^{2m}} +p^2\|s\|_{H^{2m+1}}).
\end{align}

 Thus  if $Q_1,Q_2$ are
differential operators of order $2m,2m'$  with compact support in
$U_i$, $U_j$ respectively,
by \eqref{2c21} and \eqref{2c23}, as in \eqref{c17}, we get
 for $s\in \cC^\infty_0(V_{\epsilon_0},E_p)$,
\begin{multline}\label{2c22}
\|Q_1 \phi_p(\mL_{p,D}) Q_2 s\|_{L^2}
\leq C p^{2m+2m'}\sum_{j_1=0}^m\sum_{j_2=0}^{m'}
 \|(D^2_{p,D})^{j_1}\phi_p(\mL_{p,D}) (D^2_{p,D})^{j_2}s\|_{L^2}\\
\leq C p^{2m+2m'}\sum_{j_1=0}^m\sum_{j_2=0}^{m'}
 \|(\mL_{p,D})^{j_1}\phi_p(\mL_{p,D}) (\mL_{p,D})^{j_2}s\|_{L^2}.
\end{multline}

From  (\ref{0c22}), \eqref{2c22}, as in (\ref{c19}), we get
 \begin{align}\label{0c23}
\|Q_1 \phi_p( \mL_{p,D}) Q_2
s\|_{L^2}\leq C_l p^{-l} \|s\|_{L^2}.
 \end{align}

By using Sobolev inequality as in the proof of Proposition
\ref{0t3.0}, from (\ref{1c3}), (\ref{1c22}) and (\ref{0c23}),  we
get Theorem \ref{t0.0}.
\end{proof}

\subsection{Induced operator on $U/G$}\label{0s3.2}

Let $U$ be a $G$-neighborhood of $P=\mu^{-1}(0)$ in $X$ such that
$G$ acts freely on $\ov{U}$, the closure of $U$. We will use the
notation as in Introduction and Sections \ref{s4.2}, \ref{as4.2}
with $X$ therein replaced by $U$, especially $B=U/G$.

Let $\pi : U\to B$ be the natural projection with fiber $Y$.
Let $TY$ be the sub-bundle of $TU$ generated by the $G$-action,
let $g^{TY}$, $g^{TP}$ be the metrics on $TY$, $TP$ induced by
$g^{TX}$.

Let $T^HU$, $T^HP$ be the orthogonal complements of $TY$ in $TU$,
$(TP, g^{TP})$.
Let $g^{T^HU}$ be the metric on $T^HU$ induced by $g^{TX}$,
and it induces naturally a Riemannian metric $g^{TB}$ on $B$.

Let $dv_{B}$ be the Riemannian volume form on $(B, g^{TB})$.


Recall that in \eqref{h12}, we defined the isometry $\Phi = h
\pi_G: (\cC^{\infty} (U, E_p)^G,\langle \, ,\,  \rangle)$ $ \to
(\cC^{\infty}  (B, E_{p,B}), \langle \, ,\,  \rangle)$.

By \eqref{ah7}, $\mu^{E_p}$ defines a $G$-invariant section
$\wi{\mu}^{E_p}$ of $TY\otimes \End(E_p)$ on $U$.

Remark that $\om_d$, $\tau$, $\mathbf{c}(R)$ in \eqref{c24} are
$G$-invariant.
We still denote by $\om_d$, $\tau$, $\mathbf{c}(R)$
 the induced sections on $B$.

As a direct corollary of Theorem \ref{t4.1} and \eqref{Lich},
we get the following result,
\begin{prop} \label{0t3.2} As an operator on $\cC^{\infty}  (B, E_{p,B})$,
 \begin{multline}\label{b4}
\Phi \mL_p \Phi^{-1} =\Phi D_p^2 \Phi^{-1}\\
=  \Delta^{E_{p,B}}
-  \langle\wi{\mu}^{E_p},\wi{\mu}^{E_p}\rangle_{g^{TY}}
- \frac{1}{h}\Delta_B h -2p\om_d- p\tau+\tfrac{1}{4}r^X+ \mathbf{c}(R).
\end{multline}
\end{prop}

From Theorem \ref{t0.0}, Prop. \ref{0t3.0} and \eqref{c20},
modulo $\cO(p^{-\infty})$, $P^G_p(x,x')$ depends only the restriction
of $\mL_p$ on $U$.

To get a complete picture on $P^G_p(x,x')$, we explain  now that
modulo $\cO(p^{-\infty})$, $P^G_p(x,x')$ depends only on the
restriction of $\Phi \mL_p \Phi^{-1}$ on any neighborhood of $X_G$
in $B$.

As in the proof of Theorem  \ref{t0.0}, we will fix
$\epsilon_0 >0$ small enough such that
$X_{2\epsilon_0}=\mu ^{-1}(B^{\kg^*}(0, 2\epsilon_0))$$\subset U$,
 and the constant $\var>0$ which will be fixed later, verifying
that  $d^X(x,y)>4\var$ for any $x\in  X_{\epsilon_0}$, $y\in
X\setminus U$. Set $B_{\epsilon_0}=\pi(X_{\epsilon_0})$.

First we will extend all objects from a neighborhood of $P$ to
 the total space of the normal bundle $N$ of $P$ in $X$.

Let $\pi_N : N\to P$ be the normal bundle of $P$ in $X$. We identify
$N$ to the orthogonal complement of $TP$ in $(TX, g^{TX})$. Then
$G$ acts on $N$ and the action extends naturally on
$\pi_N ^*(L|_P)$, $\pi_N ^*( E|_P)$.

By (\ref{0a0}), we have an orthogonal decomposition of $TX$,
\begin{align}\label{b11}
TX|_P = T^H P \oplus TY|_{P} \oplus N,\quad {\rm and}\, \,
TY|_{P} \simeq P\times \kg, \quad N= \bJ TY|_P\simeq P\times \kg.
\end{align}
Denote by  $P^{TY}$, $P^{TP}$, $P^{N_P}$ the orthogonal
projections from $TX$ on $TY$, $TP$ and $N|_P$ by this
identification.

From  (\ref{b11}), we have
\begin{align}\label{b12}
TN\simeq \pi_N ^* TX|_P\simeq \pi_N ^*(TP\oplus \kg).
\end{align}

For $\var>0$, we denote by
$B^N_\var= \{(y,Z)\in N, y\in P, |Z|_{g^{TX}}\leq \var\}$.

Then for $\var_0$ small enough,  the map $(y,Z)\in
N  \to \exp^X_y(Z) \in X$ is  a diffeomorphism from $B^N_{2\var_0}$
onto a tubular neighborhood $\mathcal{U}_{2 \var_0}$ of $P$ in $X$.

 From now on, we use the notation $(y,Z)$ instead of $\exp^X_y(Z)$.
  We identify $y\in P$ with $(y,0)\in N$.
From \eqref{b11}, \eqref{b12}, we may and we will
identify $TN$ to $\pi^*_N TP\oplus \kg$.

For $Z\in N_y$, $|Z|\leq 2\var_0$, we identify $L_Z, E_Z$ to
$L_y,E_y$ by using parallel transport with respect to $\nabla ^L$,
$\nabla ^E$ along the curve $[0,1]\ni u\to uZ$. In this way, we
identify the Hermitian bundles $(\pi_N ^* L|_P,\pi_N ^* h^L)$,
$(\pi_N ^* E|_P,\pi_N ^*h^ E) $ to $(L, h^L),(E,h^E)$
 on  $B^N_{2\var_0}$.

 Let $\var>0$ with $\var<\var_0/2$.
Let $\varphi: \bR\to [0,1]$ be a smooth even function such that
\begin{align}\label{b13}
\varphi (v)=1  \  \  {\rm if} \  \  |v|<2;
\quad \varphi (v)=0 \   \   {\rm if} \  |v|>4.
\end{align}
 Let $\psi_\var : N \to N$
be the map defined by $\psi_\var (Z)= \varphi(|Z|/\var) Z\in N_y$
for $Z\in N_y$.

Let $g^{TN}_Z= g^{TX}_{\psi_\var(Z)}$, $J^N_Z=  J_{\psi_\var(Z)}$
be the induced metric and almost-complex structure on $N$.

 Let $\nabla ^{\pi_N^* E}= \psi_\var ^* \nabla ^{E}$,
then  $\nabla ^{\pi_N^* E}$ is the extension of $ \nabla ^{E}$ on
$B^N_{2\var}$.

Let $\nabla ^{\pi^*_N L}$ be the Hermitian connection on $(\pi
_N^* L,\pi_N ^*  h^{L})$ defined by that for $Z\in N_y$,
\begin{align}\label{b14}
&\nabla ^{\pi_N ^* L} = \psi_\var ^* \nabla ^{L}
+(1-\varphi (\tfrac{|Z|}{\var}) )  R^{L}_{y} (Z,P^{TP} \cdot)
+\frac{1}{2}(1-\varphi ^2 (\tfrac{|Z|}{\var} ) )
R^{L}_{y} (Z, P^{N_P} \cdot).
\end{align}
Then by using the identification \eqref{b11}, and \eqref{b12},
 we calculate directly that its curvature
 $R^{\pi_N ^* L}= (\nabla ^{\pi_N ^* L})^2$ is
\begin{multline}\label{b16}
R^{\pi_N ^* L}_Z = \psi_\var ^* R^L
+ d \Big((1-\varphi  (\tfrac{ |Z|}{\var} ) )
R^{L}_{y} (Z,P^{TP} \cdot)
+\frac{1}{2}(1-\varphi ^2 (\tfrac{|Z|}{\var} ))
R^{L}_{y} (Z, P^{N_P} \cdot) \Big)  \\
= R^L_{\psi_\var(Z)} (P^{TP}\cdot, P^{TP}\cdot  )
+ R^{L}_{y} (P^{N_P}\cdot, \cdot)
+\varphi ^2(\tfrac{|Z|}{\var})
 (R^L_{\psi_\var(Z)} -  R^{L}_{y})(P^{N_P}\cdot, P^{N_P}\cdot)\\
+\varphi (\tfrac{|Z|}{\var})(R^L_{\psi_\var(Z)}-  R^{L}_{y})
 (P^{N_P}\cdot, P^{TP}\cdot)\\
- \varphi'(\tfrac{|Z|}{\var}) \frac{Z^*}{\var |Z|}\wedge
[R^{L}_{y} (Z,P^{TP}\cdot)
- R^{L}_{\psi_\var(Z)} (Z,P^{TP}\cdot)]\\
- (\varphi \varphi')(\tfrac{|Z|}{\var})\frac{Z^*}{\var |Z|}\wedge [
R^{L}_{y} (Z,P^{N_P}\cdot)- R^{L}_{\psi_\var(Z)} (Z,P^{N_P}\cdot)]\\
+ d_y \Big((1-\varphi  (\tfrac{ |Z|}{\var} ) )R^{L}_{y} (Z,P^{TP} \cdot)
+\frac{1}{2}(1-\varphi ^2 (\tfrac{|Z|}{\var} ))
R^{L}_{y} (Z, P^{N_P} \cdot) \Big).
\end{multline}
Here $Z^*\in N^*$ is the dual of $Z\in N$ with respect to the metric $g^N$.

From (\ref{b16}), one deduces that $R^{\pi_N ^* L}$ is positive in
the sense of (\ref{0a1}) when $\var$ is small enough,
 with the corresponding constant $\nu_0$ for  $R^{\pi_N ^* L}$
 being larger
than $\frac{4}{5}\nu_0$.

Note that $G$ acts naturally on  the normal bundle $N$, and under
our identification, the $G$-actions on $L,E$ on $B^N_\var$ are
exactly the $G$-actions on $L|_{P},E|_{P}$ on $P$.

Now we define the $G$-actions on  $\pi_N ^*L, \pi_N ^*E$ by their
$G$-actions on $P$, then they extend the $G$-actions on $L,E$ on
$B^N_\var$ to $N$.

By (\ref{a6}), the moment map $\mu_N: N\to \kg^*$ of the
$G$-action on $N$ is defined by
\begin{align}\label{b17}
-2 \pi \sqrt{-1} \mu_N(K) = L_K  -\nabla ^{\pi_N ^*L}_{K^N} ,  \,\,
K\in \kg.
\end{align}

Observe that $\psi_{\var*} K^N_{(y,Z)} = K^{P}_y\in TP$,
thus from  (\ref{a6}), (\ref{b11}), and (\ref{b14}),
\begin{multline}\label{b18}
2 \pi \sqrt{-1} \mu_N(K)_{(y,Z)}= (1-\varphi (|Z|/\var) )  R^{L}_{y} (Z,K^{P})
+ 2 \pi \sqrt{-1}\mu (K)_{\psi_\var(Z)}\\
= R^{L}_{y} (Z,K^{P}) + \cO( \varphi ^2(|Z|/\var)|Z|^2).
\end{multline}
Thus $\mu_N^{-1}(0)= \mu ^{-1}(0)=P$ for $\var$ small enough, and
for $|Z|\geq 4 \var$,
\begin{align}\label{b19}
2 \pi \sqrt{-1} \mu_N(K)_{(y,Z)}= R^{L}_{y} (Z,K^{P}).
\end{align}

From now on, we fix  $\var$ as above.

Let $\wi{F}(\Phi\mL_p\Phi ^{-1})(x,x')$ $(x,x'\in B_{\epsilon_0})$
be the smooth kernel of
$\wi{F}(\Phi \mL_p\Phi ^{-1})$ with respect to $dv_{B}(x')$.
We will also view  $\wi{F}(\Phi\mL_p\Phi ^{-1})$ as a $G\times G$-invariant
section of  $ {\rm pr}_1^* E_p \otimes {\rm pr}_2^* E_p^* $
 on $X_{\epsilon_0}\times X_{\epsilon_0}$.

\begin{thm}\label{0t3.3} For any $l,m\in \bN$, there exists $C_{l,m}>0$
such that for $p\geq 1$, $x,x'\in X_{\epsilon_0}$,
\begin{align}\label{b10}
| h(x)h(x') P^G_p(x,x')
- \wi{F}(\Phi \mL_p \Phi^{-1}) (\pi(x),\pi(x'))|
_{\cC^m(X_{\epsilon_0}\times X_{\epsilon_0})}
\leq C_{l,m} p^{-l}.
\end{align}
\end{thm}
\begin{proof}
Let $D_p^{N}$ be the Dirac operator on $N$ associated to the
above data  by the construction in Section \ref{s3.01}. By the
argument in \cite[p. 656-657]{MM02} and the proof of Theorem \ref{ta2},
 we know that Theorems \ref{t2.1}, \ref{ta2}
still hold for $D_p^{N}$.

Let $\mL_p^N$ be the operator on $N$ defined as in (\ref{a11}).
Then there exists $C>0$ such that for $p\geq 1$,
\begin{align}\label{b20}
&{\spec} \left(\mL_p^N\right)   \subset \{0\}\cup [p\nu-C,+\infty[.
\end{align}

Let $P^{N,G}_p$ be the orthogonal projection from
$\Omega ^{0,\bullet} (N, \pi_N ^* (L^p\otimes E))$ on $(\Ker D^N_p)^G$,
then by (\ref{b20}) and the arguments as in the proof of
Theorem \ref{0t3.0},
for any $l,m\in \bN$, $V\subset N$ a compact subset of $N$,
 there exists $C_{l,m}>0$ such that for $p\geq 1$, $x,x'\in  V$,
\begin{align}\label{b21}
|\wi{F}(\mL_p^N)(x,x') - P_p^{N,G}(x,x')|_{\cC^m(V\times V)}
\leq C_{l,m} p^{-l}.
\end{align}

 Let $P^{N/G}_p$ be the projection from
$(L^2  (N/G, (\Lambda (T^{*(0,1)}N) \otimes \pi_N ^* (L^p\otimes E))_{N/G}),
\langle \, ,\,  \rangle) $ onto $\Ker (\Phi \mL_p^N \Phi^{-1})$,
and let $P^{N/G}_p(z,z')$ be the smooth kernel of the operator
$P^{N/G}_p$ with respect to $dv_{N/G}(z')$.

We still denote by ${\rm pr}_1, {\rm pr}_2$ the projections from
$N\times N$ onto the first and second factor $N$.
We will also view $P^{N/G}_p(z,z')$ as a $G \times G$-invariant section of
$${\rm pr}_1^*(\Lambda (T^{*(0,1)}N) \otimes \pi_N^* (L^p\otimes E))
\otimes {\rm pr}_2^*(\Lambda (T^{*(0,1)}N) \otimes \pi_N^* (L^p\otimes E))^*$$
 on $N\times  N$.

 As $\Phi$ in (\ref{h12}) defines an isometry from
$(\Ker D^N_p)^G=\Ker \mL_p^N$
 onto  $\Ker (\Phi \mL_p^N \Phi^{-1})$, one has
\begin{align}\label{b23}
h(x)h(x')P^{N,G}_p(x,x') =  P^{N/G}_p(\pi(x),\pi(x')).
\end{align}

On $N/G$, by the arguments as in the proof of Theorem \ref{0t3.0}, we get
\begin{align}\label{b24}
|\wi{F}(\Phi\mL_p\Phi^{-1})(z,z') - P_p^{N/G}(z,z')|_{\cC^m(V/G\times V/G)}
\leq C_{l,m} p^{-l}.
\end{align}

By the finite propagation speed (\ref{c20}),
 we know that for $x,x'\in X_{\epsilon_0}$,
\begin{align}\label{b24a}
\wi{F}(\mL_p^N)(x,x') = \wi{F}(\mL_p)(x,x').
\end{align}

Now  we get  (\ref{b10}) from (\ref{1c3}), (\ref{b21})-(\ref{b24a}).
\end{proof}

Let $d^B(\cdot,\cdot)$ be the Riemannian distance on $B$.

 By \eqref{b4} and the finite propagation speed for solutions
of hyperbolic equations \cite[\S 7.8]{ChPi81}, \cite[\S 4.4]{Taylor81},
$\wi{F}(\Phi\mL_p\Phi^{-1})(x,x')$ only depends on the restriction of
$\Phi\mL_p\Phi^{-1}$ to $B^{B}(x,\var)$ and
\begin{align}\label{b25}
\wi{F}(\Phi\mL_p\Phi^{-1})(x, x')=0, \quad {\rm if} \, \, d^B (x, x') \geq \var .
\end{align}
Thus we have localized our problem near $X_G$.

Theorem \ref{0t3.3} helps us to understand that the asymptotic of
 $P^G_p(x,x')$ is local near $X_G$.
In the rest, we will not use directly Theorem \ref{0t3.3},
but the argument of its proof will be used in Section \ref{s3.2}.

\subsection{Rescaling and a Taylor expansion of
the operator $\Phi \mL_p \Phi ^{-1}$}\label{s3.2}

Recall that $N_G$ is the normal bundle of $X_G$ in $B$,
and we identify $N_G$ as the orthogonal complement of $TX_G$
in $(TB, g ^{TB})$.

Let $P^{TX_G}, P^{N_G}$ be the orthogonal projection
from $TB$ on $TX_G$, $N_G$ on $X_G$.

Recall that $\nabla ^{N_G}, {^0\nabla}^{TB}$ are connections on
$N_G$, $TB$ on $X_G$, and $A$ is the associated second fundamental form
defined in \eqref{a0.6}.

We fix  $x_0\in X_G$.

If $W\in T_{x_0}X_G$, let $\bR \ni t \to
x_t=\exp^{X_G}_{x_0}(tW)\in X_G$ be the geodesic in $X_G$ such
that $x_t|_{t=0}=x_0$, $\frac{dx}{dt}|_{t=0}=W$.

If $W\in T_{x_0}X_G$, $|W|\leq \var$, $V\in N_{x_0}$, let $\tau_W
V\in N_{G, \exp^{X_G}_{x_0}(W)}$ be the natural parallel transport
 of $V$ with respect to the connection
$\nabla^{N_G}$ along
the curve $[0,1] \ni t \to \exp_{x_0}^{X_G}(tW)$.

If $Z\in T_{x_0}B$, $Z=Z^0+Z^\bot$, $Z^0\in T_{x_0}X_G$,
$Z^\bot\in N_{x_0}$, $|Z^0|,|Z^\bot|\leq \var$, we identify $Z$
with $\exp^B_{\exp_{x_0}^{X_G} (Z^0)}(\tau_{Z^0} Z^\bot)$. This
identification is a diffeomorphism from
$B^{TX_G}_{x_0}(0,\var)\times B^{N_G}_{x_0}(0,\var)$ into an open
neighborhood $\cU(x_0)$ of $x_0$ in $B$. We denote it by
$\Psi$, and $\cU(x_0)\cap X_G= B^{TX_G}_{x_0}(0,\var) \times
\{0\}$.

From now on, we use indifferently the notation
$B^{TX_G}_{x_0}(0,\var)\times B^{N_G}_{x_0}(0,\var)$ or
$\cU(x_0)$, $x_0$ or $0$, $\cdots$.

 We identify $(L_B)_Z, (E_B)_Z$ and $(E_{p,B})_Z$ to $(L_B)_{x_0}, (E_B)_{x_0}$
and $(E_{p,B})_{x_0}$ by using parallel transport with respect
to $\nabla ^{L_B},\nabla ^{E_B}$ and $\nabla ^{E_{p,B}}$
 along the curve $\gamma _u: [0,1]\ni u\to uZ$.

Recall that $T^H U\subset TX$ is the horizontal bundle for $\pi:
U\to B$ defined in Section \ref{0s3.2}.

Let $P^{T^H U}$ be the orthogonal projection from $TX$ onto $T^H
U$.

For $W\in TB$, let $W^H\in T^H U$ be the lift of $W$.

For $y_0\in \pi^{-1}(x_0)$, we define the curve $\wi{\gamma} _u
:[0,1] \to X$ to be the lift of the curve $\gamma _u$ with
$\wi{\gamma} _0=y_0$
 and $\frac{\partial \wi{\gamma}_u}{\partial u} \in T^H U$.
Then on $\pi^{-1}(B^{TB}(0,\var))$, we use the  parallel transport
with respect to $\nabla ^{L},\nabla ^{E}$ and $\nabla ^{E_{p}}$
 along the curve $\wi{\gamma} _u$ to trivialized the corresponding bundles.
By \eqref{h5}, the previous trivialization is naturally induced by
  this one.

Let $\{e^0_i\}$, $\{e^\bot_j\}$ be orthonormal basis of
$T_{x_0}X_G$, $N_{G,x_0}$, then $\{e_i\}=\{e^0_i,e^\bot_j \}$ is
an orthonormal basis of $T_{x_0}B$. Let $\{e^i\}$ be its dual basis.
 We will also denote $\Psi_*(e^0_i),\Psi_*(e^\bot_j)$
by $e^0_i,e^\bot_j$. Thus in our coordinate,
\begin{align}\label{2c14}\tfrac{\partial}{\partial Z^0_i}=e^0_i,
\ \ \ \ \ \tfrac{\partial}{\partial Z^\bot_j}=e^\bot_j.
\end{align}

For $\var >0$ small enough,  we will extend the geometric objects
on $B^{TB}(x_{0},\var)$ to $\bR^{2n-n_0} \simeq T_{x_0}B$ (here we
identify $(Z_1,\cdots, Z_{2n-n_0}) \in \bR^{2n-n_0}$ to $\sum_i
Z_i e_i\in T_{x_0}B$) such that $D_p$ will become the restriction
of a spin$^c$ Dirac operator on $G\times \bR^{2n-n_0}$ associated
to a Hermitian line bundle with positive curvature. In this way,
we can replace $X$ by  $G\times \bR^{2n-n_0}$.

First of all, we denote by $L_0$, $E_0$ the trivial bundles $L|_{Gy_0},
E|_{Gy_0}$ on $X_0= G\times\bR^{2n-n_0}$, and we still denote by
 $\nabla ^L,\nabla ^E$, $h^L$ etc. the connections and metrics on  $L_0$,
$E_0$ on $\pi ^{-1}(B^{T_{x_0}B}(0,4\var))$ induced by the above
identification. Then $h^L$, $h^E$ is identified with the constant
metrics $h^{L_0}=h^{L_{y_0}}$,  $h^{E_0}=h^{E_{y_0}}$.

Set
\begin{align}\label{2c15}
\mR^\bot =\sum_j Z_j^\bot e_j^\bot= Z^\bot,\quad
 \mR^0=\sum_i Z_i^0 e_i^0= Z^0,\quad  \mR =\mR^\bot+ \mR^0=Z.
\end{align}
Then $\mR$ is the radial vector field on $\bR^{2n-n_0}$.

 Let $\varphi_\var : X_0 \to X_0$ be the map defined by
$\varphi_\var(g,Z)= (g,\varphi(|Z|/\var) Z)$
for $(g,Z)\in G\times \bR^{2n-n_0}$.

Let $g^{TX_0}(g,Z)= g^{TX}(\varphi_\var(g,Z))$,
$J_0(g,Z)=  J(\varphi_\var(g,Z))$
be the metric and almost-complex structure on $X_0$.

 Let $\nabla ^{E_0}= \varphi_\var ^* \nabla ^{E}$, then  $\nabla ^{E_0}$
is the extension of $ \nabla ^{E}$ on $\pi ^{-1}(B^{T_{x_0}B}(0,\var))$.

Let $\nabla ^{L_0}$ be the Hermitian connection on $(L_0, h^{L_0})$
on $G\times\bR^{2n-n_0}$ defined by for $Z\in \bR^{2n-n_0}$,
\begin{align}\label{1c15}
&\nabla ^{L_0}= \varphi_\var ^* \nabla ^{L}
+\Big(1-\varphi (\tfrac{|Z|}{\var}) \Big)
R^{L}_{y_0} (\mR^H,P^{TY}_{y_0} \cdot)
+\frac{1}{2}\Big(1-\varphi ^2(\tfrac{|Z|}{\var}) \Big)
 R^{L}_{y_0} (\mR^H, P^{T^H U}_{y_0} \cdot).
\end{align}
As in (\ref{b16}), its curvature  $R^{L_0}$ is positive in the
sense of (\ref{0a1}) for $\var$ small enough, and  the
corresponding constant $\nu_0$ for $R^{L_0}$ is bigger than
$\frac{4}{5}\nu_0$ uniformly for $y_0\in P$.

 From now on, we fix $\var$ as above.

Now $G$ acts naturally on $X_0$,
and under our identification, the $G$-action on $L,E$
on $G\times B^{T_{x_{0}}B}(0,\var)$ is exactly the $G$-action on
$L|_{Gy_0},E|_{Gy_0}$.

We define a $G$-action on  $L_0, E_0$ by its $G$-action on $Gy_0$,
 then it extends the $G$-action on
$L,E$ on $G\times  B^{T_{x_{0}}B}(0,\var)$ to $X_0$.

By  (\ref{a5}), for any $K\in \kg$, $W\in TP$ on $P=\mu^{-1}(0)$, we have
\begin{align}\label{2c16}
\begin{split}
& R^{L}(W, K^X)= - 2 \pi \sqrt{-1} \omega (W, K^X)
= 2 \pi \sqrt{-1} W( \mu(K) )=0, \\
&R ^L_{(1,Z^0)} (\mR^H, K^X) = R ^L_{(1,Z^0)} ( (\mR^\bot)^H, K^X).
\end{split}\end{align}

Observe that for $(1,Z)\in G\times \bR^{2n-n_0}$,
 $\varphi_{\var*} K^{X_0}_{(1,Z)} = K^X_{y_0}$ for $K\in \kg$,
by (\ref{a6}), the moment map $\mu_{X_0}: X_0\to \kg^*$ of
the $G$-action on $X_0$ is given by
\begin{align}\label{1c16}
2 \pi \sqrt{-1} \mu_{X_0}(K)_{(1,Z)}
= (1-\varphi (\tfrac{|Z|}{\var}) )  R^{L}_{y_0} (\mR^H,K^X_{y_0})
+ 2 \pi \sqrt{-1}\mu (K)_{\varphi_\var(1, Z)}.
\end{align}
Now from the choice of our coordinate, we know that
$\mu_{X_0}=0$ on $ G\times \bR^{2n-2n_0}\times \{0\}$. Moreover,
\begin{align}\label{2c17}
2 \pi \sqrt{-1}\mu (K)_{\varphi_\var(1, Z)}
= R^{L}_{(1,Z)}(\varphi (\tfrac{|Z|}{\var}) (\mR^\bot)^H, K^X)
+\cO(\varphi (\tfrac{|Z|}{\var})|Z| |Z^\bot|).
\end{align}
From our construction, \eqref{1c16} and  \eqref{2c17}, we know
that
\begin{align}\label{2c18}
\mu_{X_0}^{-1}(0)= G\times \bR^{2n-2n_0}\times \{0\}.
\end{align}
By \eqref{2c16} and  \eqref{1c16}, for $Z\in T_{x_{0}}B$, $|Z|\geq 4 \var$,
\begin{align}\label{1c17}
2 \pi \sqrt{-1} \mu_{X_0}(K)_{(1,Z)}= R^{L}_{y_0} ((\mR^\bot)^H,K^X_{y_0}).
\end{align}

Let $D_p^{X_0}$ be the Dirac operator on $X_0$ associated to the
above data  by the construction in Section \ref{s3.01}.
As in (\ref{b20}), the analogue of Theorems \ref{t2.1}, \ref{ta2}
still holds for $D_p^{X_0}$.

Let $g^{TB_0}$ be the metric on $B_0=\bR^{2n-n_0}$ induced by $g^{TX_0}$,
and let $dv_{B_0}$ be the  Riemannian volume form on
$(B_0, g^{TB_0})$.

 The operator $\Phi \mL_p^{X_0} \Phi ^{-1}$ is also well-defined
 on $T_{x_{0}}B\simeq \bR ^{2n-n_0}$.

Let $P_{x_0,p}$ be the orthogonal projection from
$L^2(\bR ^{2n-n_0},(\Lambda (T^{*(0,1)}X_0)\otimes L^p_0\otimes E_0)_{B_0})$
onto $\Ker (\Phi \mL_p^{X_0} \Phi ^{-1})$ on $\bR ^{2n-n_0}$.
Let $P_{x_0,p}(Z,Z^{\prime})$ $(Z,Z^{\prime}\in\bR ^{2n-n_0})$
be the smooth kernel of $P_{x_0,p}$ with respect to $dv_{B_0}(Z^{\prime})$.
As before, we view $P_{x_0,p}$ as a $G \times G$-invariant section of
$${\rm pr}_1^*(\Lambda (T^{*(0,1)}X_0) \otimes L^p_0\otimes E_0)
\otimes {\rm pr}_2^*(\Lambda (T^{*(0,1)}X_0) \otimes  L^p_0\otimes E_0)^*$$
 on $X_0\times  X_0$.

Let $P_{0,p}^{G}$ be the orthogonal projection from
$\Omega^{0,\bullet}(X_0,L^p_0\otimes E_0)$
onto $(\Ker D_p^{X_0})^G$, and let  $P_{0,p}^{G}(x,x')$ be the smooth kernel
of $P_{0,p}^{G}$
with respect to the volume form $dv_{X_0}(x')$.

 Note that $\Phi$ in (\ref{h12}) defines an isometry from
$(\Ker D^{X_0}_p)^G=\Ker \mL_p^{X_0}$
 onto  $\Ker (\Phi \mL_p^{X_0} \Phi^{-1})$, as in \eqref{b23}, we get
\begin{align}\label{1c19}
h(x)h(x') P_{0,p}^{G}(x,x')= P_{x_0,p} (\pi(x), \pi(x')).
\end{align}

\begin{prop} \label{p3.2} For any $l,m\in \bN$, there exists $C_{l,m}>0$
 such that for $x,x' \in G\times B^{T_{x_0}B}(0,\var)$,
\begin{align}\label{1c20}
\Big |(P_{0,p}^{G}- P_p^G)(x,x')\Big |_{\cC^m}\leq C_{l,m} p^{-l}.
\end{align}
\end{prop}
\begin{proof} By the analogue of Theorems \ref{t2.1}, \ref{ta2}, we know that
for $x,x' \in G\times B^{T_{x_0}B}(0,\var)$,
$P_{0,p}^{G}-\wi{F}(\mL_p^{X_0})$ verifies also (\ref{1c3}),
and for $x,x'\in G\times B^{T_{x_0}B}(0,\var)$,
$$\wi{F}(\mL_p^{X_0})(x,x')=\wi{F}(\mL_p) (x,x')$$
by finite propagation speed.
Thus we get (\ref{1c20}).
\end{proof}

Let $T^{*(0,1)}X_0$ be the anti-holomorphic cotangent bundle of $(X_0,J_0)$.
Since $J_0(g,Z)=  J(\varphi_\var(g, Z))$,
$ T_{Z,J_0}^{*(0,1)}X_0$ is naturally identified with
$T_{\varphi_\var(g,Z),J}^{*(0,1)}X_0$.

Let $\nabla ^{{\rm Cliff}_0}$ be the Clifford connection on
$\Lambda ( T^{*(0,1)}X_0)$ induced by the Levi-Civita connection
$\nabla ^{TX_0}$ on $(X_0, g^{TX_0})$. Let $R^{E_0}, R^{TX_0}$,
$R^{\text{Cliff}_{0}}$ be the corresponding curvatures on
$E_0,TX_0$ and $\Lambda ( T^{*(0,1)}X_0)$ (cf. (\ref{0a3})).

 We identify $\Lambda ( T^{*(0,1)}X_0)_{(g,Z)}$ with
$\Lambda ( T^{*(0,1)}_{(g,0)}X)$ by identifying first
$\Lambda ( T^{*(0,1)}X_0)_{(g,Z)}$
with $\Lambda (T^{*(0,1)}_{\varphi_\var(g,Z),J}X_0)$, which in turn is
identified with $\Lambda ( T^{*(0,1)}_{Gy_0}X)$ by using parallel transport
along $u\to u \varphi_\var(g,Z)$ with respect to
$\nabla ^{{\rm Cliff}_0}$. We also trivialize
$\Lambda ( T^{*(0,1)}X_0)$  in this way.

Let $S_L$ be a $G$-invariant unit section of $L|_{Gy_0}$.
 Using $S_L$ and  the above discussion, we
get an isometry
$$\Lambda ( T^{*(0,1)}X_0)\otimes E_0\otimes L_0^p \simeq
(\Lambda ( T^{*(0,1)}X)\otimes E)|_{\pi^{-1}(x_0)}=\bE|_{\pi^{-1}(x_0)}.$$

For any $1\leq i\leq 2n-n_0$, let $\wi{e}_i(Z)$
be the parallel transport of $e_i$ with respect to
the connection $^0\nabla^{TB}$ along $[0,1]\ni u\to u Z^0$,
and  with respect to the connection $\nabla^{TB}$ along
$[1,2]\ni u\to Z^0+(u-1) Z^\bot$.

If $\alpha = (\alpha_1,\cdots, \alpha_{2n-n_{0}})$ is a multi-index,
set $Z^\alpha = Z_1^{\alpha_1}\cdots Z_{2n-n_{0}}^{\alpha_{2n-n_{0}}}$.

Recall that $A$, $\mR^\bot$ have been defined in \eqref{a0.6}, \eqref{2c15}.

The following Lemma extends \cite[Prop. 1.28]{BeGeVe}
(cf. also \cite[Lemma 4.5]{DLM04a}).
\begin{lemma}\label{l3.2} The Taylor expansion of $\wi{e}_i(Z)$
with respect to the basis $\{e_i\}$ to order $r$ is a polynomial of
 the Taylor expansion of the curvature coefficients of $R^{TB}$
to order $r-2$ and $A$ to order $r-1$.
\end{lemma}
\begin{proof}
Let $\partial_i=\nabla _{e_i}$  be the partial derivatives along $e_i$.

Let $\Gamma ^{TB}$ be the connection form of  $\nabla^{TB}$
 with respect to
the frame $\{\wi{e}_i\}$ of $TB$.
By the definition of our fixed frame, we have
$i_{\mR^\bot}\Gamma ^{TB} =0$.
As in  \cite[(1.12)]{BeGeVe},
\begin{align}\label{0c24}
L_{\mR^\bot}  \Gamma ^{TB} = [i_{\mR^\bot}, d] \Gamma ^{TB}
= i_{\mR^\bot} (d \Gamma ^{TB} + \Gamma ^{TB} \wedge \Gamma ^{TB})
= i_{\mR^\bot} R^{TB}.
\end{align}

Let $\Theta (Z) = (\theta _j^i(Z))_{i,j=1}^{2n-n_0}$
be the $(2 n-n_0) \times (2n-n_0)$-matrix such that
\be\label{0c25}
e_i = \sum_j \theta ^j_i(Z) \wi{e}_j (Z), \quad
\wi{e}_j (Z)= (\Theta (Z)^ {-1})_j^k e_k.
\ee

Set $\theta ^j (Z) = \sum_i \theta ^j_i(Z) e^i$ and
\begin{align}\label{0c26}
&\theta = \sum_j e^j \otimes e_j
= \sum_j \theta ^j \wi{e}_j \in T^*B\otimes TB.
\end{align}

As $\nabla ^{TB}$ is torsion free, $ \nabla ^{TB}\theta=0$. Thus
the $\bR ^{2n-n_0}$-valued one-form $\theta= (\theta ^j(Z))$
 satisfies the structure equation,
\be\label{0c27}
d \theta + \Gamma ^{TB} \wedge \theta =0.
\ee

By the same proof of \cite[Prop. 1.27]{BeGeVe}, we have
\begin{align}\label{0c28}
&\mR^\bot= \sum_j Z^\bot_j \wi{e}^\bot_j (Z),\quad
i_{\mR^\bot} \theta = \sum_j Z^\bot_j e^\bot_j = Z^\bot.
\end{align}
Here under our trivialization by $\{\wi{e}_i\}$,
we consider $Z^\bot= (0,Z^\bot_1,\cdots, Z^\bot_{n_0})$
as a $\bR^{2n-n_0}$-valued function.

Substituting  (\ref{0c28}) and $(L _{\mR^\bot} -1) Z^\bot =0$
into the identity $i_{\mR^\bot} (d \theta +\Gamma ^{TB} \wedge
\theta)=0$, we obtain
\begin{align}\label{0c29}
(L _{\mR^\bot} -1) L _{\mR^\bot} \theta =
(L _{\mR^\bot} -1) ( dZ^\bot + \Gamma ^{TB} Z^\bot)
= (L _{\mR^\bot} \Gamma ^{TB}) Z^\bot
= (i_{\mR^\bot} R^{TB})Z^\bot.
\end{align}
Here we consider $R^{TB}$ as a matrix of $2$-forms, so that $R^{TB}Z^\bot$
is a vector of $2$-forms, and $\theta$ is a $\bR^{2n-n_0}$-valued $1$-form.

By \eqref{0c28} and \eqref{0c29}, we get
\begin{align}\label{0c30}
i_{e_j} (L _{\mR^\bot} -1) L _{\mR^\bot} \theta ^i (Z)=
\left \langle R^{TB} ( \mR^\bot,e_j) \mR^\bot, \wi{e}_i\right \rangle (Z) .
\end{align}

We will denote by $\partial ^\bot$, $\partial ^0$ the partial derivatives
along $N_G$, $TX_G$ respectively.
Then we have the following Taylor expansions of \eqref{0c30}:
 for $j\in \{2(n-n_0)+1, \cdots, 2n-n_0\}$, i.e. $e_j\in N_G$,
 by $L_{\mR^\bot}e^j = e^j$, we have
\begin{align}\label{0c31}
\sum_{|\alpha^\bot| \geq 1} ( |\alpha^\bot|^2 + |\alpha^\bot|)
((\partial^\bot) ^{\alpha^\bot}\theta ^i_j)({Z^0})
\frac{(Z^\bot)^{\alpha^\bot}}{\alpha^\bot !}=
\left \langle R^{TB} ( \mR^\bot,e_j) \mR^\bot, \wi{e}_i\right \rangle (Z).
\end{align}
and for $j\in \{1, \cdots, 2(n-n_0)\}$, i.e. $e_j\in TX_G$,
 by $L_{\mR^\bot}e^j = 0$, we have
\begin{align}\label{0c32}
\sum_{|\alpha^\bot| \geq 1} ( |\alpha^\bot|^2 - |\alpha^\bot|)
((\partial^\bot) ^{\alpha^\bot}\theta ^i_j)({Z^0})
\frac{(Z^\bot)^{\alpha^\bot}}{\alpha^\bot !}=
\left \langle R^{TB} ( \mR^\bot,e_j) \mR^\bot, \wi{e}_i\right \rangle (Z) .
\end{align}

From \eqref{0c31}, \eqref{0c32}, we still need to obtain the Taylor expansions
for $\theta ^i_j({Z^0})$, ($1\leq i,j\leq 2n-n_0$)
and $(\partial^\bot_k\theta ^i_j)({Z^0})$, $(1\leq j\leq 2(n-n_0))$.

By our construction, we know that
 for $i$ or $j\in \{2(n-n_0)+1, \cdots, 2n-n_0\}$,
\begin{align}\label{0c33}
\wi{e}^\bot_k(Z^0)= {e}^\bot_k (Z^0),
\quad \theta ^i_j (Z^0) = \delta_{i,j}.
\end{align}

By \cite[(1.21)]{BeGeVe} (cf. \cite[(4.35)]{DLM04a}),
we know that on $\bR^{2n-2n_0}\times \{0\}$, for
$i,j \in \{1, \cdots, 2(n-n_0)\}$,
\begin{align}\label{0c34}
\begin{split}
&\theta ^i_j(0)= \delta_{i,j} ,\\
&\sum_{|\alpha^0| \geq 1} ( |\alpha^0|^2 + |\alpha^0|)
((\partial^0) ^{\alpha^0}\theta ^i_j)({0})
\frac{(Z^0)^{\alpha^0}}{\alpha^0 !}=
\left \langle R^{TX_G} (\mR^0,e_j) \mR^0, \wi{e}_i\right \rangle (Z^0).
\end{split}
\end{align}
while by (\ref{a0.6}), \eqref{0c25}, and $[e^\bot_i, e^\bot_j]=0$,
we get
\begin{multline}\label{1c34}
(\partial^\bot_k \theta ^i_j )(Z^0)
= e^\bot_k \langle e_j^0, \wi{e}_i^0\rangle (Z^0)
=\langle \nabla^{TB} _{e^\bot_k}e_j^0, \wi{e}_i^0\rangle (Z^0)\\
= \langle \nabla^{TB} _{e_j^0}e^\bot_k, \wi{e}_i^0\rangle (Z^0)
= - \langle \nabla^{TB} _{e_j^0}\wi{e}_i^0, e^\bot_k\rangle (Z^0)
=-\langle  A(e_j^0)\wi{e}_i^0, e^\bot_k \rangle (Z^0).
\end{multline}

Let $R^{TX_G}$, $R^{N_G}$ be the curvatures of $\nabla ^{TX_G},\nabla^{N_G}$.
By \eqref{a0.6},
  \begin{align}\label{c42}
R^{TX_G} + R^{N_G}+ A^2 + {^0\nabla^{TB}_\cdot} A = R^{TB}|_{X_G}
\in \Lambda ^2(TX_G)\otimes \End(TB).
 \end{align}

For $1\leq j\leq 2(n-n_0)$, $2(n-n_0)+1\leq i\leq 2n-n_0$, $i'= i-2(n-n_0)$,
by $[e^\bot_k, e^0_j]=0$, as in \eqref{1c34}, we get
\begin{align}\label{1c35}
(\partial^\bot_k \theta ^i_j )(Z^0)
=  e^\bot_k \langle e_j^0, \wi{e}_{i'}^\bot\rangle (Z^0)
= \langle  \nabla^{TB} _{e_j^0}e^\bot_k, \wi{e}_{i'}^\bot\rangle (Z^0)
= \langle  \nabla^{N_G}_{e_j^0}e^\bot_k, e_{i'}^\bot\rangle (Z^0).
\end{align}
By \cite[Prop. 1.18]{BeGeVe} (cf. \eqref{0c39}) and \eqref{1c35},
the Taylor expansion of $(\partial^\bot_k \theta ^i_j )(Z^0)$ at $0$
to order $r$ only determines by those of $R^{N_G}$ to order $r-1$.

Now by  (\ref{0c25}), (\ref{0c31})-(\ref{1c35}) determine
 the Taylor expansion of $\theta ^i_j(Z)$
to order $m$ in terms of the Taylor expansion of the curvature coefficients
 of $R^{TB}$ to order $m-2$ and  $A$ to order $m-1$.

By \eqref{0c25}, we get Lemma  \ref{l3.2}.
\end{proof}

Let $dv_{TB}$ be the Riemannian volume form on $(T_{x_0}B,
g^{TB})$.

Let $\kappa (Z)$ $(Z\in \bR^{2n-n_0})$ be the smooth positive
function defined by the equation \be\label{u3} dv_{B_0}(Z) =
\kappa (Z) dv_{TB}(Z), \ee with $\kappa (0)=1$.

For  $s \in \cC^{\infty}(\bR^{2n-n_0}, \bE_{x_0})$ and $Z\in
\bR^{2n-n_0}$, for $t=\frac{1}{\sqrt{p}}$, set
\begin{align}\label{c27}
\begin{split}
&(S_{t} s ) (Z) =s (Z/t),\quad    \nabla_{t}=  S_t^{-1}
t\kappa ^{\frac{1}{2}}\nabla ^{E_{p,B_0}} \kappa ^{-\frac{1}{2}}S_t, \\
&  \cL^t_2= S_t^{-1}  t^2\kappa ^{\frac{1}{2}}
 \Phi D_p^{X_0,2}\Phi ^{-1} \kappa ^{-\frac{1}{2}}S_t.
\end{split}\end{align}

As in \eqref{h8}, we denote by $R^{L_B}, R^{E_B}, R^{{\rm Cliff}_B}$
the curvatures on $L_B, E_B$, $ \Lambda (T^{*(0,1)}X)_B$ induced by
$\nabla^L, \nabla^E$,$\nabla^ {{\rm Cliff}}$ on $X$.

As in \eqref{ah7}, $\wi{\mu}\in TY$,  $\wi{\mu}^E\in TY\otimes \End(E)$,
$\wi{\mu}^{{\rm Cliff}}\in TY\otimes \End(\Lambda (T^{*(0,1)}X))$
are sections induced by $\mu, \mu^E, \mu^{{\rm Cliff}}$
 in \eqref{a5}, \eqref{0a9}.

 Denote by  $\nabla_V$ the ordinary differentiation
 operator on $T_{x_0}B$ in the direction $V$.

Denote by $(\partial ^\alpha R^{L_B})_{x_0}$ the tensor
$(\partial^\alpha R^{L_B})_{x_0}(e_i,e_j) :=\partial ^\alpha(
R^{L_B}(e_i,e_j))_{x_0}$.

\begin{thm}\label{t3.3}
There exist $\mA_{i,j,r}$ {\rm (} resp. $\mB_{i,r}$, $\mC_{r}${\rm)}
{\rm(}$r\in \bN, i,j\in \{1,\cdots, 2n-n_0\}${\rm)}  polynomials
in $Z$, and $\mA_{i,j,r}$ is a monomial in $Z$ with degree $r$,
the degree on $Z$ of $\mB_{i,r}$ {\rm(}resp. $\mC_{r}${\rm)}
has the same parity with $r-1$ {\rm(}resp. $r${\rm)},
 with the following properties:

-- the coefficients of $\mA_{i,j,r}$ are polynomials in $R^{TB}$
 {\rm (}resp. A{\rm )} and their derivatives at $x_0$ to order $r-2$
{\rm (}resp. $r-1${\rm )};

-- the coefficients of $\mB_{i,r}$ are polynomials in
$R^{TB}$, $A$, $R^{{\rm Cliff}_B}$, $R^{E_B}$,
 {\rm (}resp. $R^{L_B}${\rm )} and their derivatives at $x_0$ to order $r-1$
{\rm (}resp. $r${\rm )};

--  the coefficients of $\mC_{r}$ are polynomials in
$R^{TB}$, $A$, $R^{{\rm Cliff}_B}$, $R^{E_B}$, $\wi{\mu}^E$,
$\wi{\mu}^{{\rm Cliff}}$ {\rm (}resp. $r^X$, $\tr[R^{T^{(1,0)}X}]$, $R^E$;
resp. $h$, $R^L$, $R^{L_B}$; resp. $\mu${\rm )}
and their derivatives at $x_0$ to order $r-1$
 {\rm (}resp. $r-2$; resp. $r$; resp. $r+1${\rm )}.

--  if we denote by
 \begin{align}\label{0c35}
\begin{split}&\mO_r =  \mA_{i,j,r}\nabla_{e_i} \nabla_{e_j}
+ \mB_{i,r}\nabla_{e_i}+ \mC_{r},\\
&\cL^0_2 =-\sum_{j=1}^{2n-n_{0}} \Big(\nabla_{e_{j}}
+\frac{1}{2}R^{L_B}_{x_0}(\mR, e_j)\Big)^2
-2 \omega_{d,x_0} -\tau_{x_0}+ 4\pi ^2 |P^{TY}\bJ_{x_{0}}\mR|^2,
\end{split}\end{align}
 then
\begin{align}\label{c30}
\cL^t_2=\cL^0_2 + \sum_{r=1}^m t^r \mO_r + \cO(t^{m+1}).
\end{align}
Moreover, there exists $m'\in \bN$ such that for any $k\in \bN$, $t\leq 1$,
$|tZ|\leq \var$,
the derivatives of order $\leq k$ of the coefficients of the operator
 $\cO(t^{m+1})$ are dominated by $C t^{m+1} (1+|Z|)^{m'}$.
\end{thm}
\begin{proof}
Let $\Gamma ^{E_B}$, $\Gamma ^{L_B}$ and $\Gamma ^{\text{Cliff}_B}$
be the connection forms of $\nabla^{E_B}$,
$\nabla^{L_B}$ and $\nabla^{\text{Cliff}_B}$
 with respect to any fixed frames for
$E_B$, $L_B$ and $\Lambda (T^{*(1,0)}X)_B$
 which are parallel along the curve $\gamma_u: [0,1]\ni u\to uZ$
under our trivialization on $B^{T_{x_0}B}(0,\var)$.
Then $\Gamma^{E_B}$ is $\End(\bC^{\dim E})$-valued $1$-form on $\bR^{2n-n_0}$
 and $\Gamma^{L_B}$ is $1$-form on $\bR^{2n-n_0}$.

Now for $\Gamma ^\bullet= \Gamma ^{E_B}, \Gamma ^{L_B}$
or $\Gamma ^{\text{Cliff}_B}$ and
 $R^\bullet= R^{E_B}, R^{L_B}$ or $R^{\text{Cliff}_B}$ respectively,
by the definition of our fixed frame and
\cite[Proposition 1.18]{BeGeVe} (cf. also \cite[(4.45)]{DLM04a}),
the Taylor coefficients of
$\Gamma ^\bullet (e_j) (Z)$ at $x_0$
to order $r$ only determines by those of $R^\bullet$ to order $r-1$, and
\begin{align}\label{0c39}
\sum_{|\alpha|=r}  (\partial^\alpha
 \Gamma ^\bullet ) _{x_0} (e_j) \frac{Z^\alpha}{\alpha !}
=\frac{1}{r+1} \sum_{|\alpha|=r-1}
(\partial^\alpha R^\bullet ) _{x_0}(\mR, e_j)
  \frac{Z^\alpha}{\alpha !}.
\end{align}
Especially,
\begin{align}\label{0c40}
\Gamma ^\bullet _Z(e_j)= \frac{1}{2}R^\bullet_{x_0}(\mR, e_j) + \cO(|Z|^2).
\end{align}

 By (\ref{c27}),  for $t= 1/\sqrt{p}$,
if $|Z|\leq \sqrt{p} \var$, then
\begin{align}\label{c37}
&\nabla_t  = \kappa ^{\frac{1}{2}}(tZ)\Big(\nabla + (t \Gamma ^{{\rm Cliff}_B}
+ t \Gamma ^{E_B}  + \frac{1}{t} \Gamma ^{L_B})  (tZ)\Big)
\kappa ^{-\frac{1}{2}}(tZ).
\end{align}
Moreover, set
\begin{align}\label{ac37}
(\nabla ^{TB}_{e_{i}}e_{j})(Z)=\Gamma_{ij}^k(Z) e_{k}, \quad
g_{ij}(Z)=g^{TB}(e_{i},e_{j})(Z) = \theta^k_i \theta^k_j(Z),
\end{align}
then $\Gamma_{ij}^k$ is the connection form of $\nabla ^{TB}$
with respect to the frame $\{e_i\}$.

Let $(g^{ij})$ be the inverse matrix of $(g_{ij})$, then
\begin{align}\label{c38}
\Delta ^{E_{p,B}} = -\sum_{ij} g^{ij}
\Big(\nabla ^{E_{p,B}}_{e_{i}} \nabla ^{E_{p,B}}_{e_{j}}
- \Gamma_{ij}^k\nabla ^{E_{p,B}}_{e_{k}}\Big),
\end{align}
and by \eqref{ah1}, \eqref{u3},
\begin{align}\label{ac38}
\begin{split}
&\kappa(Z)= (\det g_{ij})^{1/2}(Z),\\
&\Gamma_{ij}^k= \frac{1}{2} g^{kl} (\partial_i g_{jl}
+\partial_j g_{il} -\partial_l g_{ij}).
\end{split}\end{align}

By (\ref{b4}), (\ref{c27}) and (\ref{c38}),
\begin{multline}\label{1ue1}
 \cL^t_2(Z) = - g^{ij}(tZ) (\nabla_{t,e_i}\nabla_{t,e_j}
 -t \Gamma_{ij}^k(tZ) \nabla_{t,e_k})
- \langle t \wi{\mu}^{E_{p}},t \wi{\mu}^{E_{p}}\rangle_{g^{TY}} (tZ) \\
-2  \om_d(tZ)- \tau(tZ)+ t^2 \Big(\frac{1}{4}r^X+ \mathbf{c}(R)
-\frac{1}{h}\Delta_{B_0}h\Big)(tZ).
\end{multline}

By \eqref{0a9},
 \begin{align}\label{c39}
 \langle t \wi{\mu}^{E_{p}},t \wi{\mu}^{E_{p}}\rangle_{g^{TY}}
= -4\pi ^2 |\frac{1}{t}\wi{\mu}|^2_{g^{TY}}
+   \langle  4\pi \sqrt{-1} \wi{\mu}
+ t^2(\wi{\mu}^{{\rm Cliff}} +\wi{\mu}^E),
\wi{\mu}^{{\rm Cliff}} +\wi{\mu}^E\rangle_{g^{TY}}.
\end{align}

By \eqref{0a0}, \eqref{a5}, and $\wi{\mu}_{y_0}=0$,
for  $y_0\in P, \pi(y_0)=x_0$,
we get for $K\in \kg$,
  \begin{align}\label{c40}
- \langle \bJ  e_i^H, K^X\rangle_{y_0} = \omega(K^X, e_i^H)
= \nabla_{e_i^H}(\mu(K))
= \langle \nabla_{e_i^H}^{TY} \wi{\mu}, K^X\rangle_{y_0},
 \end{align}
thus
  \begin{align}\label{c41}
|\wi{\mu}|^2_{g^{TY}} (Z)=  |\nabla^{TY}_\mR \wi{\mu}|^2_{g^{TY}}
+ \cO(|Z|^3)= |P^{TY}\bJ_{x_{0}}\mR|^2+ \cO(|Z|^3).
 \end{align}

By  Lemma \ref{l3.2},
 \eqref{0c39}, \eqref{c37}, \eqref{1ue1} and \eqref{c41},
we know that $\cL^t_2$ has the expansion \eqref{c30},
in particular, we get the formula $\cL^0_2$ in \eqref{0c35}.

By \eqref{c42}, \eqref{0c39} and \eqref{1ue1}, we get the properties on
$\mA_{i,j,r}$, $\mB_{i,r}$.

By  \eqref{c42},  \eqref{1ue1} and \eqref{c39},
we get the properties on $\mC_r$.

The proof of Theorem \ref{t3.3} is complete.
\end{proof}

\subsection{Uniform estimate on the $G$-invariant Bergman kernel}
\label{s3.3}

Recall that the operators $\cL^t_2$, $\nabla_{t}$ were defined in (\ref{c27}),
 and ${\bf E}_0= \Lambda(T^{*(0,1)}X_0)\otimes E_0$.
We have trivialized
 the bundle ${\bf E}_{0,B_0}$ to ${\bf E}_{B,x_0}$ in Section \ref{s3.2}.
We still denote by $h^{{\bf E}_{0,B_0}}$ the metric on the trivial
bundle ${\bf E}_{B,x_0}$ on $\bR^{2n-n_0 }$ induced by the corresponding
metric on ${\bf E}_{0,B_0}$. Note that $h^{{\bf E}_{0,B_0}}$ is not
a constant metric on  $\bR^{2n-n_0 }$.

We also denote by $\left\langle \ ,\ \right\rangle_{0,L^2}$
 and $\|\ \|_{0,L^2}$
the scalar product and the $L^2$ norm on $\cC^\infty (T_{x_0}B, \bE_{B,x_0})$
induced by $g^{T_{x_0}B}, h^{{\bf E}_{0,B_0}}$ as in (\ref{h10}).

Let $\wi{\mu}_{X_0}, \wi{\mu}^{E_{0,p}}$ be
the $G$-invariant sections of $TY$, $TY\otimes \End(E_{0,p})$ on $X_0$
induced by $\mu_{X_0}$, $\mu^{E_{0,p}}$ as in \eqref{ah7}.

Let $\{f_l\}$ be a $G$-invariant orthonormal frame of $TY$
on $\pi^{-1}(B^B(x_0,\var))$, then
$(f_{0,l})_Z=(f_l) _{\varphi_\var(Z)}$ is a
 $G$-invariant orthonormal frame of $TY_0$ on $X_0$.

\begin{defn} \label{tu0} Set
\begin{align}\label{u2}
\mD_t = \{ \nabla_{t,e_i}, 1\leq i\leq 2n-n_0;
 \frac{1}{t}\langle \wi{\mu}_{X_0}, f_{0,l}\rangle (tZ), 1\leq j\leq n_0\}.
\end{align}
 For $k\in \bN^*$, let $\mD_t^k$ be the family of operators acting on
$\cC^\infty (T_{x_0}B, \bE_{B,x_0})$ which can be written in the form
$Q= Q_1\cdots Q_k$, $Q_i\in \mD_t$.
\end{defn}


For $s\in \cC^{\infty}( T_{x_0}B, \bE_{B,x_0}) $, $k\geq 1$, set
\begin{align}\label{u4}
\begin{split}
&\|s\|_{t,0}^2 = \int_{\bR^{2n-n_0}} |s(Z)|^2_{h^{{\bf E}_{0,B_0}}(tZ)}dv_{T_{x_0}B}(Z)= t^{-2n+n_0}\|S_t s\|_{0, L^2}^2,\\
&\| s \|_{t,k}^2 =  \| s\|_{t,0}^2 +\sum_{l=1}^k \sum_{Q\in \mD_t^l}
\| Q s\|_{t,0}^2 .
\end{split}\end{align}

We denote by $\left\langle s ', s \right\rangle_{t,0}$
 the inner product on $\cC^\infty (T_{x_0}B, \bE_{B,x_0})$
 corresponding to $\|\quad\|^2_{t,0}$.

Let $H^m_t$ be the Sobolev space of order $m$ with norm $\|\quad\|_{t,m}$.
Let $H^{-1}_t$ be the Sobolev space of order $-1$ and let $\|\quad\|_{t,-1}$
 be the norm on  $H^{-1}_t$ defined by
$\|s\|_{t,-1} = \sup_{0\neq s'\in  H^1_t }$
$|\left\langle s,s'\right\rangle_{t,0}|/\|s'\|_{t,1}$.

If $A\in \cL (H^{m}_t, H^{m'}_t)$ $(m,m' \in \bZ)$,
 we denote by  $\|A\|^{m,m'}_t$ the norm of $A$ with respect to the norms
$\| \quad \|_{t,m}$ and $\| \quad \|_{t,m'}$.

Then $\cL^t_2$ is a formally self-adjoint elliptic operator with respect to
 $\|\quad\|^2_{t,0}$, and is a smooth family of
operators with respect to the parameter $x_0\in X_G$.

\begin{thm}  \label{tu1} There exist constants $ C_1, C_2, C_3>0$
such that for $t\in ]0,1]$ and any
$s,s'\in C ^{\infty}_0(\bR^{2n-n_0}, \bE_{B,x_0})$,
\begin{align}\label{u5}
\begin{split}& \left\langle \cL^t_2 s,s\right\rangle_{t,0}
\geq C_1\|s\|_{t,1}^2 -C_2 \|s\|_{t,0}^2  , \\
& |\left\langle \cL^t_2 s, s'\right\rangle_{t,0}|
 \leq C_3 \|s\|_{t,1}\|s'\|_{t,1}.
\end{split}\end{align}
\end{thm}

\begin{proof}
By (\ref{1c17}) and our construction
for $L_0, E_0$ on $X_0$, we know for $Z\in T_{x_0} B$, $|Z|>4\var$,
\begin{align}\label{0u5}
\mu ^{E_{0,p}}(K)_{(1,Z)}= p\,   R^L_{y_0}((\mR^\bot)^H, K^X_{y_0}).
\end{align}
Thus from \eqref{1ue1} and (\ref{u4}),
\begin{multline}\label{u7}
\left\langle  \cL^t_2 s,s\right\rangle_{t,0} = \|\nabla_ts\|_{t,0}^2
- t^2 \left\langle \langle \wi{\mu}^{E_{0,p}}, \wi{\mu}^{E_{0,p}}\rangle_{g^{TY}}(tZ) s,s\right\rangle_{t,0}\\
 + \left\langle\left  (-2 S_t^{-1}\om_d
-S_t^{-1}\tau + t^2S_t^{-1}(\tfrac{1}{4}r^X+\mathbf{c}(R)
-\frac{1}{h}\Delta_{B_0}h)\right )s,
s\right\rangle_{t,0}.
\end{multline}
From  (\ref{1c16}), (\ref{c39}), (\ref{0u5}),
and our construction on $\nabla ^{E_0}$,
\begin{align}\label{u8}
-t^2 \left\langle \langle \wi{\mu}^{E_{0,p}},
 \wi{\mu}^{E_{0,p}}\rangle_{g^{TY}}(tZ) s,s\right\rangle_{t,0}\
\geq 2\pi ^2 \sum_{l=1}^{n_0} \Big\|\frac{1}{t}
\langle \wi{\mu}_{X_0}, f_{0,l}\rangle (tZ)s\Big\|_{t,0}^2 - Ct\| s\|_{t,0}^2.
\end{align}
From (\ref{u7}) and (\ref{u8}), we get (\ref{u5}).
\end{proof}

Recall that $\nu$ is the constant in  \eqref{a12}.

Let $\delta$ be the counterclockwise oriented circle in $\bC$
of center $0$ and radius $\nu/4$,
and let $\Delta$ be the oriented path in $\bC$ which goes parallel to the real
axis from $+\infty +i$ to $\frac{\nu}{2}+i$ then parallel to the imaginary
 axis to $\frac{\nu}{2}-i$ and the parallel to the real
axis to $+\infty -i$.
\begin{figure}[ht]
\centerline{\epsfxsize=3.1in\epsfbox{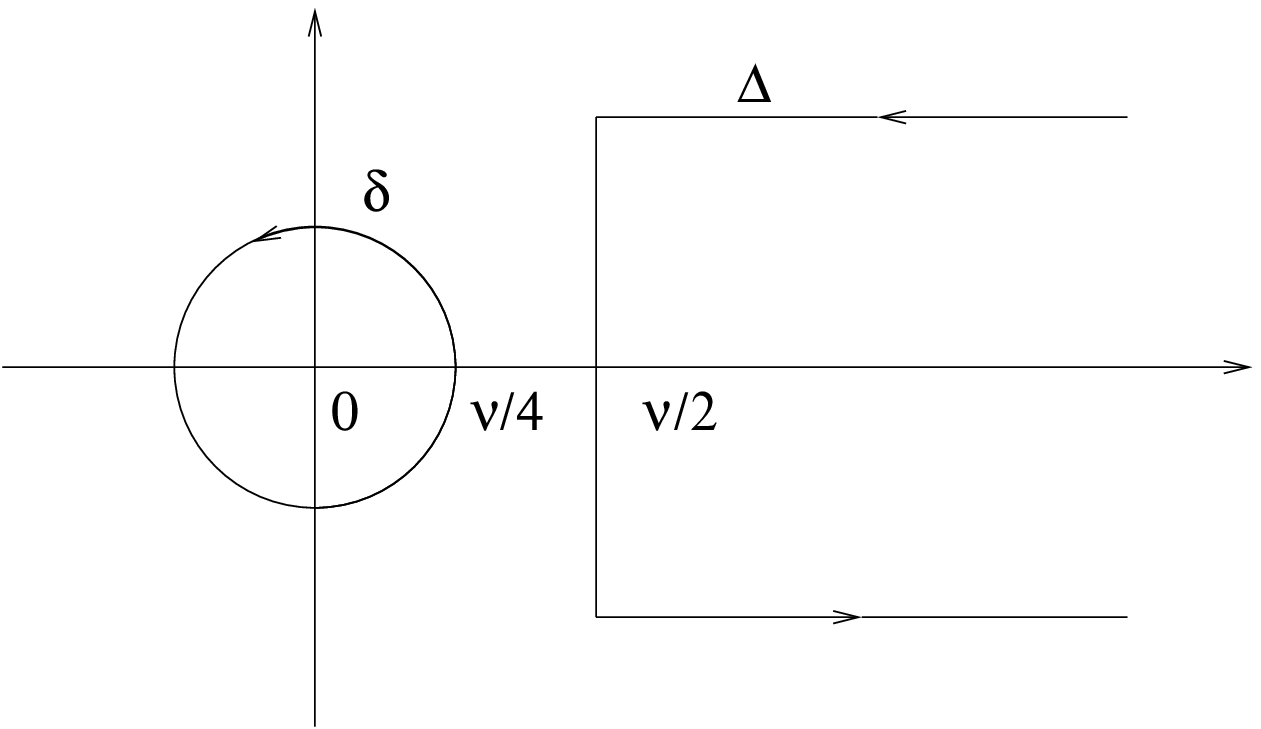}}
\end{figure}

\begin{thm}\label{tu4}  There exist $t_0>0$, $C>0$ such that for
$t\in ]0,t_0]$, $\lambda \in  \delta\cup \Delta$
and $x_0\in X_G$, $(\lambda-\cL^t_2)^{-1}$ exists and
\begin{align}\label{ue2}
\begin{split}& \| (\lambda- \cL^t_2)^{-1}\|^{0,0}_t \leq C,\\
& \| (\lambda- \cL^t_2)^{-1}\|^{-1,1}_t
\leq C (1+|\lambda |^2).
\end{split}\end{align}
\end{thm}
\begin{proof}
By (\ref{a12}), (\ref{b4}) for $D^{X_0}_p$, and (\ref{c27}),
there exists $t_0>0$ such that for $t\in ]0,t_0]$,
\begin{align}\label{u9}
\spec \left( \cL^t_2\right)\subset \{0\}\cup
\left[\nu,+ \infty \right[.
\end{align}
Thus $(\lambda- \cL^t_2)^{-1}$ exists for $\lambda \in \delta\cup \Delta$.

 The first inequality of (\ref{ue2}) is from (\ref{u9}).

By (\ref{u5}), for $\lambda_0\in\bR$, $\lambda_0\leq -2C_2$,
 $(\lambda_0- \cL^t_2)^{-1}$ exists, and we have
$\|(\lambda_0- \cL^t_2)^{-1}\|^{-1,1}_t \leq \frac{1}{C_1}$. Now,
\be\label{ue7}
(\lambda- \cL^t_2)^{-1}= (\lambda_0- \cL^t_2)^{-1}
- (\lambda-\lambda_0) (\lambda- \cL^t_2)^{-1}(\lambda_0- \cL^t_2)^{-1}.
\ee
Thus for  $\lambda\in \delta\cup \Delta$, from (\ref{ue7}), we get
\be\label{ue8}
\|(\lambda-\cL^t_2)^{-1}\|^{-1,0}_t  \leq
\frac{1}{C_1} \Big(1+\frac{4}{\nu}|\lambda-\lambda_0|\Big).
\ee

Now we  change the last two factors in (\ref{ue7}), and apply (\ref{ue8}),
we get
\begin{align}\label{ue9}
\|(\lambda-\cL^t_2)^{-1}\|^{-1,1}_t & \leq  \frac{1}{C_1}+
 \frac{|\lambda-\lambda_0|}{{C_1}^2} \Big(1+\frac{4}{\nu}|\lambda
 -\lambda_0|\Big)\\
&\leq C(1+|\lambda|^2).\nonumber
\end{align}

The proof of our Theorem is complete.
\end{proof}

\begin{prop} \label{tu5} Take $m \in \bN^*$. There exists $C_m>0$ such that
 for $t\in ]0,1]$,
$Q_1, \cdots, Q_m$ $\in \mD_t$ $\cup \{Z_i\}_{i=1}^{2n-n_0}$
and  $s,s'\in \cC^{\infty}_{0}(\bR^{2n-n_0}, \bE_{B,x_0})$,
\be\label{ue11}
\left |\left\langle [Q_1,
[Q_2,\ldots, [Q_m,  \cL^t_2]] \ldots ]s,
s'\right\rangle_{t,0} \right |
\leq C_m  \|s\|_{t,1} \|s'\|_{t,1}.
\ee
\end{prop}
\begin{proof} Note that $[\nabla_{t,e_i}, Z_j]=\delta_{ij}$.
By (\ref{1ue1}), we know that $[Z_j, \cL^t_2]$ verifies (\ref{ue11}).

Recall that by (\ref{1c16}) and (\ref{1c17}),
$(\nabla_{e_i}\langle \wi{\mu}_{X_0}, f_{0,l}\rangle)(tZ)$
is uniformly bounded with its derivatives for $t\in [0,1]$ and
 \begin{align}\label{1ue0}
\nabla_{e_i}\langle \wi{\mu}_{X_0}, f_{0,l}\rangle = (e_i \langle
\wi{\mu}_{X_0}, f_{0,l}\rangle)_{x_0} =\om(f_{0,l}, e_i)_{x_0}
\end{align}
for $|Z|\geq 4\var$. Thus $[\frac{1}{t}\langle \wi{\mu}_{X_0},
f_{0,l}\rangle(tZ), \cL^t_2]$ also verifies  (\ref{ue11}).

Note that by  (\ref{c27}),
\begin{align}\label{1ue2}
[\nabla_{t,e_i},\nabla_{t,e_j}]= \left(R^{L_{0,B_0}}(tZ)
+ t^2  R^{{\bf E}_{0,B_0}}(tZ)\right)(e_i,e_j).
\end{align}
Thus from (\ref{1ue1}), (\ref{1ue0}) and (\ref{1ue2}), we know that
$[\nabla_{t,e_k},  \cL^t_2]$
has the same structure as  $\cL^t_2$ for $t \in ]0,1]$,
i.e. $[\nabla_{t,e_k},  \cL^t_2]$ has the type as
\begin{multline}\label{1ue3}
\sum_{ij} a_{ij}(t,tZ) \nabla_{t,e_i}\nabla_{t,e_j}
+\sum_{i}c_{i}(t,tZ) \nabla_{t,e_i} \\
 + \sum_l \Big [c'_{l}(t,tZ)
\frac{1}{t}\langle \wi{\mu}_{X_0}, f_{0,l}\rangle(tZ)
+ a | \frac{1}{t}\wi{\mu}_{X_0}|^2_{g^{TY}}(tZ)\Big ]+ c(t,tZ),
\end{multline}
where  $a\in \bC$;
$a_{ij}(t,Z),c_{i}(t,Z), c'_{j}(t,Z), c(t,Z)$
 and their derivatives on $Z$
are uniformly bounded for $Z\in \bR ^{2n-n_0}, t\in [0,1]$;
moreover, they are polynomials in $t$.
In fact, for $[\nabla_{t,e_k},  \cL^t_2]$,  $a=0$ in \eqref{1ue3}.

Let $(\nabla_{t,e_i})^*$ be the adjoint of $\nabla_{t,e_i}$ with respect to
$\left\langle \ , \ \right\rangle_{t,0}$, then by (\ref{u4}),
\begin{align}\label{1ue4}
(\nabla_{t,e_i})^* =- \nabla_{t,e_i}
-t (k^ {-1}\nabla_{e_i}k)(tZ),
\end{align}
the last term of (\ref{1ue4}) and its derivatives in  $Z$
 are uniformly bounded in $Z\in \bR ^{2n-n_0}, t\in [0,1]$.

By (\ref{1ue3}) and (\ref{1ue4}), (\ref{ue11}) is verified for $m=1$.

By iteration, we know that $[Q_1,[Q_2,\ldots, [Q_m,  \cL^t_2]] \ldots ]$
has the same structure (\ref{1ue3})  as $\cL^t_2$.
By (\ref{1ue4}), we get Proposition \ref{tu5}.
\end{proof}

\begin{thm}\label{tu6} For any $t\in ]0,t_0]$, $\lambda \in \delta\cup \Delta$,
 $m \in \bN$, the resolvent $(\lambda-\cL^t_2)^{-1}$ maps $H^m_t$
into $H^{m+1}_t$. Moreover for any $\alpha\in \bZ^{2n-n_0}$,
 there exist $N \in \bN$, $C_{\alpha, m}>0$
such that for  $t\in ]0,t_0]$, $\lambda \in \delta\cup \Delta$,
 $s\in \cC ^\infty_0 (\bR^{2n-n_0},\bE_{B,x_0})$,
\be\label{ue12}
\|Z^\alpha (\lambda-\cL^t_2)^{-1} s\|_{t,m+1}
\leq C_{\alpha, m}  (1+|\lambda|^2)^N \sum_{\alpha' \leq \alpha}
\|Z^{\alpha'} s\|_{t,m} .
\ee
\end{thm}
\begin{proof} For $Q_1, \cdots, Q_m\in \mD_t$,
$Q_{m+1},\cdots, Q_{m+|\alpha|}\in\{Z_i\}_{i=1}^{2n-n_0}$,  we can
express $Q_1\cdots$ $Q_{m+|\alpha|}(\lambda-\cL^t_2)^{-1}$ as a
linear combination of operators of the type \be\label{ue13} [Q_1,
[ Q_2,\ldots [Q_{m'},(\lambda-\cL^t_2)^{-1}]]\ldots ]Q_{m'+1}
\cdots Q_{m+|\alpha|}, \quad m'\leq m+|\alpha|. \ee

Let $\cR_{t}$ be the family of operators $$\cR_{t} = \{ [ Q_{j_1},
[Q_{j_2},\ldots [Q_{j_l},\cL^t_2]]\ldots ] \}. $$

Clearly, any commutator $[Q_1, [ Q_2,\ldots
[Q_{m'},(\lambda-\cL^t_2)^{-1}]]\ldots ]$ is a linear combination
of operators of the form \be\label{ue14}
(\lambda-\cL^t_2)^{-1}R_1(\lambda-\cL^t_2)^{-1}R_2 \cdots
R_{m'}(\lambda-\cL^t_2)^{-1} \ee with  $R_1, \cdots, R_{m'} \in
\cR_{t}$.

By Proposition \ref{tu5}, the norm $\|\quad\|_t^{1,-1}$ of the
 operators $R_j\in \cR_{t}$ is
uniformly bound by $C$.

By Theorem \ref{tu4}, we find that there exist $C>0$, $N \in \bN$
such that the norm $\|\quad\|_t^{0,1}$
of operators (\ref{ue14}) is dominated by
$C  (1+|\lambda|^2)^N$.
\end{proof}

Let $\pi_B : TB \times_{B} TB \to B$ be the natural projection from the
fiberwise product of $TB$ on $B$.

Let $e ^{-u\cL^t_2}(Z,Z^{\prime})$, $(\cL^t_2e ^{-u\cL^t_2})(Z,Z^{\prime})$
be the smooth kernels of
 the operators  $e ^{-u\cL^t_2}$, $\cL^t_2e ^{-u\cL^t_2}$
with respect to $dv_{T_{x_0}B}(Z^{\prime})$.

Note that $\cL^t_2$ are families of differential operators with coefficients
in $\End (\bE_{B,x_0}) =\End (\Lambda (T^{*(0,1)}X)\otimes E)_{B,x_0}$.
Thus we can view  $e ^{-u\cL^t_2}(Z,Z^{\prime})$, $(\cL^t_2e
^{-u\cL^t_2})(Z,Z^{\prime})$ as smooth sections of $\pi_B ^* (\End
(\Lambda (T^{*(0,1)}X)\otimes E)_{B})$ on  $TB\times_{B} TB$.

 Let $\nabla ^{\End (\bE_{B})}$ be the connection on
$\End (\Lambda (T^{*(0,1)}X)\otimes E)_{B}$ induced by
$\nabla ^{\text{Cliff}_{B}}$ and $\nabla ^{E_{B}}$.
And $\nabla ^{\End (\bE_{B})}$, $h^E$ and $g^{TX}$
induce naturally a $\cC^m$-norm for the parameter $x_0\in X_G$.

As in Introduction, for $Z\in T_{x_0}B$, we will write
 $Z=Z^0+Z^\bot$, with $Z^0\in T_{x_0}X_G$, $Z^\bot\in N_{G,x_0}$.

\begin{thm}\label{tue8} There exists $C'' >0$ such that
 for any $m,m' , m'', r \in \bN$, $u_0>0$, there exists $C>0$
 such that for $t \in ]0,t_0]$, $u\geq u_0$, $Z,Z^{\prime}\in T_{x_0}B$,
 \begin{align}\label{ue15}
\begin{split}&\sup_{|\alpha|+|\alpha'|\leq m}
 (1+|Z^\bot|+|Z^{\prime\bot}|)^{m''}\Big |\frac{\partial^{|\alpha|+|\alpha'|}}
{{\partial Z}^{\alpha} {\partial Z^{\prime}}^{\alpha'}}
\frac{\partial^{r}}{\partial t^{r}}
e ^{-u\cL^t_2} \left (Z, Z^{\prime}\right )\Big |_{\cC^{m'}(X_G)} \\
&\hspace{10mm} \qquad   \leq C (1+|Z^0|+|Z^{\prime 0}|)^{2(n+r+m^\prime+1)+m}
\exp \Big(\frac{1}{2}\nu u- \frac{2C''}{u} |Z-Z^{\prime}|^2\Big), \\
&\sup_{|\alpha|+|\alpha'|\leq m}
(1+|Z^{\bot}|+|Z^{\prime\bot}|)^{m''}
\Big |\frac{\partial^{|\alpha|+|\alpha'|}}
{{\partial Z}^{\alpha} {\partial Z^{\prime}}^{\alpha'}}
\frac{\partial^{r}}{\partial t^{r}}
(\cL^t_2e ^{-u\cL^t_2})  \left (Z, Z^{\prime}\right )
\Big |_{\cC^{m'}(X_G)}  \\
&\hspace{10mm} \qquad   \leq C (1+|Z^0|+|Z^{\prime
0}|)^{2(n+r+m^\prime+1)+m} \exp \Big( -\frac{1}{4}\nu u- \frac{2C''}{u}
|Z-Z^{\prime}|^2\Big),
\end{split}\end{align}
where $\cC^{m'}(X_G)$ is the $\cC^{m'}$ norm for the parameter
$x_0\in X_G$.
\end{thm}
\begin{proof} By (\ref{u9}), for any $k\in \bN^*$,
\begin{align}\label{1ue15}
\begin{split}
&e ^{-u\cL^t_2}= \frac{(-1)^{k-1} (k-1)!}{2\pi i u^{k-1}}
\int_{\delta\cup\Delta} e^{-u\lambda} (\lambda - \cL^t_2)^{-k} d \lambda,\\
&\cL^t_2e ^{-u\cL^t_2}= \frac{(-1)^{k-1} (k-1)!}{2\pi i u^{k-1}}
\int_{\Delta} e^{-u\lambda}\Big[ \lambda (\lambda - \cL^t_2)^{-k}
-  (\lambda - \cL^t_2)^{-k+1} \Big]  d \lambda.
\end{split}\end{align}

From Theorem \ref{tu6}, we deduce that
if $Q\in \cup_{l=1}^m\mD_t^l$, there are $N\in \bN$, $C_m>0$ such that
for any $\lambda \in \delta\cup \Delta$,
\be\label{ue18}
\| Q (\lambda-\cL^t_2)^{-m}\|_t^{0,0}
 \leq C_m  (1+|\lambda |^2) ^N.
\ee

Recall that $\cL_t^2$ is self-adjoint with respect to $\|\, \,
\|_{t,0}$. After taking the adjoint of (\ref{ue18}),  we get
\be\label{ue20} \| (\lambda-\cL^t_2)^{-m}Q\|_t^{0,0} \leq C_m
(1+|\lambda |^2)^N. \ee

From  (\ref{1ue15}), (\ref{ue18}) and (\ref{ue20}), we get
  if $Q,Q' \in \cup_{l=1}^m\mD_t^l$,
\begin{align}\label{ue21}
\begin{split}
&\|Q e ^{-u\cL^t_2}Q' \|^{0,0}_t \leq C_{m} e^{\frac{1}{4}\nu u},\\
&\|Q (\cL^t_2e ^{-u\cL^t_2})Q' \|^{0,0}_t \leq C_{m} e^{-\frac{1}{2}\nu u}.
\end{split}\end{align}

Let $|\quad|_m$ be the usual Sobolev norm on
$\cC^\infty(\bR^{2n-n_0}, \bE_{B,x_0})$ induced by
$h^{\bE_{B,x_0}}= h^{(\Lambda (T^{*(0,1)}X)\otimes E)_{B,x_0}}$
and the volume form $dv_{T_{x_0}B}(Z)$ as in (\ref{u4}).

Observe that by (\ref{c37}), (\ref{u4}), there exists $C>0$ such
that for $s\in \cC^\infty (T_{x_0}B, \bE_{B,x_0})$, $\supp (s)
\subset B^{T_{x_0}B}(0,q)$, $m\geq 0$,
\begin{align}\label{1ue21}
&\frac{1}{C} (1+q)^{-m} \|s\|_{t,m}\leq |s|_m \leq C(1+q)^m \|s\|_{t,m}.
\end{align}

Now (\ref{ue21}), (\ref{1ue21})  together with Sobolev's
inequalities imply  that if $Q,Q' \in \cup_{l=1}^m\mD_t^l$, for
$\mathcal{K}_u (\cL^t_2) =e ^{-\frac{1}{4}\nu u} e ^{-u\cL^t_2}$
or $e ^{\frac{1}{2}\nu u} \cL^t_2 e ^{-u\cL^t_2}$, we have
\begin{align}\label{ue23}
&\sup_{|Z|,|Z^{\prime}|\leq q}| Q_Z Q'_{Z^{\prime}}
\mathcal{K}_u (\cL^t_2) (Z,Z^{\prime})  | \leq C(1+q)^{2n+2}.
\end{align}

By  (\ref{1c16}), (\ref{2c17}) and \eqref{1c17},
\begin{align}\label{0ue23}
&\sum_{l=1}^{n_0} \Big|\frac{1}{t}\langle \wi{\mu}_{X_0},
f_{0,l}\rangle (tZ)\Big|^2
= |\frac{1}{t}\wi{\mu}_{X_0}|^2 _{g^{TY}} (tZ)
\geq C |Z^\bot|^2 .
\end{align}

Thus by (\ref{c37}), (\ref{ue23}), (\ref{0ue23}), we derive
(\ref{ue15}) with the exponentials $e^{\frac{1}{4}\nu u}$,
$e^{-\frac{1}{2}\nu u}$ for the case when $r =m'=0$ and $C''=0$, i.e.
\begin{align}\label{0ue24}
&\sup_{|\alpha|+|\alpha'|\leq m}
 (1+|Z^\bot|+|Z^{\prime\bot}|)^{m''}\Big |\frac{\partial^{|\alpha|+|\alpha'|}}
{{\partial Z}^{\alpha} {\partial Z^{\prime}}^{\alpha'}}
\mathcal{K}_u (\cL^t_2) \left (Z, Z^{\prime}\right )\Big |\\
&\hspace{10mm} \qquad   \leq C (1+|Z^0|+|Z^{\prime 0}|)^{2n+m+2}.\nonumber
\end{align}

To obtain (\ref{ue15}) in general, we proceed as in the proof of
\cite[Theorem 11.14]{B95}.

Note that the function $f$ is defined in (\ref{c2}).
For $\varrho >1$, put
\begin{align}  \label{ue24}
&K_{u,\varrho}(a)= \int_{-\infty}^{+\infty} \exp(i v\sqrt{2u} a)
 \exp(-\frac{v^2}{2})
\Big (1-f (\frac{1}{\varrho}\sqrt{2u} v) \Big ) \frac{dv}{\sqrt{2\pi}}.
\end{align}
Then there exist $C',C_1>0$ such that for any $c>0$, $m,m'\in \bN$,
there is $C>0$ such that for $ u\geq u_0$,
$a\in \bC, |{\rm Im} (a)|\leq c$,  we have
\begin{align} \label{1ue24}
|a|^m |K_{u,\varrho}^{(m')}(a)| \leq C \exp \Big( C'c^2 u- \frac{C_1}{u} \varrho^2\Big).
\end{align}

For any $c>0$, let  $V_c$ be the images of $\{\lambda \in \bC,
|\Im (\lambda) |\leq c\}$ by the map $\lambda \to \lambda ^2$.
Then $$ V_c=\{\lambda \in \bC, {\rm Re} (\lambda) \geq
\frac{1}{4c^2} {\rm Im} (\lambda)^2 -c^2\},$$ and $\delta \cup
\Delta\subset V_c$ for $c$ large enough.

Let $\wi{K}_{u,\varrho}$ be the holomorphic function such that
$\wi{K}_{u,\varrho}(a^2)=K_{u,\varrho}(a)$.

By (\ref{1ue24}), for $\lambda \in V_c$, \be\label{ue25}
|\lambda|^m |\wi{K}_{u,\varrho}^{(m')}(\lambda)| \leq C
\exp\Big(C'c^2 u - \frac{C_1}{u} \varrho^2\Big). \ee

Using finite propagation speed of solutions of hyperbolic
equations
 and (\ref{ue24}), we find that there exists a fixed constant
 (which depends on $\varepsilon$) $c'>0$ such that
\begin{align}\label{1ue25}
\wi{K}_{u,\varrho}(\cL^t_2)(Z,Z^{\prime}) = e ^{-u \cL^t_2}(Z,Z^{\prime})\quad {\rm if} \,\,
|Z-Z^{\prime}|\geq c' \varrho.
\end{align}

By (\ref{ue25}), we see that given $k\in \bN$, there is a unique
holomorphic function $\wi{K}_{u,\varrho,k} (\lambda)$ defined on a neighborhood
of $V_c$ such that  it verifies the same estimates as $\wi{K}_{u,\varrho}$ in (\ref{ue25}) and
$\wi{K}_{u,\varrho,k} (\lambda) \to 0$ as $\lambda \to +\infty$; moreover
\be\label{ue16}
\wi{K}_{u,\varrho,k}^{(k-1)} (\lambda)/(k-1)! = \wi{K}_{u,\varrho} (\lambda).
\ee
Thus as in (\ref{1ue15}),
\begin{align} \label{1ue26}
\begin{split}
&\wi{K}_{u,\varrho}(\cL^t_2)= \frac{1}{2\pi i} \int_{\delta\cup\Delta}\wi{K}_{u,\varrho,k}
(\lambda)(\lambda-\cL^t_2)^{-k} d \lambda,\\
&\cL^t_2\wi{K}_{u,\varrho}(\cL^t_2)= \frac{1}{2\pi i} \int_{\Delta}
\wi{K}_{u,\varrho,k} (\lambda)\Big[\lambda (\lambda-\cL^t_2)^{-k}
- (\lambda-\cL^t_2)^{-k+1}\Big]  d \lambda.
\end{split}\end{align}

By (\ref{ue18}), (\ref{ue20}) and by proceeding
as in (\ref{ue21})-(\ref{ue23}),
we find that for ${\bf K}_u(a)= \wi{K}_{u,\varrho}(a)$
or $a\wi{K}_{u,\varrho}(a)$, for  $|Z|, |Z^{\prime}|\leq q$,
 \begin{multline}  \label{ue27}
\sup_{|\alpha|+|\alpha'|\leq m}
 (1+|Z^{\bot}|+|Z^{\prime\bot}|)^{2n+m+m''+2} \Big |\frac{\partial^{|\alpha|+|\alpha'|}}
{\partial Z^{\alpha} {\partial Z^{\prime}}^{\alpha'}}
{\bf K}_u(\cL^t_2) (Z, Z^{\prime})\Big |\\
 \leq C (1+q)^{2n+2+m}
\exp (C'c^2 u- \frac{C_1}{u} \varrho^2).
\end{multline}

Setting $\varrho\in \bN^*$, $|\varrho-
\frac{1}{c'}|Z-Z^{\prime}||<1$ in (\ref{ue27}), we get for
$\alpha, \alpha'$ verifying  $|\alpha|+|\alpha'|\leq m$,
\begin{multline}  \label{1ue27}
(1+|Z^{\bot}|+|Z^{\prime\bot}|)^{m''}\Big |\frac{\partial^{|\alpha|+|\alpha'|}}
{\partial Z^{\alpha} {\partial Z^{\prime}}^{\alpha'}}
{\bf K}_u(\cL^t_2) (Z, Z^{\prime})\Big |\\
 \leq C (1+|Z^0|+|Z^{\prime 0}|)^{2n+m+2}
\exp (C'c^2  u- \frac{C_1}{2{c'}^2u} |Z-Z^{\prime}|^2).
\end{multline}

Take $\delta_1= \frac{C'c^2 +\frac{1}{4} \nu}{C'c^2+\frac{1}{2} \nu}$,
from \eqref{0ue24}$^{\delta_1}\times$ \eqref{1ue27}$^{1-\delta_1}$ and \eqref{1ue25},
we get (\ref{ue15}) for $r=m'=0$.

To get (\ref{ue15}) for $r \geq 1$, note that from (\ref{1ue15}), for $k\geq 1$
 \begin{align}\label{ue28}
\frac{\partial^{r}}{\partial t^{r}}
e ^{-u\cL^t_2} =& \frac{(-1)^{k-1} (k-1)!}{2\pi i u^{k-1}}
\int_{\delta\cup \Delta}  e ^{-u\lambda}
\frac{\partial^{r}}{\partial t^{r}}(\lambda-\cL^t_2)^{-k} d \lambda.
\end{align}

We have the similar equation for
$\frac{\partial^{r}}{\partial t^{r}}(\cL^t_2e ^{-u\cL^t_2})$.

Set \be\label{ue29} I_{k,r} = \Big \{ (\bk,\br)=(k_i,r_i)|
\sum_{i=0}^j k_i =k+j, \sum_{i=1}^j r_i =r,\, \,   k_i, r_i \in
\bN^*\Big \}. \ee Then  there exist $a ^{\bk}_{\br} \in \bR$ such
that
\begin{align}\label{ue30}
\begin{split}
& A^{\bk}_{\br}  (\lambda,t) = (\lambda-\cL^t_2)^{-k_0}
\frac{\partial^{r_1}\cL^t_2}{\partial t^{r_1}}  (\lambda-\cL^t_2)^{-k_1}
\cdots\frac{\partial^{r_j}\cL^t_2}{\partial t^{r_j}}
(\lambda-\cL^t_2)^{-k_j},\\
& \frac{\partial^{r}}{\partial t^{r}}
(\lambda-\cL^t_2)^{-k}=
\sum_{(\bk,\br)\in I_{k,r} } a ^{\bk}_{\br}  A ^{\bk}_{\br}  (\lambda,t).
\end{split}
\end{align}

We claim that $A ^{\bk}_{\br}(\lambda,t)$ is well defined and
for any $m\in \bN$, $k>2(m+r+1)$, $Q,Q'\in \cup_{l=1}^m \mD^l_t$,
there exist $C>0$, $N\in \bN$ such that for $\lambda\in \delta\cup  \Delta$,
\begin{align}\label{1ue30}
\|Q  A ^{\bk}_{\br}(\lambda,t)Q' s\|_{t,0}
\leq C (1+|\lambda|)^N \sum_{|\beta|\leq 2r} \|Z^\beta s\|_{t,0}.
\end{align}

In fact, by (\ref{1ue1}), $\frac{\partial^{r}}{\partial
t^{r}}\cL^t_2$ is a combination of $\frac{\partial^{r_1}}{\partial
t^{r_1}}(g^{ij}(tZ))$, $(\frac{\partial^{r_2}}{\partial
t^{r_2}}\nabla_{t,e_i})$,
 $\frac{\partial^{r_3}}{\partial t^{r_3}}(q(tZ))$,
$\frac{\partial^{r_4}}{\partial t^{r_4}} (t \langle
\wi{\mu}^{E_{0,p}}, f_{0,l} (tZ)\rangle)$, where $q$ runs over the
functions $r^X$, etc.,  appearing in (\ref{1ue1}).

Now $\frac{\partial^{r_1}}{\partial t^{r_1}}(q(tZ))$ (resp.
$\frac{\partial^{r_1}}{\partial t^{r_1}}$
 $(t \langle \wi{\mu}^{E_{0,p}}, f_{0,l} \rangle (tZ))$,
$\frac{\partial^{r_1}}{\partial t^{r_1}}\nabla_{t,e_i}$) ($r_1\geq 1$)
 are functions of  the type as
 $q'(tZ)Z^\beta$, $|\beta|\leq r_1$ (resp. $r_1+1$) (where $q'$, as $q$, runs over the
functions $r^X$, etc.,  appearing in (\ref{1ue1})), with $q'(Z)$
and its derivatives on $Z$ being  bounded smooth functions on $Z$.

Let $\cR'_t$ be the family of operators of the type $$\cR'_{t} =
\{ [f_{j_1} Q_{j_1}, [f_{j_2} Q_{j_2},\ldots [f_{j_l}
Q_{j_l},\cL^t_2]]\ldots ] \}$$ with $f_{j_i}$ smooth bounded (with
its derivatives) functions and  $Q_{j_i}\in \mD_t\cup \{
Z_j\}_{j=1}^{2n-n_0}$.

Now for the operator $ A ^{\bk}_{\br}(\lambda,t)Q'$,
we will move first all the term $Z^\beta$ in $d'(tZ)Z^\beta$
as above to the right hand side of this operator, to do so,
we always use the commutator trick, i.e., each time,
we consider only the commutation for $Z_i$,
 not for $Z^\beta$ with $|\beta|>1$.

Then $A ^{\bk}_{\br}(\lambda,t)Q'$ is as the form
$\sum_{|\beta|\leq 2r} L^t_\beta Q''_\beta  Z^\beta$,
and $Q''_\beta$ is obtained from
 $Q'$ and its commutation with $Z^\beta$.

Now we move all the terms $\nabla_{t,e_i}$,
$\langle \frac{1}{t}\wi{\mu},  f_{0,l} \rangle (tZ)$ in
$\frac{\partial^{r_j}\cL^t_2}{\partial t^{r_j}}$
to the right hand side of the operator $L^t_\beta$.

Then as in the proof of Theorem \ref{tu6},
we get finally  that $Q  A ^{\bk}_{\br}(\lambda,t)Q'$
is as the form  $\sum_{\beta} \cL^t_\beta Z^\beta$ where $\cL^t_\beta$
is a linear combination of operators of the form
\be\label{2ue30}
Q (\lambda-\cL^t_2)^{-k'_0}R_1(\lambda-\cL^t_2)^{-k'_1}R_2
\cdots R_{l'}(\lambda-\cL^t_2)^{-k'_{l'}} Q''' Q'',
\ee
with  $R_1, \cdots, R_{l'} \in \cR'_{t}$,  $Q'''\in \cup_{l=1}^{2r} \mD^l_t$,
 $Q''\in  \cup_{l=1}^{m} \mD^l_t$, $|\beta|\leq 2 r$,
 and $Q''$ is obtained from $Q'$ and its commutation with $Z^\beta$.

By the argument as in \eqref{ue18} and  \eqref{ue20}, as $k>2(m+r+1)$,
we can split the above operator to two parts
\begin{align*}
&Q (\lambda-\cL^t_2)^{-k'_0}R_1(\lambda-\cL^t_2)^{-k'_1}R_2
\cdots R_{i}(\lambda-\cL^t_2)^{-k''_{i}};\\
&(\lambda-\cL^t_2)^{-(k'_{i}-k''_{i}) }\cdots
R_{l'}(\lambda-\cL^t_2)^{-k'_{l'}} Q''' Q'',
\end{align*}
and the $\|\quad \|^{0,0}_t$-norm of each part is bounded
by $C(1+|\lambda|^2)^N$.

 Thus
the proof of  (\ref{1ue30}) is complete.

By (\ref{ue28}), (\ref{ue30}) and  (\ref{1ue30}), we get the similar
estimate \eqref{0ue24},
 (\ref{1ue27}) for
$\frac{\partial^{r}}{\partial t^{r}}e ^{-u\cL^t_2}$,
$\frac{\partial^{r}}{\partial t^{r}}
(\cL^t_2e ^{-u\cL^t_2})$ with the exponential $2n+m+2r+2$
instead of  $2n+m+2$ therein.

Thus  we get (\ref{ue15}) for $m'=0$.

Finally, for $U\in TX_G$ a vector on $X_G$,
\begin{align}\label{0ue30}
& \nabla ^{\pi ^* \End (\bE_B)}_U e ^{-u\cL^t_2}
= \frac{(-1)^{k-1} (k-1)!}{2\pi i u^{k-1}}
\int_{\delta\cup \Delta}  e ^{-u\lambda}
\nabla ^{\pi ^* \End (\bE_B)}_U (\lambda-\cL^t_2)^{-k} d \lambda .
\end{align}

Now, by using the similar formula (\ref{ue30})  for
$\nabla ^{\pi ^* \End (\bE_B)} _U(\lambda-\cL^t_2)^{-k}$ by replacing
$\frac{\partial^{r_1}\cL^t_2}{\partial t^{r_1}}$ by
$\nabla ^{\pi ^* \End (\bE_B)}_U \cL^t_2$, and remark that
$\nabla ^{\pi ^* \End (\bE_B)}_U \cL^t_2$ is a differential operator
on $T_{x_0} B$ with the same structure  as $\cL^t_2$.

Then by the above argument, we get (\ref{ue15}) for $m'\geq 1$.
\end{proof}

Let $P_{0,t}$ be the orthogonal projection from
$\cC^\infty (T_{x_0}B, \bE_{B,x_0})$ to the kernel of $\cL^t_2$ with respect to
$\left\langle \, ,\,\right\rangle_{t,0}$. Set
\begin{align}\label{1ue31}
F_u(\cL^t_2) = \frac{1}{2\pi i} \int_{\Delta}
 e ^{-u\lambda}(\lambda-\cL^t_2)^{-1} d \lambda .
\end{align}
 By (\ref{u9}),
\begin{align}\label{1ue32}
F_u(\cL^t_2)=e ^{-u\cL^t_2}-P_{0,t}
= \int_u^{+\infty} \cL^t_2e ^{-u_1\cL^t_2}du_1.
\end{align}

Let $P_{0,t}(Z,Z^{\prime})$, $F_u(\cL^t_2)(Z,Z^{\prime})$
be the smooth kernels of
 $P_{0,t},F_u(\cL^t_2)$ with respect to $dv_{T_{x_0}B}(Z^{\prime})$.

\begin{cor}\label{0tue8} With the notation in Theorem \ref{tue8},
\begin{multline}\label{1ue33}
\sup_{|\alpha|+|\alpha'|\leq m}
(1+|Z^{\bot}|+|Z^{\prime\bot}|)^{m''}\Big |\frac{\partial^{|\alpha|+|\alpha'|}}
{\partial Z^{\alpha} {\partial Z^{\prime}}^{\alpha'}}
\frac{\partial^{r}}{\partial t^{r}} F_u(\cL^t_2)
\left (Z, Z^{\prime}\right )\Big |_{\cC^{m'}(P)}  \\
   \leq C (1+|Z^0|+|Z^{\prime 0}|)^{2n+m+2m^\prime+2r+2}
\exp ( -\frac{1}{8}\nu u- \sqrt{C''\nu} |Z-Z^{\prime}|).
\end{multline}
\end{cor}
\begin{proof} Note that $\frac{1}{8}\nu u
+\frac{2C''}{u} |Z-Z^{\prime}|^2\geq \sqrt{C''\nu} |Z-Z^{\prime}|$, thus
\begin{multline}\label{1ue34}
\int_u^{+\infty}
e ^{ -\frac{1}{4}\nu u_1- \frac{2C''}{u_1} |Z-Z^{\prime}|^2}du_1
\leq e ^{-\sqrt{C''\nu}|Z-Z^{\prime}|}
\int_u^{+\infty} e ^{-\frac{1}{8}\nu u_1}du_1 \\
 = \frac{8}{\nu}e ^{ -\frac{1}{8}\nu u- \sqrt{C''\nu} |Z-Z^{\prime}|}.
\end{multline}
By  (\ref{ue15}), (\ref{1ue32}) and (\ref{1ue34}), we get (\ref{1ue33}).
\end{proof}

For $k$ large enough, set
\begin{align}\label{ue31}
\begin{split}
& F_{r,u}
= \frac{(-1)^{k-1} (k-1)!}{2\pi i \, r! u^{k-1}}\int_{\Delta}
e ^{-u\lambda}   \sum_{(\bk,\br)\in I_{k,r} }
 a ^{\bk}_{\br}  A ^{\bk}_{\br}  (\lambda,0)d \lambda ,\\
&J_{r,u}= \frac{(-1)^{k-1} (k-1)!}{2\pi i \, r! u^{k-1}}
\int_{\delta\cup \Delta}  e ^{-u\lambda}   \sum_{(\bk,\br)\in I_{k,r} }
 a ^{\bk}_{\br}  A ^{\bk}_{\br}  (\lambda,0)d \lambda ,\\
&F_{r,u,t} = \frac{1}{r!}\frac{\partial^{r}}{\partial t^{r}}
F_u (\cL^t_2)- F_{r,u},
\quad J_{r,u,t} = \frac{1}{r!}\frac{\partial^{r}}{\partial t^{r}}
e ^{-u\cL^t_2}- J_{r,u}.
\end{split}
\end{align}

Certainly, as $t\to 0$, the limit of  $\|\quad\|_{t,m}$ exists,
and we denote it by  $\|\quad\|_{0,m}$.
\begin{thm} \label{tue9} For any $r\geq 0$, $k>0$, there exist $C>0$, $N\in \bN$ such that for
$t \in [0,t_0], \lambda \in \delta\cup \Delta$,
\begin{align}\label{ue32}
&  \left \|\Big (\frac{\partial^{r}\cL^t_2}{\partial t^{r}} -
\frac{\partial^{r}\cL^t_2}{\partial t^{r}}|_{t=0} \Big ) s \right \|_{t,-1}
\leq Ct \sum_{|\alpha|\leq r+3} \|Z^\alpha s\|_{0,1},\\
& \Big \|\Big (\frac{\partial^{r}}{\partial t^{r}} (\lambda-\cL^t_2)^{-k}
-\sum_{(\bk,\br)\in I_{k,r} }
 a ^{\bk}_{\br}  A ^{\bk}_{\br} (\lambda,0)\Big )s\Big \|_{0,0}
\leq C t  (1+|\lambda|^2)^N
\sum_{|\alpha|\leq 4r+3} \|Z^\alpha s\|_{0,0}.\nonumber
\end{align}
\end{thm}
\begin{proof}  Note that by (\ref{c37}), (\ref{u4}),
for $t\in [0,1]$, $k\geq 1$
\begin{align}\label{1ue35}
\|s\|_{t,0}\leq C\|s \|_{0,0},\quad \|s\|_{t,k}
\leq C \sum_{|\alpha|\leq k} \|Z^\alpha s\|_{0,k}.
\end{align}
 An application of Taylor expansion for
(\ref{1ue1}) leads to the following equation, if $s,s'$ have
compact support, \be\label{ue33} \Big | \left\langle  \Big
(\frac{\partial^{r}\cL^t_2}{\partial t^{r}} -
\frac{\partial^{r}\cL^t_2}{\partial t^{r}} |_{t=0} \Big )s,s'
\right\rangle_{0,0}\Big | \leq C t \|s'\|_{t,1}\sum_{|\alpha|\leq
r+3} \|Z^\alpha s\|_{0,1}. \ee Thus we get the first inequality of
(\ref{ue32}).

Note that
\begin{align}\label{ue34}
& (\lambda-\cL^t_2)^{-1}- (\lambda-\cL^0_2)^{-1}
=(\lambda-\cL^t_2)^{-1}(\cL^t_2-\cL^0_2) (\lambda-\cL^0_2)^{-1}.
\end{align}
Now from  (\ref{ue2}),  (\ref{ue33}) and (\ref{ue34}),
\begin{align}\label{0ue34}
& \left \|\left ((\lambda-\cL^t_2)^{-1}- (\lambda-\cL^0_2)^{-1} \right)s
 \right \|_{0,0}
\leq Ct  (1+|\lambda |^4)   \sum_{|\alpha|\leq 3} \|Z^\alpha s\|_{0,0}.
\end{align}
After taking the limit, we know that Theorems \ref{tu4}-\ref{tu6} still hold
for $t=0$.

 Note that
$\nabla_{0,e_j} = \nabla_{e_j}+\frac{1}{2} R^{L_B}_{x_0}(\mR, e_j)$
by (\ref{c37}).

If we denote by $\cL_{\lambda,t}=\lambda -\cL_2^{t}$, then
\begin{multline}\label{0ue35}
 A^{\bk}_{\br}  (\lambda,t)-  A^{\bk}_{\br}  (\lambda,0)
 = \sum_{i=1}^j \cL_{\lambda,t}^{-k_0}\cdots
\left ( \frac{\partial^{r_i}\cL_2^t}{\partial t^{r_i}}
- \frac{\partial^{r_i}\cL_2^t}{\partial t^{r_i}} |_{t=0}\right )
\cL_{\lambda,0}^{-k_{i}}\cdots \cL_{\lambda,0}^{-k_j}\\
+ \sum_{i=0}^j \cL_{\lambda,t}^{-k_0}\cdots
\left (\cL_{\lambda,t}^{-k_i}- \cL_{\lambda,0}^{-k_i} \right)
\left (\frac{\partial^{r_{i+1}}\cL_2^t}{\partial t^{r_{i+1}}}|_{t=0}\right)
\cdots \cL_{\lambda,0}^{-k_j}.
\end{multline}

Now from the first inequality  of (\ref{ue32}),
(\ref{ue2}),  (\ref{ue30}), (\ref{0ue34}) and (\ref{0ue35}),
we get (\ref{ue32}).
\end{proof}

\begin{thm} \label{tue12}
There  exist $C>0$, $N\in \bN$
such that for $t \in ]0,t_0]$, $u\geq u_0$,
$q\in \bN$, $Z,Z^{\prime}\in T_{x_0}B$, $|Z|, |Z^{\prime}|\leq q$,
\begin{align}\label{ue42}
\Big |F_{r,u,t}(Z,Z^{\prime})\Big |
\leq & C t^{1\over 2(2n-n_0+1)} (1+q)^N   e ^{-\frac{1}{8}\nu u},\\
\Big |J_{r,u,t}(Z,Z^{\prime})\Big | \leq & C t^{1\over
2(2n-n_0+1)} (1+q)^N e ^{\frac{1}{2}\nu u}.\nonumber
\end{align}
\end{thm}
\begin{proof}
Let $J^0_{x_0, q}$ be the vector space of square integrable sections of
 $\bE_{B,x_0}$ over $\{ Z\in T_{x_0}B, |Z|\leq q+1\}$.

 If $s\in J^0_{x_0, q}$,
 put $\|s\|^2_{(q)} = \int_{|Z|\leq q+1} |s|^2_{\bE_{B,x_0}} dv_{TB}(Z)$.
Let $\|A\|_{(q)}$ be the  operator norm of $A\in \cL(J^0_{x_0, q})$
with respect to $\|\quad\|_{(q)}$.

By (\ref{ue28}),
 (\ref{ue31}) and (\ref{ue32}), we get:
 there exist  $C>0, N\in \bN$
such that for $t \in ]0,t_0],  u\geq u_0$,
\begin{align}\label{ue39}
& \|F_{r,u,t}\|_{(q)}
\leq C t (1+q)^N   e ^{-\frac{1}{2}\nu u},\\
& \|J_{r,u,t}\|_{(q)}
\leq C t (1+q)^N \ e ^{\frac{1}{4}\nu u}.\nonumber
\end{align}

 Let $\phi : \bR\to [0,1]$ be a smooth function with compact support,
 equal $1$ near $0$, such that  $\int_{T_{x_0}B} \phi
(Z) dv_{T_{x_0}B}(Z)=1$.

Take $\varsigma \in ]0,1]$.

By the proof of Theorem \ref{tue8}, $F_{r,u}$ verifies the similar
inequality as in (\ref{1ue33}).
Thus by (\ref{1ue33}), there exists $C>0$ such that if
$|Z|,|Z^{\prime}|\leq q$, $U,U'\in \bE_{B,x_0}$,
\begin{multline}\label{ue43}
\Big |  \left\langle
F_{r,u,t} (Z,Z^{\prime}) U,U'     \right\rangle
-\int_{T_{x_0}B \times T_{x_0}B}  \left\langle
F_{r,u,t}(Z-W,Z^{\prime}-W') U,U'     \right\rangle\\
  \frac{1}{\varsigma ^{4n-2n_0}} \phi (W/\varsigma) \phi (W'/\varsigma)
dv_{T_{x_0}B}(W)dv_{T_{x_0}B}(W')\Big |
\leq C \varsigma (1+q)^N e ^{ -\frac{1}{8}\nu u} |U||U'|.
\end{multline}

On the other hand, by (\ref{ue39}),
\begin{multline}\label{ue44}
\Big |\int_{T_{x_0}B \times T_{x_0}B}  \left\langle
F_{r,u,t}  (Z-W,Z^{\prime}-W') U,U'     \right\rangle
  \frac{1}{\varsigma ^{4n-2n_0}} \phi (W/\varsigma) \phi (W'/\varsigma)\\
  dv_{T_{x_0}B}(W)dv_{T_{x_0}B}(W')\Big | \leq C
t\frac{1}{\varsigma ^{2n-n_0}}  (1+q)^N e ^{-\frac{1}{2}\nu u} |U||U'|.
\end{multline}

By taking $\varsigma = t^{1/2(2n-n_0+1)}$, we get (\ref{ue42}).

In the same way, we get  (\ref{ue42}) for $J_{r,u,t}$.
\end{proof}

\begin{thm} \label{tue14} There exists $C''>0$ such that for
any $k,m,m',m''\in \bN$, there exist $N\in \bN$, $C>0$  such that if
$t\in ]0,t_0], u\geq u_0$, $Z,Z^{\prime}\in T_{x_0}^HU$,
$\alpha,\alpha \in \bZ^{2n-n_0}$, $|\alpha|+|\alpha'|\leq m$,
 \begin{align}\label{0ue45}
\begin{split}
&(1+|Z^\bot|+|Z^{\prime\bot}|)^{m''}\Big |\frac{\partial^{|\alpha|+|\alpha'|}}
{\partial Z^{\alpha} {\partial Z^{\prime}}^{\alpha'}}\Big (F_u (\cL^t_2)
- \sum_{r=0}^k F_{r,u}t^r\Big ) (Z,Z^{\prime})
\Big |_{\cC^{m'}(X_G)} \\
 &\hspace*{10mm}  \leq C  t^{k+1} 
(1+|Z^0|+|Z^{\prime 0}|)^{2(n+k+m^\prime+2)+m}
\exp ( -\frac{1}{8} \nu  u- \sqrt{C''\nu} |Z-Z^{\prime}|),\\
&
(1+|Z^\bot|+|Z^{\prime\bot}|)^{m''}\Big |\frac{\partial^{|\alpha|+|\alpha'|}}
{\partial Z^{\alpha} {\partial Z^{\prime}}^{\alpha'}}\Big (e ^{-u\cL^t_2}
- \sum_{r=0}^k J_{r,u} t^r\Big ) (Z,Z^{\prime})
\Big |_{\cC^{m'}(X_G)} \\
&\hspace*{20mm}    \leq C  t^{k+1} 
 (1+|Z^0|+|Z^{\prime 0}|)^{2(n+k+m^\prime+2)+m}
\exp (\frac{1}{2} \nu u- \frac{2C''}{u} |Z-Z^{\prime}|^2).
\end{split}\end{align}
\end{thm}
\begin{proof} By (\ref{ue31}) and (\ref{ue42}),
 \begin{align}\label{0ue47}
&\frac{1}{r!}  \frac{\partial^{r}}{\partial t^{r}}
F_u (\cL^t_2) |_{t=0} = F_{r,u},\quad
 \frac{1}{r!}\frac{\partial^{r}}{\partial t^{r}} e ^{-u\cL^t_2}|_{t=0}
= J_{r,u}.
\end{align}

Now by Theorem \ref{tue8} and (\ref{ue31}), $ J_{r,u}, F_{r,u}$
have the same estimates as $\frac{\partial^{r}}{\partial t^{r}} e
^{-u\cL^t_2}$, $\frac{\partial^{r}}{\partial t^{r}}F_u (\cL^t_2)$
in (\ref{ue15}), (\ref{1ue33}).

 Again from (\ref{ue15}), (\ref{1ue33}), (\ref{ue31}), (\ref{ue42}),
and the Taylor expansion
 \begin{align}\label{0ue48}
G(t)- \sum_{r=0}^k \frac{1}{r !} \frac{\partial ^r G}{\partial t^r}(0) t^r
= \frac{1}{k!}\int_0^t (t-t_0)^k \frac{\partial ^{k+1} G}{\partial
t^{k+1} }(t_0) dt_0,
\end{align} we get (\ref{0ue45}).
\end{proof}

\subsection{Evaluation of $J_{r,u}$} \label{s3.5}

For $u >0$, we will write $u\Delta_j$ for the rescaled simplex
 $\{ (u_1, \cdots, u_j)|$ $0\leq u_1\leq u_2\leq$ $\cdots\leq u_j\leq u\}$.

 Let $e ^{-u\cL^0_2} (Z,Z^{\prime})$ be the smooth kernel of $e ^{-u\cL^0_2}$
with respect to $dv_{T_{x_0}B}(Z^{\prime})$.

Recall that the $\mO _r$'s have been defined in (\ref{0c35}).

\begin{thm}  \label{tue15}  For $r\geq 0$, we have
\begin{multline}\label{ue49}
J_{r,u}=\sum_{\sum_{i=1}^j r_i =r,\, r_i\geq 1} (-1)^j
\int_{u\Delta_j} e ^ {-(u-u_j)\cL^0_2} \mO _{r_j}
e ^ {-(u_j-u_{j-1})\cL^0_2}\\
\cdots   \mO _{r_1} e ^ {-u_1\cL^0_2} du_1\cdots du_j,
\end{multline}
where the product in the integrand is the convolution product. Moreover,
\be\label{ue52}
J_{r,u}(Z,Z^{\prime}) = (-1)^r J_{r,u}(-Z,-Z^{\prime}).
\ee
\end{thm}
\begin{proof}
We introduce an even extra-variable $\sigma$ such that $\sigma ^{r+1}=0$.

Set $[\quad]^{[r]}$  the coefficient of $\sigma ^r$,
 $\cL_\sigma= \cL^0_2 + \sum_{j=1}^r  \mO _j \sigma ^j$.

From (\ref{ue31}), (\ref{0ue47}), we  know
\begin{align}\label{0ue51}
 J_{r,u}(Z,Z^{\prime})= \frac{1}{r!} \frac{\partial ^r }{\partial t ^r}
 e ^{-u\cL^t_2}(Z,Z^{\prime}) |_{t=0}
=  [ e ^{-u\cL_\sigma}]^{[r]}(Z,Z^{\prime}).
\end{align}

Now from (\ref{0ue51}) and the Volterra expansion of
$e ^{-u\cL_\sigma}$ (cf. \cite[\S 2.4]{BeGeVe}),
we get (\ref{ue49}).

We  prove (\ref{ue52}) by iteration.

From (\ref{0c35})
\begin{multline}\label{ue54}
\cL^0_2 = -\sum_{j=1}^{2n-n_{0}} (\nabla_{e_{j}})^2
-\pi ^2 \langle ((P^{T^HU}\bJ P^{T^HU})^2
+4 P^{T^HU}\bJ P^{TY}\bJ P^{T^HU})_{x_0} \mR,\mR\rangle \\
 + 2\pi \sqrt{-1}\nabla_{P^{T^HU}\bJ P^{T^HU} \mR}
-2 \omega_{d,x_0} -\tau_{x_0}.
\end{multline}
Here the matrix $((P^{T^HU}\bJ P^{T^HU})^2 +4 P^{T^HU}\bJ
P^{TY}\bJ P^{T^HU})_{x_0}$ need not   commute with $P^{T^HU}\bJ
P^{T^HU}$.

Thus
 \cite[(6.37), (6.38)]{B90d} does not apply directly here, and we could not get
a precise formula for $ e ^{-u\cL^0_2}$ as in  \cite[(4.106)]{DLM04a}.

By the uniqueness of the solution of heat equations and
(\ref{ue54}), we know
\begin{align}\label{0ue52}
e ^{-u\cL^0_2} (Z,Z^{\prime})= e ^{-u\cL^0_2} (-Z,-Z^{\prime}).
\end{align}

By (\ref{ue49}),
\begin{align}\label{0ue53}
J_{0,u}(Z,Z^{\prime})=e ^{-u\cL^0_2} (Z,Z^{\prime}).
\end{align}

Thus we get (\ref{ue52}) for $r=0$.

If (\ref{ue52}) holds for $r\leq k$,
then by  (\ref{ue49}),  (\ref{0ue53}),
\begin{align}\label{0ue54}
J_{k+1,u}= -\sum_{j=1}^{k+1}
\int_0^u  e ^{-(u-u_1)\cL^0_2} \mO_j J_{k+1-j,u_1}du_1.
\end{align}

By the iteration, Theorem \ref{t3.3} and   (\ref{0ue53}),
and note that  $\nabla_{e_i}$ in $\mO_j$ will change
the parity of the polynomials we obtained,
we get (\ref{ue52}) for $r=k+1$.
\end{proof}

\subsection{Proof of Theorem \ref{t0.1}}
\label{s3.6}

By (\ref{1ue32}) and (\ref{0ue45}),  for any $u>0$ fixed,
 there exists $C_u>0$ such that for
$t=\frac{1}{\sqrt{p}}$, $Z,Z^{\prime}\in T_{x_0}B$, $x_0\in P$,
$\alpha,\alpha'\in \bZ^{2n-n_0}$, $|\alpha|+|\alpha'|\leq m$, we have
\begin{multline}\label{1ue54}
 (1+|Z^\bot|+|Z^{\prime\bot}|)^{m''}
\left |\frac{\partial^{|\alpha|+|\alpha'|}}
{\partial Z^{\alpha} {\partial Z^{\prime}}^{\alpha'}} \Big (P_{0,t}
- \sum_{r=0}^k t^r(J_{r,u}- F_{r,u})\Big)(Z,Z^{\prime})\right |_{\cC ^{m'}(X_G)}\\
\leq C_u  t^{k+1}  (1+|Z^0|+|Z^{\prime 0}|)^{2(n+k+m^\prime+2)+m}
\exp (- \sqrt{C''\nu } |Z-Z^{\prime}|).
\end{multline}

Set
\begin{align}\label{1ue53}
P^{(r)}=J_{r,u}- F_{r,u}.
\end{align}

Then  $P^{(r)}$ does not depend on $u>0$ by (\ref{1ue54}), as
$P_{0,t}$ does not depend on $u$.

Moreover, by taking the limit of (\ref{1ue33}) as $t\to 0$,
\begin{multline}\label{1ue56}
 (1+|Z^\bot|+|Z^{\prime\bot}|)^{m''}\Big | F_{r,u}(Z,Z^{\prime})\Big |_{\cC^{m'}(X_G)}\\
\leq C (1+|Z^0|+|Z^{\prime 0}|)^{2n+2r+2m^\prime+2} \exp ( -\frac{1}{8} \nu  u- \sqrt{C''\nu} |Z-Z^{\prime}|).
\end{multline}
Thus
\begin{align}\label{1ue57}
J_{r,u}(Z,Z^{\prime}) = P^{(r)}(Z,Z^{\prime})+F_{r,u}(Z,Z^{\prime})
= P^{(r)}(Z,Z^{\prime})
+\mathcal{O}( e^ {-\frac{1}{8} \nu  u}),
\end{align}
uniformly on any compact set of $T_{x_0}B\times T_{x_0}B$.

Especially, from \eqref{ue52}, \eqref{1ue57}, we get
\begin{align}\label{ue64}
& P^{(r)}(Z,Z^{\prime})= (-1)^r P^{(r)}(-Z,-Z^{\prime}).
\end{align}

By (\ref{c27}), for $Z,Z^{\prime}\in T_{x_0}B$,
\begin{align}\label{1ue60}
&P_{x_0,p}(Z,Z^{\prime})= p^{n-\frac{n_0}{2}}
\kappa ^{-\frac{1}{2}}(Z) P_{0,t}(Z/t,Z^{\prime}/t) \kappa ^{-\frac{1}{2}}(Z^{\prime}).
\end{align}

We note in passing that, as a consequence of (\ref{1ue54})
and (\ref{1ue60}), we obtain the following estimate.
\begin{thm} \label{tue17}
For any $k,m,m', m''\in \bN$, there exists $C>0$ such that for
$Z,Z^{\prime}\in T_{x_0}B$, $|Z|, |Z^{\prime}|\leq  \var$, $x_0\in X_G$,
\begin{multline}\label{ue66}
\sup_{|\alpha|+|\alpha'|\leq m}
(1+\sqrt{p}|Z^{\bot}|+\sqrt{p}|Z^{\prime\bot}|)^{m''}
\left |\frac{\partial^{|\alpha|+|\alpha'|}}
{\partial Z^{\alpha} {\partial Z^{\prime}}^{\alpha'}}\right. \\
\left. \left (p^{-n+\frac{n_0}{2}}  P_{x_0,p}(Z,Z^{\prime})
-\sum_{r=0}^k  P^{(r)} (\sqrt{p} Z,\sqrt{p} Z^{\prime})
\kappa ^{-\frac{1}{2}}(Z)\kappa ^{-\frac{1}{2}}(Z^{\prime})
p^{-r/2}\right )\right |_{\cC ^{m'}(X_G)}\\
\leq C  p^{-(k+1-m)/2}
(1+\sqrt{p} |Z^{0}|+\sqrt{p} | Z^{\prime 0}|)^{2(n+k+m^\prime+2)+m}
\exp (- \sqrt{C''\nu } \sqrt{p} |Z-Z^{\prime}|).
\end{multline}
\end{thm}

\comment{
By \eqref{ue66} and the Taylor expansion of $\kappa (Z)$,
with the notation in \eqref{ue66},
\begin{multline}\label{ue67}
\sup_{|\alpha|+|\alpha'|\leq m}
(1+\sqrt{p}|Z^{\bot}|+\sqrt{p}|Z^{\prime\bot}|)^{m''}
\left |\frac{\partial^{|\alpha|+|\alpha'|}}
{\partial Z^{\alpha} {\partial Z^{\prime}}^{\alpha'}}\right. \\
\left. \left (p^{-n+\frac{n_0}{2}}  P_{x_0,p}(Z,Z^{\prime})
-\sum_{r=0}^k  P^{(r)} (\sqrt{p} Z,\sqrt{p} Z^{\prime})
p^{-r/2}\right )\right |_{\cC ^{m'}(X_G)}\\
\leq C  p^{-(k+1-m)/2}
(1+\sqrt{p} |Z^{0}|+\sqrt{p} | Z^{\prime 0}|)^{2(n+k+1)+m}
\exp (- \sqrt{C''\nu } \sqrt{p} |Z-Z^{\prime}|).
\end{multline}
}
 From (\ref{1c19}), (\ref{1c20}), (\ref{ac38}) and (\ref{ue66}),  we get
Theorem \ref{t0.1} without knowing the properties \eqref{0.7},
\eqref{a0.7} for $ P^{(r)}$.

To prove the uniformity part of Theorem \ref{t0.1}, we notice that
in the proof of Theorem \ref{tue8},
we only use the derivatives of the coefficients of
$\cL^t_2$ with order $\leqslant 2n+m+m'+r+2$.
Thus the constants in Theorems \ref{tue8} and \ref{tue12},
(resp. Theorem  \ref{tue14})
are uniformly bounded, if with respect to a fixed metric $g^{TX}_0$,
the $\cC^{2n+m+m'+r+2}$ (resp. $\cC^{2n+m+m'+k+3}$)\,--\,norms on $X$ of the data
{\rm(}$g^{TX}$, $h^L$, $\nabla ^L$, $h^E$, $\nabla ^E$, $J${\rm)}
 are bounded, and $g^{TX}$ is bounded below.

Moreover, taking derivatives with respect to the parameters we obtain
a similar equation as \eqref{0ue30}, where
$x_0\in X_G$ plays now a role of a parameter. Thus
the $\cC^{m'}$--\,norm in \eqref{ue66} can also include the parameters
if the $\cC^{m'}$--\,norms (with respect to the parameter $x_0\in X_G$)
 of the derivatives of above data with order
$\leqslant 2n+k+m+3$ are bounded.

Thus we can take $C_{k,\, l}$ in \eqref{0.6} independent of $g^{TX}$
under our condition.

This achieves the proof of Theorem \ref{t0.1}
except \eqref{0.7} and \eqref{a0.7}
which will be proved in Theorem \ref{0t3.6}
under the condition in Theorem \ref{t0.1}.


\comment{
Moreover, from Theorems \ref{t3.3}, \ref{tue15}, (\ref{c30}),
we get the property on $b_{r}$ in  Theorem \ref{t0.1}.
To get the last part of Theorem \ref{t0.1}, we notice that the
constants in Theorems \ref{tue8} and \ref{tue12} will be uniform bounded under
our condition, thus we can take $C_{k,\, l}$
in (\ref{0.8}) independent of $g^{TX}$.
Thus we have proved Theorem \ref{t0.1}.
}
\newpage

\section{Evaluation of $P^{(r)}$}\label{s30}

In this Section, inspired by the method  in \cite[\S 1.4,
1.5]{MM04a}, we develop  a direct and effective method to compute
$P^{(r)}$. In particular, we get \eqref{0.7} and \eqref{a0.7}
under the condition in Theorem \ref{t0.1}.

This section is organized as follows.
In Section \ref{s3.10}, we study the spectrum of the limited operator
$\cL^0_2$. In  Section \ref{s3.11},
we get a direct method to evaluate $P^{(r)}$ in \eqref{0.7},
especially, we prove \eqref{0.7} and \eqref{a0.7}.
 In  Section \ref{s3.111}, we compute explicitly $\mO_1$ in \eqref{c30},
and get a general formula for $P^{(2)}$ by using the operators
$\mO_1$, $\mO_2$.
 In  Section \ref{s3.12}, we compute  explicitly
 an interesting example:
the line bundle $\mO(2)$ on $(\bC P^1, 2\, \om_{FS})$.
We verify that Theorem \ref{t0.1} coincides with
our computation here
if $0$ is a regular value of the moment map $\mu$,
but it does not hold if $0$ is a singular value.

We use the notations in  Section \ref{s3.2},
and we suppose that \eqref{g2} is verified.

\subsection{Spectrum of $\cL^0_2$}\label{s3.10}
Recall that $T^HP$ is the orthogonal complement of $TY$ in $(TP,
g^{TP})$. Note that by (\ref{0a0}) and (\ref{a5}), we have the
following orthogonal splitting of vector bundles on
$P=\mu^{-1}(0)$,
\begin{align}  \label{g1}
&TX= T^HP\oplus TY\oplus \bJ TY,& TP=  T^HP\oplus TY.
\end{align}
In the rest of this Section, we suppose that on $P$
\begin{align}  \label{g2}
&\bJ T^HP = T^HP, &\bJ^2 TY =TY.
\end{align}
\eqref{0a00} and \eqref{g2} imply that $-J\bJ$ preserves $TY$ and $\bJ TY$.

As $g^{TX}$ is $J$-invariant, we get
\begin{align}  \label{0g2}
\bJ\, TY =J\,  TY, \quad J\,  T^HP = T^HP.
\end{align}
Thus $(\bJ TY)_B|_{X_G}$ is the orthogonal complement of $TX_G$ in $TB$,
and $\bJ$ induces naturally $ \bJ_G\in \End(TX_G)$.
We will identify $(\bJ TY)_B|_{X_G}$ to the normal bundle of $X_G$ in $B$.

For $U,V\in T_{x_0}B$, $x_0\in X_G$, we have
\begin{align}  \label{g3}
 \om(U^H,V^H)= \om_G(P^{TX_G}U, P^{TX_G}V).
\end{align}

From the above discussion,  for $x_0\in X_G$,  we can choose
$\{w^0_j\}_{j=1}^{n-n_0}$, $\{e ^\bot_j\}_{j=1}^{n_0}$ orthonormal
basis of $T^{(1,0)}_{x_0}X_G$, $(\bJ TY)_{B,x_0}\subset TB$ such
that
\begin{align}  \label{g4}
\begin{split}
&\bJ|_{T^{(1,0)}_{x_0}X_G}= \frac{\sqrt{-1}}{2\pi}
{\rm diag} (a_1, \cdots, a_{n-n_0})\in \End(T^{(1,0)}_{x_0}X_G),\\
&\bJ^2|_{(\bJ TY)_B} =  \frac{-1}{4\pi ^2}
{\rm diag} (a^{\bot,2}_{1}, \cdots, a^{\bot,2}_{n_0})
\in \End((\bJ TY)_{B,x_0}),
\end{split}\end{align}
with $a_j, a^\bot_{j}>0$, and let $\{w^{0,j}\}_{j=1}^{n-n_0}$,
$\{e ^{\bot j}\}_{j=1}^{n_0}$ be their dual basis, then
$$e^0_{2j-1}=\tfrac{1}{\sqrt{2}}(w^0_j+\overline{w}^0_j) \ \ \mbox{
and}\ \ \
e^0_{2j}=\tfrac{\sqrt{-1}}{\sqrt{2}}(w^0_j-\overline{w}^0_j)\,,$$
 $j=1,\dotsc,n-n_0\, ,$
forms an orthonormal basis of $T_{x_0}X_G$.

 From now on, we use the coordinate in Section \ref{s3.2}
induced by the above basis.

Denote by $Z^0= (Z^0_1,\cdots, Z^0_{2n-2n_0})$, $Z^\bot=
(Z^\bot_{1},\cdots, Z^\bot_{n_0})$, then $Z=(Z^0,Z^\bot)$.

In what follows we will use the complex coordinates
$z^0=(z^0_1,\cdots,z^0_{n-n_0})$, thus
$Z^0=z^0+\overline{z}^0$, and $w^0_i=\sqrt{2}\tfrac{\partial}{\partial z^0_i}$,
$\overline{w}^0_i=\sqrt{2}\tfrac{\partial}{\partial\overline{z}^0_i}$, and
\begin{align}\label{0g4}
&e^0_{2i-1}=\tfrac{\partial}{\partial z^0_i}
+ \tfrac{\partial}{\partial\overline{z}^0_i},\quad
e^0_{2i}=\sqrt{-1}(\tfrac{\partial}{\partial z^0_i}
-\tfrac{\partial}{\partial\overline{z}^0_i}).
\end{align}
We will also identify $z^0$ to $\sum_i z^0_i\tfrac{\partial}{\partial z^0_i}$
and $\overline{z}^0$ to
$\sum_i\overline{z}^0_i\tfrac{\partial}{\partial\overline{z}^0_i}$ when
we consider $z^0$ and $\overline{z}^0$ as vector fields. Remark that
\begin{equation}\label{g5}
\Big\lvert\tfrac{\partial}{\partial z^0_i}\Big\rvert^2=
\Big\lvert\tfrac{\partial}{\partial\overline{z}^0_i}\Big\rvert^2
=\dfrac{1}{2}\,,
\quad\text{so that $|z^0|^2=|\overline{z}^0|^2=\dfrac{1}{2} |Z^0|^2$\,.}
\end{equation}

It is very useful to rewrite $\cL^0_2$ by using
the creation and annihilation operators.
Set
\begin{align}\label{g6}
\begin{split}
&b_i=-2{\tfrac{\partial}{\partial z^0_i}}
+\frac{1}{2}a_i\overline{z}^0_i\,,\quad
b^{+}_i=2{\tfrac{\partial}{\partial\overline{z}^0_i}}
+\frac{1}{2}a_i z^0_i\,,\quad b=(b_1,\cdots,b_{n-n_0})\, ; \\
&b^\bot_j = -{\tfrac{\partial}{\partial Z^\bot_j}} + a^\bot_jZ^\bot_j,\quad
b^{\bot +}_j = {\tfrac{\partial}{\partial Z^\bot_j}} + a^\bot_jZ^\bot_j,
\quad b^{\bot}=(b^\bot_1,\cdots,b^\bot_{n_0}).
\end{split}\end{align}
Then for any polynomial  $g(Z^0,Z^\bot)$ on $Z^0$ and $Z^\bot$,
\begin{align}\label{g7}
&[b_i,b^{+}_j]=b_i b^{+}_j-b^{+}_j b_i =-2a_i \delta_{i\,j},
& [b_i,b_j]=[b^{+}_i,b^{+}_j]=0\, ,\nonumber\\
& [g,b_j] =  2 \tfrac{\partial}{\partial z_j^0}g,
& [g,b_j^+]
= - 2\tfrac{\partial}{\partial \overline{z}_j^0}g\,,  \\
&[b^\bot_i, b^{\bot +}_j]=-2a^\bot_i \delta_{i\,j},\, &
[b^\bot_j, b^{\bot}_k]=[b^{\bot +}_j, b^{\bot +}_k]=0,\nonumber\\
& [g,b^{\bot}_j]=  - [g,b^{\bot +}_j]
= \tfrac{\partial}{\partial Z^\bot_j}g\,. &  \nonumber
\end{align}

Set
\begin{align}\label{g8}
\cL=\sum_{j=1}^{n-n_{0}} b_j b^{+}_j ,\quad
\cL^{\bot}= \sum_{j=1}^{n_{0}} b^{\bot}_j b^{\bot +}_j,\quad
\nabla_{0,\cdot}= \nabla_{\cdot}
+ \frac{1}{2} R^{L_B}_{x_0}(\mR,\cdot).
\end{align}

From  \eqref{h8} and \eqref{g3}, for $U,V\in T_{x_0}B$, we get
\begin{align}\label{0g5}
&R^{L_B}_{x_0}(U,V) = -2\pi \sqrt{-1} \left\langle \bJ P^{TX_G}U,
P^{TX_G}V\right\rangle.
\end{align}

By  (\ref{c24}),  \eqref{g4}, \eqref{g6}, \eqref{g8} and \eqref{0g5}, we have
\begin{align}\label{0g6}
\begin{split}
&b_i= - 2 \nabla_{0,\tfrac{\partial}{\partial z ^0_i}},\quad
b^{+}_i=2\nabla_{0,\tfrac{\partial}{\partial\overline{z}^0_i}},\quad
\nabla_{0,e^\bot_j} = \nabla_{e^\bot_j},\\
&\tau_{x_0}= \sum_j a_j + \sum_j a_j^\bot.
\end{split}\end{align}

From (\ref{0c35}), (\ref{g8}) and \eqref{0g6}, we get
\begin{align}\label{g9}
&\cL^0_2 =-\sum_{j=1}^{2n-2n_{0}} (\nabla_{0, e^0_{j}} )^2
-\sum_{j=1}^{n_{0}}\Big( (\nabla_{e^{\bot}_{j}})^2
-|a^\bot_jZ^\bot_j|^2\Big) -2 \omega_{d,x_0} -\tau_{x_0}\\
&\hspace*{5mm}= \cL+ \cL^{\bot} -2 \omega_{d,x_0}.\nonumber
\end{align}

By \cite[\S 8.6]{Taylor96}, \cite[Theorem 1.15]{MM04a}, we know
\begin{thm}\label{t3.4}
The spectrum of the restriction of $\cL$ on $L^2(\bR^{2n-2n_0})$
 is given by
\begin{equation}\label{g10}
{\spec}\left({{\cL}|_{{L^2(\bR^{2n-2n_0})}}}\right)= \Big\lbrace
2\sum_{i=1}^{n-n_0}\alpha_i^0 a_i\,:\,
 \alpha^0 =(\alpha_1^0 ,\cdots,\alpha_{n-n_0}^0)\in\bN ^{n-n_0}\Big\rbrace,
\end{equation}
and an orthogonal basis of the eigenspace of
$2\sum_{i=1}^{n-n_0}\alpha_i^0  a_i $ is given by
\begin{equation}\label{g11}
b^{\alpha^0 }\Big((z^0)^{\beta}\exp\Big({- \frac{1}{4}\sum_i
a_i|z^0_i|^2}\Big)\Big)    \,,\quad\text{with $\beta\in\bN^{n-n_0}$}\,.
\end{equation}

The spectrum of the restriction of $\cL^{\bot}$ on $L^2(\bR^{n_0})$
 is given by
\begin{equation}\label{g12}
{\spec}\left({{\cL^{\bot}}|_{{L^2(\bR^{n_0})}}}\right)=
\Big\lbrace 2 \sum_{i=1}^{n_0}\alpha_i^\bot a^\bot _i\,:\, \alpha
^\bot =(\alpha ^\bot_1,\cdots,\alpha ^\bot_{n_0}) \in\bN ^{n_0}
\Big\rbrace,
\end{equation}
and  the eigenspace of
$2 \sum_{i=1}^{n_0}\alpha_i^\bot a_i^\bot$ is one dimensional and
an orthonormal basis is given by
\begin{equation}\label{g14}
\Big(\prod_{i=1}^{n_0} \sqrt{\frac{\pi}{a^\bot_i}} (2a^\bot_i)^{\alpha^\bot_i}
(\alpha^\bot_i !)\Big)^{-1/2}
(b^\bot)^{\alpha ^\bot} \exp\Big({- \frac{1}{2}\sum_i
a^\bot_i|Z^\bot_i|^2}\Big)   \,.
\end{equation}
\end{thm}

Especially,  the orthonormal basis of $\Ker(
\cL|_{L^2(\bR^{2n-2n_0})})$; $\Ker(
{\cL^{\bot}}|_{{L^2(\bR^{n_0})}})$ are
\begin{align}\label{g15}
\begin{split}
&\Big (\frac{a ^\beta}{2 ^{|\beta|} \beta!}
\prod_{i=1}^{n-n_0} \frac{a_i}{2\pi}\Big)^{\frac{1}{2}}
\Big((z^0)^\beta
\exp\Big (-\frac{1}{4} \sum_{j=1}^{n-n_0} a_j |z^0_j|^2\Big),
\, \, \beta \in \bN^{n-n_0};\\
&G^\bot(Z^\bot)=\Big(\prod_{i=1}^{n_0} \frac{a^\bot_i}{\pi}\Big)^{\frac{1}{4}}
\exp\Big (-\frac{1}{2}\sum_{i=1}^{n_0} a ^\bot_i|Z^\bot_i|^2\Big).
\end{split}
\end{align}

Let $P_{\cL}(Z^0,Z^{\prime 0})$,
$P_{\cL^\bot}(Z^\bot,Z^{\prime\bot})$ (resp. $P(Z,Z^{\prime})$)
 be the kernels of the orthogonal projections $P_{\cL}$, $P_{\cL^\bot}$
(resp. $P$) from $L^2(\bR^{2n-2n_0})$ onto $\Ker( \cL) $,
$L^2(\bR^{n_0})$ onto $\Ker (\cL^{\bot})$ (resp.
$L^2(\bR^{2n-n_0})$ onto $\Ker (\cL +\cL^{\bot})$).

From \eqref{g15}, we get
\begin{equation}\label{g16}
\begin{split}
P_{\cL}(Z^0,Z^{\prime 0}) =&\Big (\prod_{i=1}^{n-n_0}\frac{a_i}{2\pi}\Big)
\:\:\exp\Big(-\frac{1}{4}\sum_{i=1}^{n-n_0}
a_i\big(|z^0_i|^2+|z^{\prime 0}_i|^2
-2z^0_i\overline{z}^{\prime 0}_i\big)\Big),\\
P_{\cL^\bot}(Z^\bot,Z^{\prime\bot})=&
\Big (\prod_{i=1}^{n_0}\sqrt{\frac{a ^\bot_i}{\pi}}\Big)
\exp\Big(-\frac{1}{2}\sum_{i=1}^{n_0}
a ^\bot_i(|Z^\bot_i|^2+|Z^{\prime\bot}_i|^2)\Big),\\
P(Z,Z^{\prime}) =&P_{\cL}(Z^0,Z^{\prime 0})
P_{\cL^\bot}(Z^\bot,Z^{\prime\bot}).
\end{split}
\end{equation}

Let $P^N $  be the orthogonal projection from $L^2(\bR^{2n-n_0},
(\Lambda(T^{*(0,1)}X)\otimes E)_{x_0})$ onto $N=\Ker (\cL^0_2)$.
Let $P^N(Z,Z^{\prime})$ be the associated kernel.

Recall that the projection $I_{\bC \otimes E_B} $ from $(\Lambda
(T^{*(0,1)}X)\otimes E)_B$ onto $\bC\otimes E_B$ is defined in
Introduction.

By \eqref{0a00},  \eqref{0a1}, (\ref{c24}) and \eqref{g4},
  \begin{align}\label{0g17}
-\om_{d,x_0}\geq \nu_0 \quad \mbox{on} \,\,  \Lambda ^{>0} (T^{*(0,1)}X),
\end{align}
 thus
\begin{align}\label{g17}
&P^N(Z,Z^{\prime}) =P(Z,Z^{\prime}) I_{\bC\otimes E_B}.
\end{align}

If $\bJ =J$ on $P$, then by \eqref{g16} and \eqref{g17},
\begin{align}\label{g18}
\begin{split}
&P^N(Z,Z^{\prime}) =\exp\Big(-\frac{\pi}{2}\sum_{i=1}^{n-n_0}
\big(|z^0_i|^2+|z^{\prime 0}_i|^2 -2z^0_i\overline{z}^{\prime 0}_i\big)\Big)\\
&\hspace*{30mm}\times 2 ^{\frac{n_0}{2}}
\exp\Big(- {\pi}\big(|Z^\bot|^2+|Z^{\prime\bot}|^2\big) \Big)
I_{\bC\otimes E_B},
\end{split}\\
&P^N((0,Z^\bot),(0,Z^{\bot})) =2 ^{\frac{n_0}{2}}
\exp\Big(- 2{\pi}|Z^\bot|^2 \Big)I_{\bC\otimes E_B}.\nonumber
\end{align}

\subsection{Evaluation of  $P^{(r)}$:
a proof of \eqref{0.7} and \eqref{a0.7}}\label{s3.11}

Recall that $\delta$ is the counterclockwise oriented circle in $\bC$
of center $0$ and radius $\nu/4$.

By  (\ref{u9}),
\begin{align}\label{0c55}
P_{0,t} = \frac{1}{2 \pi i} \int_{\delta}(\lambda
-\cL^t_2)^{-1}d\lambda .
\end{align}

Let $f(\lambda, t)$ be a formal power series with values  in $
\End(L^2(\bR^{2n-n_0},(\Lambda(T^{*(0,1)}X)\otimes E)_{B,x_0}))$
\be\label{c54} f(\lambda, t)= \sum_{r=0}^\infty t^r f_r(\lambda),
\quad f_r(\lambda)\in
\End(L^2(\bR^{2n-n_0},(\Lambda(T^{*(0,1)}X)\otimes E)_{B,x_0})).
\ee

 By (\ref{c30}),  consider
the equation of formal power series for  $\lambda\in \delta$,
\begin{align}\label{c55}
(\lambda-\cL^0_2 - \sum_{r=1}^\infty t^r \mO_r) f(\lambda, t) =
\Id_{L^2(\bR^{2n-n_0},(\Lambda(T^{*(0,1)}X)\otimes E)_{B,x_0})}.
\end{align}

Let $N^\bot$ be the orthogonal space of $N$ in
$L^2(\bR^{2n-n_0},(\Lambda(T^{*(0,1)}X)\otimes E)_{B,x_0})$, and
$P^{N^\bot}$ be the orthogonal projection from
$L^2(\bR^{2n-n_0},(\Lambda(T^{*(0,1)}X)\otimes E)_{B,x_0})$ onto
$N^\bot$.

We decompose $f(\lambda, t)$ according the splitting
$L^2(\bR^{2n-n_0}, (\Lambda(T^{*(0,1)}X)\otimes E)_{B,x_0})=
N\oplus N^\bot$,
\begin{align}\label{c56}
 g_r(\lambda) = P^{N}f_r(\lambda),\quad
 f^\bot_r(\lambda)=P^{N^\bot}f_r(\lambda).
\end{align}

Using Theorem \ref{t3.4}, \eqref{g9}, \eqref{0g17},  (\ref{c56})
and identifying the powers of $t$ in (\ref{c55}), we find that
\begin{align}\label{c57}
\begin{split}
&g_0(\lambda) = \frac{1}{\lambda}P^{N},
\quad  f^\bot_0(\lambda)= (\lambda -\cL^0_2)^{-1}P^{N^\bot} ,\\
&f^\bot_r(\lambda)=(\lambda -\cL^0_2)^{-1}
 \sum_{j=1}^{r} P^{N^\bot}\mO_j   f_{r-j}(\lambda),\\
&  g_r(\lambda)  = \frac{1}{\lambda } \sum_{j=1}^{r} P^{N}\mO_j
f_{r-j}(\lambda).
\end{split}\end{align}

Recall that $P^{(r)}$ $(r\in \bN)$
is defined in \eqref{1ue53} and \eqref{ue66}.

\begin{thm} \label{0t3.6} There exist  $J_{r}(Z,Z^{\prime})$  polynomials
in $Z,Z^{\prime}$ with the same parity as $r$,  whose coefficients
are polynomials in $A$, $R^{TB}$, $R^{{\rm Cliff}_B}$, $R^{E_B}$,
$\mu^E$, $\mu ^{\rm Cliff}$
{\rm (}resp. $r^X$, $\tr[R^{T^{(1,0)}X}]$, $R^E$;
resp. $h$, $R^L$, $R^{L_B}$; resp. $\mu${\rm )} and their
derivatives at $x_0$ up to order $r-1$
 {\rm (}resp. $r-2$; resp. $r$; resp. $r+1${\rm )},
and in the inverses of the linear combination of the eigenvalues of $\bJ$
at $x_0$\,, such that
\begin{align}\label{c86}
&P^{(r)}(Z,Z^{\prime})= J_{r}(Z,Z^{\prime})P(Z,Z^{\prime}).
\end{align}
Moreover,
\begin{align}\label{1c52}
&P^{(0)}(Z,Z^{\prime})
=P^N(Z,Z^{\prime})=P(Z,Z^{\prime})I_{\bC\otimes E_B}.
\end{align}
\end{thm}
\begin{proof} By (\ref{0c55}), for $M>0$,
by combining Theorems \ref{tu1}-\ref{tu6}
and the argument as in \cite[\S 1.3]{MM04a},
we get another proof of the existence of
the asymptotic expansion of $P_{0,t}(Z,Z^{\prime})$
for $|Z|,|Z^{\prime}|\leq M$ when $t\to 0$.

By \eqref{1c19}, \eqref{1c20} and \eqref{1ue60}, this gives another proof
of Theorems \ref{t0.1}, \ref{tue17} for $|Z|,|Z^{\prime}|\leq M/\sqrt{p}$.
Moreover, by (\ref{ue28}), (\ref{ue31}) and (\ref{c56}),
\begin{align}\label{c90}
P^{(r)}= \frac{1}{2\pi i} \int_{\delta}  g_r (\lambda)d
\lambda +  \frac{1}{2\pi i} \int_{\delta} f_r^\bot
(\lambda)d \lambda.
\end{align}
From (\ref{c57}), (\ref{c90}), we get  (\ref{1c52}).

Generally, from Theorems \ref{t3.3}, \ref{t3.4}, (\ref{g7}),
(\ref{c57}), (\ref{c90}) and the residue formula, we conclude
 Theorem \ref{0t3.6}.
\end{proof}

\begin{proof}[Proof of \eqref{0.7} and \eqref{a0.7}]
As $\bJ =J$ on $\mu^{-1}(0)$, the condition \eqref{g2}
is verified.

From Theorem \ref{0t3.6}, \eqref{g18}, we get \eqref{0.7} and
\eqref{a0.7}.
\end{proof}

From Theorem \ref{t3.4}, (\ref{c57}), (\ref{c90}),
and the residue formula, we can get $P^{(r)}$ by using the operators
 $(\cL^0_2)^{-1}$, $P^N, P^{N^\bot}, \mO_k\ (k\leq r)$.

This gives us a direct method to compute $P^{(r)}$
in view of  Theorem \ref{t3.4}. In particular,
\begin{equation}\label{g51}
\begin{split}
 P^{(1)}=& -P^N\mO_1P^{N^\bot}(\cL^0_2)^{-1}P^{N^\bot}
-P^{N^\bot}(\cL^0_2)^{-1}P^{N^\bot}\mO_1P^N,
\end{split}
\end{equation}
and
\begin{equation}\label{g19}
\begin{split}
 P^{(2)}=& \frac{1}{2 \pi i} \int_{\delta}
  \Big[(\lambda -\cL^0_2)^{-1} P^{N^\bot} (\mO_1 f_1 + \mO_2 f_0)(\lambda)
+ \frac{1}{\lambda}  P^N  (\mO_1 f_1+ \mO_2 f_0)(\lambda)\Big]d\lambda \\
=& \frac{1}{2 \pi i} \int_{\delta}
\Big\{(\lambda -\cL^0_2)^{-1} P^{N^\bot}
\Big[\mO_1 \Big((\lambda -\cL^0_2)^{-1} P^{N^\bot} \mO_1
+ \frac{1}{\lambda}  P^N \mO_1 \Big) + \mO_2 \Big]\\
&+ \frac{1}{\lambda}  P^N
\Big[\mO_1 \Big((\lambda -\cL^0_2)^{-1} P^{N^\bot} \mO_1
+ \frac{1}{\lambda}  P^N \mO_1 \Big) + \mO_2 \Big]
\Big\}(\lambda -\cL^0_2)^{-1} d\lambda\\
=&(\cL^0_2)^{-1}P^{N^\bot}\mO_1 (\cL^0_2)^{-1}P^{N^\bot}\mO_1 P^N
- P^{N^\bot}  (\cL^0_2)^{-2} \mO_1 P^N\mO_1 P^N\\
&+ (\cL^0_2)^{-1}P^{N^\bot}\mO_1 P^N \mO_1 (\cL^0_2)^{-1}P^{N^\bot}
- (\cL^0_2)^{-1}P^{N^\bot}\mO_2 P^N\\
&+ P^N\mO_1(\cL^0_2)^{-1}P^{N^\bot}\mO_1(\cL^0_2)^{-1}P^{N^\bot}
-  P^N \mO_1 (\cL^0_2)^{-2} P^{N^\bot} \mO_1 P^N\\
&- P^N\mO_1 P^N\mO_1(\cL^0_2)^{-2}P^{N^\bot}
- P^N \mO_2(\cL^0_2)^{-1}P^{N^\bot}
 .
\end{split}
\end{equation}

In the next Subsection we will prove $P^N\mO_1 P^N =0$, thus
the second and seventh terms in \eqref{g19} are zero.

\subsection{A formula for $\mO_1$}\label{s3.111}
We will use the notation in Section \ref{s1}. All tensors in this
Subsection will be evaluated at the base point $x_0\in X_G$.

For $\psi$ a tensor on $X$, we denote by $\nabla^X \psi$ its
covariant derivative induced by $\nabla ^{TX}$.

If $\psi_1$ is a $G$-equivariant tensor, then we can consider it
as a tensor on $B=U/G$ with the covariant derivative $\nabla
^B\psi_1$, we will denote  by $$(\nabla ^B\nabla ^B
\psi_1)_{(c_je_j,c'_k e_k)} := c_jc'_k (\nabla ^B_{e_j}\nabla
^B_{e_k} \psi_1)_{x_0},$$ etc.

We denote by $\{e_a\}$ an orthonormal basis of $(TX, g^{TX})$.

To simplify the notation, we often denote by $U$ the lift $U^H$
of $U\in TB$.

 Recall that $\wi{\mu}\in TY$ is defined by \eqref{ah7}
 and the moment map $\mu$ \eqref{a6},
and that $A$ is the second fundamental form of $X_G$
defined by \eqref{0.6}.

\begin{lemma} \label{0t3.7} The following identities hold,
\begin{align}\label{g3.43}
&(\nabla ^{TY}_{\mR}\wi{\mu})_{x_0}= -\bJ \mR^\bot, \nonumber\\
&(\nabla ^{TY}_{\cdot}\nabla ^{TY}_{\cdot}\wi{\mu})_{(\mR,\mR)}
:= (\nabla ^{TY}_{e^H_j}\nabla ^{TY}_{e^H_i} \wi{\mu})_{x_0} Z_j Z_i\\
&\hspace{10mm}=  -P^{TY}\left( (\nabla ^{X}_{\mR^0}\bJ) (\mR^0+2\mR^\bot)
+ (\nabla ^{X}_{\mR^\bot}\bJ) \mR^\bot\right) \nonumber\\
&\hspace{15mm}-\bJ A(\mR^0)\mR^0- \frac{1}{2} T(\mR^0, \bJ\mR^0)
+  T(\mR^\bot, \bJ \mR^\bot)
. \nonumber
\end{align}
\end{lemma}
\begin{proof}
Recall that $P^{TY}, P^{T^HX}$ are the orthogonal projections from $TX$
onto $TY, T^HX$ defined in Section \ref{s4.2}.
Note that on $P$, by \eqref{0g2},
\begin{align}\label{g3.20}
\bJ e^{\bot,H}_i \in TY,\quad
\bJ e^{0,H}_i= (\bJ_G e^{0}_i)^H\in T^H P .
\end{align}

By (\ref{ah7}) and (\ref{a5}), for $K\in \kg$,
\begin{align}\label{g3.25}
-\left\langle  \bJ e^H_i,K^X \right\rangle =\nabla _{e^H_i} \mu(K)
= \left\langle \nabla ^{TY}_{e^H_i} \wi{\mu},K^X \right\rangle
+  \left\langle \wi{\mu}, \nabla ^{TY}_{e^H_i} K^X \right\rangle.
\end{align}
From (\ref{h1}), (\ref{h2b}), (\ref{h3}) and (\ref{g3.25}),
\begin{align}\label{g3.26}
\nabla ^{TY}_{e^H_i} \wi{\mu}
= -P^{TY}\bJ e^H_i - \frac{1}{2}\dot{g}^{TY}_{e_i^H} \wi{\mu}
= -P^{TY}\bJ e^H_i - T(e^H_i, \wi{\mu}) .
\end{align}

From (\ref{g3.26}) and the fact that $\widetilde{\mu}=0$ on $P$,
one gets the first equation in (\ref{g3.43}).

Now for $W$ (resp. $Y$) a smooth section of $TX$ (resp. $TY$), by
\eqref{ah4},
\begin{multline}\label{g3.40}
 \left\langle  \nabla ^{TY}_{e^H_j}P^{TY} W, Y \right\rangle
= e^H_j \left\langle   W, Y \right\rangle
- \left\langle  P^{TY}  W, \nabla ^{TY}_{e^H_j}Y \right\rangle\\
= \left\langle  \nabla ^{TX}_{e^H_j}W, Y \right\rangle
+  \frac{1}{2}\left\langle T(e^H_j,P^{T^HX} W), Y \right\rangle .
\end{multline}

By \eqref{g3.40},
\begin{align}\label{g3.41}
\nabla ^{TY}_{e^H_j}P^{TY} W= P^{TY} \nabla ^{TX}_{e^H_j}W
+\frac{1}{2}T(e^H_j,P^{T^HX} W).
\end{align}

By  \eqref{g3.26} and \eqref{g3.41},
\begin{multline}\label{g3.42}
\nabla ^{TY}_{e^H_j}\nabla ^{TY}_{e^H_i} \wi{\mu}
=-P^{TY} (\nabla ^{X}_{e^H_j}\bJ) e^H_i
  -P^{TY} \bJ \nabla ^{TX}_{e^H_j} e^H_i\\
- \frac{1}{2}T(e^H_j,P^{T^HX} \bJ e^H_i )
- \frac{1}{2}(\nabla ^{TY}_{e^H_j} \dot{g}^{TY}_{e^H_i}) \wi{\mu}
- \frac{1}{2}\dot{g}^{TY}_{e^H_i}(\nabla ^{TY}_{e^H_j}\wi{\mu}).
\end{multline}

 By (\ref{h0}) and (\ref{h4}), for $U_1,U_2$ sections of $TB$ on $B$,
\begin{align}\label{g3.29}
&\nabla ^{TX}_{U^H_2}U^H_1 = (\nabla ^{TB}_{U_2}U_1)^H
-\frac{1}{2}T(U^H_2,U^H_1).
\end{align}

By the definition of our basis $\{e^0_i, e^\bot_j\}$ in Section
\ref{s3.2},
\begin{align}\label{g3.31}
(\nabla ^{TB}_{e^0_i}e^0_j)_{x_0}=A(e^0_i)e^0_j,
\quad (\nabla ^{TB}_{e^0_i}e^\bot_j)_{x_0} =
 (\nabla ^{TB}_{e^\bot_j}e^0_i)_{x_0}=A(e^0_i)e^\bot_j,
\quad
(\nabla ^{TB}_{e^\bot_j}e^\bot_i)_{x_0}=0.
\end{align}
Thus by \eqref{h3}, \eqref{g2},  \eqref{g3.26}, \eqref{g3.42},
\eqref{g3.29},  \eqref{g3.31} and the facts that $A$ exchanges
$N_G$ and $TX_G$ on $X_G$,  and that $\wi{\mu}=0$ on $P$, we get
\begin{align}\label{g3.30}
&(\nabla ^{TY}_{\cdot}\nabla ^{TY}_{\cdot}\wi{\mu})_{(\mR,\mR)}
=  -P^{TY} (\nabla ^{X}_{\mR}\bJ) \mR    -\bJ A(\mR^0)\mR^0
- \frac{1}{2} T(\mR, \bJ\mR^0) +  T(\mR, \bJ \mR^\bot).
\end{align}

We use the closeness of $\om$ to complete the proof of \eqref{g3.43}.

 From \eqref{0.1}, for $U,V,W\in TX$,
\begin{align}\label{g3.23}
\langle(\nabla ^{X}_U\bJ )V,W\rangle=(\nabla ^{X}_U\omega)(V,W),
\end{align}
 thus
\begin{equation} \label{g24}
\langle(\nabla ^{X}_U\bJ )V,W\rangle
+\langle(\nabla ^{X}_V\bJ )W,U\rangle+\langle
(\nabla ^{X}_W\bJ )U,V\rangle=d\omega(U,V,W)=0.
\end{equation}

 By \eqref{h0}, \eqref{h4},  \eqref{g24} and \eqref{g3.20},
for $Y$ a smooth section of $TY$,
$$\langle \bJ \nabla ^{TX}_Y e^0_j, e^\bot_i \rangle
= -\langle \nabla ^{TX}_Y e^0_j,\bJ e^\bot_i \rangle=-\langle
T(e^0_j,\bJ e^\bot_i), Y \rangle$$ and
\begin{multline}\label{g25}
\langle T(e^\bot_i, \bJ e^0_j), Y\rangle
= -2 \langle \nabla ^{TX}_Y (\bJ e^0_j), e^\bot_i \rangle
= -2 \langle (\nabla ^{X}_Y \bJ)e^0_j, e^\bot_i \rangle
+ 2  \langle  T(e^0_j,\bJ e^\bot_i), Y \rangle\\
= 2 \langle (\nabla ^{X}_{e^0_j}\bJ) e^\bot_i,Y \rangle
-2   \langle (\nabla ^{X}_{e^\bot_i}\bJ) e^0_j,Y \rangle
+ 2  \langle  T(e^0_j,\bJ e^\bot_i), Y \rangle.
\end{multline}

From  \eqref{g3.30}, \eqref{g25}, we get the second equation of
\eqref{g3.43}.
\end{proof}

\begin{thm}\label{0t3.8} The following identity holds,
\begin{align}\label{g20}
\mO_1= & -\frac{2}{3} (\partial_j R^{L_B})_{x_0} (\mR,e_i)Z_j  \nabla_{0,e_i}
 -\frac{1}{3} (\partial_i R^{L_B})_{x_0} (\mR, e_i)\\
&- 2 \left \langle A(e^0_i)e^0_j,\mR^\bot \right \rangle
\nabla_{0,e^0_i} \nabla_{0,e^0_j}
   -\pi\sqrt{-1}\left \langle(\nabla_{\mR}^X \bJ) e_a, e_b\right \rangle
\,c(e_a)\,c(e_b) \nonumber\\
&+4 \pi ^2 \left \langle(\nabla_{\mR^0}^X  \bJ) (\mR^0  +2 \mR^\bot)
+ (\nabla_{\mR^\bot}^X  \bJ)\mR^\bot- T(\mR^\bot,\bJ \mR^\bot),
\bJ \mR^\bot\right \rangle\nonumber\\
 &+4 \pi ^2 \left \langle  \bJ A(\mR^0)\mR^0
+ \frac{1}{2} T(\mR^0, \bJ\mR^0),
 \bJ \mR^\bot\right \rangle
+ 4\pi \sqrt{-1}\left \langle \wi{\mu}^{\mathrm{Cliff}}+  \wi{\mu}^{E},
 \bJ \mR^\bot\right \rangle. \nonumber
\end{align}
\end{thm}

\begin{proof} For $\psi\in (T^*X \otimes \End(\Lambda (T^{*(0,1)}X)))_B
\simeq (T^*X \otimes (C(TX)\otimes_\bR \bC))_B$, where $C(TX)$ is
the Clifford
 bundle of $TX$, we denote by $\nabla^X\psi$ the covariant
derivative of $\psi$ induced by $\nabla^{TX}$.

From $[\nabla^{\text{Cliff}}_W, c(e_a)]= c(\nabla^{TX}_W e_a)$,
   we observe that for $W\in TB$,
\begin{align}\label{0c38}
\nabla^X _W (\psi({e}_a)c({e}_a))
&= (\nabla^X _W \psi)({e}_a)c({e}_a)
+ \psi(\nabla^{TX}_{W^H}{e}_a)c({e}_a)
+ \psi({e}_a)c(\nabla^{TX}_{W^H}{e}_a)\\
&= (\nabla^X _W \psi)({e}_a)c({e}_a) .\nonumber
\end{align}
Thus by \eqref{c24} and \eqref{0c38}, for $k\geq 2$,
\begin{multline}\label{0c38a}
-(2\om_d+\tau)(tZ)
=\frac{1}{2}\left(R^L({e}_a,{e}_b) \,c({e}_a)\,c({e}_b)\right)(tZ)\\
= \frac{1}{2}\sum_{r=0}^k \tfrac{\partial^r}{\partial t^r}
\left[(R^L({e}_a,{e}_b) \,c({e}_a)\,c({e}_b))(tZ)\right] |_{t=0} \frac{t^r}{r !} + \cO(t^{k+1})\\
=\frac{1}{2}\Big (R^L_{x_0}+ t (\nabla^X _\mR R^L)_{x_0} \Big )
(e_a,e_b)  \,c(e_a)\,c(e_b) + \cO(t^{2}).
\end{multline}

By Lemma \ref{0t3.7} and \eqref{c39}, we have
\begin{align}\label{g26}
-t^2 \langle \wi{\mu}^{E_p},  \wi{\mu}^{E_p}\rangle (tZ)
=& 4 \pi ^2\sum_{k=2}^3  \frac{1}{k!}
\frac{\partial ^k}{\partial t^k}\Big(|\wi{\mu}|_{g^{TY}}^2 (tZ)\Big)
\Big|_{t=0}\,  t^{k-2}\\
&+ 4\pi \sqrt{-1} t \left \langle \wi{\mu}^{\text{Cliff}}+  \wi{\mu}^{E},
 \bJ \mR^\bot\right \rangle_{x_0}
+\cO(t^2). \nonumber
\end{align}

The following two formulas are clear,
\begin{align}
\begin{split}
\frac{1}{2}\left. \frac{\partial^2}{\partial t^2}
|\widetilde{\mu}|^2_{g^{TY}}(tZ)\right|_{t=0}
&=\frac{1}{2}
\left.\left(\nabla\nabla|\widetilde{\mu}|^2_{g^{TY}}(Z)\right)_{(\mR,\mR)}
\right|_{Z=0}
=|\nabla_\mR ^{TY}\widetilde{\mu} |^2,\\
\frac{1}{3!}\left. \frac{\partial^3}{\partial t^3}
|\widetilde{\mu}|^2_{g^{TY}}(tZ)\right|_{t=0}
&=\frac{1}{6}\left.\left(\nabla\nabla\nabla|\widetilde{\mu}|^2_{g^{TY}}(Z)
\right)_{(\mR,\mR,\mR)}\right|_{Z=0} \label{g3.481}\\
&=\langle(\nabla^{TY}\nabla^{TY}\widetilde{\mu})_{(\mR,\mR)},
\nabla^{TY}_\mR\widetilde{\mu} \rangle.
\end{split}\end{align}

From Lemma \ref{0t3.7} and \eqref{g26}-\eqref{g3.481},
 we see that the contribution from
$-t^2 \langle \wi{\mu}^{E_p},  \wi{\mu}^{E_p}\rangle (tZ)$
is the last three terms of \eqref{g20}.

By \eqref{0c39}, \eqref{c37} and \eqref{g8}, we have
\begin{align}\label{g21}
\nabla_{t,e_i} = \nabla_{0, e_i}
+ \frac{t}{3} (\partial_j R^{L_B})_{x_0}Z_j (\mR,e_i)
-\frac{t}{2} (\frac{1}{\kappa}\nabla_{e_i}\kappa)(tZ) + \cO(t^2).
\end{align}

By $g_{ij}(Z)=\theta_i^k(Z) \theta_j^k(Z)$ and
\eqref{0c33}-\eqref{1c34}, we know
\begin{align}\label{g22}
\begin{split}
g_{ij}(Z)=& \delta_{ij}-2 \left \langle A(e^0_i)e^0_j,\mR^\bot \right \rangle
+ \cO(|Z|^2) \quad \mbox{for} \, \, 1\leq i,j\leq 2(n-n_0),\\
& \delta_{ij}+ \cO(|Z|^2) \quad \mbox{otherwise};\\
\kappa(Z)=& \det (g_{ij}(Z))^{1/2} = 1- \langle A(e^0_i)e^0_i, \mR^\bot\rangle
+  \cO(|Z|^2).
\end{split}\end{align}

From \eqref{g3.31}, \eqref{g21} and \eqref{g22},
 the first three terms of the right hand side of \eqref{g20} is
the coefficient
$t^1$ of the Taylor expansion of
 $-g^{ij}(tZ)(\nabla_{t,e_i}\nabla_{t,e_j}$ $-t
 \nabla_{t,\nabla^{TB}_{e_i}e_j(tZ)})$.

By \eqref{1ue1}, \eqref{g3.23} and the above argument,
 the proof of Theorem \ref{0t3.8} is complete.
\end{proof}

\begin{thm} \label{0t3.9} We have the relation
\begin{align}\label{g22.1}
P^N \mO_1 P^N=0.
\end{align}
\end{thm}
\begin{proof} By \eqref{g6} and \eqref{g16},
\begin{align} \label{g3.57}
\begin{split}
&b^+_iP^N =b^{\bot+}_iP^N =0\,,\quad \,
(b^{\bot}_i P^N)(Z,Z^{\prime})=2a^{\bot}_i Z^\bot_i P^N(Z,Z^{\prime}),\\
&(b_iP^N)(Z,Z^{\prime})=a_i(\ov{z}^0_i-\ov{z}^{\prime 0}_i)P^N(Z,Z^{\prime}).
\end{split}\end{align}

We learn from \eqref{g3.57} that for any polynomial $g(Z^\bot)$ in
$Z^\bot$, we can write $g(Z^\bot)P^N(Z,Z^{\prime})$ as sums of
$g_{\beta^\bot}(b^\bot)^{\beta^\bot} P^N(Z,Z^{\prime})$ with
constants $g_{\beta^\bot}$.
By Theorem \ref{t3.4},
\begin{equation}\label{g3.59}
 P_{\cL^\bot} (b^\bot)^{\alpha^\bot}   g(Z^\bot)P^N=0\,,
\quad\text{for $|\alpha^\bot|>0$}.
\end{equation}

Let $\{w_a\}$ be an orthonormal basis of $(T^{(1,0)}X,g^{TX})$.

Note that if $f,g$ are two $\bC$-linear forms, then
$$f(e_a)g(e_a)= f(w_a)g(\ov{w}_a)+f(\ov{w}_a)g(w_a).$$
Thus by
Theorem \ref{t3.4}, \eqref{0a2}, \eqref{g17} and \eqref{g3.57},
\begin{align}\label{g22.2}
P^N \left \langle(\nabla_{\mR}^X \bJ) e_a, e_b\right \rangle
\,c(e_a)\,c(e_b)P^N  = -2  P^N \left \langle(\nabla_{\mR}^X \bJ)w_a,
\ov{w}_a \right \rangle  P^N \\
=-2 \left \langle(\nabla_{z^0 +\ov{z}^{\prime 0}}^X \bJ)w_a,
\ov{w}_a \right \rangle P^N
= \sqrt{-1}  \tr|_{TX} [J(\nabla_{z^0 +\ov{z}^{\prime 0}}^X \bJ) ]P^N.\nonumber
\end{align}

By \eqref{g6}, \eqref{0g6}, \eqref{g17}, \eqref{g20},
\eqref{g3.57}-\eqref{g22.2}, we get
\begin{multline} \label{g22.3}
P^N \mO_1 P^N= P^N
\Big\{\frac{2}{3} (\partial_\mR R^{L_B})_{x_0}
(\mR,\tfrac{\partial}{\partial \ov{z}^0_i}) b_i
-\frac{1}{3}(\partial_{e^0_i} R^{L_B})_{x_0}(\mR,e^0_i) \\
+ \frac{1}{3} (\partial_\mR R^{L_B})_{x_0}
(\mR,e^\bot_j) b^\bot_j
-\frac{1}{3} (\partial_{e^\bot_j} R^{L_B})_{x_0}(\mR^0,e^\bot_j)\\
+ \pi \tr|_{TX} [J(\nabla_{\mR^0}^X \bJ)]
+ 8\pi^2 \left\langle (\nabla_{\mR^0}^X \bJ)\mR^\bot,\bJ \mR^\bot \right\rangle
\Big\}P^N.
\end{multline}

By  \eqref{g7}, \eqref{g3.57} and \eqref{g3.59},
\begin{align} \label{g22.4}
P^N Z^\bot_j  Z^\bot_k    P^N
=\frac{1}{2 a^\bot_k} P^N Z^\bot_j b ^\bot_k   P^N
= \frac{1}{2 a^\bot_k} \delta_{jk} P^N.
\end{align}

For $\psi$ a tensor on $X_G$, let $\nabla^{X_G}\psi$ be the covariant
derivative of $\psi$ induced by the Levi-Civita connection $\nabla ^{TX_G}$.

For $U,V,W\in T_{x_0}X_G$, by \eqref{g2}, \eqref{0g2} and
\eqref{0g5}, we have
\begin{align}\label{g22.5}
(\partial_{U} R^{L_B})_{x_0}  (V,W)
= -2\pi \sqrt{-1} \left\langle(\nabla^{X_G}_U \bJ_G) V,W\right\rangle
= -2\pi \sqrt{-1}\left\langle (\nabla^{X}_U \bJ) V,W\right\rangle.
\end{align}

From \eqref{g2},  \eqref{g4}, we know that
\begin{align}\label{ag22.5}
 \bJ e^\bot_j = \frac{a_j}{2\pi} Je^\bot_j.
\end{align}

 By Theorem \ref{t3.4}, \eqref{h8}, \eqref{0a00}, \eqref{g7}, \eqref{g24}
and  \eqref{g3.57}-\eqref{ag22.5}, we get
\begin{multline} \label{g22.6}
P^N \mO_1 P^N= P^N \Big\{-\frac{4\pi\sqrt{-1}}{3} \Big[
2 \left\langle (\nabla^{X}_{\mR^0} \bJ)\tfrac{\partial}{\partial z^0_i},
\tfrac{\partial}{\partial \ov{z}^0_i}\right\rangle
+ \Big\langle (\nabla^{X}_{\tfrac{\partial}{\partial z^0_i}} \bJ)
 \mR^0,\tfrac{\partial}{\partial \ov{z}^0_i} \Big\rangle \\
- \Big\langle (\nabla^{X}_{\tfrac{\partial}{\partial \ov{z}^0_i}}  \bJ)
 \mR^0,\tfrac{\partial}{\partial z^0_i} \Big\rangle \Big]
+ \pi \tr|_{TX} [J(\nabla_{\mR^0}^X \bJ)]
+2\pi \left\langle (\nabla_{\mR^0}^X \bJ)e^\bot_j,
Je^\bot_j\right\rangle \Big\}  P^N\\
=\pi \Big[ -4 \sqrt{-1}\left\langle (\nabla_{z^0+\ov{z}^{\prime 0}}^{X}\bJ)
 \tfrac{\partial}{\partial z^0_i},
\tfrac{\partial}{\partial \ov{z}^0_i}\right\rangle
  + \tr|_{TX} [J(\nabla_{z^0+\ov{z}^{\prime 0}}^X \bJ)]
-  2 \left\langle J (\nabla_{z^0+\ov{z}^{\prime 0}}^X \bJ)e^\bot_j,
e^\bot_j\right\rangle \Big]P^N =0 .
\end{multline}

The proof of Theorem \ref{0t3.9} is complete.
\end{proof}

From \eqref{g19} and Theorem \ref{0t3.9}, we get the following
general formula which will be used in Section \ref{s4},
\begin{equation}\label{g22.7}
\begin{split}
 P^{(2)}=&(\cL^0_2)^{-1}P^{N^\bot}\mO_1 (\cL^0_2)^{-1}P^{N^\bot}\mO_1 P^N
- (\cL^0_2)^{-1}P^{N^\bot}\mO_2 P^N\\
&+ P^N\mO_1(\cL^0_2)^{-1}P^{N^\bot}\mO_1(\cL^0_2)^{-1}P^{N^\bot}
- P^N \mO_2(\cL^0_2)^{-1}P^{N^\bot}\\
&+ (\cL^0_2)^{-1}P^{N^\bot}\mO_1 P^N \mO_1 (\cL^0_2)^{-1}P^{N^\bot}
-  P^N \mO_1 (\cL^0_2)^{-2} P^{N^\bot} \mO_1 P^N.
\end{split}
\end{equation}

\subsection{Example $(\bC P^1, 2\, \om_{FS})$}\label{s3.12}
Let $\om_{FS}$ be the K\"ahler form
associated to the Fubini-Study metric $g^{T\bC P^1}_{FS}$ on $\bC P^1$.
We will use the metric $g^{T\bC P^1}=2 \, g^{T\bC P^1}_{FS}$ on $\bC P^1$
in this Subsection.

Let $L$ be the holomorphic line bundle $\mO(2)$ on $\bC P^1$.
Recall that $\mO(-1)$ is the tautological line bundle of  $\bC P^1$.

We  will use the homogeneous coordinate $(z_0,z_1)\in \bC^2$
 for $\bC P^1\simeq  (\bC^2\setminus \{0\})/\bC^*$.

Denote by $U_i=\{[z_0,z_1]\in \bC P^1; z_i\neq 0\}$, $(i=0,1)$,
the open subsets of  $\bC P^1$, and the two coordinate charts are
defined by $\phi_i : U_i\simeq \bC$, $\phi_i([z_0,z_1])=
\frac{z_j}{z_i}$, $j\neq i$.

For any $i_0$, $i_1 \in \bN$,
 $ z_0^{i_0}z_1^{i_1}$ is naturally identified
to a holomorphic section of $\mO(-i_0 - i_1)^*$ on $\bC P^1$.
For any $k\in \bN$, we have
\begin{align}\label{g30}
H^0(\bC P^1, \mO(k)) = \bC \{s_{k,i_0}:= z_0^{i_0}z_1^{i_1},\, \,
i_0+i_1 =k,\,  {\rm and} \,\,  i_0, i_1 \in \bN\}.
\end{align}

On $U_i$, the trivialization of the line bundle $L$ is defined by
$L\ni s \to s/z_i^2$, here $z_i^2$ is considered as a holomorphic section of
$\mO(2)$.

In the following, we will work on $\bC$ by using
 $\phi_1: U_1\to \bC$. Then for $z\in \bC$,
 \begin{align}\label{g31}
\om_{FS}(z)= \frac{\sqrt{-1}}{2\pi}\overline{\partial}\partial
 \log ((1+|z|^2)^{-1}) = \frac{\sqrt{-1}}{2\pi}
\frac{dz\wedge d\ov{z}}{(1+|z|^2)^2}.
\end{align}

Let $h^L$ be the smooth Hermitian  metric on $L$ on $\bC P^1$
defined by for $z\in \bC$,
\begin{align}\label{g32}
|s_{2,0}|_{h^L}^2 (z)= (1+|z|^2)^{-2}.
\end{align}
Let $\nabla ^L$ be the holomorphic Hermitian connection
of $(L, h^L)$ with its curvature $R^L$.

By \eqref{g31} and \eqref{g32}, under our trivialization on  $\bC$
\begin{align}\label{g33}
&\nabla ^L= \ov{\partial} +\partial + \partial \log( |s_{2,0}|_{h^L}^2),\\
&\frac{\sqrt{-1}}{2\pi}R^L =
 \frac{\sqrt{-1}}{2\pi}\overline{\partial}\partial \log |s_{2,0}|_{h^L}^2
=  2 \, \om_{FS}=: \om. \nonumber
\end{align}

Let $K$ be the canonical basis of ${\rm Lie}\,  S^1= \bR$, i.e.
for $t\in \bR, \exp(tK)= e ^{2\pi \sqrt{-1} t}\in S^1$.

We define an  $S^1$-action on $\bC P^1$ by $g\cdot
[z_0,z_1]=[gz_0,z_1]$ for $g\in S^1$.

On our local coordinate, $g\cdot z= g z$, and the vector field
$K^{\bC P^1}$ on $\bC P^1$ induced by $K$ is
\begin{align}\label{g34}
K^{\bC P^1}(z):= \tfrac{\partial}{\partial t} \exp(-tK)\cdot z |_{t=0}
= -2\pi \sqrt{-1}\Big(z \tfrac{\partial}{\partial z}
- \ov{z} \tfrac{\partial}{\partial \ov{z}}\Big).
\end{align}

Set $$\mu(K)([z_0, z_1])= \frac{2|z_0|^2}{|z_0|^2+|z_1|^2} -1.$$

 Then, on $\bC$,
\begin{align}\label{g35}
\mu(K)= 2|z|^2\,  (1+|z|^2)^{-1} -1.
\end{align}

By \eqref{g31}, \eqref{g34} and \eqref{g35}, we verify easily
that $\mu$ is a moment map associated to the $S^1$-action on $(\bC
P^1, \om)$ in the sense of \eqref{a5}.

The  ${\rm Lie}\, S^1$-action on the sections of $L$ defined by
\eqref{a6} induces a holomorphical $S^1$-action on $L$. In
particular, from \eqref{g33}-\eqref{g35},
\begin{align}\label{g36}
\tfrac{\partial}{\partial t} \exp(-tK)\cdot s_{2,j} |_{t=0}
=: L_K s_{2,j} = 2\pi \sqrt{-1}(1-j)\, s_{2,j}.
\end{align}

By \eqref{g36}, the $S^1$-invariant sub-space of $H^0(\bC P^1,
L^p)$ and $\mu^{-1}(0)$ are
\begin{align}\label{g37}
H^0(\bC P^1, L^p)^{S^1} = \bC\,   s_{2p,p}, \quad
\mu^{-1}(0)=\{z\in \bC, |z|=1\},
\end{align}
and $S^1$ acts freely on $\mu^{-1}(0)$, thus $(\bC
P^1)_{S^1}=\{{\rm pt}\}$.

Under our trivialization of $L$, $s_{2p,j}\in H^0(\bC P^1, L^p)$
is the function $z ^j$, and from \eqref{g32},
\begin{align}\label{g38}
\| s_{2p,j}\|^2_{L^2} = \int_{\bC} \frac{|z|^{2j}}{(1+|z|^2)^{2p}} 2\, \om_{FS}
= \int_0^\infty  \frac{2 t^j\, dt}{(1+t)^{2p+2}}
= \frac{2\, j!\, (2p-j)!}{(2p+1)!}.
\end{align}
Thus $(\frac{(2p+1)!}{2\, (p!)^2})^{1/2}s_{2p,p}$  is an orthonormal
basis of $H^0(\bC P^1, L^p)^{S^1}$.

Let $\ov{\partial}^{L^p*}$ be the formal adjoint of the Dolbeault
operator $\ov{\partial}^{L^p}$. For $p\geq 1$,
the spin$^c$ Dirac operator $D_p$ in \eqref{defDirac}
and its kernel are given by
\begin{align}\label{ag38}
D_p=\sqrt{2}\left(\ov{\partial}^{L^p}+ \ov{\partial}^{L^p *}\right),
\quad \Ker D_p = H^0(\bC P^1, L^p).
\end{align}

 Finally, by Def. \ref{ta0}, for $p\geq 1$, we get
\begin{align}\label{g39}
\begin{split}
&P^G_p(z,z')=\frac{(2p+1)!}{2\, (p!)^2} s_{2p,p}(z)\otimes s_{2p,p}(z')^*,\\
&P^G_p(z,z)=\frac{(2p+1)!}{2\, (p!)^2} |s_{2p,p}|_{h^{L^p}}^2(z)
= \frac{(2p+1)!}{2\, (p!)^2}\frac{|z|^{2p}}{(1+|z|^2)^{2p}}.
\end{split}\end{align}

Note that our trivialization by $s_{2,0}$ is not unitary, thus we
do not see directly the off-diagonal decay \eqref{0.8} from
\eqref{g39}.

Here we will only verify that \eqref{g39} is compatible with
\eqref{a0.7}, \eqref{aa0.8} and \eqref{a0.8}.

Recall that Stirling's formula \cite[(3.A.40)]{Taylor96}
 tells us that as $p\to +\infty$,
\begin{align}\label{g40}
p! = (2\pi p)^{1/2} p^p\,  e^{-p}
\left(1+\cO\left(\frac{1}{p}\right)\right).
\end{align}

By \eqref{g40},
\begin{align}\label{g41}
\frac{(2p+1)!}{2\, (p!)^2}= \frac{\sqrt{p}}{\sqrt{\pi}\, e} \, 2^{2p}
\Big(1+\frac{1}{2p}\Big)^{2p}\Big(1+\cO(\frac{1}{p})\Big)
= \sqrt{\frac{p}{\pi}}\,  2^{2p} \Big(1+\cO(\frac{1}{p})\Big) .
\end{align}

Now, $\bC^*$ is an open neighborhood of $\mu^{-1}(0)$
and $B= \bC^*/S^1\simeq \bR^+$ by mapping $z\in \bC^*$ to $r=|z|\in \bR^+$.

By \eqref{g31}, the metrics on
 $\{|z|=r\}= \{r e ^{2\pi \sqrt{-1} \theta}; \theta\in \bR/\bZ \}$,
 $B\simeq \bR^+$ induced by $\om=2\, \om_{FS}$ is
\begin{align}\label{g42}
8\pi \, r^2 \, (1+r^2)^{-2}\, d\theta \otimes d\theta,
\quad \quad  g^{TB}= \frac{2}{\pi}(1+r^2)^{-2}\, dr \otimes dr.
\end{align}

From \eqref{g42}, the fiberwise volume function $h^2(r)$ in
\eqref{0.6} on $\bR^+$ is
\begin{align}\label{g43}
h^2(r)= \sqrt{8\pi} \, r \, (1+r^2)^{-1}.
\end{align}

From \eqref{g39}, \eqref{g41} and \eqref{g43}, we get for
$|z|=r$,
\begin{align}\label{g44}
h^2(r) P^G_p(z,z)
= \sqrt{8\pi}\frac{(2p+1)!}{2\, (p!)^2}
\Big(\frac{r}{1+r^2}\Big)^{2p+1}=
\sqrt{2 \, p}  \Big(\frac{2 r}{1+r^2}\Big)^{2p+1}\Big(1+\cO(\frac{1}{p})\Big).
\end{align}

When $|z|=1$, from \eqref{g44},
we re-find \eqref{aa0.8} and \eqref{a0.8}.

From \eqref{g42}, $\sqrt{2 \pi}\frac{\partial}{\partial r}$ is
an orthonormal basis of $(B, g^{TB})$ at $r=1$,
thus the normal coordinate $Z^\bot$ has the form
$r-1= \sqrt{2\pi} (Z^\bot + \cO(|Z^\bot|^2)$. Thus
\begin{align}\label{g45}
(2 r \, (1+r^2)^{-1})^{2p+1}
=  e ^{ (2p+1) \log (1-  \pi (Z^\bot )^2+\cO(|Z^\bot|^3) )}
= e ^{-2 \pi p (Z^\bot )^2}+\cdots.
\end{align}
This means that \eqref{g44}, \eqref{g45} are compatible with
\eqref{a0.7} and \eqref{g18}.

If we consider the sub-space $H^0(\bC P^1, L^p)_{-p}$ of $H^0(\bC
P^1, L^p)$ with the weight $-p$ of $S^1$-action,
 then by \eqref{a6} as in \eqref{g36}, and \eqref{g38},
$\sqrt{p+\frac{1}{2}}\,  s_{2p,0}$
is an orthonormal basis of  $H^0(\bC P^1, L^p)_{-p}$.

Thus the smooth kernel $P_p^{-p}(z,z')$ of the orthogonal projection from
$\cC^\infty ( \bC P^1, L^p)$ onto $H^0(\bC P^1, L^p)_{-p}$ is
\begin{align}\label{g46}
P_p^{-p}(z,z') = (p+\frac{1}{2})\,  s_{2p,0}(z)\otimes s_{2p,0}(z')^*,\quad
P_p^{-p}(z,z) = (p+\frac{1}{2}) (1+|z|^2)^{-2p}.
\end{align}

Note that $\mu^{-1}(-1)= \{0\}$, i.e. $-1$ is a singular value of
$\mu$.

Let $\mu_1$ be the moment map defined by $\mu_1(K)=\mu(K) +1$,
then  $H^0(\bC P^1, L^p)_{-p}$ is the corresponding
$S^1$-invariant holomorphic sections of $L^p$ with respect to the
corresponding $S^1$-action.

Thus $0$ is a singular value  of $\mu_1$ and this explains why we
have a factor $p$ in \eqref{g46} instead of
 $p^{1/2}$ in \eqref{g44}.


\section{Applications}\label{s6}

This Section is organized as follows. In Section \ref{s6.1}, we
explain Theorem \ref{t6.1}, the version of Theorem \ref{t0.1} when
we only assume that $\mu$ is regular at $0$. In Section
\ref{s6.2}, We explain how to handle the $\vartheta$-weight
Bergman kernel. In Section \ref{s6.3}, we deduce \eqref{aa0.8},
and \eqref{a0.8}
 from \cite[Theorem 4.18$^\prime$]{DLM04a}.
In Section \ref{s6.4}, we explain Theorem \ref{t0.1}
implies Toeplitz operator type properties on $X_G$.
In Section \ref{s6.5}, we extend our results
for non-compact manifolds and for covering spaces.
In Section \ref{s6.6},
we explain the relation on
the $G$-invariant Bergman kernel on $X$ and the Bergman kernel on $X_G$.

We use the notation in Introduction.

\subsection{Orbifold case}\label{s6.1}

In this Subsection, we only suppose that $0\in \kg^*$ is a regular
value of $\mu$, then $G$ acts only infinitesimal freely on
$P=\mu^{-1}(0)$, thus $X_G = P/G$ is a compact symplectic
orbifold.

Let $G^0= \{g\in G, g\cdot x =x \, \, \mbox{\rm for any}\, x\in
P\}$, then $G^0$ is a finite normal sub-group of $G$ and $G/G_0$
acts effectively on $P$.

We  will use the notation for the orbifold as in \cite[\S 1]{M05},
\cite[\S 4.2]{DLM04a}.

Let $U$ be a $G$-neighborhood of $P=\mu^{-1}(0)$ in $X$ such that
$G$ acts infinitesimal freely on $\ov{U}$, the closure of $U$.
From the construction in Section \ref{as4.2}, any $G$-equivariant
vector bundle $F$ on $U$ induces an orbifold vector bundle $F_{B}$
on the orbifold $B=U/G$.

The function $h$ in (\ref{0.6}) is only $\cC ^\infty$ on the
regular part of the orbifold $B$, and we extend continuously $h$
to $U/G$ from its regular part, which is $\cC ^\infty$ and we
denote it by $\widehat{h}$, then $\widehat{h}$ is also $\cC
^\infty$ on $U$.

As we work on $P$ in Section \ref{0s3.2}, we need not to modify
this part.

We need to modify Section \ref{s3.2} as follows.

Observe first that the construction in Section \ref{s4.2} works
well if we only assume that $G$ acts locally freely on $X$
therein.

Denote by $\nabla ^{T^H U}$ the connection on $T^H U$ as in
Section \ref{s4.2}, and on $P$, let $\nabla ^N$, $\nabla ^{T^H
P}$, ${^0\nabla} ^{T^H U}$ be the connections on $N,T^H U$ in
Section \ref{0s3.2} as in \eqref{a0.6}.

For $y_0\in P$, $W\in T^H U$ (resp.$T^H P$), we define $ \bR \ni t
\to x_t=\exp^{T^H U}_{y_0}(tW)\in U$  (resp. $\exp^{T^H
P}_{y_0}(tW)\in P$) the curve such that $x_t|_{t=0}=y_0$,
$\frac{dx}{dt}|_{t=0}=W$, $\frac{dx}{dt}\in T^H U$, $\nabla ^{T^H
U}_{\frac{dx}{dt}}{ \frac{dx}{dt}}=0$ (resp. $\frac{dx}{dt}\in T^H
P$, $\nabla ^{T^H P}_{ \frac{dx}{dt}}\frac{dx}{dt}=0$).

 By proceeding as in Section \ref{s3.2},
 we identify $B^{T^HU}(y_0,\var)$ to a subset of $U$ as following,
for $Z\in B^{T^HU}(y_0,\var)$,
$Z=Z^0+Z^\bot$, $Z^0\in T_{y_0}^H P$, $Z^\bot\in N_{x_0}$,
we identify $Z$ with
$\exp^{T^H U}_{\exp_{x_0}^{T^H P} (Z^0)}(\tau_{Z^0} Z^\bot)$.

 Set $G_{y_0}= \{g\in G, gy_0=y_0\}$,
then $G\cdot B^{T^H U}(y_0,\var) = G\times_{G_{y_0}}B^{T^H U}(y_0,\var)$
 is a $G$-neighborhood of $Gy_0$,
 and $(G_{y_{0}},B^{T^H U}(y_0,\var)) $ is a local coordinate of $B$.

As the construction in Section \ref{s3.2} is
$G_{y_0}$-equivariant, we  extend the geometric objects on
$G\times_{G_{y_0}} B^{T^H U}(y_{0},\var)$ to $G\times_{G_{y_0}}
\bR^{2n-n_0}=X_0$.

Thus we get the corresponding geometric objects on $G\times
\bR^{2n-n_0}$ by using the covering $G\times \bR^{2n-n_0} \to
G\times_{G_{y_0}} \bR^{2n-n_0}$, especially,
$\widehat{\mL}_p^{X_0}$ (where we use the\ $\widehat{\cdot}$\
notation to indicate the modification) is defined similarly on
$G\times \bR^{2n-n_0}$, and Theorem \ref{ta2} holds for
$\widehat{\mL}_p^{X_0}$.

Let $\widehat{\pi}_G: G\times \bR^{2n-n_0} \to  \bR^{2n-n_0}$ be
the natural projection and as in \eqref{h12}, we define
$\widehat{\Phi}=\widehat{h}\widehat{\pi}_G$, then the operator
 $\widehat{\Phi} \widehat{\mL}_p^{X_0} \widehat{\Phi} ^{-1}$
is well-defined on $T^H_{y_{0}}U\simeq \bR ^{2n-n_0}$.

Let $g^{T^HX_0}$ be the metric on $\bR^{2n-n_0}$ induced by $g^{TX_0}$,
and let $dv_{T^H X_0}$ be the  Riemannian volume form on
$(\bR^{2n-n_0}, g^{T^HX_0})$.

Let $P_{y_0,p}$ be the orthogonal projection from $L^2(\bR
^{2n-n_0}, (\Lambda (T^{*(0,1)}X)\otimes L^p\otimes E)_{y_0})$
onto $\Ker (\widehat{\Phi} \widehat{\mL}_p^{X_0} \widehat{\Phi}
^{-1})$ on $\bR ^{2n-n_0}$. Let $P_{y_0,p}(Z,Z^{\prime})$
$(Z,Z^{\prime}\in\bR ^{2n-n_0})$ be the smooth kernel of
$P_{y_0,p}$ with respect to $dv_{T^HX_0}(Z^{\prime})$.

Let $P_{0,p}^{G}$ be the orthogonal projection from
$\Omega^{0,\bullet}(X_0,L^p_0\otimes E_0)$
on $(\Ker D_p^{X_0})^G$, and let  $P_{0,p}^{G}(x,x')$ be the smooth kernel
of $P_{0,p}^{G}$
with respect to the volume form $dv_{X_0}(x')$.

Let $P_p^{X_0/G}(y,y')$ $(y,y'\in X_0/G)$ be the smooth kernel
associated to the operator on $X_0/G$ induced by
$\widehat{\Phi} \widehat{\mL}_p^{X_0} \widehat{\Phi} ^{-1}$ as in (\ref{b23}).

Note that our trivialization
 of the restriction of  $L$ on $B^{T^HU}(y_{0},\var)$ as in Section \ref{s3.2}
is not $G_{y_0}$-invariant, except that $G_{y_0}$ acts trivially on $L_{y_0}$.

For $x,x'\in X_0$, with their representatives $\wi{x},
\wi{x}'\in\bR ^{2n-n_0}$, we have
\begin{align}\label{6.1}
\widehat{h}(x) \widehat{h}(x')P_{0,p}^{G}(x,x')= P_p^{X_0/G}(\pi(x),\pi(x'))
= \frac{1}{|G^0|}
\sum_{g\in G_{y_0}} (g,1)\cdot P_{y_0,p} (g^{-1}\wi{x}, \wi{x}').
\end{align}
The second equation of \eqref{6.1} is from direct computation
(cf. \cite[(5.19)]{DLM04a}).

As we work on $G\times \bR^{2n-n_0}$,   for the operator
$\widehat{\Phi} \widehat{\mL}_p^{X_0} \widehat{\Phi} ^{-1}$, Prop.
\ref{p3.2} and Sections \ref{s3.3}-\ref{s3.6} still holds.

From Theorem \ref{tue17} for $P_{y_0,p}$ and \eqref{6.1}, we get
\begin{thm} \label{t6.1} Under the same notation in
Theorems \ref{t0.1}, \ref{tue17},
for $\alpha, \alpha'\in \bN^{2n-n_0}$, $|\alpha|+|\alpha'|\leq m$, we have
\begin{multline}\label{6.2}
(1+\sqrt{p}|Z^\bot|+\sqrt{p}|Z^{\prime\bot}|)^{m''}
\left |\frac{\partial^{|\alpha|+|\alpha'|}}
{\partial Z^{\alpha} {\partial Z^{\prime}}^{\alpha'}}
 \left (p^{-n+\frac{n_0}{2}} (\widehat{h}\kappa^{\frac{1}{2}} )(Z)
(\widehat{h}\kappa^{\frac{1}{2}})(Z^{\prime})
P^G_p\circ\Psi(Z,Z^{\prime}) \right.\right.\\
- \frac{1}{|G^0|} \sum_{r=0}^k  \sum_{g\in G_{y_0}}
(g,1)\cdot  P^{(r)}_{y_0} (g^{-1}\sqrt{p} Z,\sqrt{p} Z^{\prime})
p^{-\frac{r}{2}}\Big)\Big|_{\cC ^{m'}(P)}\\
\leq C  p^{-(k+1-m)/2}  (1+\sqrt{p} |Z^0|+\sqrt{p} |Z^{\prime 0}|)^{2(n+k+m'+2)+m}
\exp (- \sqrt{C''\nu p} \inf_{g\in G_{y_0}}$ $|g^{-1}Z-Z^{\prime}|)\\
+ \cO(p^{-\infty}).
\end{multline}
\end{thm}

If $Z=Z^{\prime}= Z^0$, then for $g\in G_{y_0}$, $gZ^0=Z^0$,
we use Theorem \ref{tue17} for $Z= Z^{\prime}= 0$ with the base point $ Z^0$,
and for the rest element in $G_{y_0}$,
we use  Theorem \ref{tue17} for $Z=Z^{\prime}= Z^0$ with the base point $y_0$,
then we get
\begin{multline}\label{a6.2}
\left |p^{-n+\frac{n_0}{2}} (\widehat{h}^2\kappa)(Z^0)P^G_p\circ\Psi(Z^0,Z^0)
- \frac{1}{|G^0|} \sum_{r=0}^k
\sum_{g\in G_{y_0}, gZ^0=Z^0}
(g,1)\cdot P^{(2r)}_{Z^0} (0,0)p^{-r}\right.\\
\left. - \frac{1}{|G^0|} \sum_{r=0}^{2k}  \sum_{g\in G_{y_0}, gZ^0\neq Z^0}
(g,1)\cdot P^{(r)}_{y_0}  (g^{-1}\sqrt{p} Z^0,\sqrt{p} Z^0)
p^{-\frac{r}{2}}\right |\\
\leq C  p^{-(2k+1)/2} \Big(1+ (1+\sqrt{p} |Z^0|)^{2(n+2k+2)}
\exp (- \sqrt{C'''\nu  p}|Z^0|)\Big).
\end{multline}

Note that if $g\in G_{y_0}$ acts as the multiplication by
$e^{i \theta}$ on $L_{y_0}$, then $(g,1)\cdot P^{(r)}_{y_0}$, $(g,1)\cdot P^{(r)}_{Z^0}$
in \eqref{a6.2} have a factor
$e^{i \theta p}$ which depends on $p$.

Of course, after replacing  $L$ by some power of $L$, we can
assume that $G_{y_0}$ acts as identity on $L$ for any $y_0\in P$,
in this case, $(g,1)\cdot P^{(r)}_{y_0}$, $(g,1)\cdot
P^{(r)}_{Z^0}$ do not depend on $p$.

From Theorem \ref{0t3.6} and \eqref{a6.2},
 if the singular set of $X_G$ is not empty,
analogous to the usual orbifold case \cite[(5.27)]{DLM04a},
$p^{-n+\frac{n_0}{2}}$ $P^G_p(y_0, y_0)$, $(y_0\in P)$ does not
have a uniform asymptotic expansion in the form
$\sum_{r=0}^{\infty} c_r(y_0)  p^{-r}$.

\subsection{$\vartheta$-weight Bergman kernel on $X$}\label{s6.2}


In this section, we assume that $G$ acts on $P=\mu^{-1}(0)$
freely.

Let $\mV$ be a finite dimensional irreducible representation of
$G$, we denote it by $\rho ^\mV: G\to \End(\mV)$.
Let $\vartheta$ be the highest weight
of the representation $\mV$.
Let $\mV^*$ be the trivial bundle on $X$ with $G$-action 
$\rho^{\mV^*}$ induced by $\rho^{\mV}$.

Let $P^\mV_p$ be the orthogonal projection from
$\Omega^{0,{\scriptscriptstyle{\bullet}}}(X,L^p\otimes E)$ on
$\Hom_G(\mV, \Ker D_{p})\otimes \mV\subset \Ker D_{p}$. Let
$P^\mV_p(x,x')$, $(x,x'\in X)$, be the smooth kernel of $P^\mV_p$
with respect to $dv_X(x')$.

We call $P^\mV_p(x,x')$ the $\vartheta$-weight Bergman kernel of
$D_p$.

We explain now the asymptotic expansion of $P^\mV_p(x,x')$ as
$p\to \infty$.


We will consider the corresponding objects in Sections 
\ref{s1}-\ref{s30} by replacing $E$ by $E\otimes \mV^*$. Especially, 
we denote by $D^{\mV^*}_p$ the corresponding spin$^c$ Dirac
operator associated to the bundle $L^p\otimes E\otimes \mV^*$.

Certainly, all results in Sections \ref{s1}-\ref{s30} still hold
for  the bundle $E\otimes \mV^*$.

Let $P^{\vartheta}_p$ be the orthogonal projection from
$\cC^\infty(X, E_p\otimes \mV^*)$ onto $(\Ker D^{\mV^*}_p )^G$,
and $P^{\vartheta}_p(x,x')$, $(x,x'\in X)$  the smooth kernel of
$P^{\vartheta}_p$ with respect to $dv_X(x')$.

As $\mV$ is an irreducible representation of $G$, we get
\begin{align}\label{a6.3}
\Ker D^{\mV^*}_p = (\Ker D_p) \otimes \mV^*,\quad
(\Ker D^{\mV^*}_p )^G= \Hom_G (\mV, \Ker D_p).
\end{align}

Let $\{v_i\}$ be an orthonormal basis of $\mV$ with respect to a
$G$-invariant metric on  $\mV$ and $\{v_i^*\}$ the corresponding
dual basis.

Let $dg$ be a Haar measure on $G$.
By Schur Lemma,
\begin{align}\label{a6.4}
\int_G g\cdot (v_j \otimes v^*_i)dg
= \frac{1}{\dim_\bC \mV} \delta_{ij}\, {\rm Id}_{\mV}.
\end{align}
Thus if $W$ is a finite dimensional representation of $G$ with
the highest weight $\vartheta$, then for any $s\in W$, we have
\begin{align}\label{a6.5}
s= (\dim_\bC \mV)\, \Big(\int_G g\cdot (s\otimes v^*_i)dg\Big) \otimes v_i
\in \Hom_G (\mV, W) \otimes \mV=W .
\end{align}

From \eqref{a6.5} and the $G\times G$-invariance of the kernel
$P^{\vartheta}_p(x,x')$, we get
\begin{align}\label{6.14}
\begin{split}
&P^{\mV}_p (x,x') =(\dim_\bC \mV)\, \sum_i(P^{\vartheta}_p(x,x')v_i^*, v_i),\\
&P^{\mV}_p (x,x) =(\dim_\bC \mV)\,  \tr_{\mV^*}P^{\vartheta}_p(x,x)
\in \End (\Lambda (T^{*(0,1)}X)\otimes E)_x.
 \end{split}
\end{align}
In fact, let $\{\psi_j\}$ be an orthonormal basis of
$\Ker(D^{\mV^*}_p)^G$, then for any $j$ fixed, in view of the
second equality in (\ref{a6.3}), one sees that
\begin{align}\label{a6.6}
 \psi_j^* \psi_j\in
\End_G(\mV)\ \ \ \mbox{and}\ \ \ \tr_{\mV}[ \psi_j^*
\psi_j]=\|\psi_j\|_{L^2}^2=1.
\end{align}
 Thus by Schur Lemma,
 \begin{align}\label{a6.7}
\psi_j^* \psi_j= \frac{1}{\dim_\bC \mV}\, {\rm Id}_{\mV}
\end{align}
and
$\{(\dim_\bC \mV)^{\frac{1}{2}}\psi_j v_i\}$ is an orthonormal basis of
$\Hom_G(\mV, \Ker D_{p})\otimes \mV\subset \Ker D_{p}$.

Let $U$ be a $G$-neighborhood of $P=\mu^{-1}(0)$ as in Theorem \ref{t0.1},
 $P^{\vartheta}_p$ is  viewed as a smooth section of
${\rm pr}^*_1(E_p\otimes \mV^*)_B\otimes {\rm pr}^*_2(E_p\otimes \mV^*)_B^*$
on $B\times B$, or as a $G\times G$-invariant smooth section of
${\rm pr}^*_1(E_p\otimes \mV^*)\otimes {\rm pr}^*_2(E_p\otimes \mV^*)^*$
on $U\times U$.

Moreover,  $v_i$, $v_i^*$ are
 smooth (not $G$-invariant) sections of $U\times \mV$, $U\times \mV^*$ on $ U$.
Thus from \eqref{6.14},
 $P^{\mV}_p$ is not a $G\times G$-invariant section of
${\rm pr}^*_1(E_p)\otimes {\rm pr}^*_2(E_p^*)$ on  $U\times U$.

Now Theorem \ref{t0.1} applies well to the bundle $E\otimes  \mV^*$,
 thus we get the asymptotic expansion of $P^{\vartheta}_p(x,x')$
as $p\to +\infty$,
and the leading term in the expansion of \\
$p^{-n+\frac{n_0}{2}}
(h\kappa^{\frac{1}{2}} )(x)(h\kappa^{\frac{1}{2}} )(x^{\prime})
P^{\vartheta}_p(x,x')$ is
$P(\sqrt{p}Z,\sqrt{p}Z^{\prime})I_{\bC\otimes (E\otimes  \mV^*)_B}$.

By \eqref{6.14},
 the leading term of the asymptotic expansion of 
$p^{-n+\frac{n_0}{2}}
(h\kappa^{\frac{1}{2}} )(x)(h\kappa^{\frac{1}{2}} )(x^{\prime})
P^{\mV}_p (x,x')$ is
\begin{align}\label{6.15}
(\dim_\bC \mV)^2\,  P(\sqrt{p}Z,\sqrt{p}Z^{\prime})I_{\bC\otimes E_B},\quad P(0,0)= 2^{n_0/2}.
\end{align}

Let $\Theta$ be the curvature of $P \to X_G$ as in Section
\ref{s4.2}.
Let $\rho ^{\mV^*}_*$ denote the differential of $\rho^{\mV^*}$.
By \eqref{h8},
\begin{align}\label{a6.8}
R^{(E\otimes \mV^*)_G}= R^{E_G}+ \rho^{\mV^*}_*(\Theta).
\end{align}

In the same way, we can define $\mathscr{I}^\mV_p$ a section of
$\End (\Lambda (T^{*(0,1)}X)\otimes E)_B$ on $X_G$ by \eqref{0.10}
for $P^{\mV}_p$.
From \eqref{0.15} (which will be proved in Section \ref{s4}),
\eqref{6.14}, \eqref{6.15} and \eqref{a6.8}, we get

\begin{thm} \label{t6.3} Under the condition of
Theorem \ref{t0.6}, the first coefficients of the asymptotic
expansion of $\mathscr{I}^\mV_p\in\End (E_G)$ in \eqref{0.14} is
\begin{align}\label{6.16}
&\Phi_0=(\dim_\bC \mV)^2 ,\\
&\Phi_1=\frac{1}{8\pi}  (\dim_\bC \mV)^2\,\Big( r^{X_G}_{x_0}
+ 6 \Delta_{X_G}\log h
+ 4R^{E_G}_{x_0}(w^0_j, \ov{w}^0_j)\Big)\nonumber\\
&\hspace*{15mm}+ \frac{1}{2\pi} (\dim_\bC \mV) \,
\tr_{\mV^*}\left[ \rho^{\mV^*}_* (\Theta)(w^0_j,
\ov{w}^0_j)\right]. \nonumber
\end{align}
\end{thm}

\subsection{Averaging the Bergman kernel:
a direct proof of \eqref{aa0.8} and \eqref{a0.8}}\label{s6.3}
We use the same assumption and notation in Theorem \ref{t0.1}.

 Let $P_p(x,x')$ be the smooth kernel of the orthogonal projection $P_p$
from $\Omega^{0,{\scriptscriptstyle{\bullet}}}(X,L^p\otimes E)$
onto $\Ker D_p$ with respect to $dv_X(x')$. Then $P_p(x,x')$ is the usual
Bergman kernel associated to $D_p$.

Let $dg$ be a Haar measure on $G$.
By Schur Lemma,
\begin{align}\label{6.17}
P^G_p(x,x')=\int_G ((g,1)\cdot P_p)(x,x')dg
=\int_G (g,1)\cdot P_p(g^{-1}x,x')dg.
\end{align}

One possible way to get Theorem \ref{t0.1} is to apply the full
off-diagonal expansion \cite[Theorem 4.18$^\prime$]{DLM04a}
 to \eqref{6.17}.

Unfortunately, we do not know how to get  the full off-diagonal expansion,
especially the fast decay along $N_G$ in \eqref{0.8} in this way.

However, it is easy to get
 \eqref{aa0.8} and \eqref{a0.8} as  direct consequences of
\cite[Theorem 4.18$^\prime$]{DLM04a} and \eqref{6.17}.

As  in Section \ref{0s3.2}, we denote by $TY$ the sub-bundle of $TX$
on a neighborhood of $P=\mu^{-1}(0)$
generated by the $G$-action and by $T^HP$ the orthogonal
complement of $TY$ in $(TP, g^{TP})$.

Take $y_0\in P$.
Let $\{e_i\}_{i=1}^{2(n-n_0)}$, $\{f_l\}_{l=1}^{n_0}$ be
orthonormal basis of $T^H_{y_0}P$, $T_{y_0}Y$. Then
$\{e_i\}_{i=1}^{2(n-n_0)}\cup \{f_l, J_{y_0}f_l\}_{l=1}^{n_0}$ is
an orthonormal basis of $T_{y_0}X$.
We use this orthonormal basis to get a local coordinate of $X$ by
using the exponential map $\exp ^X_{y_0}$.

We identify $B^{T_{y_0}X} (0,\var)$ to $B^X(y_0, \var)$
by the exponential map $Z\to \exp^X_{y_0} (uZ)$.

Let $\nabla ^{{\rm Cliff}\otimes E}$ be the connection on
$\Lambda (T^{*(0,1)}X) \otimes E$ induced by
$\nabla ^{{\rm Cliff}}$ and $\nabla^E$.

For $Z\in B^{T_{y_0}X} (0,\var)$, we identify
$L_Z$, $(\Lambda (T^{*(0,1)}X) \otimes E)_Z$, $(E_p)_Z$
 to $L_{y_0}$, $(\Lambda (T^{*(0,1)}X) \otimes E)_{y_0}$, $(E_p)_{y_0}$
by parallel transport with respect to the connections
$\nabla ^{L}$, $\nabla ^{{\rm Cliff}\otimes E}$, $\nabla ^{E_p}$
along the curve $\gamma_Z: [0,1]\ni u\to uZ$.

Under this identification, for $Z,Z'\in B^{T_{y_0}X} (0,\var)$,
one has $$P_p(Z,Z^{\prime}) \in
\End(\Lambda (T^{*(0,1)}X)\otimes E)_{y_0}.$$

Let $\mathbf{\kappa}_1(Z)$ be the function on $B^{T_{y_0}X} (0,\var)$ defined by
\begin{align}\label{b6.14}
dv_X(Z)= \mathbf{\kappa}_1(Z) dv_{T_{x_0}X}.
\end{align}

By \cite[Theorem 4.18$^\prime$]{DLM04a} (i.e. Theorem \ref{t0.1}
for $G=\{1\}$), there exist
$J_{r}(Z^{\prime})\in \End (\Lambda (T^{*(0,1)}X)$ $\otimes E)_{y_0}$,
polynomials in $Z^{\prime}$ with
the same parity as $r$, such that for any $k, m'\in \bN$, there
exist $C,M>0$ such that for $Z^{\prime}\in T_{y_0}X,
|Z^{\prime}|\leq \var$,
\begin{multline}\label{6.18}
\left |\frac{1}{p^{n}}
P_p(Z^{\prime},0)- \sum_{r=0}^k J_r(\sqrt{p} Z^{\prime})
\mathbf{\kappa}_1^{-1}(Z^{\prime}) e ^{-\frac{\pi}{2}p |Z^{\prime}|^2}
p^{-\frac{r}{2}}\right |_{\cC ^{m'}(P)}\\
\leq C  p^{-(k+1)/2}  (1+\sqrt{p} |Z^{\prime}|)^M
\exp (- \sqrt{C''\nu } \sqrt{p} |Z^{\prime}|) +\cO(p^{-\infty}),
\end{multline}
and
\begin{align}\label{b6.15}
J_0(Z)= I_{\bC\otimes E}.
\end{align}

For $K\in \kg$, $|K|$ small, $e ^K$ maps
 $(\Lambda (T^{*(0,1)}X) \otimes E)_{e ^{-K}y_0}$, $L_{e ^{-K}y_0}$
to $(\Lambda (T^{*(0,1)}X) \otimes E)_{y_0}$, $L_{y_0}$,
and under our identification, we denote these maps by
\begin{align}\label{a6.16}
f^E(K)\in \End(\Lambda(T^{*(0,1)}X) \otimes E)_{y_0},
\quad f^L(K)\in \End(L_{y_0})\simeq \bC.
\end{align}

By \cite[Prop. 5.1]{BeGeVe}, if we denote by
\begin{align}\label{b6.16}
 j_{\kg}(K) = {\det}_{\kg} (\frac{1-e ^{-ad\, K}}{ad\, K})
\end{align}
 for $K\in \kg$, then in exponential coordinates of $G$,
\begin{align}\label{b6.17}
d(e^K) = j_{\kg}(K) dK.
\end{align}

By \cite[Prop. 4.1]{DLM04a} (i.e. Theorem \ref{t0.0} for
$G=\{1\}$), \eqref{6.17}, as $G$ acts freely on $P$, we know
\begin{align}\label{a6.15}
P^G_p(y_0,y_0) = \int_{K\in \kg, |K|\leq \var}
f^E(K) (f^L (K) )^p  P_p(e ^{-K} y_0, y_0)j_{\kg}(K) dK
+ \cO(p^{-\infty}).
\end{align}

Let $S^L$ be the section of $L$ on $B^{T_{y_0}X}(0,\var)$ obtained
by parallel transport of a unit vector of $L_{y_0}$ with respect
to the connection $\nabla^L$ along the curve $\gamma_Z$.
Let $\Gamma^L$ be the connection form of $L$ with respect to this
trivialization.

Recall that for $K\in \kg$, the corresponding vector field $K^X$
on $X$ is defined in Section \ref{s4.2}.
Recall that $\{ K_i\}$ is a basis of $\kg$.

By \eqref{0c40}, for $K\in \kg$,
\begin{align}\label{a6.17}
\begin{split}
&(e^K \cdot S^L)(0) = e^K \cdot S^L(e^{-K} y_0) = f^L (K) S^L(0),\, \,
 {\rm with}\, \, f^L (0)=1,\\
&\Gamma^L_Z(K^X) = \frac{1}{2} R^L_{y_0}(Z, K^X)+ \cO(|Z|^2).
\end{split}\end{align}

By \eqref{a6}, \eqref{a5}, \eqref{a6.17} and $\mu=0$ on $P$, we get
\begin{multline}\label{b6.18}
 (L_{ K_j} (L_{ K_i} S^L))(0)
= (\nabla^L_{K_j^X} (\nabla^L_{K_i^X}  S^L- 2\pi \sqrt{-1} \mu(K_i) S^L))(0) \\
= \frac{1}{2} R^L_{y_0}(K_j^X, K^X_i)S^L(0)
=\pi \sqrt{-1}\langle  d \mu(K_i),K^X_j \rangle S^L(0)=0.
\end{multline}

By \eqref{a6}, \eqref{a6.17}, \eqref{b6.18} and $\mu=0$ on $P$,
we get
\begin{align}\label{a6.18}
\begin{split}\frac{\partial f^L}{\partial K_i} (0) S^L(0) &= (L_{ K_i} S^L)(0)
= (\nabla^L_{K_i^X}S^L)(0)=0,\\
\frac{\partial^2 f^L}{\partial K_i\partial K_j} (0) S^L(0)
&= \frac{\partial^2 }{\partial t_1\partial t_2}
(e^{t_1K_i+t_2 K_j} \cdot S^L)(0) |_{t_1=t_2=0} \\
&=  (L_{ K_j} (L_{ K_i} S^L) + L_{ K_i} (L_{ K_j} S^L))(0)
=0.
\end{split}\end{align}
Thus from \eqref{a6.18},
\begin{align}\label{a6.19}
(f^L (K) )^p = (1+ \cO(|K|^3))^p.
\end{align}
Moreover,
\begin{align}\label{a6.20}
f^E(K)= {\rm Id}_{(\Lambda (T^{*(0,1)}X) \otimes E)_{x_0}} + \cO(|K|),\quad
\mathbf{\kappa}_1 (Z)= 1+ \cO(|Z|^2).
\end{align}

Let $dv_Y$ be the Riemannian volume form on $(TY, g^{TY})$.
Observe also that if we denote by $i_{y_0}: G\to Gy_0$
the map defined by $i_{y_0}(g)=gy_0$, then
\begin{align}\label{a6.21}
 \frac{1}{h^2(y)}dv_Y(y) =i^{-1}_{y_0 *} dg,
\end{align}
 this gives us a factor
$\frac{1}{h^2(y_0)}$
when we take the integral on $\kg$ instead on the normal coordinate on $X$.

By \eqref{6.17}, \eqref{6.18}, \eqref{a6.15}, \eqref{a6.19}-\eqref{a6.21}
and the Taylor expansion for $\mathbf{\kappa}_1$, $f^E$, $f^L$,
  as in \cite[Theorems 5.8, 5.9]{BeGeVe},
we know that there exist $J'_r(Z)$ polynomials in $Z$ with same parity on $r$,
and $J'_0= I_{\bC \otimes E}$, such that
\begin{align}\label{6.19}
P^G_p(y_0,y_0) \sim  p^{n} \frac{1}{h^2(y_0)}
\int_{K\in \kg, |K|\leq \var} e ^{-\frac{\pi}{2}p |K|^2}
\sum_{r=0}^\infty  J'_r(\sqrt{p} K) p^{-r/2} dK.
\end{align}
Moreover,
\begin{align}\label{6.20}
\int_{K\in \kg} e ^{-\frac{\pi}{2}p |K|^2} dK
= 2 ^{\frac{n_0}{2}} p^{-\frac{n_0}{2}}.
\end{align}

After taking the integral on $\kg$, from \eqref{6.19} and
\eqref{6.20}, we get \eqref{aa0.8} and \eqref{a0.8}.

By  \eqref{6.14}, \eqref{6.19} and \eqref{6.20}, we get also
the asymptotic expansion for $P^\mV_p(y_0,y_0)$, $y_0\in P$.

\subsection{Toeplitz operators on $X_G$}\label{s6.4}

In this Subsection, we suppose that  $(X,\om)$ is a K\"ahler manifold,
$\bJ=J$, and  $L,E$ are holomorphic
vector bundles with holomorphic Hermitian connections $\nabla ^L,\nabla ^E$.
Let $G$ be a compact connected Lie group acting holomorphically
on $X,L,E$ which preserves
$h^L$ and $h^E$.

We suppose that $G$ acts freely on $P=\mu^{-1}(0)$.
 Then $(X_G,\om_G)$ is K\"ahler and $L_G,E_G$ are holomorphic on $X_G$.

In this case, there exists a natural isomorphism from
$ (\Ker D_p)^G$ onto $\Ker D_{G,p}$.


At the end of this Subsection, we will explain the corresponding
result in the  symplectic case,
 especially, for $p\gg 1$, we construct a natural isomorphism
from $(\Ker D_p)^G$ onto $\Ker D_{G,p}$.

\comment{
Let $\ov{\partial}^{L^p\otimes E,*}$
be the formal adjoint of the Dolbeault operator
$\ov{\partial}^{L^p\otimes E}$, then
\begin{align}\label{6.21}
D_p=\sqrt{2}(\ov{\partial}^{L^p\otimes E}+ \ov{\partial}^{L^p\otimes E,*}),
\end{align}
}

In the current situation, the spin$^c$ Dirac operator $D_p$
 was given by \eqref{60.21} and $D_p^2$
preserves the $\bZ$-grading of $\Omega^{0, \bullet}(X,L^p\otimes E)$.
Similar properties hold for $D_{G,p}$.

As in Section \ref{s3.0}, let $P_{G,p}$ be the orthogonal projection from
$\Omega^{0,\bullet}(X_G, L^p_G,\otimes E_G)$ onto $\Ker D_{G,p}$,
and let $P_{G,p}(x,x')$ be the corresponding smooth kernel.

 By the Kodaira vanishing theorem, for $p$ large enough,
\begin{align}\label{6.22}
(\Ker D_p)^G = H^0(X, L^p\otimes E)^G, \quad
\Ker D_{G,p}=H^0(X_G, L^p_G\otimes E_G).
\end{align}

\comment{
Let $i:P\hookrightarrow X$ be the natural injection.
Let $\pi_G:  \cC^{\infty}(P, L^p\otimes E)^G
 \to \cC^{\infty} (X_G,L^p_G\otimes E_G)$ be the natural identification.
Then by a result of Zhang \cite[Theorem 1.1 and Proposition
1.2]{Z99}, the map
$$ \pi_G\circ i^*: \cC^{\infty}(X, L^p \otimes E)^G
 \to \cC^{\infty} (X_G,L^p_G\otimes E_G) $$
induces a natural isomorphism
\begin{align}\label{6.23}
\sigma_p=  \pi_G \circ i^* :  H^0(X, L^p\otimes E)^G
\to H^0(X_G, L^p_G\otimes E_G).
\end{align}
(When $E=\bC$, this result was first proved in \cite[Theorem
3.8]{GuSt82}.)
}

As $D^2_p, D^2_{G,p}$ preserve the $\bZ$-gradings of
$\Omega^{0,\bullet}(X, L^p\otimes E)$, $\Omega^{0,\bullet}(X_G,
L^p_G\otimes E_G)$ respectively,  we only need to take care of
their restrictions on $\cC^{\infty}(X, L^p\otimes E)$ and
$\cC^{\infty}(X_G, L^p_G\otimes E_G)$.
In this way,
\begin{align}\label{6.23}
\begin{split}
&P^G_p(x,x')\in \cC^{\infty}(X\times X,
{\rm pr}_1^* (L^p\otimes E)\otimes {\rm pr}_2^* (L^p\otimes E)^*),\\
&P_{G,p}(x_0,x_0')\in \cC^{\infty}(X_G\times X_G,
{\rm pr}_1^* (L^p_G\otimes E_G)\otimes {\rm pr}_2^* (L^p_G\otimes E_G)^*).
\end{split}
\end{align}

 Recall that the morphism $\sigma_p:  H^0(X, L^p\otimes E)^G
\to H^0(X_G, L^p_G\otimes E_G)$ was defined in \eqref{60.23}.
Set
\begin{align}\label{a6.22}
\sigma_{p}^G=\sigma_{p}\circ P^G_p:
\cC^{\infty}(X, L^p\otimes E) \to H^0(X_G, L^p_G\otimes E_G).
\end{align}

Let $\sigma_{p}^{G*}$ be the adjoint of $\sigma_{p}^G$ with
respect to the natural inner products (cf. \eqref{h10}) on $\cC^{\infty}(X,
L^p\otimes E)$, $\cC^{\infty}(X_G, L^p_G\otimes E_G)$.
Set
\begin{align}\label{a6.23}
 \mP ^{X_G}_p := p^{-\frac{n_0}{2}} \sigma_{p}^G \circ \sigma_{p}^{G*}.
\end{align}

Let $\{s_{p,i}\}_{i=1}^{d_p}$ be an orthonormal basis of $H^0(X,
L^p\otimes E)^G$.
For $y_0 \in X_G,\ x,\ x'\in X$, one verifies
\begin{align}\label{a6.24}
\begin{split}
&P^{G}_p(x,x')= \sum_{i=1}^{d_p} s_{p,i}(x)\otimes s_{p,i}(x')^*,\\
&\sigma_{p}^G(y_0,x)=P^G_p (y_0,x),\quad
 \sigma_{p}^{G*}(x, y_0)=  P^G_p (x, y_0),
\end{split}\end{align}
where by $P^G_p (y_0,x)$ (resp. $P^G_p (x, y_0)$) we
mean $P^G_p (y,x)$ (resp. $P^G_p (x, y)$) for any $y\in \pi_G^{-1}(y_0)$,
which is well-defined by the $G$-invariance of $P_p^G$.

From \eqref{60.23}, we know that $\mP ^{X_G}_p$ commutes
 with the operator $P_{G,p}$ and
\begin{align}\label{6.24}
\mP ^{X_G}_p= P_{G,p}\mP ^{X_G}_pP_{G,p}.
\end{align}

Let $P^{G}_p|_P$ be the restriction of the smooth kernel
$P^{G}_p(x,x')$ on $P\times P$. Then $$P^{G}_p|_P(x,x')\in
\cC^\infty (P\times P, {\rm pr}_1^* (L^p\otimes E)\otimes {\rm
pr}_2^* (L^p\otimes E)^*)$$ is $G\times G$-invariant.
By composing with $\pi_G$, $$(\pi_G\circ P^{G}_p|_P)(x_0,x'_0)\in
\cC^{\infty}(X_G\times X_G, {\rm pr}_1^* (L^p_G\otimes E_G)\otimes
{\rm pr}_2^* (L^p_G\otimes E_G)^*).$$
We denote by
$\pi_G \circ P^{G}_p|_P$ the operator defined by the smooth kernel
$(\pi_G\circ P^{G}_p|_P)(x_0,x'_0)$ and the volume form $dv_{X_G}(x'_0)$.
Then from \eqref{a6.24}, we verify that
\begin{align}\label{add}
\mP ^{X_G}_p(x_0,x'_0) =p^{-\frac{n_0}{2}} P^G_p(x_0,x'_0)
= p^{-\frac{n_0}{2}}\pi_G \circ P^{G}_p|_P(x_0,x'_0).
\end{align}

\begin{defn}\label{d6.6}
 A family of operators $T_p:H^0(X_G, L^p_G\otimes E_G)\to
H^0(X_G, L^p_G\otimes E_G)$ is a Toeplitz operator if there exists
a sequence of smooth sections
 $g_l\in\cC^\infty (X_G, \End(E_G))$ with an asymptotic expansion
$g(\cdot,p)$ of the form $\sum_{l=0}^\infty p^{-l}g_l(x)$
such that for any $k\in \bN^*$, there exists $C>0$ 
such that for any $p\in \bN$,
\begin{align}\label{6.26}
\|T_p-  P_{G,p}\sum_{l=0}^k p^{-l}g_l(x) P_{G,p}\|^{0,0}\leq C p^{-(k+1)}.
\end{align}
Here $\|\quad \|^{0,0}$ is the operator norm with respect 
to the norm $\|\quad\|_{L^2}$.
We call $g_0(x)$ the principal symbol of $T_p$.
If $T_p$ is self-adjoint,
then we call $T_p$ is a self-adjoint Toeplitz operator.
\end{defn}

Recall that $h$ is the fiberwise volume function defined by \eqref{0.6}.

Let $dg$ be a Haar measure on $G$.

The main result of this Subsection is the following result.

\begin{thm}\label{t6.7} Let $f$ be a smooth  section of $\End(E)$ on $X$.
Let $f^G\in \cC^\infty (X_G,$ $\End(E_G))$ be the $G$-invariant part
of $f$ on $P$ defined by $f^G(x)= \int_G g f(g^{-1}x)dg$.
Then $\mT_{f,p}= p^{-\frac{n_0}{2}}
\sigma_{p}^{G}f \sigma_{p}^{G*}$ is a Toeplitz
operator with principal symbol $2^{\frac{n_0}{2}} \frac{f^G}{h^2}(x)$.
In particular $\mP ^{X_G}_p$ is a Toeplitz operator with principal symbol
$2^{\frac{n_0}{2}} / h^2(x)$.
\end{thm}
\begin{proof} Let $f^*$ be the adjoint of $f$. By writing $$f={f+f^*\over 2}+\sqrt{-1}{f-f^*\over
2\sqrt{-1}},$$ we may and we will assume from now on that $f$ is
self-adjoint.

 We need to find a family of  sections
$g_l\in \cC^\infty(X_G, \End(E_G))$ such that for any $m\geq 1$,
\begin{align}\label{a6.26}
\mT_{f,p} = \sum_{l=0}^m P_{G,p}\, g_l \, p^{-l}\, P_{G,p}
+\cO(p^{-m-1}).
\end{align} Moreover, we can make these $g_l$'s to be self-adjoint.

Let $U$ be a $G$-neighborhood of $P=\mu^{-1}(0)$ as in Theorem \ref{t0.1}.

Let $\psi$ be a $G$-invariant function on $X$ such that $\psi=1$
on a neighborhood of $P$ and $\supp (\psi )\subset \{y\in X,
d(y,P)<\var_0/2\}\cap U$.

Write
\begin{align}\label{a6.28}
\sigma_{p}^{G} f \sigma_{p}^{G*}
= \sigma_{p}^{G} \psi f \sigma_{p}^{G*}
+ \sigma_{p}^{G}  (1-\psi) f \sigma_{p}^{G*}.
\end{align}

For $x_0,x'_0\in X_G$, let $x, x'\in P$
such that $\pi(x)=x_0$, $\pi(x')=x'_0$. By \eqref{a6.24},
\begin{align}\label{a6.29}
(\sigma_{p}^{G} ((1-\psi) f) \sigma_{p}^{G*}) (x_0,x'_0)
= \int_X P^G_p(x,y)  ((1-\psi) f) (y) P^G_p(y,x') dv_X(y).
\end{align}

From Theorem \ref{t0.0}, \eqref{a6.29} 
and $\supp ((1-\psi) f )\cap P=\emptyset$,
we know that for any $l,m\in \bN$, there exists $C_{l,m}>0$ such
that for any $p\in \bN, x_0,x'_0\in X_G$,
\begin{align}\label{a6.30}
|(\sigma_{p}^{G}  ((1-\psi) f) \sigma_{p}^{G*})
(x_0,x'_0)|_{\cC^m(X_G\times X_G)}  \leq C_{l,m} p^{-l}.
\end{align}

We define $f_B\in \cC^\infty (B, \End(E_B))$ by
\begin{align}\label{b6.30}
 f_B(x_0)=  \int_G g (\psi f)(g^{-1}x)dg
\end{align}
for $x_0\in B, x\in U$ such that $\pi(x)=x_0$. 
Clearly, if $x_0\in P$, as $\psi|_P=1$, one gets
\begin{align}\label{b6.30a}
 f_B(x_0)=  f^G(x_0).
\end{align}

From (\ref{b6.30}), for $x_0,\ x_0'\in B,\ x,\ x'\in U$ such that
$\pi(x)=x_0,\ \pi(x')=x_0'$, one gets
\begin{multline}\label{a6.31}
\sigma_{p}^{G} \psi f \sigma_{p}^{G*}(x_0,x'_0)
= \int_U  P^G_p(x,y)  (\psi f) (y)P^G_p(y,x') dv_X(y)\\
= \int_B  P^G_p(x_0, y_0)f_B(y_0) P^G_p(y_0, x'_0) h^2(y_0) dv_B(y_0).
\end{multline}

For $x_0\in X_G$, we will  work on the normal coordinate of $X_G$
 with center $x_0$ as in Theorem \ref{t0.1}.

Recall that $P_{\cL}( Z^0,Z^{\prime 0})$ was defined by
\eqref{g16} with $a_i=a^\bot_i=2 \pi$ therein.

By  \eqref{a6.29},  \eqref{a6.30} and \eqref{a6.31}, 
for $|Z^0|, |Z^{\prime 0}|\leq \var_0/2$, 
\begin{align}\label{b6.29}
\mT_{f,p}(Z^{0},Z^{\prime 0}) 
-p^{-n_0/2} \int_{\stackrel{|W|\leq \var_0,}{W\in T_{x_0}B}} 
P^G_p(Z^{0},W) (f_{B}h^2)(W)  
 P^G_p(W, Z^{\prime 0})dv_{B}(W)
=\cO(p^{-\infty}).
\end{align}

By Theorem \ref{t0.1}, \eqref{b6.29} and the Taylor expansion of $f_B$,  
there exist $Q_{0,r}\in \End(E_{G,x_0})$
polynomials on $Z^0, Z^{\prime 0}$ with same parity on  $r$ such
that the following formula, obtained through compositions, holds,
\begin{multline}\label{a6.32}
\left | p^{-n+n_0} \mT_{f,p}(Z^{0},Z^{\prime 0})
 - \sum_{r=0}^k  (Q_{0,r} P_{\cL})  (\sqrt{p} Z^0,\sqrt{p} Z^{\prime 0})
p^{-\frac{r}{2}}\right |_{\cC ^{m'}(X_G)}\\
\leq C p^{-(k+1)/2}  (1+\sqrt{p} |Z^0|+\sqrt{p} |Z^{\prime 0}|)^M
\exp (- \sqrt{C''\nu } \sqrt{p} |Z^0-Z^{\prime 0}|)
+\cO(p^{-\infty}).
\end{multline}

On the normal coordinate in $X_G$, under the trivialization
induced by the parallel transport of $\nabla ^{(L^p\otimes E)_G}$
along the geodesic,
 by \cite[Theorem 4.18$^\prime$]{DLM04a}
(i.e. Theorem \ref{t0.1} for $G=1$), we get :
 there exist  $J_r (Z^0,Z^{\prime 0})\in \End (E_G)_{x_0}$,
  polynomials in $Z^0,Z^{\prime 0}$ with the same parity as $r$, such that
for any $k,m'\in \bN$, there exist $M\in \bN, C>0$ such that for
 $x_0\in X_G$, $Z^0,Z^{\prime 0}\in T_{x_0}X_G$, $|Z^0|, |Z^{\prime 0}|\leq  \var$,
\begin{multline}\label{6.25}
\left | p^{-n+n_0}
P_{G,p}(Z^{0},Z^{\prime 0})
 - \sum_{r=0}^k  (J_r P_{\cL})  (\sqrt{p} Z^0,\sqrt{p} Z^{\prime 0})
\kappa ^{-\frac{1}{2}}(Z^0)\kappa ^{-\frac{1}{2}}(Z^{\prime 0})
p^{-\frac{r}{2}}\right |_{\cC ^{m'}(X_G)}\\
\leq C p^{-(k+1)/2}  (1+\sqrt{p} |Z^0|+\sqrt{p} |Z^{\prime 0}|)^M
\exp (- \sqrt{C''\nu } \sqrt{p} |Z^{0}-Z^{\prime 0}|)+\cO(p^{-\infty}).
\end{multline}

By using the Taylor expansion of $\kappa^{-1/2}$, from \eqref{6.25},
there exist $J_{0,r}\in \End(E_{G,x_0})$, polynomials on $Z^0,
Z^{\prime 0}$ with same parity as $r$, such that
\begin{align}\label{6.27}
&\left | p^{-n+n_0}  P_{G,p}(Z^{0},Z^{\prime 0})
 - \sum_{r=0}^k  (J_{0,r} P_{\cL})  (\sqrt{p} Z^0,\sqrt{p} Z^{\prime 0})
p^{-\frac{r}{2}}\right |_{\cC ^{m'}(X_G)}\\
&\leq C p^{-(k+1)/2}  (1+\sqrt{p} |Z^0|+\sqrt{p} |Z^{\prime 0}|)^M
\exp (- \sqrt{C''\nu } \sqrt{p} |Z^0-Z^{\prime 0}|)+\cO(p^{-\infty}).\nonumber
\end{align}
Moreover, by \eqref{a0.7} and \eqref{b6.30a}, for $Q_{0,0},
J_{0,0}$ in \eqref{a6.32} and  \eqref{6.27}, we have
\begin{align}\label{6.28}
Q_{0,0} =2^{\frac{n_0}{2}} \frac{f^G}{h^2}(x_0),
\quad  J_{0,0} = \Id _{E_G}.
\end{align}

In what follows, all operators will be defined by their kernels
with respect to $dv_{T_{x_0}X_G}$. We will add a subscript $z^0$ or
$z^{\prime 0}$ when we need to specify the operator acting on the
variables $Z^0$ or $Z^{\prime 0}$.

By Theorem \ref{t3.4}, we know that
\begin{align}\label{6.29}
b^+_{j,z^0}  P_{\cL}=0,\quad (b_{j} P_{\cL})( Z^0, Z^{\prime 0})
=b_{j,z^0} P_{\cL}( Z^0, Z^{\prime 0})
= 2\pi (\ov{z}^0_j- \ov{z}_j^{\prime 0}) P_{\cL}( Z^0, Z^{\prime
0}).
\end{align}
Thus for $F(Z^0, Z^{\prime 0})$ a polynomial on $Z^0, Z^{\prime 0}$,
by \eqref{g7},  Theorem \ref{t3.4}, \eqref{6.29},
we can replace $\ov{z}^0_j$ in $F$ of $(FP_{\cL})(Z^0, Z^{\prime 0})$
by the combination of $b_{j,z^0}$ and $\ov{z}_j^{\prime 0}$,
thus there exist polynomials $F_{\alpha^0}$ $(\alpha^0 \in \bN ^{n-n_0})$
 on $z^0, Z^{\prime 0}$
 (resp. $F_{\alpha^0,0}$ on $z^0$, $\overline{z}^{\prime 0}$) such that
\begin{align}\label{6.30}
\begin{split}
&(FP_{\cL})(Z^0, Z^{\prime 0}) =\sum_{\alpha^0} b^{\alpha^0}_{z^0}
(F_{\alpha^0} P_{\cL})( Z^0, Z^{\prime 0}),\\
&((FP_{\cL}) \circ P_{\cL})(Z^0, Z^{\prime 0}) =
\sum_{\alpha^0} b^{\alpha^0}_{z^0} F_{\alpha^0,0}(z^0, \overline{z}^{\prime 0})
 P_{\cL}( Z^0, Z^{\prime 0}).
\end{split}\end{align}
In fact, by  Theorem \ref{t3.4}, the coefficient of 
$P_{\cL}( Z^0, Z^{\prime 0})$ in the right hand side of
 the second equation of \eqref{6.30} is anti-holomorphic on ${z}^{\prime 0}$.
Moreover, by Theorem \ref{t3.4}, $|\alpha^0|+ \deg F_{\alpha^0}$,
$|\alpha^0|+ \deg F_{\alpha^0,0}$ have the same parity with the
degree of $F$ on $Z^0, Z^{\prime 0}$. In particular, $F_{0,0}(z^0,
\overline{z}^{\prime 0})$ is a polynomial on $z^0$,
$\overline{z}^{\prime 0}$ and its degree has the same parity with
$\deg F$.

 We will denote by
\begin{align}\label{a6.33}
F_\cL :=F_0, \quad   F_{\cL,0}= F_{0,0}.
\end{align}

Let $(FP_\cL)_p$ be the operator defined by the kernel
$p^{n-n_0}  (FP_\cL)(\sqrt{p}Z^0, \sqrt{p}Z^{\prime 0})$.

By Theorem \ref{t3.4}, \eqref{6.27}, \eqref{6.29},
there exist polynomials $H_r(F)$ on $Z^0,Z^{\prime 0}\in T_{x_0}X_G$,
with the same parity with  $\deg F+r$, such that 
we have the following asymptotic at center $x_0$,
\begin{align}\label{6.31}
P_{G,p}(FP_{\cL})_p P_{G,p} \sim
\sum_{r=0}^\infty (H_r(F)P_{\cL})_p p^{-r/2};
\quad H_{0}(F)= F_{0,0} = F_{\cL,0},
\end{align}
with the reminder term estimated in the sense of
\eqref{a6.32} and \eqref{6.27}.

By Theorem \ref{t3.4}, \eqref{6.27}, \eqref{6.30} and
\eqref{6.31}, the coefficient of $p^{-k/2}$ in the expansion \eqref{6.31} of
$\sum_{r=0}^{k}P_{G,p}(Q_{0,r} P_\cL)_p p^{-\frac{r}{2}}$$
P_{G,p}$ is
 \begin{align}\label{6.34}
((Q_{0,k})_{\cL,0}P_\cL)_p + \sum_{r=1}^{k}(
H_{r}(Q_{0,{k-r}})P_\cL)_p.
\end{align}

Now, by \eqref{a6.22},
\begin{align}\label{b6.31}
\mT_{f,p}= p^{-\frac{n_0}{2}} \sigma^G_p f\sigma^{G*}_p =
P_{G,p}\mT_{f,p}P_{G,p}.
\end{align}
Thus by \eqref{a6.32}, \eqref{6.34} and \eqref{b6.31}, we get
\begin{align}\label{a6.34}
Q_{0,k} = (Q_{0,k})_{\cL,0}
+ \sum_{r=0}^{k-1} H_{k-r}(Q_{0,{r}}) .
\end{align}

By \eqref{6.27} and \eqref{6.28}, for ${\bf f}\in \cC^\infty(X_G,
\End(E_G))$, there exist polynomials $G_r({\bf f})(Z^0,Z^{\prime
0})$ with the same parity as $r$ such that in the normal
coordinates as above,
\begin{align}\label{6.33}
P_{G,p} {\bf f} P_{G,p} \sim \sum_{r=0}^{\infty}
(G_r({\bf f})P_\cL)_p \, p^{-\frac{r}{2}};
 \quad {\rm with } \, \, G_0({\bf f}) (Z^0,Z^{\prime 0})= {\bf f}(x_0).
\end{align}

By \eqref{6.27}, \eqref{6.33}, $(P_{G,p})^2= P_{G,p}$ and by
proceeding as in \eqref{a6.34}, we get
\begin{align}\label{6.35}
&G_r({\bf f}) = \sum_{i=0}^r H_i(G_{r-i}({\bf f}))
= (G_r({\bf f}))_{\cL,0} +\sum_{i=1}^r H_i(G_{r-i}({\bf f})).
\end{align}

From \eqref{6.28}, we define
\begin{align}\label{6.32}
g_0(x)=2^{\frac{n_0}{2}} \frac{f^G}{h^2}(x).
\end{align}

Assume that we have found $g_l \in \cC^\infty(X_G, \End(E_G))$,
$(l\leq k_0)$,
 self-adjoint sections such that \eqref{a6.26} holds for $m=k_0$.

We claim that $Q_{0,2k_0+1}$ is determined by $g_l, (l\leq k_0)$,
and there exists $g_{k_0+1}\in \cC^\infty(X_G, \End(E_G))$
self-adjoint such that $Q_{0,2k_0+2}$ is determined by $g_l\
(l\leq k_0+1)$.

By \eqref{a6.26}, \eqref{a6.32} and  \eqref{6.33}, for $0< k\leq 2k_0$,
\begin{align}\label{a6.35}
Q_{0,k}= \sum_{2l+j=k} G_j(g_l).
\end{align}
Then by \eqref{a6.34}, \eqref{6.35} and \eqref{a6.35}, for
$m=2k_0+1$,
\begin{multline}\label{6.36}
Q_{0,m} = (Q_{0,m})_{\cL,0}
+ \sum_{r=0}^{m-1} H_{m-r}\Big ( \sum_{l=0}^{[r/2]}G_{r-2l}(g_l)\Big)\\
=(Q_{0,m})_{\cL,0} + \sum_{l=0}^{[\frac{m-1}{2}]} \sum_{r=2l}^{m-1}
 H_{m-r} (G_{r-2l}(g_l))\\
=(Q_{0,m})_{\cL,0} +\sum_{l=0}^{[\frac{m-1}{2}]}\Big(G_{m-2l}(g_l)
-(G_{m-2l}(g_l))_{\cL,0}\Big).
\end{multline}

Set
\begin{align}\label{6.37}
\mathcal{F}_m= (Q_{0,m})_{\cL,0}-\sum_{l=0}^{[\frac{m-1}{2}]}(G_{m-2l}(g_l))_{\cL,0}.
\end{align}
Then by \eqref{a6.32}, \eqref{6.31}, $\mathcal{F}_m$ is a polynomial on
$z^0, \ov{z}^{\prime 0}$ with the same parity as $m$.
Moreover, as $T_{f,p}$ and $g_l$ are self-adjoint, we know that 
\begin{align}\label{c6.31}
\mathcal{F}^{(i)}_{m,x_0}(z^0, \ov{z}^{\prime 0})
= (\mathcal{F}^{(i)}_{m,x_0}(z^{\prime 0}, \ov{z}^{0}))^*.
\end{align}
Let $\mathcal{F}^{(i)}_m$ be the degree $i$ part of the polynomial
 $\mathcal{F}_m$ on $z^0, \ov{z}^{\prime 0}$.

We need to prove that for $m=2k_0+1$,
\begin{align}\label{6.38}
\mathcal{F}^{(i)}_m=0 \quad {\rm for}\, \, i>0.
\end{align}

 Set
\begin{align}\label{b6.32}
\begin{split}
&F^{(i)}_m(x_0,y_0)= \mathcal{F}^{(i)}_{m,x_0}(0,  \ov{z}^{\prime 0})
\in \End(E_{G,x_0}), \\
&\wi{F}^{(i)}_m(x_0,y_0)= (F^{(i)}_m(y_0,x_0))^* 
\in  \End(E_{G,y_0}),
\end{split} \end{align}
with $y_0= \exp_{x_0}^{X_G}(Z^{\prime 0})$, they define smooth
sections on a neighborhood of the diagonal of $X_G\times X_G$.
Clearly, $\wi{F}^{(i)}_m(x_0,y_0)$'s need not be polynomials of $z^0$ and
$\overline{z}'^0$.

Let $\psi:\bR\to [0,1]$ be an even function such that 
$\psi(u)=1$ for $|u|\leq \var_0/4$ and $0$ for $|u|>\var_0/2$.

Let $d^{X_0}(x_0,y_0)$ be the Riemannian distance on $X_G$.

We denote by $(\psi F^{(i)}_mP_{G,p})$,
$(P_{G,p}\psi \wi{F}^{(i)}_m)$ the operators defined by the kernel
$(\psi(d^{X_0}) F^{(i)}_mP_{G,p})(x_0,y_0)$, 
$(P_{G,p}\psi(d^{X_0}) \wi{F}^{(i)}_m)(x_0,y_0)$
with respect to $dv_{X_G}(y_0)$.
Set
 \begin{align}\label{b6.33}
P_{p,k_0}= \mT_{f,p}- P_{G,p}\sum_{l=0}^{k_0}g_lp^{-l}  P_{G,p}
-\sum_i (\psi F^{(i)}_{2k_0+1}P_{G,p})\, p^{(-2k_0+i-1)/2}.
\end{align}

By \eqref{a6.32}, \eqref{6.33},  \eqref{a6.35} and \eqref{6.36},
\begin{multline}\label{a6.36}
\left | p^{-n+n_0} \Big(\mT_{f,p}
- P_{G,p}\sum_{l=0}^{k_0}g_lp^{-l}  P_{G,p}\Big)(0,Z^{\prime 0})
- p^{-(2k_0+1)/2}(\mathcal{F}_{2k_0+1}P_\cL)
(0, \sqrt{p} Z^{\prime 0})\right |\\
\leq C p^{-k_0-1}  (1+\sqrt{p} |Z^{\prime 0}|)^M
\exp (- \sqrt{C''\nu } \sqrt{p} |Z^{\prime 0}|)+\cO(p^{-\infty}).
\end{multline}
Then by  \eqref{6.34} and \eqref{a6.36}, there exist polynomials
$Q_{0,r,k_0}$ on $Z^0, Z^{\prime 0}$ with the same parity as $r$
such that for $k>2k_0+2$, the kernel of the operator $P_{p,k_0}$
has the expansion at the normal coordinate of $x_0$, as
 \begin{multline}\label{b6.34}
\Big| p^{-n+n_0} P_{p,k_0} (0,Z^{\prime 0})
 - \sum_{r=2k_0+2}^k  (Q_{0,r,k_0} P_{\cL})  (0,\sqrt{p} Z^{\prime 0})
p^{-\frac{r}{2}}\Big|\\
\leq C p^{-(k+1)/2}  (1+\sqrt{p} |Z^{\prime 0}|)^M
\exp (- \sqrt{C''\nu } \sqrt{p} |Z^{\prime 0}|)+\cO(p^{-\infty}).
\end{multline}
We denote by $Q_{0,r,k_0}^p$
the operator defined as in \eqref{b6.33} by the kernel
$Q_{0,r,k_0}^p(x_0,y_0) = p^{n-n_0}\psi(d^{X_0}(x_0,y_0)) 
 (Q_{0,r,k_0} P_{\cL})(0,\sqrt{p} Z^{\prime 0})$.

Set 
 \begin{align}\label{a6.37}
K_{p,k}(x_0,y)= \psi(d^{X_0}(x_0,y))  P_{p,k_0}(x_0,y)
-  \sum_{r=2k_0+2}^k  Q_{0,r,k_0}^p p^{-\frac{r}{2}}(x_0,y).
\end{align}
Then by \eqref{b6.34},
\begin{multline}\label{b6.37}
|K_{p,k}(x_0,y)|\\
\leq C  p^{n-n_0-(k+1)/2}  (1+\sqrt{p} d^{X_0}(x_0,y))^M
\exp (- \sqrt{C''\nu } \sqrt{p} d^{X_0}(x_0,y))+\cO(p^{-\infty}).
\end{multline}
Thus for any $s\in \cC^\infty(X_G, L^p_G\otimes E_G)$,
\begin{align}\label{a6.38}
\begin{split}
\|K_{p,k}s \|_{L^2} ^2
&\leq \int_{x_0\in X_G} \Big(\int_{y_0\in X_G}
|K_{p,k}(x_0,y_0)| dv_{X_G}(y_0)\Big)\\
&\hspace*{5mm} \times \Big(\int_{y_0\in X_G}
|K_{p,k}(x_0,y_0)| |s|^2(y_0)
dv_{X_G}(y_0)\Big)dv_{X_G}(x_0) \\
&\leq C p^{-(k+1)}\|s\|_{L^2}^2.
\end{split}\end{align}
In the same way as in \eqref{a6.38}, 
\begin{align}\label{b6.38}
\|Q_{0,r,k_0}^p s \|_{L^2} ^2\leq C\|s\|_{L^2}^2.
\end{align}
Moreover, by Theorem \ref{t0.0},  \eqref{b6.33}, we get
\begin{align}\label{b6.36}
|((1-\psi(d^{X_0}))P_{p,k_0})(x_0,y_0)|= \cO(p^{-\infty}).
\end{align}

\comment{
Note that if we take $Z^0=0$ in \eqref{6.27}, we get that
uniformly on $x_0,y_0\in X_G$, $p\in \bN$,
 \begin{align}\label{a6.37}
\Big| P_{G,p} (x_0,y_0)\Big|\leq C p^{n-n_0}.
\end{align}

Take $k$ large enough in \eqref{b6.34},
by \eqref{a6.37}, uniformly on $y_0\in X_G$,
\begin{align}\label{b6.37}
\int_{x_0\in X_G}\Big|\Big(\Big( P_{p,k_0}
-  \sum_{r=2k_0+2}^k  Q_{0,r,k_0}^p p^{-\frac{r}{2}}\Big)
\circ P_{G,p}\Big)(x_0, y_0)\Big| dv_{X_G}(x_0) \leq C p^{-(k_0+1)}.
\end{align}
But for each operator $Q_{0,r,k_0}^p$, by \eqref{6.27},
\begin{align}\label{a6.38}
\Big|(Q_{0,r,k_0}^p\circ P_{G,p})(x_0, y_0)\Big|
\leq C p^{n-n_0} (1+\sqrt{p} \, d^{X_0}(x_0,y_0) )^M e^{-c\sqrt{p} \, d^{X_0}(x_0,y_0)}
+ \cO(p^{-\infty}).
\end{align}

Thus for any $s\in \cC^\infty(X_G, L^p_G\otimes E_G)$,
\begin{align}\label{b6.38}
\begin{split}
\|(Q_{0,r,k_0}^p P_{G,p})s \|_{L^2} ^2
&\leq \int_{x_0\in X_G} \Big(\int_{y_0\in X_G}
| Q_{0,r,k_0}^p P_{G,p}|(x_0,y_0) dv_{X_G}(y_0)\Big)\\
&\hspace*{5mm} \times \Big(\int_{y_0\in X_G}
| Q_{0,r,k_0}^p P_{G,p}|(x_0,y_0) |s|^2(y_0)
dv_{X_G}(y_0)\Big)dv_{X_G}(x_0) \\
&\leq C\|s\|_{L^2}^2.
\end{split}\end{align}

But for any
$s\in H^0(X_G, L_G^p\otimes E_G)$, we have
\begin{align}\label{b6.36}
s= P_{G,p}s.
\end{align}
}

 From \eqref{b6.37}, \eqref{a6.38}, \eqref{b6.38} and \eqref{b6.36},
we know that there exists $C>0$ such that for any
$s\in \cC^\infty(X_G, L_G^p\otimes E_G)$,
 \begin{align}\label{b6.35}
\| P_{p,k_0} s\|_{L^2} \leq C p^{-(k_0+1)} \|s\|_{L^2}.
\end{align}

Let $P_{p,k_0}^*$ be the adjoint of $ P_{p,k_0}$.
By \eqref{b6.35},
\begin{align}\label{b6.39}
\|P_{p,k_0}^* s  \|_{L^2} \leq C p^{-(k_0+1)} \|s\|_{L^2}  .
\end{align}
But
\begin{align}\label{b6.40}
P_{p,k_0}^* = \mT_{f,p}- P_{G,p}\sum_{l=0}^{k_0}g_lp^{-l}
P_{G,p} - \sum_i (P_{G,p}\psi \wi{F}^{(i)}_{2k_0+1})
p^{(-2k_0+i-1)/2}.
\end{align}

By \eqref{6.27} and the Taylor expansion of
$\wi{F}^{(i)}_{2k_0+1}$ under our trivialization of $E_G$ 
by using parallel transport along the path $[0,1]\ni u\to uZ^{\prime 0}$, 
we have that in the sense of
\eqref{a6.36},
\begin{multline}\label{a6.39}
p^{-n+n_0} \sum_i (P_{G,p}\wi{F}^{(i)}_{2k_0+1})(0,Z^{\prime 0})
p^{(-2k_0+i-1)/2}\\
\sim \sum_{\alpha, i,r}\left(J_{0,r} P_\cL\right) (0, \sqrt{p} Z^{\prime 0})
\frac{\partial ^\alpha \wi{F}^{(i)}_{2k_0+1}}
{\partial (Z^{\prime 0})^\alpha}(x_0,0)
\frac{(\sqrt{p}Z^{\prime 0})^\alpha}{\alpha !}
p^{(-2k_0+i -|\alpha|-r-1)/2}.
\end{multline}

By \eqref{a6.36}, \eqref{b6.39}, \eqref{b6.40} and \eqref{a6.39},
we know that all coefficients of $p^{(-2k_0-1+j)/2}$ for $j>0$
 of the right hand side of \eqref{a6.39} should be zero.
Thus we get for any $j>0$,
\begin{align}\label{a6.40}
\sum_{i=j}^{\deg F_m} \sum_{|\alpha|+r=0}^{i-j}
J_{0,r} (0, \sqrt{p} Z^{\prime 0})
\frac{\partial ^\alpha \wi{F}^{(i)}_{2k_0+1}}
{\partial (Z^{\prime 0})^\alpha}(x_0,0)
\frac{(\sqrt{p}Z^{\prime 0})^\alpha}{\alpha !}=0.
\end{align}

From \eqref{a6.40}, we will prove by recurrence
 that for any $j>0$
\begin{align}\label{a6.41}
\frac{\partial ^\alpha \wi{F}^{(i)}_{2k_0+1}}
{\partial (Z^{\prime 0})^\alpha}(x_0,0)=0 \quad \mbox{for}\,\,
 i-|\alpha|\geq j>0.
\end{align}

In fact, for $j=\deg F_{2k_0+1}$ in \eqref{a6.40}, by \eqref{6.28},
we get  $\wi{F}^{(\deg F_{2k_0+1})}_{2k_0+1}(x_0,0)=0$,
thus \eqref{a6.41} holds for  $j=\deg F_{2k_0+1}$.

Assume that  for $j> j_0>0$,  \eqref{a6.41} holds.
Then for $j=j_0$, the coefficient with $r>0$ in \eqref{a6.40} is zero,
thus by \eqref{6.28}, \eqref{a6.40} reads as
\begin{align}\label{a6.42}
 \sum_{\alpha}
\frac{\partial ^\alpha \wi{F}^{(j_0+|\alpha|)}_{2k_0+1}}
{\partial (Z^{\prime 0})^\alpha}(x_0,0)
\frac{(\sqrt{p}Z^{\prime 0})^\alpha}{\alpha !}=0.
\end{align}
From \eqref{a6.42}, we get \eqref{a6.41} for $j=j_0$.
The proof of \eqref{a6.41} is complete.

By \eqref{b6.39}, \eqref{a6.39}
and \eqref{a6.41},  by comparing the coefficient of $p^{-(2k_0+1)/2}$ in
\eqref{a6.36} and \eqref{b6.40}, we get
\begin{align}\label{a6.43}
\wi{F}^{(i)}_{2k_0+1}(x_0,Z^{\prime 0})
=\mathcal{F}^{(i)}_{{2k_0+1},x_0} (0, \ov{z}^{\prime 0}) 
+ \cO(|Z^{\prime 0}|^{i+1}).
\end{align}
Thus from \eqref{b6.32} and \eqref{a6.43}
\begin{align}\label{b6.41}
F^{(i)}_{2k_0+1}(Z^{\prime 0}, x_0)
= (\mathcal{F}^{(i)}_{2k_0+1,x_0}(0, \ov{z}^{\prime 0}))^*
+ \cO(|Z^{\prime 0}|^{i+1}).
\end{align}

Let $j_{X_G}: X_G\to X_G\times X_G$ be the diagonal injection.
By \eqref{b6.32},
\begin{align}\label{b6.410}
\frac{\partial}{\partial z^{\prime 0}_j}
F^{(i)}_{2k_0+1}=0 \quad \mbox{near}\, \,  j_{X_G} (X_G).
\end{align}

By \eqref{b6.32} again and recurrence, for $\alpha \in
\bN^{n-n_0}$, if $\alpha_j>0$, by taking $\alpha^\prime=(\alpha_0,
\cdots, \alpha_j-1,\cdots, \alpha_{n-n_0})$, one has
\begin{align}\label{b6.42}
\frac{\partial ^\alpha}{\partial z^{0,\alpha}} F^{(i)}_{2k_0+1}(\cdot, x_0)
|_{x_0}
=\frac{\partial }{\partial z^{0}_j}
 j_{X_G}^* \Big (\frac{\partial ^{\alpha^\prime}}{\partial z^{0,\alpha^\prime}}
 F^{(i)}_{2k_0+1}\Big)\Big|_{x_0}
- \frac{\partial ^{\alpha^\prime}}{\partial z^{0,\alpha^\prime}}
\frac{\partial }{\partial z^{\prime 0}_j}F^{(i)}_{2k_0+1}(\cdot, \cdot)\Big|_{0,0}
=0.
\end{align}
But by \eqref{b6.41}, for $|\alpha|\leq i$,
\begin{align}\label{b6.43}
\frac{\partial ^\alpha}{\partial z^{0,\alpha}} F^{(i)}_{2k_0+1}(\cdot, x_0)
|_{x_0}
= \Big(\Big(\frac{\partial ^\alpha}{\partial \ov{z}^{\prime 0,\alpha}}
\mathcal{F}^{(i)}_{2k_0+1,x_0}\Big)(0, \ov{z}^{\prime 0})\Big)^*.
\end{align}

From \eqref{b6.42} and \eqref{b6.43}, the $\alpha$-derivative for
$|\alpha|\leq i$ of $F^{(i)}_{m,x_0}(x_0, \cdot)$ is zero at
$x_0$. Thus
\begin{align}\label{c6.38}
\mathcal{F}^{(i)}_{m,x_0}(0, \ov{z}^{\prime 0})
= \mathcal{F}^{(i)}_{m,x_0}(z^0,0)=0.
\end{align}
Now, we consider the operator
\begin{align}\label{c6.39}
\frac{1}{\sqrt{p}} P_{G,p} \left(\nabla^{L^p\otimes E}
_{\varphi(d^{X_0}(x_0,y_0))
(\frac{\partial}{\partial z^0_j}+ \frac{\partial}{\partial \ov{z}^0_j})_{x_0}} 
\left(\mT_{f,p}- \sum_{l=0}^{k_0} P_{G,p} g_l p^{-l} 
P_{G,p}\right)\right)P_{G,p},
\end{align}
then the leading term of its asymptotic expansion as in \eqref{a6.32} is
\begin{align}\label{c6.40}
(\tfrac{\partial}{\partial z^0_j}\mathcal{F}_{2k_0+1,x_0})
(\sqrt{p}z^0, \sqrt{p}\ov{z}^{\prime 0})
P_{\cL}(\sqrt{p}z^0, \sqrt{p}\ov{z}^{\prime 0}) p^{-(2k_0+1)/2}.
\end{align}
Here by \eqref{c6.38}), 
$(\frac{\partial}{\partial z^0_j}\mathcal{F}_{2k_0+1,x_0})
(z^0,\ov{z}^{\prime 0})$
is an even degree polynomial on $z^0,\ov{z}^{\prime 0}$, and its constant
 term is zero. Now, by proceeding as \eqref{b6.33}--\eqref{b6.43} 
for the operator \eqref{c6.39}, by \eqref{c6.31}, we get
\begin{align}\label{c6.41}
\frac{\partial}{\partial z^0_j}\mathcal{F}_{2k_0+1,x_0}^{(i)}(z^0,\ov{z}^{\prime 0})
= \frac{\partial}{\partial z^{\prime 0}_j}
\mathcal{F}_{2k_0+1,x_0}^{(i)}(z^0,\ov{z}^{\prime 0})=0.
\end{align}
By continous this processus, we get \eqref{6.38}.

This means that $Q_{0, 2k_0+1}$ verifies also \eqref{a6.35}.

By the same argument, \eqref{6.36} still holds for $m=2k_0+2$.
Thus we can define
\begin{align}\label{b6.44}
g_{k_0+1}(x_0)= \mathcal{F}^{(0)}_{2k_0+2,x_0}
= \mathcal{F}^{(0)}_{2k_0+2,x_0}(0,0).
\end{align}
By proceeding exactly the same proof as before, we get
\eqref{6.38} for $m=2k_0+2$.
Thus for $k=2k_0+2$, \eqref{a6.35} still holds.

As $\pi_{G,p}, {\mT}_{f,p}, g_l$ $(1\leq l\leq k_0)$ are self-adjoint,
$g_{k_0+1}$ is also self-adjoint.

 By recurrence, we know that there exist $g_l$'s
such that \eqref{a6.26} holds for any $m$.

The proof of Theorem \ref{t6.7} is complete.
\end{proof}

\begin{cor}\label{aat6.8} 
For $f_1,f_2\in \cC^\infty (X)$,
we have
\begin{align}\label{c6.43}
[\mT_{f_1,p}, \mT_{f_2,p}]=
\frac{2^{n_0}}{\sqrt{-1}\, p}  P_{G,p}
\left\{\frac{f_1^G}{h^2},\frac{f_2^G}{h^2} \right\} P_{G,p}
+ \cO(p^{-2}).
\end{align}
Here $\{ \, , \, \}$ is the Poisson bracket
on $(X_G,2\pi \om_G)$:  for $g_1, g_2\in \cC^\infty (X_G)$,
if $\xi_{g_2}$ is the Hamiltonian vector field generated by $g_2$
which is defined by $2 \pi i_{\xi_{g_2}}\om_G=dg_2$, then 
\begin{align}\label{c6.44}
\{g_1, g_2\}= \xi_{g_2}(dg_1).
\end{align}
\end{cor}
\begin{proof} By applying \cite{MM04c} or \cite[\S 5.5]{MM05b},
(cf. \cite{BMS94} for another approach where they worked for $E=\bC$),
from Theorem \ref{t6.7}, we get immediately \eqref{c6.43}.
\end{proof}

\begin{lemma}\label{at6.8}  Let
$$T_p = \sum_{l=0}^\infty P_{G,p} g_l p^{-l} P_{G,p} +\cO(p^{-\infty})
: H^0(X_G, L^p_G\otimes E_G)\to
H^0(X_G, L^p_G\otimes E_G)$$
 be a Toeplitz operator with principal symbol
$g_0\in \cC^{\infty}(X_G, \End(E_G))$
.
Then

i) If $g_0$ is invertible, then $T_p^{-1}$ is a
Toeplitz operator with principal symbol
$g^{-1}\Id_{E_G}$.

ii) If $g_0= g \Id_{E_G}$ with $g\in \cC^{\infty}(X_G)$, $g>0$,
and $T_p$ is self-adjoint,
  then $T_p^{1/2}$ is a self-adjoint Toeplitz operator
with principal symbol $g^{1/2}\Id_{E_G}$.
\end{lemma}
\begin{proof} We only prove ii), the proof of i) is similar and simpler.

As $g>0$, there exist $C_0, C_1>0$ such that $C_0<g<C_1$. Thus for
any $s\in H^0(X_G, L^p_G\otimes E_G)$,
\begin{align}\label{b6.46}
\left\langle T_p s,s\right\rangle = \left\langle g_0
s,s\right\rangle +\cO\left(\frac{1}{p}\right) \|s\|_{L^2}^2 \geq
\left(C_0+\cO\left(\frac{1}{p}\right)\right) \|s\|_{L^2}^2.
\end{align}
Thus for $p$ large enough, $T^{1/2}_p : H^0(X_G, L^p_G\otimes
E_G)\to  H^0(X_G, L^p_G\otimes E_G)$ is well defined.

Let $\delta_1$ be the smooth bounded closed contour on $\{\lambda
\in \bC, {\rm Re}(\lambda )>0\}$ such that $[\frac{1}{2} C_0, 2 C_1]$
is in the interior domain got by $\delta_1$.

As in the proof of Theorem \ref{t6.7}, by recurrence, we will find
$f_l\in \cC^\infty (X_G, \End(E_G))$ such that
\begin{align}\label{b6.47}
T_p = (T_{m,p})^2 + \cO(p^{-m-1}) \quad {\rm with}
\quad T_{m,p} = \sum_{l=0}^m P_{G,p}f_l p^{-l} P_{G,p}.
\end{align}
Then for $p$ large enough,
\begin{multline}\label{b6.48}
T^{1/2}_p - T_{m,p} = \frac{1}{2 \pi i} \int_{\lambda\in \delta_1}
\lambda ^{1/2}
\Big [ (\lambda-T_p)^{-1} -(\lambda -(T_{m,p})^2)^{-1}\Big] d\lambda\\
= \frac{1}{2 \pi i} \int_{\lambda\in \delta_1}\lambda ^{1/2}
(\lambda-T_p)^{-1} (T_p-(T_{m,p})^2)(\lambda -(T_{m,p})^2)^{-1} d\lambda.
\end{multline}

If \eqref{b6.47} holds, then by \eqref{b6.48} we know that in the
sense of the operator norm,
\begin{align}\label{b6.49}
T^{1/2}_p - T_{m,p} = \cO(p^{-m-1}).
\end{align}

To complete the proof of Lemma \ref{at6.8},
 it remains to establish \eqref{b6.47}.

By \eqref{6.27}, there exist $Q_{0,r}\in  \End(E_G)_{x_0}$ such
that in the sense of \eqref{a6.32}, \eqref{6.27} and \eqref{6.31},
\begin{align}\label{b6.50}
T_p \sim \sum_{r=0}^\infty (Q_{0,r}P_\cL)_p p^{-r/2}.
\end{align}

We will prove by recurrence that there exist
$f_l\in \cC^\infty (X_G, \End(E_G))$ self-adjoint such that for any $k\in \bN$,
\begin{multline}\label{b6.51}
\Big| p^{-n+n_0} (T_p-(T_{k,p})^2) (\sqrt{p} Z^0, \sqrt{p} Z^{\prime 0})\Big|\\
\leq p^{-(2k+1)/2} (1+\sqrt{p} |Z^0|+\sqrt{p} |Z^{\prime 0}|)^M
\exp (- \sqrt{C''\nu } \sqrt{p} |Z^0-Z^{\prime 0}|)+\cO(p^{-\infty}).
\end{multline}

Set $f_0= g^{1/2} \Id_{E_G}$. Then \eqref{b6.47} is verified for
$m=0$.

Assume that for $k\leq m$, we have found $f_l$ such that
\eqref{b6.47} holds. If we denote the expansion of $(T_{m,p})^2$
in the sense of \eqref{6.27},
\begin{align}\label{b6.52}
(T_{m,p})^2 \sim \sum_{r=0}^\infty (\wi{Q}_{0,r}^m P_\cL)_p p^{-r/2}.
\end{align}
Then by the proof of \eqref{a6.35} for $2k_0+1$,
\begin{align}\label{b6.53}
\wi{Q}_{0,2m+1}^m = \sum_{2l+r=2m+1} G_r(g_l)=Q_{0,2m+1}.
\end{align}
Thus by \eqref{b6.53},
 \eqref{b6.51} stills holds when we replace the factor $p^{-(2m+1)/2}$ by
$p^{-m-1}$ at the right hand side of \eqref{b6.51}.
Thus
\begin{align}\label{b6.54}
T_p-(T_{k,p})^2 \sim \sum_{r=2m+2}^\infty
((Q_{0,r}-\wi{Q}_{0,r}^m)P_\cL)_p p^{-r/2}.
\end{align}

By \eqref{6.31}, \eqref{b6.54},  we know that
\begin{align}\label{b6.55}
Q_{0,2m+2}-\wi{Q}_{0,2m+2}^m =(Q_{0,2m+2}-\wi{Q}_{0,2m+2}^m)_{\cL,0}.
\end{align}
This means that $Q_{0,2m+2}-\wi{Q}_{0,2m+2}^m$ is a polynomial on
 $z^0,\ov{z}^{\prime 0}$ with even degree.

 Set
\begin{align}\label{b6.56}
f_{m+1}(x_0) = - \frac{1}{2} g^{-1/2} (Q_{0,2m+2}-\wi{Q}_{0,2m+2}^m)(0,0).
\end{align}
Then by the proof of \eqref{a6.35}, for $2k_0+2$, we know that
 the polynomial $Q_{0,2m+2}-\wi{Q}_{0,2m+2}^m$ equals to the constant $-2g^{1/2}f_{m+1}$.
Thus we prove \eqref{b6.51} for $k=m+1$.

By the above argument, we have established \eqref{b6.47}, thus
Lemma \ref{at6.8}.
\end{proof}

Since the isomorphism $\sigma_p : H^0(X, L^p\otimes E)^G$ $ \to
H^0(X_G, L^p_G\otimes E_G)$ is not an isometry, we define the
associated unitary operator,
\begin{align}\label{6.39}
\Sigma_p = \sigma_p^{G*} (\sigma_p^{G}\circ \sigma_p^{G*})^{-1/2} :
H^0(X_G, L^p_G\otimes E_G) \to H^0(X, L^p\otimes E)^G.
\end{align}
\comment{We extend $\Sigma_p$ to $L^2(X_G,L^p_G\otimes E_G)$ in
such a way that it vanishes over $(H^0(X_G, L^p_G\otimes
E_G))^\bot$. Then by \eqref{a6.24},
\begin{align}\label{6.40}
\Sigma_p\Sigma_p^* = P^G_{p},\quad
\Sigma_p^*\Sigma_p = P_{G,p},  \quad
  P^G_{p}  \Sigma_p P_{G,p}= \Sigma_p.
\end{align}
}

\begin{thm} \label{t6.8} Let $f$ be a $\cC^\infty$   section
 of $\End(E)$ on $X$. Then
\begin{align}\label{6.41}
T^G_{f,p} = \Sigma_p^* f \Sigma_p :
H^0(X_G, L^p_G\otimes E_G)\to  H^0(X_G, L^p_G\otimes E_G)
\end{align}
is a Toeplitz operator on $X_G$. Its principal symbol is
$f^G\in \cC^\infty (X_G, \End(E_G))$.
\end{thm}
\begin{proof} By \eqref{6.24} and \eqref{6.39},
\begin{align}\label{6.42}
T^G_{f,p} = (\mP^{X_G}_{p})^{-\frac{1}{2}}
\mT_{f,p}(\mP^{X_G}_{p})^{-\frac{1}{2}} .
\end{align}

By Theorem \ref{t6.7}, \eqref{a6.23},
 $\mP^{X_G}_{p}=
p^{-\frac{n_0}{2}}\sigma_p^{G}\circ \sigma_p^{G*}$, $\mT_{f,p}$
are Toeplitz operators on $X_G$ with principal symbols $
2^{n_0/2}/h^2(x)$, $2^{n_0/2}\frac{f^G}{h^2}(x)$ respectively.

By Lemma \ref{at6.8}, we know that
$(\mP^{X_G}_{p})^{-\frac{1}{2}}$ is a Toeplitz operator on $X_G$.

By \eqref{6.42}, $T^G_{f,p}$ has the expansion as \eqref{a6.32}.
By the proof of Theorem \ref{t6.7}, we then know that $T^G_{f,p}$
is a Toeplitz operator. \comment{ Now, we can use the formal power
series as in the proof
}
\end{proof}

\begin{rem}\label{at6.10} i) Certainly, by combining the argument
here and Section \ref{s6.1}, we can get the corresponding version
when $X_G$ is an orbifold.

ii) When $E=\bC$, and $f=1$, from Theorem \ref{t6.7},
$\mP^{X_G}_{p}$ is an elliptic (i.e. its principal symbol is invertible)
 Toeplitz operator.
 This is the analytic core result claimed in  \cite[\S 8]{Pao04}.

iii) When $E=\bC$ and $G$ is the torus $\bT^{n_0}$,   Theorem
\ref{t6.8} is one of the main results of Charles \cite[Theorem
1.2]{Charles04}, and in \cite[\S 5.6]{Charles04}, he knew also
that $\mP^{X_G}_{p}$ is an elliptic  Toeplitz operator. Moreover,
he established the corresponding version when $X_G$ is an
orbifold.
\end{rem}

If $X$ is only symplectic and $\bJ=J$, then as the argument
in \cite[\S 3e)]{TZ98}, $J$ induces an almost complex structure $J_G$
on $(TX)_B$, and $J_G$ preserves $N_{G,J}= N_G \oplus J_G N_G$ and $TX_G$.
Thus one can construct canonically the Hermitian vector bundles
$N^{(1,0)}_{G,J}$ etc, which further gives the canonical identification
of Hermitian vector bundles
\begin{align} \label{a6.44}
\Lambda (T^{*(0,1)}X)_B|_{X_G} = \Lambda (N^{*(0,1)}_{G,J})
\widehat{\otimes} \Lambda (T^{*(0,1)}X_G).
\end{align}

Let $q$ be the canonical orthogonal projection
\begin{align} \label{a6.45}
q: \Lambda (N^{*(0,1)}_{G,J} )
\widehat{\otimes} \Lambda (T^{*(0,1)}X_G)\otimes L_G^p\otimes E_G
\to \Lambda (T^{*(0,1)}X_G)\otimes L_G^p\otimes E_G
\end{align}
which acts as identity on $\Lambda (T^{*(0,1)}X_G)\otimes
L_G^p\otimes E_G$ and maps each $ \Lambda ^i(N^{*(0,1)}_{G,J})
\widehat{\otimes} \Lambda (T^{*(0,1)}X_G) $ $\otimes L_G^p\otimes
E_G$, $i\geq 1$, to zero.

We define
\begin{align} \label{6.45}
\sigma_p:= P_{G,p} q \pi_G i^*  P^G_p : (\Ker D_p)^G \to \Ker D_{G,p}.
\end{align}

Certainly in the K\"ahler case, $\sigma_p$ coincides with
\eqref{60.23}.

By using Theorems \ref{t0.0}, \ref{t0.1} as in the proof of 
Theorem \ref{t6.7} 
(cf. \cite{MM04c}, \cite[\S 5.5]{MM05b} for more details on 
the Toeplitz operators in the symplectic setting), we get
\begin{thm} \label{t6.9} Let $f$ be a smooth section of $\End(E)$ on $X$,
then $\mT_{f,p}= \sigma_p f \sigma_p^*: \Ker D_{G,p}  \to \Ker D_{G,p}$
is a Toeplitz operator with principal symbol
$2^{n_0/2}\frac{f^G}{h^2}(x) I_{\bC\otimes E_G}
\in \End(\Lambda (T^{*(0,1)}X_G)\otimes E_G)$.
\end{thm}
\begin{cor} \label{at6.11}
For $p$ large enough, $\sigma_p$ in \eqref{6.45} is an
isomorphism. Thus $\sigma_p$ defines a natural identification  for
`quantization commutes with reduction' in the (asymptotic)
symplectic case.
\end{cor}
\begin{proof} From Theorem \ref{t6.9} for $f=1$, we get
\begin{align} \label{a6.46}
\sigma_p  \sigma_p^* = 2^{n_0/2} P_{G,p}h^{-2} I_{\bC\otimes E_G}P_{G,p}
+ \cO(\frac{1}{p}).
\end{align}

But from the argument as \eqref{b6.38} and Theorem \ref{t0.1} for $G=1$,
we get for any $s\in \Omega^{0,\bullet}(X_G, L^p_G\otimes E_G)$, we have
\begin{align} \label{a6.47}
\|(I_{\bC\otimes E_G}P_{G,p}-P_{G,p})s\|_{L^2}
\leq \frac{C}{\sqrt{p}}\|s\|_{L^2}.
\end{align}
Thus for $p$ large enough, $\sigma_p\sigma_p^* $ is an isomorphism.
Thus $\sigma_p$ is surjective.

In view of \eqref{a0.3}, $\sigma_p$ in \eqref{6.45} is an isomorphism.
\end{proof}

\begin{rem} \label{at6.14} If we replace the condition $\bJ=J$
by \eqref{g2}, then the canonical map $\sigma_p$ in \eqref{6.45}
is still well defined. From the argument here, we  still 
know that $\sigma_p$ is an
isomorphism for $p$ large enough,.
\end{rem}

\subsection{Generalization to non-compact manifolds}\label{s6.5}

In this Subsection, let $(X,\omega)$ be a symplectic manifold, and
$(L,\nabla ^L)$ (resp. $(E,\nabla ^E)$) be Hermitian line (vector)
bundle on $X$, and the compact connected Lie group $G$ acts on $X$
as in Introduction, especially,
 $\omega=\frac{\sqrt{-1}}{2\pi} R^L $. But we only suppose
that $(X, g^{TX})$ is a complete manifold.

If $G=1$, these kind results were studied in \cite[\S 3]{MM04a}.

By the argument in Section \ref{s3.0}, if the square of the spin$^c$
Dirac operator $D^2_p$ has a spectral gap as in \eqref{main1},
then we can localize our problem and get a version of Theorems \ref{t0.0},
 \ref{t0.1} from Section \ref{s3.2}. In particular, if the geometric
data on $X$ verify the bounded geometry, then $D^2_p$ verify
the spectral gap \eqref{main1}.

We explain in more detail now.

We suppose

i) The tensors $R^E, r^X, \tr[R^{T^{(1,0)}X}]$
are uniformly bounded with respect on  $(X,g^{TX})$.

ii) There exists $c>0$ such that
\begin{align}\label{6.46}
\sqrt{-1} R^L(\cdot, J\cdot)\geq c g^{TX}(\cdot, \cdot).
\end{align}

\comment{
i) The tensors $R^L, R^E, J R^{TX}, \mu ^E, \mu^{TX}$ verify
the condition of the bounded geometry, i.e. the norms of these tensors
and their derivatives are bounded with respect to $g^{TX}$.

ii) There exists $c>0$ such that
\begin{align}\label{6.46}
\sqrt{-1} R^L(\cdot, J\cdot)\geq c g^{TX}(\cdot, \cdot).
\end{align}

iii) Set $X_\var = \mu^{-1} (B^g(0,\var))$. There exists
$\delta_0>0$ such that for any $0<\var \leq \delta_0$, there
exists $\var'>0$ such that $d^X(x,y)>4\var '$ for any $x\in
X_{\var}, y\in X\setminus X_{2\var}$, and $\exp^X_y: \{Z\in T_y X,
|Z|<\delta_0\}\to X$ defines a diffeomorphism onto a neighborhood
of $y$ in $X$.
}

\begin{rem} \label{at6.12} For the operator 
$D_p=\sqrt{2}(\ov{\partial}^{L^p\otimes E}+ \ov{\partial}^{L^p\otimes E,*}) $
in the holomorphic case, the above condition i) can be replaced by 

i')  The tensors $R^E, R^{T^{(1,0)}X}, \partial T$ is uniformly bounded 
 with respect to $(X,g^{TX})$, here $T$ is the torsion of $(X,\omega)$ 
as in \cite[\S 3.5]{MM04a}.
\end{rem}

Then by the argument in \cite[p. 656]{MM02} (cf.  \cite[\S 3]{MM04a}),
we know that Theorem \ref{t2.1}
still holds. Thus Theorem \ref{ta2} still holds.

Let $P^G_p$ be the orthogonal projection from $L^2(X, E_p)$ onto
$(\Ker D_p)^G$, and $P^G_p(x,x')$ $(x,x\in X)$ be its kernel as in
Def. \ref{ta0}.

Note that $\Ker D_p$ and $(\Ker D_p)^G$ need not be finite dimensional.

By the proof of Prop. \ref{0t3.0}, we know that for any compact
set $K\subset X$, $l,m\in \bN$, there exist $C_{l,m}(K)>0$ such
that for $p\geq C_L/\nu$,
\begin{align}\label{6.47}
|\wi{F}(\mL_p)(x,x') - P_p^G(x,x')|_{\cC^m(K\times K)} \leq C_{l,m}(K) p^{-l}.
\end{align}

By the proof of Theorem \ref{t0.0}, we get
\begin{thm} \label{t6.10} For any compact set $K \subset X$,
 $0<\var_0 \leq \delta_0$, $l,m \in \bN$,
there exists $C_{l,m}>0$ (depend on $K$, $\var$)
such that for $p\geq 1$,  $x,x'\in K, d^X(Gx,x')\geq \var_0$
or   $x,x'\in (X\setminus X_{2\var_0})\cap K $,
\begin{align}\label{6.48}
|P^G_p(x, x')|_{\cC^m} \leq C_{l,m} p^{-l}.
\end{align}
\end{thm}


From Section \ref{s3.2}, we get Theorem \ref{t0.1}, but now the
norm $\cC^{m'}(X_G)$ in \eqref{0.8} should be replaced by $\cC
^{m'}({ K})$ for any compact set ${  K}\subset X_G$.

One interesting case of the above discussion is when $P=\mu^{-1}(0)$ is
 compact,
 by the same argument as in
Theorems \ref{t6.7}, \ref{t6.9}, we can prove a version of Section \ref{s6.4}.

In fact, when $X=\bC^n, G= \bT^{n_0}$,
$L$ is the trivial line bundle with the metric $|1|_{h^L}(Z)= e^{- |z|^2}$,
 the Toeplitz operator type properties was studied
by Charles \cite{Charles04}.

Another interesting case  is a version of Theorem \ref{t0.1} for
covering manifolds.

Let $\tx$ be a para-compact smooth manifold, such that there is a
discrete group $\Gamma$ acting freely on $\tx$ with a compact
quotient $X=\tx/\Gamma$.

Let $\pi_{\Gamma}:\tx\longrightarrow X$ be the projection. Assume
that all above geometric data on $X$ can be lift on $\tx$. We
denote by $\wi{\bJ}$, $g^{T \tx}$, $\tom$, $\tj$, $\tl$, $\te$
 the pull-back of the corresponding
objects in Section \ref{s0} by the projection $\pi_{\Gamma}:\tx
\to X$, moreover, we assume that the $G$-action and the
$\Gamma$-action on them commute.

By the above arguments (cf. \cite[Theorems 4.4 and 4.6]{MM02}),
there exists a spectral gap for the square of the spin$^c$ Dirac
operator $\wi{D}_p$ on $\tx$.

By the finite propagation speed of solutions of hyperbolic
equations
 \eqref{b25}, we get an extension of
\cite[Theorem 3.13]{MM04a} where $G=1$.

\begin{thm}\label{t6.11} We fix
$0<\var_0 < \inf_{x\in X}\{\text{injectivity radius of $x$}\}$.
For  any compact set $K\subset\tx$ and $k,l\in \bN$, there exists
$C_{k,\,l,\,K}>0$ such that for $x,x'\in K$, $p\in \bN$,
\begin{equation}\label{6.49}
\begin{split}
&\left|\widetilde{P}^G_{p}(x,x') -
P^G_{p}(\pi_{\Gamma}(x),\pi_{\Gamma}(x'))\right|_{\cC ^l(K\times
K)}
\leqslant C_{k,\,l,\,K}\: p^{-k-1}\, ,\quad  {\rm if}\,\, d^X(x,x')< \var_0,\\
&\left|\widetilde{P}^G_{p}(x,x')\right|_{\cC ^l(K\times K)}
\leqslant C_{k,\,l,\,K}\: p^{-k-1}\, ,\quad
{\rm if}\,\,  d^X(x,x')\geqslant\var_0.
\end{split}
\end{equation}
Especially, $\widetilde{P}^G_{p}(x,x)$ has the same asymptotic expansion as
$P^G_{p}(\pi_{\Gamma} (x),\pi_{\Gamma} (x))$ in Corollary \ref{t0.5}
 on any compact set $K\subset\tx$.
\end{thm}

\subsection{Relation on the Bergman kernel on $X_G$}\label{s6.6}

From \eqref{b4}, if the operator $\Phi \mL_p \Phi^{-1}$
has the form $D^2_{G,p} + \Delta_N + 4\pi |\mu|^2 - 2\pi n_0$
under the splitting \eqref{a6.44}, then we will find the full
asymptotic expansion of the Bergman kernel on $X_G$ from $P^G_p(x,x')$.

In this Subsection, we suppose that $X$ is compact
 and $G$ is a torus $\bT^{n_0}=\bR ^{n_0}/ \bZ^{n_0}$.

Let $\theta: TP \to \kg$ be a connection form for the $G$-principal bundle
$\pi: P=\mu^{-1}(0)\to X_G$ with curvature $\Theta$.
Let $T^H P=\Ker \theta \subset TP$.

Set $M=P\times \kg^*$,
${\bf q}: M\to  \kg^*$ be the natural projection
and
\begin{align}\label{6.50}
\om ^M= \pi^* \om_G +d  \langle {\bf q}, \theta\rangle
= \pi^* \om_G +\langle {\bf q}, \Theta\rangle
+ \langle d{\bf q}, \theta\rangle .
\end{align}

By the normal crossing formula \cite[Prop. 40.1]{GuSt84}, we know there exists
a symplectic diffeomorphism such that on a neighborhood $U$ of $P$,
\begin{align}\label{6.51}
\Psi_{sym} : (X,\om)\simeq (M,\om^M),
\end{align}
and under this identification, the moment map $\mu$ (cf.
\eqref{a6}) is defined by $-{\bf q}$.

From now on, we use this neighborhood of $P$ and we will choose
metrics and connections.

Let $g^{\kg}$ be the metric on $\kg$ induced by the canonical flat
metric on $\bR^{n_0}$, and $\{K_i\}$ be the canonical unitary
basis of $\bR^{n_0}$.

Now we choose $J$ an almost -complex structure on $TX$ compatible with $\om$
such that on $T^HP$ on $U$, $J$ is induced by an almost-complex structure
on $TX_G$ which is compatible with $\om_G$, and on $\kg\oplus \kg^*$,
for $K\in\kg$, $JK\in \kg^*$ is defined by
$(JK, K')= \langle K,K' \rangle_\kg $ for $K' \in\kg$.

We also suppose $\Theta$ is  $J$-invariant.

 Let $g^{TX}$ be a $J$-invariant metric on $TX$ such that
\begin{align}\label{6.52}
g^{TX}= \pi^* g^{TX_G} \oplus g^\kg \oplus g^{\kg^*} \quad {\rm on} \, \, U.
\end{align}
As $g^\kg$ is a constant metric on $TY= \kg$, $\nabla ^{TY}$ is
the trivial connection on $TY$.  By \eqref{h4}, on $U$,
\begin{align}\label{6.53}
\nabla^{TP}_{U_1^H} = \nabla^{TX_G}_{U_1^H}
+\nabla ^{TY}_{U_1^H}+ S(U_1^H).
\end{align}

Let $\nabla ^{\Lambda (N^{*(0,1)}_{G,J} )}$ be the trivial
connection on the trivial bundle $\Lambda (N^{*(0,1)}_{G,J})$ (cf.
\eqref{a6.44}) on $U$, and $\nabla ^{\text{Cliff}_{X_G}}$ be the
Clifford connection on $\Lambda (T^{*(0,1)}X_G)$.

By \eqref{6.53}, under the identification \eqref{a6.44}, on $U$,
we have
\begin{align} \label{6.54}
\nabla ^{ \text{Cliff}}_{e^H_i}
&= \nabla ^{ \text{Cliff}_{X_G}}_{e^H_i}\otimes {\rm Id}
+ {\rm Id}\otimes \nabla ^{\Lambda (N^{*(0,1)}_{G,J} )} _{e^H_i}
+ \frac{1}{2} \langle S(e^H_i)e^H_j, K_l  \rangle c(e^H_j)c(K_l )\\
& = \nabla ^{ \text{Cliff}_{X_G}}_{e^H_i}\otimes {\rm Id}
+ {\rm Id}\otimes \nabla ^{\Lambda (N^{*(0,1)}_{G,J} )} _{e^H_i}
+ \frac{1}{4} \langle \Theta (e_i,e_j), K_l\rangle c(e^H_j)c(K_l ).\nonumber
\end{align}
However, the last term does not preserve
 $\Lambda (T^{*(0,1)}X_G)$ and $\Lambda (N^{*(0,1)}_{G,J} )$.

 From \eqref{b4} and \eqref{6.54}, in general, $\Phi \mL_p \Phi^{-1}$
will not preserve  $\Lambda (T^{*(0,1)}X_G)$ and $\Lambda (N^{*(0,1)}_{G,J})$
if $ \Theta$ is not null.

Now, we suppose that $ \Theta=0$ on $X_G$.

In this situation, on $B=U/G\subset X_G\times \kg^*$, by
\eqref{b4}, we have
\begin{align} \label{6.55}
\Phi \mL_p \Phi^{-1} = D^2_{G,p}
-\sum_l (\nabla ^{\Lambda (N^{*(0,1)}_{G,J} )}_{K_l})^2
+ 4\pi ^2 |{\bf q}|^2 -2n_0 \pi.
\end{align}

By Theorem \ref{t0.1},  Section \ref{s3.11} and \eqref{g16}, we
know that the asymptotic expansion of the Bergman kernel
 has the following relation
for $(x,Z^\bot)\in N_{G,x},\ (x',Z^{\prime\bot})\in N_{G,x'}$,
\begin{align} \label{6.56}
P^G_p((x,Z^\bot), (x',Z^{\prime\bot}))
= P_{G,p} (x,x') P_{\cL}(Z^\bot,Z^{\prime\bot}) +\mO (p^{-\infty}).
\end{align}

\newpage


\section{Computing the coefficient $\Phi_1$ and $P^{(2)}(0,0)$}\label{s4}

In this Section, $(X,\om,J)$ is a compact
 K\"ahler manifold,
$g^{TX}$ is a $G$-invariant Riemannian metric on $TX$ which is
compatible with $J$.
 $(E,h^E)$, $(L,h^L)$ are holomorphic
 Hermitian vector bundles on $X$,
and $\nabla ^E,\nabla ^L$ are the holomorphic  Hermitian
connections on $(E,h^E)$, $(L,h^L)$. Moreover,
$$\frac{\sqrt{-1}}{2 \pi} R^L=\om.$$

The action of $G$ is holomorphic and $G$ acts freely on
$P=\mu^{-1}(0)$. Thus $(X_G,\om_G,J_G)$ is a compact K\"ahler
manifold.

In Sections \ref{s4.3}-\ref{s4.7}, we suppose that in \eqref{0.1},
$\bJ=J$ on a $G$-neighborhood $U$ of $P=\mu ^{-1}(0)$.

The main purpose here is to compute the coefficient $\Phi_1$ in (\ref{0.14})
and $P^{(2)}(0,0)$ in \eqref{a0.8}.

By \eqref{0.13} (cf. also Theorem \ref{tue17}),
\begin{align} \label{6.560}
\Phi_1(x_0)= \int_{Z\in N_{G,x_0}} P^{(2)}_{x_0}(Z,Z) dv_{N_G}(Z).
\end{align}

We will first compute explicitly the terms $\mO_1$ and $\mO_2$
involved in $P^{(2)}$ in \eqref{g19}, \eqref{g22.7},
 and then compute
the integration of $P^{(2)}$ along the normal spaces to $X_G$.

Sometimes the computations seem to be long  and tedious, involving
many subtle relations between metrics, connections and curvatures
near $X_G$, but fortunately the final result on $\Phi_1$ is still of a simple
form, as expected.

Throughout the  computations below, a key idea is to rewrite all
operators by using the creation and annihilation operators
$b_i,b_i^+$, $b_j^\bot, b_j^{\bot+}$, then under the help of
\eqref{g7} and Theorem \ref{t3.4}, we can do the operations and to
obtain the crucial Lemmas \ref{t4.8}, \ref{t4.9}.

To get the final simple formula \eqref{0.15}, we still need to
prove a highly non-trivial identity \eqref{4.16}.

The formula for $P^{(2)}(0,0)$ in Theorem \ref{t4.15} is quite complicate,
it involves $h$, the volume function of the orbit
 and the curvature for the principal bundle $P\to X_G$.

This Section is organized as follows. In  Section \ref{s4.3}, we
explain various relations of the curvature of the fibration $P\to
X_G$ and the second fundamental form of $P$. In Section
\ref{s4.4}, we obtain the explicit formulas for
 the operators $\mO_1$, $\mO_2$.
In  Section \ref{s4.5}, we apply the formulas in Section
\ref{s4.4} and \eqref{6.560} to \eqref{g51}, and we get a formula
for the coefficient $\Phi_1$. In  Section \ref{s4.6}, we compute
finally $\Phi_1$, thus prove Theorem \ref{t0.6}. In Section
\ref{s4.7}, we  compute $P^{(2)}(0,0)$ in Theorem \ref{t4.15}. In
Section \ref{s8.1}, we explain how to reduce the general case to
the case $\bJ=J$ which has been worked out in Sections
\ref{s4.3}-\ref{s4.6}.

In the whole Section, if there is no other specific notification,
 when we meet the operation $|\quad|^2$,
 we will first do this operation, then take the sum of the indices.

\subsection{The second fundamental form of $P$}\label{s4.3}

We use the notations in Sections \ref{s3.01}, \ref{s3.0}.
Then the normal
bundle $N_{G}$ of $X_G$ in $U/G$ is $(J TY)_G$.

Let $\iota: X_G\to U/G$ be the natural embedding.

We will apply the notation in Section \ref{s4.2} to $B=U/G$.

 Let $\nabla ^{TX_G}$, $\nabla ^{N_{G}}$ be connections on
$TX_G,N_{G}$ on $X_G$ induced by projections of $\nabla
^{TB}|_{X_G}$. Then $\nabla ^{TX_G}$ is the Levi-Civita connection
on $(TX_G, g^{TX_G})$.

Let
\begin{align}\label{3.009}{^0\nabla} ^{TB}= \nabla
^{TX_G} \oplus\nabla ^{N_{G}}
\end{align}
 be the connection on $TB$ on $X_G$ induced by
$\nabla ^{TX_G}$, $\nabla ^{N_{G}}$ with curvature $^0R^{TB}$.

Set
\begin{align}\label{3.10}
A=\nabla ^{TB}|_{X_G}- {^0\nabla} ^{TB}.
\end{align}
Then $A$ is a $1$-form on  $X_G$ taking values in the skew-adjoint
endomorphisms of $(TB)|_{X_G}$ which exchange $TX_G$ and $N_{G}$.

We recall the following properties of $R^{TB}$: for
$U,V,W,W_2\in TB$,
\begin{equation}\label{a3.15}
\begin{split}
&\left \langle R^{TB}(U,V)W,W_2 \right \rangle
=\left \langle R^{TB}(W,W_2) U,V\right \rangle,\\
&R^{TB}(U,V)W+ R^{TB}(V,W)U+ R^{TB}(W,U)V=0.
\end{split}
\end{equation}

On $X_G$, let  $\{e^0_i\}$ be an orthonormal frame of $TX_G$,
let $\{e ^\bot_j\} $ be an orthonormal frame  of $N_{G}$,
then $\{e_i\}= \{e^0_i,e ^\bot_j\}$ is an orthonormal frame of $TB$.

The following result gives detail informations on the torsion $T$
of the fibration, as well as  the second fundamental form $A$.

\begin{thm} \label{t4.2} On $P$, the restriction of the tensor
$\left\langle JT(\cdot, J \cdot), \cdot \right\rangle$
on $(N_{G})^{\otimes \, 3}$ is symmetric,
 and
\begin{subequations}
\begin{align}
&(A(e^0_i)e^0_j)^H= \frac{1}{2} J T(e^{0,H}_i, Je^{0,H}_j), \label{3.11a}\\
&T(e^{0,H}_i, e^{0,H}_j)=T((J_Ge^0_i)^H, (J_Ge^0_j)^H),\label{3.11b}\\
&T(e^{0,H}_i, e^{\bot,H}_j)= 2 T((J_Ge^0_i)^H, Je^{\bot,H}_j),\label{3.11c}\\
&\left\langle T(e^{0,H}_i, e^{\bot,H}_j), Je^{\bot,H}_k\right\rangle
= \left\langle T(e^{0,H}_i, e^{\bot,H}_k),
Je^{\bot,H}_j\right\rangle ,\label{3.11d}\\
&\sum_k \left\langle T(e^{\bot,H}_k, e^{\bot,H}_j),
Je^{\bot,H}_k\right\rangle =0.\label{3.11e}
\end{align}
\end{subequations}
\end{thm}
\begin{proof} Observe first that we have
\begin{subequations}
\begin{align}\label{3.12a}
 &\nabla ^{TX}J=0;\\
&(J_G e^0_i)^H= J e^{0,H}_i \quad {\rm on } \quad  P.\label{3.12b}
\end{align}
\end{subequations}

Let $Z$ be a smooth section of $TY$, then $JZ\in  JTY\subset T^HX$
on $P$, by  (\ref{h4}), (\ref{g1}) and (\ref{3.12a}), on $P$, we
have
\begin{multline}\label{3.13}
\left\langle  J (A(e^0_i)e^0_j)^H, Z\right\rangle =-\left\langle
\nabla ^{T^HX} _{e^{0,H}_i}e^{0,H}_j, JZ\right\rangle
=- \left\langle \nabla ^{TX} _{e^{0,H}_i} e^{0,H}_j, J Z\right\rangle\\
=\left\langle \nabla ^{TX} _{e^{0,H}_i}(Je^{0,H}_j), Z\right\rangle
=\left\langle S(e^{0,H}_i) Je^{0,H}_j,  Z \right\rangle
=-\frac{1}{2} \left\langle T(e^{0,H}_i,Je^{0,H}_j),
Z\right\rangle.
\end{multline}
Thus we get (\ref{3.11a}), as $ A(e^0_i)e^0_j\in N_G$.

Note that $[Z,e^H_i]\in TY$, by (\ref{h2b}), (\ref{h4}) and
(\ref{3.12a}),
 \begin{align}\label{3.14}
\left\langle T(e^{H}_i,e^{H}_j), Z\right\rangle =2 \left\langle
\nabla ^{TX} _{e^{H}_i} Z, e^H_j\right\rangle =2 \left\langle
\nabla ^{TX}_{Z}e^{H}_i, e^H_j\right\rangle = 2 \left\langle
\nabla ^{TX}_{Z}(Je^{H}_i), Je^H_j\right\rangle .
\end{align}

From (\ref{3.12b}) and (\ref{3.14}), we get (\ref{3.11b}).

From (\ref{h4}), (\ref{3.14})
and $J e^{\bot,H}_j, J e^{\bot,H}_k\in TY$ on $P$, we get
\begin{align}\label{3.15}
\left\langle T(e^{0,H}_i,e^{\bot,H}_j), Z\right\rangle = 2
\left\langle S(Z) (Je^{0,H}_i),Je^{\bot,H}_j \right\rangle = 2
\left\langle T(Je^{0,H}_i,Je^{\bot,H}_j), Z\right\rangle .
\end{align}
Thus we get  (\ref{3.11c}).
By  (\ref{h3}), (\ref{3.15}), we get
 \begin{align}\label{3.16}
\left\langle T(e^{0,H}_i,e^{\bot,H}_j), J e^{\bot,H}_k\right\rangle
=2 \left\langle T(Je^{0,H}_i,Je^{\bot,H}_j), J e^{\bot,H}_k\right\rangle
= \left\langle T(e^{0,H}_i,e^{\bot,H}_k), J e^{\bot,H}_j\right\rangle.
\end{align}
Thus we get  (\ref{3.11d}).
 By (\ref{h4}), (\ref{3.12a}) and $Je^{\bot,H}_j\in TY$ on $P$,
\begin{multline}\label{3.17}
\left\langle T(e^{\bot,H}_i, Je^{\bot,H}_j), Je^{\bot,H}_k\right\rangle
=\left\langle  \nabla ^{TX}_{Je^{\bot,H}_k}e^{\bot,H}_i,
J e^{\bot,H}_j\right\rangle\\
=-\left\langle \nabla ^{TX}_{Je^{\bot,H}_k}(J e^{\bot,H}_i),
e^{\bot,H}_j\right\rangle
=\left\langle  \nabla ^{TX}_{Je^{\bot,H}_k}e^{\bot,H}_j,
J e^{\bot,H}_i\right\rangle
= \left\langle T(e^{\bot,H}_j, Je^{\bot,H}_i), Je^{\bot,H}_k\right\rangle .
\end{multline}

By (\ref{h3}) and  (\ref{3.17}),
 $\left\langle JT(\cdot, J \cdot), \cdot \right\rangle$ is symmetric
on the horizontal lift of $N_{G}^{\otimes \, 3}$.

Note that $\{Je^{\bot,H}_k\}$ is a  $G$-invariant orthonormal frame
 of $TY$ on $P$, by \eqref{3.14},
 \begin{align}\label{a3.14}
\left\langle T(e^{\bot,H}_i,e^{\bot,H}_j), Je^{\bot,H}_k \right\rangle
= 2 \left\langle \nabla ^{TY}_{Je^{\bot,H}_k}(Je^{\bot,H}_i),
 Je^{\bot,H}_j\right\rangle .
\end{align}

By \eqref{bh4} and \eqref{a3.14}, we get \eqref{3.11e}. The proof
of Theorem \ref{t4.2} is complete.
\end{proof}

\begin{rem} \label{r4.3}
From (\ref{h3}) and (\ref{3.11b}), $\Theta|_{X_G}$ is a
$(1,1)$-form on $X_G$. Especially, for any complex representation
$V$ of $G$, $P\times_G V$ is a holomorphic vector bundle on $X_G$.
Moreover, by \eqref{3.11a}, for $U\in TX_G, V\in N_G$, we have at $x_0$,
\begin{align} \label{33.36}
A(U)V= \left\langle A(U)V, e^0_j\right\rangle e^0_j
=-\left\langle V, A(U)e^0_j\right\rangle e^0_j
= \frac{1}{2} \left\langle T(U,Je^0_j), JV\right\rangle e^0_j.
\end{align}
For $x_0\in X_G$, if $\{e^\bot_j\}$ is a fixed orthonormal basis of
$N_{G,x_0}$ as above, then for $U\in T_{x_0}X_G$, we will denote by
 \begin{align}\label{34.21}
\mT_{ijk}=  \left\langle JT(e^\bot_i, J e^\bot_j), e^\bot_k \right\rangle,
\, \,
\wi{\mT}_{ijk}=  \left\langle JT(e^\bot_i, e^\bot_j), e^\bot_k \right\rangle,
 \, \,  \mT_{jk}(U)= \left\langle JT(U,e^\bot_j), e^\bot_k \right\rangle .
\end{align}
By Theorem \ref{t4.2}, $\mT_{ijk}$ is symmetric on $i,j,k$ and
 $\mT_{jk} \in T_{x_0}^*X_G$ is  symmetric on $j,k$,
$\wi{\mT}_{ijk}$ is anti-symmetric on $i,j$.
\end{rem}

\begin{rem}\label{r4.4} From Remark \ref{p4.02} and \eqref{a3.14},
we know that
$\left\langle JT(\cdot,  \cdot), \cdot \right\rangle$
is anti-symmetric on $(N_{G})^{\otimes \, 3}$ if $g^{TY}$ is induced by
a family of $Ad$-invariant metric on $\kg$.
If $G$ is abelian, then by \eqref{bh7}, \eqref{a3.14},
$T(\cdot,  \cdot)=0$ on $(N_{G})^{\otimes \, 2}$, thus $\wi{\mT}_{ijk}=0$.
\end{rem}

\comment{Finally, we recall the following properties of  $R^{TB}$: for
$U,V,W,W_2\in TB$,
\begin{equation}\label{a3.15}
\begin{split}
&\left \langle R^{TB}(U,V)W,W_2 \right \rangle
=\left \langle R^{TB}(W,W_2) U,V\right \rangle,\\
&R^{TB}(U,V)W+ R^{TB}(V,W)U+ R^{TB}(W,U)V=0.
\end{split}
\end{equation}

By \eqref{a3.10}, and  \eqref{3.11a},
\begin{multline}\label{a3.13}
\left \langle R^{TB}(e^0_j,e^0_i)e^0_k,e^0_l\right \rangle
= \left \langle R^{TX_G}(e^0_j,e^0_i)e^0_k,e^0_l\right \rangle
-  \left \langle A(e^0_i)e^0_k,A(e^0_j)e^0_l \right \rangle
+ \left \langle A(e^0_j)e^0_k,A(e^0_i)e^0_l \right \rangle\\
=  \left \langle R^{TX_G}(e^0_j,e^0_i)e^0_k,e^0_l\right \rangle
- \frac{1}{4}\left \langle JT(e^0_i, J e^0_k),
JT(e^0_j, J e^0_l)\right \rangle
+ \frac{1}{4}\left \langle JT(e^0_j, J e^0_k),
JT(e^0_i, J e^0_l)\right \rangle.
\end{multline}
}

\subsection{Operators $\mO_1$,  $\mO_2$ in (\ref{c30})}\label{s4.4}

We use the notation in Sections \ref{s3.2}, \ref{s3.10}, and all
tensors will be evaluated at $x_0\in X_G$.

Recall that $(X,\omega)$ is K\"ahler and ${\bJ}=J$ on a
$G$-neighborhood $U$ of $P=\mu^{-1}(0)$, then in (\ref{g4})
\begin{align} \label{b3.11}
a_i=a^\bot_j=2\pi.
 \end{align}

Clearly, on $U$,  the Levi-Civita connection $\nabla ^{TX}$ preserves
$T^{(1,0)}X$ and $T^{(0,1)}X$, and $\nabla ^{T^{(1,0)}X} =
P^{T^{(1,0)}X}\nabla ^{TX} P^{T^{(1,0)}X}$ is the holomorphic
Hermitian connection on $T^{(1,0)}X$, while the Clifford
connection $\nabla ^{{\rm Cliff}}$ on $\Lambda (T^{*(0,1)}X)$ is
$\nabla ^{ \Lambda (T^{*(0,1)}X)}$, the natural  connection
induced by $\nabla ^{T^{(1,0)}X}$.

Let $\overline{\partial}^{L^p\otimes E,*}$  be the canonical formal adjoint
of the Dolbeault operator $\overline{\partial} ^{L^p\otimes E}$
 on  $\Omega ^{0,\bullet}(U, L^p\otimes E)$.
Then  the operator $D_p$ in (\ref{defDirac}) is
\begin{align} \label{b3.12}
D_p = \sqrt{2}\left( \overline{\partial} ^{L^p\otimes E}+
\overline{\partial} ^{L^p\otimes E,*}\right).
\end{align}

 Note that $D_p^2$ preserves the $\bZ$-grading of
  $\Omega ^{0,\bullet}(U, L^p\otimes E)$.

Set
\begin{align}\label{b3.13}
D^2_{p,i}= D^2_p|_{\Omega ^{0,i}(U, L^p\otimes E)}.
\end{align}

Since $\nabla ^{{\rm Cliff}}$ preserves the $\bZ$-grading of
 $\Lambda (T^{*(0,1)}X)$, the operator $\cL^t_2$ in \eqref{c27}
also preserves the $\bZ$-grading on $\Lambda (T^{*(0,1)}X_0)$.
Moreover, $\cL^t_2$ is invertible on $\oplus_{q=1}^n
\Omega^{0,q}(X_0, L^p_0\otimes E_0)$ for $t$ small enough.

From Section \ref{s3.11}, for $P^{(r)}$ in \eqref{0.7},
 \be \label{3.19}
P^{(r)}= I_{\bC\otimes E_G}P^{(r)} I_{\bC\otimes E_G}.
 \ee
Thus we only need
to do the computation for $D^2_{p,0}$.

In what follows, we compute everything on $\cC ^\infty (U,
L^p\otimes E)$.

Take $x_0\in X_G$.

If $Z\in T_{x_0}B$, $Z=Z^0+Z^\bot$, $Z^0\in T_{x_0}X_G$,
$Z^\bot\in N_{x_0}$, $|Z^0|,|Z^\bot|\leq \var$, in Section
\ref{s3.2}, we identify $Z$ with $\exp^B_{\exp_{x_0}^{X_G}
(Z^0)}\tau_{Z^0} (Z^\bot)$. This identification is a
diffeomorphism from $B^{TX_G}_{x_0}(0,\var)\times
B^N_{x_0}(0,\var)$ into an open neighborhood $\cU(x_0)$ of $x_0$
in $B$, we denote it by $\Psi$. Then $\cU(x_0)\cap X_G=
B^{TX_G}_{x_0}(0,\var) \times \{0\}$.

In what follows, we use indifferently the notation
$B^{TX_G}_{x_0}(0,\var)\times B^N_{x_0}(0,\var)$ or $\cU(x_0)$,
$x_0$ or $0$, $\cdots$.

From now on, we replace $U/G$ by $\bR^{2n-n_0} \simeq T_{x_0}B$ as
in Section \ref{s3.2}, and we use the notation therein.
Especially,
\be \label{3.20}
\nabla_t= t S^{-1}_t \kappa ^{1/2}
\nabla ^{(L^p\otimes E)_B} \kappa ^{-1/2} S_t,
\ee
and $\mO_r$ in (\ref{c30}) takes value in $\End (E_B)$.

\comment{
If $U\in T_{x_0}X_G$, let $t\in \bR\to
x_t=\exp^{X_G}_{x_0}(tU)\in X_G$ be the geodesic in $X_G$ such
that $x_t|_{t=0}=x_0$, $\frac{dx}{dt}|_{t=0}=U$. If $U\in
T_{x_0}X_G$, $|U|\leq \var$, $V\in N_{x_0}$, let $\tau_U V\in
N_{\exp^{X_G}_{x_0}(U)}$ be the natural parallel transport
 of $V$ with respect to the connection $\nabla^{N}$ along
the curve $t\in [0,1]\to \exp_{x_0}^{X_G}(tU)$.
}

Let $\{e^0_i\}$, $\{e^\bot_j\}$ be orthonormal basis of
$T_{x_0}X_G$, $N_{G,x_0}$ respectively. We will also denote
$\Psi_*(e^0_i),\, \Psi_*(e^\bot_j)$ by
$e^0_i,\ e^\bot_j$.

Let $\{e_i\}$ denote the basis $\{e_i^0,\, e_j^\perp\}$. Thus in
our coordinates,
\begin{align}\label{a33.16}
\tfrac{\partial}{\partial Z^0_i}=e^0_i,\ \ \ \
\tfrac{\partial}{\partial Z^\bot_j}=e^\bot_j.
\end{align}

We denote by $(g^{ij}(Z))$ the inverse of the matrix
$(g_{ij}(Z))=(g^{TB}_{ij}(Z))$.

Recall $\Gamma _{ij}^l$ is the connection form of $\nabla ^{TB}$,
with respect to the frame $\{e_i\}$, defined in \eqref{ac37}.
Also recall that $\mR$, $\mR^0$ and $\mR^\perp$ are defined in \eqref{2c15}.

As in \eqref{ah7}, the moment map $\mu$ induces a $G$-invariant
$\cC^\infty$ section $\wi{\mu}$ of $TY$ on $U$.

Note also that  by \eqref{c24}, $R^{E}_\tau \in \End(E)$ defines a
section of $\End(E_B)$ on $B=U/G$.

 Set
\begin{align}\label{a3.16}
 \cL^t_3(Z) &= - g^{ij}(tZ) \Big(\nabla_{t,e_i}\nabla_{t,e_j}
 -t \Gamma_{ij}^k(tZ) \nabla_{t,e_k}\Big)\\
& +t^2 \Big(\frac{1}{h}g^{ij} (\nabla_{e_i}\nabla_{e_j}h
 - \Gamma_{ij}^k \nabla_{e_k}h)\Big)(tZ)
- t^2 R^{E}_\tau(tZ) -2\pi n .\nonumber
\end{align}

By (\ref{b4}), (\ref{c27}) and (\ref{a3.16}), we can reformulate
(\ref{0c35}), (\ref{1ue1}), in using the notations in (\ref{g8}),
as follows,
\begin{align}\label{a3.17}
\begin{split}
&\nabla_{0,\cdot}= \nabla_\cdot + \frac{1}{2} R^{L_{B}}_{x_0}  (\mR,\cdot)
= \nabla_\cdot
 -\pi \sqrt{-1}\left\langle J_{x_0}Z^0,\cdot\right\rangle_{x_0},\\
&\cL^0_2=\sum_{j=1}^{n-n_0} b_jb^+_j+\sum_{j=1}^{n_0}b^{\bot}_jb^{\bot +}_j
= -\sum_j (\nabla_{0,e_j})^2 +4\pi ^2 |Z^\bot|^2 -2\pi n, \\
 &\cL^t_2(Z) = \cL^t_3(Z)
+ 4\pi ^2 \Big|\frac{1}{t}\wi{\mu}\Big|^2_{g^{TY}}(tZ)
- \left \langle 4\pi\sqrt{-1} \wi{\mu}
+t^2 \wi{\mu}^E,\wi{\mu}^E\right \rangle_{g^{TY}}(tZ). 
\end{split}\end{align}

If there is no another specification, we will evaluate our tensors
at $x_0$, and most of time, we will omit the subscript $x_0$.

Set $h_0=h_{x_0}$, and for $U\in T_{x_0}B$, set
\begin{align}\label{a3.18}
\begin{split}
&B(Z,U)= \frac{1}{2}\sum_{|\alpha|=2}(\partial ^{\alpha}R^{L_B})_{x_0}
\frac{Z^\alpha}{\alpha !}(\mR,U),\\
&I_1= - B(Z,e_i)\nabla_{0,e_i} -\frac{1}{2} \nabla_{e_i}(B(Z,e_i)),\\
&I_2= \Big(\left \langle \frac{1}{3}R^{TX_G} (\mR^0,e^0_i)\mR^0
+ \nabla_{\mR^0}^{TX_G}(A(e^0_i)\mR^\bot),e^0_j\right \rangle
+\left \langle e^0_i, \nabla_{\mR^0}^{TX_G}(A(e^0_j)\mR^\bot)\right \rangle\\
&\hspace*{5mm}
-3 \left \langle A(e^0_i)\mR^\bot, A(e^0_j)\mR^\bot\right \rangle
+ \left \langle R^{TB}(\mR^\bot,e^0_i)\mR^\bot, e^0_j\right \rangle\Big)
\nabla_{0,e^0_i} \nabla_{0,e^0_j} \\
&\hspace*{5mm}
+ \Big(\left \langle R^{N_G}(\mR^0,e^0_j)\mR^\bot, e^\bot_i\right \rangle
+ \frac{4}{3}\left \langle R^{TB}(\mR^\bot,e^0_j)\mR^\bot,
e^\bot_i\right \rangle\Big) \nabla_{0,e^\bot_i} \nabla_{0,e^0_j}\\
&\hspace*{5mm}+  \frac{1}{3}\left \langle R^{TB}(\mR^\bot,e^\bot_i)\mR^\bot,
e^\bot_j\right \rangle \nabla_{0,e^\bot_i} \nabla_{0,e^\bot_j}.
\end{split}
\end{align}

Recall that the operator $\cL$ has been defined in (\ref{g8}).

Set also
\begin{align}\label{b3.15}
\begin{split}
\Gamma_{i i}(\mR) =& \frac{2}{3} R^{TX_G}_{x_0} (\mR^0, e^0_i) e^0_i
+ \nabla^{TB}_{\mR^0}(A(e^0_i) e^0_i) + R^{TB}( \mR^\bot, e^0_i)e^0_i\\
 &+ A(e^0_i) A(e^0_i)\mR^\bot + \nabla^{TX_G}_{ e^0_i}(A(e^0_i) \mR^\bot)
-  A(\mR^0)A(e^0_i) e^0_i,\\
K_2(\mR)= &
 \frac{1}{3} \left \langle R^{TX_G}(\mR^0, e^0_i)\mR^0, e^0_i\right \rangle
+ \left \langle R^{TB}(\mR^\bot, e^0_i)\mR^\bot, e^0_i\right \rangle\\
&+ \frac{1}{3} \left \langle R^{TB} (\mR^\bot, e^\bot_i)\mR^\bot,
e^\bot_i\right \rangle
+2 \Big(\sum_i \left \langle A(e^0_i)e^0_i, \mR^\bot\right \rangle\Big)^2\\
&-  | A(e^0_i)\mR^\bot|^2
+ 2\left\langle \nabla^{TX_G}_{\mR^0}(A(e^0_i)\mR^\bot), e^0_i\right\rangle.
\end{split}\end{align}

\begin{lemma}\label{t4.4} There exist second order
differential operators $\mO^\prime_r$
as in Theorem \ref{t3.3} such that for $|t|\leq 1$,
\begin{align}\label{a3.19}
\cL^t_3 =\cL^0_3 + \sum_{r=1}^m t^r \mO^\prime_r +\cO(t^{m+1}),
\end{align}
with
\begin{align}\label{a3.20}
\cL^0_3=& \cL -\sum_{j=1}^{n_0} (\nabla_{e^\bot_j})^2 -2\pi n_0
= \cL^0_2 -4\pi ^2 |Z^\bot|^2,\\
\mO^\prime_1 =& -\frac{2}{3} ( \partial_j R^{L_B})_{x_0} (\mR,e_i)Z_j  \nabla_{0,e_i}
 -\frac{1}{3} (\partial_i R^{L_B})_{x_0} (\mR, e_i)
- 2 \left \langle A(e^0_i)e^0_j,\mR^\bot \right \rangle
\nabla_{0,e^0_i} \nabla_{0,e^0_j}, \nonumber\\
\mO^\prime_2 =&I_1+I_2+\Big[\frac{1}{4} K_2(\mR)-\frac{3}{8}\Big(\sum_l \left \langle A(e^0_l)e^0_l,
\mR^\bot\right \rangle\Big)^2, \cL_2^0\Big]
\nonumber\\
&- 2 \left \langle A(e^0_i)e^0_j,\mR^\bot \right \rangle
 \Big(\frac{2}{3} ( \partial_k R^{L_B})_{x_0} (\mR,e^0_j)Z_k  \nabla_{0,e^0_i}
 +\frac{1}{3} (\partial^0_j R^{L_B})_{x_0} (\mR, e^0_i) \Big)\nonumber\\
&+ \left \langle \Gamma_{i i}(\mR), e_j\right \rangle  \nabla_{0,e_j}
- \frac{1}{2}\left \langle A(e^0_l)e^0_l, \mR^\bot\right \rangle
\nabla_{A(e^0_k)e^0_k}
+2 \left \langle A(e^0_i)e^0_j, \mR^\bot\right \rangle \nabla_{A(e^0_i)e^0_j}
\nonumber\\
&+ \frac{2}{3} \left \langle R^{TB}(\mR^\bot, e^\bot_i)e^\bot_i,
 e_j\right \rangle  \nabla_{0,e_j}
- R^{E_B}_{x_0}(\mR,e_i)\nabla_{0,e_i}
-R^{E_B}_{\tau,x_0}\nonumber\\
&-\frac{1}{9}\sum_{i}\Big[\sum_j(\partial_j R^{L_B})_{x_0} (\mR,e_i)Z_j\Big]^2
+\frac{1}{h_0} (\nabla_{e_j}\nabla_{e_j} h
-\nabla_{A(e^0_i)e^0_i}h )_{x_0} .\nonumber
\end{align}
\end{lemma}
\begin{proof} By (\ref{0c39}) and (\ref{3.20}),
\begin{multline}\label{3.28}
\nabla_{t,e_i}= \kappa ^{1/2}(tZ)\Big(\nabla_{e_i}
+ \Big(\frac{1}{2} R^{L_{B}}_{x_0}
+ \frac{t}{3} (\partial _k R^{L_{B}})_{x_0}  Z_k \\
+ \frac{t^2}{4} \sum_{|\alpha|=2} (\partial ^\alpha R^{L_{B}})_{x_0}
\frac{Z^\alpha}{\alpha !}+ \frac{t^2}{2} R^{E_{B}}_{x_0} \Big) (\mR, e_i)
+ \cO(t^3)\Big) \kappa ^{-1/2}(tZ).
\end{multline}

To get \eqref{a3.20}, we could use \eqref{0c31}-\eqref{1c34},
while here we will get it directly from the local computation.

By \cite[Prop. 1.28]{BeGeVe} (cf. \cite[(1.31)]{MM04a}) and
\eqref{0c39},
\begin{align}\label{3.21}
\begin{split}
&\left \langle e^0_i,e^0_j\right \rangle_{Z^0} = \delta_{ij}
+\frac{1}{3}\left \langle R^{TX_G}_{x_0} (\mR^0,e^0_i) \mR^0,
e^0_j\right \rangle_{x_0} + \cO(|Z^0|^3),\\
&(\nabla^{N_G}_{e^0_k}\nabla^{N_G}_{e^0_i}e^\bot_j)_{x_0}
= \frac{1}{2} R^{N_G}_{x_0}(e^0_k,e^0_i)e^\bot_j.
\end{split}\end{align}
Moreover, for $W,V\in N_{x_0}$, $\gamma_s(t)=(Z^0, t(W+sV))$ is a
family of geodesics  from $(Z^0,0)$.

Set $Y=\frac{\partial}{\partial t}\gamma_s(t)$,
$X(\gamma_s(t))=\frac{\partial}{\partial s}\gamma_s(t)=tV$.

Since $\nabla^{TB}_Y Y=0,\ \nabla^{TB}_Y X- \nabla^{TB}_X Y= [Y,X]
= \gamma_*[\frac{\partial}{\partial t},\frac{\partial}{\partial
s}]=0$, we get
\begin{align}\label{a3.22}
0=\nabla^{TB}_X\nabla^{TB}_Y Y
=\nabla^{TB}_Y\nabla^{TB}_Y X-R^{TB}(Y,X)Y.
\end{align}

Take $V=e^\bot_i$, we get at $s=t=0$,
\begin{align}\label{a3.21}
(\nabla^{TB}_W\nabla^{TB}_W e^\bot_i )_{Z^0}
= \frac{1}{3}\nabla^{TB}_Y\nabla^{TB}_Y\nabla^{TB}_Y X
=\frac{1}{3}R^{TB}( W,e^\bot_i)W.
\end{align}Under our coordinate, we have
\begin{align}\label{b3.16}
&(\nabla^{TB}_{e^\bot_j}e^\bot_i)_{x_0}
=(\nabla^{TX_G}_{e^0_i}e^0_j)_{x_0} =
(\nabla^{N_G}_{e^0_j}e^\bot_i)_{x_0}=0, \quad
(\nabla^{TB}_{e^0_i}e^0_j)_{x_0} = A_{x_0}(e^0_i)e^0_j,\nonumber\\
&(\nabla^{TB}_{e^\bot_j}e^0_i)_{x_0}
=(\nabla^{TB}_{e^0_i}e^\bot_j)_{x_0}
= A_{x_0}(e^0_i)e^\bot_j,\nonumber\\
&(\nabla^{TB}_{\mR^\bot}e^\bot_i)_Z=0,\\
&(\nabla^{TB}_{e^\bot_j}e^\bot_i)_{Z^0}
=(\nabla^{N_G}_{\mR^0}e^\bot_i)_{Z^0}
=0.\nonumber
\end{align}
Moreover, by  (\ref{a3.15}), \eqref{3.21}, \eqref{a3.21} and
\eqref{b3.16}
 (comparing with \cite[(1.31)]{MM04a}),
we have at $x_0$ that
\begin{align}\label{b3.17}
\begin{split}
&\nabla^{TB}_{e^\bot_k}\nabla^{TB}_{e^\bot_j}e^\bot_i
=\frac{1}{3} R^{TB}(e^\bot_k,e^\bot_j)e^\bot_i + \frac{1}{3} R^{TB}(e^\bot_k,e^\bot_i)e^\bot_j,\\
&\nabla^{TB}_{e^0_k}\nabla^{TB}_{e^\bot_j}e^\bot_i =0, \\
& \nabla^{TB}_{e^\bot_k}\nabla^{TB}_{e^\bot_j}e^0_i
=\nabla^{TB}_{e^\bot_k}\nabla^{TB}_{e^0_i}e^\bot_j
= R^{TB}(e^\bot_k,e^0_i)e^\bot_j,\\
&\nabla^{TB}_{e^0_k}\nabla^{TB}_{e^\bot_j}e^0_i
= \nabla^{TB}_{e^0_k}\nabla^{TB}_{e^0_i}e^\bot_j
= \nabla^{N_G}_{e^0_k}\nabla^{N_G}_{e^0_i}e^\bot_j
+ A(e^0_k)A(e^0_i)e^\bot_j + \nabla^{TX_G}_{e^0_k}(A(e^0_i)e^\bot_j)\\
&\hspace*{5mm}= \frac{1}{2} R^{N_G}(e^0_k,e^0_i)e^\bot_j
+  A(e^0_k)A(e^0_i)e^\bot_j + \nabla^{TX_G}_{e^0_k}(A(e^0_i)e^\bot_j),\\
&\nabla^{TB}_{e^\bot_j}\nabla^{TB}_{e^0_k}e^0_i
=R^{TB}(e^\bot_j, e^0_k)e^0_i
+ \nabla^{TB}_{e^0_k}\nabla^{TB}_{e^\bot_j}e^0_i, \\
&\nabla^{TB}_{e^0_k}\nabla^{TB}_{e^0_j}e^0_i
= \nabla^{TX_G}_{e^0_k}\nabla^{TX_G}_{e^0_j}e^0_i
+  \nabla^{TB}_{e^0_k}(A(e^0_j)e^0_i)\\
&\hspace*{5mm}= \frac{1}{3}R^{TX_G}(e^0_k,e^0_j)e^0_i
+ \frac{1}{3} R^{TX_G}(e^0_k,e^0_i)e^0_j
+ \nabla^{TB}_{e^0_k}(A(e^0_j)e^0_i).
\end{split}\end{align}

In the following, for a tensor $\psi$ and the covariant derivative
$\nabla^{B}$ acting on $\psi$ induced by $\nabla ^{TB}$, we denote by
$$(\nabla^{B} \nabla^{B}  \psi)_{(c_je_j,c'_ke_k)}
= c_jc'_k (\nabla^{B} _{e_j}\nabla^{B} _{e_k}\psi)_{x_0}.$$

From \eqref{b3.17}, we get at $x_0$ the following formula which
will be used in (\ref{b3.24}), (\ref{b3.25}), (\ref{b3.31}),
 (\ref{b3.32}) and (\ref{4.41}),
\begin{align}\label{b3.18}
\begin{split}
&(\nabla^{TB}\nabla^{TB} e^0_i)_{(\mR^0,\mR^0)}
= \frac{1}{3}R^{TX_G} (\mR^0,e^0_i)\mR^0
+ \nabla_{\mR^0}^{TB}(A(e^0_j)e^0_i)Z^0_j,\\
&(\nabla^{TB}\nabla^{TB} e^0_i)_{(\mR^0,\mR^\bot)}
=\frac{1}{2} R^{N_G}(\mR^0, e^0_i)\mR^\bot
+ A(\mR^0)A(e^0_i)\mR^\bot + \nabla_{\mR^0}^{TX_G}(A(e^0_i)\mR^\bot),\\
&(\nabla^{TB}\nabla^{TB} e^0_i)_{(\mR^\bot,\mR^\bot)}
=R^{TB}(\mR^\bot,e^0_i)\mR^\bot,\\
&(\nabla^{TB}\nabla^{TB} e^\bot_j)_{(\mR^0,\mR^0)}
=A(\mR^0)A(\mR^0)e^\bot_j+ \nabla_{\mR^0}^{TX_G}(A(e^0_k)e^\bot_j)Z^0_k,\\
&(\nabla^{TB}\nabla^{TB} e^\bot_j)_{(\mR^0,\mR^\bot)}=0,\\
&(\nabla^{TB}\nabla^{TB} e^\bot_j)_{(\mR^\bot,\mR^\bot)}
=\frac{1}{3} R^{TB}(\mR^\bot,e^\bot_j)\mR^\bot,\\
&(\nabla^{TB}\nabla^{TB} e_j)_{(\mR^\bot,\mR^0)}
=(\nabla^{TB}\nabla^{TB} e_j)_{(\mR^0,\mR^\bot)}
+ R^{TB}(\mR^\perp,\mR^0)e_j.
\end{split}\end{align}

\comment{\begin{align}\label{b3.18}
&(\nabla^{TB}\nabla^{TB} e^0_i)_{(\mR,\mR)}
= \frac{1}{3}R^{TX_G} (\mR^0,e^0_i)\mR^0
+\nabla_{\mR^0}^{TB}(A(e^0_j)e^0_i)Z^0_j\\
&\hspace*{10mm}+ R^{N_G}(\mR^0, e^0_i)\mR^\bot
 + 2 A(\mR^0)A(e^0_i)\mR^\bot +2
\nabla_{\mR^0}^{TX_G}(A(e^0_i)\mR^\bot)\nonumber\\
&\hspace*{30mm} +
R^{TB}(\mR^\bot,e^0_i)\mR^\bot+R^{TB}(\mR^\perp,\mR^0)e_i^0
,\nonumber\\
&(\nabla^{TB}\nabla^{TB} e^\bot_j)_{(\mR,\mR)}
= \frac{1}{3} R^{TB}(\mR^\bot,e^\bot_j)\mR^\bot
+ R^{TB}(\mR^\bot,\mR^0)e^\bot_j\nonumber\\
&\hspace*{30mm}
+ A(\mR^0)A(\mR^0)e^\bot_j + \nabla_{\mR^0}^{TX_G}(A(e^0_k)e^\bot_j)Z^0_k.
\nonumber
\end{align}
} 

Note that  by \eqref{b3.16}, $\nabla^{TB}_\mR (A_{x_0}
(e^0_i)e^0_i) = A(\mR^0)A_{x_0} (e^0_i)e^0_i$.

From \eqref{b3.16}, (\ref{b3.17}), we get
\begin{align}\label{b3.21}
\begin{split}
&(\nabla^{TB}_{e^\bot_i}e^\bot_i)_Z
= \frac{2}{3} R^{TB}(\mR^\bot,e^\bot_i)e^\bot_i +\cO(|Z|^2),\\
&(\nabla^{TB}_{e^0_i}e^0_i)_Z  = A_{x_0} (e^0_i)e^0_i-
\nabla^{TB}_{\mR}(A_{x_0} (e^0_i)e^0_i) +\frac{2}{3}
R^{TX_G}(\mR^0, e^0_i)e^0_i + \nabla^{TB}_{\mR^0}(A(e^0_i) e^0_i)\\
&\hspace*{10mm} + A(e^0_i) A(e^0_i)\mR^\bot
+ \nabla^{TX_G}_{ e^0_i}(A(e^0_i) \mR^\bot)
+ R^{TB}(\mR^\bot, e^0_i)e^0_i +\cO(|Z|^2)\\
&\hspace*{5mm}= A_{x_0} (e^0_i)e^0_i+ \Gamma_{i i}(\mR)+\cO(|Z|^2)
,
\end{split}\end{align}

Thus  by \eqref{b3.16}, \eqref{b3.17} and \eqref{b3.18}, at $x_0$,
\begin{align}\label{b3.22}
\nabla_{\mR^0}\nabla_{\mR^\bot}\left \langle e^\bot_j, e^0_i\right \rangle
&= \left \langle \nabla^{TB}_{\mR^0} e^\bot_j,
\nabla^{TB}_{\mR^\bot} e^0_i\right \rangle
+ \left \langle e^\bot_j,
\nabla^{TB}_{\mR^0}\nabla^{TB}_{\mR^\bot} e^0_i\right \rangle\\
&= \frac{1}{2} \left \langle R^{N_G}(\mR^0, e^0_i)\mR^\bot,
e^\bot_j\right \rangle.\nonumber
\end{align}

On the other hand, we have the following expansion for $\left
\langle e_j, e_i\right \rangle_Z$,
\begin{multline}\label{b3.23}
\left \langle e_i, e_j\right \rangle_{Z} = \left \langle e_i,
e_j\right \rangle _{Z^0} + (\nabla_{\mR^\bot}\left \langle e_i,
e_j\right \rangle) _{Z^0} + \frac{1}{2} (\nabla\nabla\left \langle
e_i, e_j\right \rangle)
_{(\mR^\bot,\mR^\bot), x_0} +\cO(|Z|^3)\\
= \left \langle e_i, e_j\right \rangle _{Z^0}
+ (\nabla_{\mR^\bot}\left \langle e_i, e_j\right \rangle)_{x_0}
+ (\nabla_{\mR^0}\nabla_{\mR^\bot}\left \langle e_i, e_j\right \rangle)_{x_0}
+ \left \langle \nabla^{TB}_{\mR^\bot}e_i,
\nabla^{TB}_{\mR^\bot} e_j\right \rangle _{x_0}\\
+ \frac{1}{2}\left \langle (\nabla^{TB} \nabla^{TB} e_i)
_{(\mR^\bot,\mR^\bot)},  e_j\right \rangle
+  \frac{1}{2}\left \langle e_i,  (\nabla^{TB} \nabla^{TB} e_j)
_{(\mR^\bot,\mR^\bot)}\right \rangle +\cO(|Z|^3).
\end{multline}

Thus by \eqref{3.21}, \eqref{b3.16},
\eqref{b3.18} and \eqref{b3.22}-\eqref{b3.23},
\begin{multline}\label{b3.24}
\left \langle e^0_i, e^0_j\right \rangle_Z
= \delta_{ij} -2 \left \langle A_{x_0}(e^0_i) e^0_j, \mR^\bot \right \rangle
+\frac{1}{3} \left \langle R^{TX_G}(\mR^0, e^0_i)\mR^0, e^0_j\right \rangle\\
+\left \langle \nabla ^{TX_G}_{\mR^0}(A(e^0_i)\mR^\bot), e^0_j\right \rangle
+\left \langle e^0_i, \nabla ^{TX_G}_{\mR^0}(A(e^0_j)\mR^\bot)\right \rangle\\
+ \left \langle A(e^0_i)\mR^\bot, A(e^0_j)\mR^\bot\right \rangle
+ \left \langle R^{TB}(\mR^\bot, e^0_i)\mR^\bot, e^0_j\right \rangle
 +\cO(|Z|^3),
\end{multline}
and
\begin{align}\label{b3.25}
\begin{split}
&\left \langle e^0_i, e^\bot_j\right \rangle_Z
= \frac{1}{2}\left \langle R^{N_G} (\mR^0, e^0_i)\mR^\bot,
e^\bot_j\right \rangle
+ \frac{2}{3}\left \langle R^{TB} (\mR^\bot, e^0_i)\mR^\bot,
e^\bot_j\right \rangle +\cO(|Z|^3),\\
&\left \langle e^\bot_i, e^\bot_j\right \rangle_Z
= \delta_{ij}+ \frac{1}{3} \left \langle R^{TB} (\mR^\bot, e^\bot_i)\mR^\bot,
e^\bot_j\right \rangle +\cO(|Z|^3). 
\end{split}\end{align}

From \eqref{b3.15},  \eqref{b3.24} and \eqref{b3.25}, we get
\begin{align}\label{b3.26}
&\det g_{ij}(Z) = 1-2 \left\langle A_{x_0}(e^0_i)e^0_i,\mR^\bot \right\rangle
 +K_2(\mR)+ \cO(|Z|^3),\\
&\kappa ^{\frac{1}{2}}(tZ) = (\det g_{ij})^{1/4} (tZ)\nonumber\\
&\hspace*{5mm}
= 1-\frac{t}{2}\left \langle A(e^0_i)e^0_i, \mR^\bot\right \rangle
- \frac{3t^2}{8} \Big(\sum_i \left \langle A(e^0_i)e^0_i,
\mR^\bot\right \rangle\Big)^2
+\frac{t^2}{4} K_2(\mR) + \cO(t^3), \nonumber\\
&\kappa ^{-\frac{1}{2}}(tZ) = 1+\frac{t}{2}\left \langle A(e^0_i)e^0_i,
\mR^\bot\right \rangle
+ \frac{5t^2}{8} \Big(\sum_i \left \langle A(e^0_i)e^0_i,
\mR^\bot\right \rangle\Big)^2
-\frac{t^2}{4} K_2(\mR) + \cO(t^3). \nonumber
\end{align}

\comment{
By \cite[Prop. 1.28]{BeGeVe}, the Taylor expansion of $g_{ij}(Z)$
with respect to the basis $\{e_i\}$ to order $r$ is a polynomial of
 the Taylor expansion of $R^{TB}$
to order $r-2$, moreover
\begin{equation}\label{3.21}
\begin{split}
&g_{ij}(Z) = \delta_{ij} +  \frac{1}{3}
\left \langle R^{TB}_{x_0} (\mR,e_i) \mR, e_j\right \rangle_{x_0}
 + \cO(|Z|^3),\\
&\kappa(Z)= |\det (g_{ij}(Z))|^{1/2}  = 1 +
\frac{1}{6} \left \langle R^{TB}_{x_0} (\mR,e_i) \mR, e_i\right \rangle_{x_0}
 + \cO(|Z|^3).
\end{split}
\end{equation}
Owing to \eqref{3.21},
\begin{equation}\label{3.22}
\begin{split}
\Gamma _{ij}^l (Z)& =  \frac{1}{2} g^{lk} (\partial_i g_{jk}
+ \partial_j g_{ik}-\partial_k g_{ij})(Z)\\
&= \frac{1}{3}\Big  [
\left \langle R^{TB}_{x_0} (\mR, e_j) e_i, e_l\right \rangle _{x_0}
+ \left \langle R^{TB}_{x_0} (\mR, e_i) e_j, e_l\right \rangle_{x_0}\Big ]
 + \cO(|Z|^2).
\end{split}
\end{equation}
}

Moreover, as a $2(n-n_0)\times 2(n-n_0)$-matrix, we have
\begin{multline}\label{b3.260}
\Big((\delta_{ij}-2  \left \langle A_{x_0}(e^0_i) e^0_j, \mR^\bot
\right \rangle)\Big)^{-1} = \left(\delta_{ij}+ 2  \left \langle
A_{x_0}(e^0_i) e^0_j,
\mR^\bot \right \rangle\right)\\
+4 \Big(\left\langle A_{x_0}(e^0_i)\mR^\bot,
A_{x_0}(e^0_j)\mR^\bot \right\rangle\Big)
+ \cO(|Z|^3) .
\end{multline}

Note that from \eqref{g7}, \eqref{a3.17},
\begin{align}\label{b3.27}
&[\left \langle A(e^0_i)e^0_i, \mR^\bot\right \rangle, \cL^0_2]
=2 \left \langle A(e^0_i)e^0_i,
e^\bot_k\right \rangle \nabla_{0,e^\bot_k}.
\end{align}
Thus from \eqref{b3.15}, \eqref{3.28}, \eqref{b3.21},
\eqref{b3.24}-\eqref{b3.26},
 the coefficients of $t$, $t^2$ in the expansion of
$g^{ij}(tZ) t \Gamma^k_{ij}(tZ)\nabla_{t,e_k}=tg^{ij}(tZ)
\nabla_{t,(\nabla^{TB}_{e_i}e_j)(tZ)}$ are
\begin{align}\label{b3.28}
&\left \langle A(e^0_i)e^0_i,e^\bot_k\right \rangle \nabla_{0,e^\bot_k};
\\
&2 \left \langle A(e^0_i)e^0_j, \mR^\bot\right \rangle \nabla_{A(e^0_i)e^0_j}
+\left \langle \Gamma_{i i}(\mR), e_j \right \rangle\nabla_{0,e_j}
+\frac{2}{3} \left \langle R^{TB} (\mR^\bot, e^\bot_i)
e^\bot_i, e_j\right \rangle \nabla_{0,e_j}\nonumber\\
&\hspace*{5mm}
- \Big[\frac{1}{2} \left \langle A(e^0_l)e^0_l, \mR^\bot\right \rangle,
 \nabla_{A(e^0_i)e^0_i}\Big]+\frac{1}{3}(\partial_k
 R^{L_B})_{x_0}Z_k(\mR,A(e_i^0)e_i^0) .\nonumber
\end{align}

By \eqref{a3.16}, \eqref{3.28} and \eqref{b3.24}-\eqref{b3.28},
the coefficient of $t$ in the expansion of $\cL^t_3$ is
$\mO^\prime_1$ in \eqref{a3.20}.

We denote by $[A,B]_+=AB+BA$.

By \eqref{a3.16}, \eqref{3.28}, \eqref{b3.21} and
\eqref{b3.24}-\eqref{b3.260}, the coefficient of $t^2$ in the
expansion of $\cL^t_3-(g^{ij}t \Gamma^k_{ij})(tZ)\nabla_{t,e_k}$
is
\begin{multline}\label{b3.29}
I_2 - 2\left \langle A(e^0_i)e^0_j,\mR^\bot \right \rangle
\Big[\frac{1}{3}\nabla_{0,e^0_i} ( \partial_k R^{L_B})_{x_0} (\mR,e^0_j)Z_k\\
+ \frac{1}{3} ( \partial_k R^{L_B})_{x_0} (\mR,e^0_i)Z_k\nabla_{0,e^0_j}
-\frac{1}{2} [\left \langle A(e^0_l)e^0_l,\mR^\bot \right \rangle ,
\nabla_{0,e^0_i}\nabla_{0,e^0_j}]\Big]\\
+I_1+ \Big[\frac{1}{2}\left \langle A(e^0_l)e^0_l,\mR^\bot \right \rangle,
 [\frac{1}{3} ( \partial_k R^{L_B})_{x_0} (\mR,e_i)Z_k,\nabla_{0,e_i}]_+
  \Big]\\
+ \Big[\frac{1}{4}K_2(\mR)
-\frac{3}{8} \Big(\sum_l \left \langle A(e^0_l)e^0_l, \mR^\bot
\right \rangle \Big)^2 , \cL^0_2\Big]\\
 -\frac{1}{4}[\left \langle A(e^0_l)e^0_l, \mR^\bot \right \rangle,
\cL^0_2] \left \langle A(e^0_k)e^0_k,
\mR^\bot \right \rangle
- R^{E_B}_{x_0}(\mR,e_i)\nabla_{0,e_i}\\
 -\frac{1}{9} \sum_i \Big[\sum_j(\partial_j R^{L_B})_{x_0} (\mR,e_i)Z_j\Big]^2
- R^E_{\tau,x_0}
+\frac{1}{h_0} (\nabla_{e_j}\nabla_{e_j} h
-\nabla_{A(e^0_i)e^0_i}h )_{x_0} .
\end{multline}
Here $I_2$ is from the coefficient of $t^2$ in the expansion of $g^{ij}$,
 the second term is the product of the coefficients of $t^1$
in the expansion of $g^{ij}$ and $\nabla_{t,e_i}\nabla_{t,e_j}$;
$I_1$ is from the coefficient of $t^2$ in the expansion of $R^{L_B}$,
the fourth term is from the product of the coefficients of $t^1$
in $\kappa^{1/2}, \kappa^{-1/2}$ and in
$\kappa^{-1/2}\nabla_{t,e_i}\nabla_{t,e_i} \kappa^{1/2}$ (cf. \eqref{3.28}),
the fifth and sixth terms is from the coefficients
 of $t^2$ in the expansion of
$\kappa^{1/2}, \kappa^{-1/2}$ and the product of the coefficients of $t^1$
in the expansion of $\kappa^{1/2}$ and $\kappa^{-1/2}$;
the seventh term is from $R^{E_B}$, and the eighth term is from the product
 of the coefficients of $t^1$ in the expansion of $R^{L_B}$.

Certainly,
\begin{align}\label{b3.290}
\frac{1}{6}\Big[\left \langle A(e^0_l)e^0_l,\mR^\bot \right \rangle,
 [ ( \partial_k R^{L_B})_{x_0} (\mR,e_i)Z_k,\nabla_{0,e_i}]_+
  \Big] = -\frac{1}{3} ( \partial_k R^{L_B})_{x_0} (\mR,A(e^0_l)e^0_l)Z_k.
\end{align}

 By \eqref{b3.27}, \eqref{b3.28}, \eqref{b3.290}
and by the fact that $A(e^0_i)e^0_j$ is symmetric on $i,j$, we see
that the coefficient of $t^2$ in the expansion of $\cL^t_3$ is
$\mO^\prime_2$ in \eqref{a3.20}.
\end{proof}

To simplify the notation, we will often denote by $e_i$
the lift $e^H_i$ of $e_i$.
\begin{lemma} \label{at4.4} The following identities hold,
\begin{subequations}
\begin{align}
&(\partial _i R^{L_{B}})_{x_0}(\mR,e_l)Z_i  = -3\sqrt{-1}\pi
\left\langle J T(\mR, e_l) - J T(\mR^0, P^{TX_G}e_l),
\mR^\bot \right\rangle,\label{3.32}\\
&\frac{\sqrt{-1}}{\pi}B(Z, e^0_l)
= \frac{1}{6}\left\langle R^{TX_G}(\mR^0, J \mR^0)\mR^0, e^0_l\right\rangle
-\frac{5}{4}\left\langle  J\mR^\bot,
\nabla ^{TY}_{\mR} (T(e_i,e^0_l))Z_i\right\rangle\label{a3.32}\\
&\hspace*{5mm}+\frac{1}{2}
\left\langle 2 \nabla^{TX_G}_{\mR^0} (A(e^0_l)e^\bot_j) Z^\bot_j
+ R^{TB}(\mR^\bot, e^0_l ) \mR^\bot
+R^{TB}(\mR^\bot, \mR^0 )e^0_l , J\mR^0\right\rangle\nonumber\\
&\hspace*{5mm}
-\frac{1}{2}  \left\langle 3  \nabla^{TX_G}_{\mR^0}
(A(e^0_i)e^\bot_j) Z^0_iZ^\bot_j
+2  R^{TB}(\mR^\bot, \mR^0)\mR^\bot
+ R^{TB}(\mR^\bot, \mR^0)\mR^0, Je^0_l\right\rangle
\nonumber\\
&\hspace*{5mm}
+ \frac{1}{2}\left\langle J\mR^\bot, T(\mR^0-\frac{1}{4}\mR, e^0_i)
\right\rangle
 \left\langle J\mR^\bot, T(e^0_i, J e^0_l)\right\rangle
\nonumber\\
&\hspace*{5mm}
+\frac{1}{8}\left\langle  T(\mR^0, \mR^\bot), T(e^0_l, J \mR^0)\right\rangle
+ \frac{1}{8}\left\langle  T(\mR^0, J\mR^0),
T(\mR^\bot, e^0_l)\right\rangle\nonumber\\
&\hspace*{5mm}
- \frac{1}{8}\left\langle  T(\mR^\bot, J\mR^0), T(\mR,e^0_l)\right\rangle
+\frac{1}{2} \left\langle
  T(\mR^\bot, J \mR^\bot), T(\mR,e^0_l)\right\rangle\nonumber\\
&\hspace*{5mm}
- \frac{1}{8}\left\langle T(e^0_l, J\mR^0), e^\bot_j \right\rangle
\left\langle J\mR^\bot,  T(\mR^\bot,e^\bot_j)\right\rangle.\nonumber
\end{align}
\end{subequations}
\end{lemma}
\begin{proof}
By  (\ref{h3}), (\ref{ah7}), (\ref{h8}) and (\ref{a6}),
\begin{align}\label{a3.28}
\frac{\sqrt{-1}}{2\pi} R^{L_{B}}(e_k,e_l)
&= \left\langle Je^H_k, e^H_l \right\rangle  + \mu(\Theta)(e_k,e_l)\\
&= \left\langle Je^H_k, e^H_l \right\rangle
+ \left\langle \wi{\mu}, T(e_k,e_l)\right\rangle.\nonumber
\end{align}

Thus by \eqref{g3.43}, \eqref{3.11a} and \eqref{3.12a},
we get at $x_0$ the following formulas which will be used in \eqref{3.38},
\begin{align}\label{a3.43}
\wi{\mu}_{x_0}=0, \quad
(\nabla ^{TY}_{\mR}\wi{\mu})_{x_0}= -J \mR^\bot, \quad
(\nabla ^{TY}_{\cdot}\nabla ^{TY}_{\cdot}\wi{\mu})_{(\mR,\mR)}
=  T(\mR^\bot, J\mR^\bot).
\end{align}

By (\ref{g3.26}) and $\mu=0$ on $P$, we have at $x_0$,
\begin{align}\label{3.30}
(\nabla _{e_i} \left\langle \wi{\mu}, T(e_k,e_l)\right\rangle)_{x_0}
&= \left\langle \nabla ^{TY}_{e^H_i}\wi{\mu}, T(e_k,e_l)\right\rangle
+ \left\langle \wi{\mu}, \nabla ^{TY}_{e^H_i}(T(e_k,e_l))\right\rangle\\
&= \left\langle JT(e_k, e_l), e_i\right\rangle . \nonumber
\end{align}

By (\ref{g3.29}), (\ref{3.12a}) and \eqref{b3.16}, we have
\begin{multline}\label{3.29}
(\nabla _{e^H_i}\left\langle Je^H_k, e^H_l \right\rangle )_{x_0}
= \left\langle J\nabla ^{TX}_{e^H_i} e^H_k, e^H_l \right\rangle_{x_0}
+\left\langle Je^H_k,\nabla ^{TX}_{e^H_i} e^H_l \right\rangle_{x_0}\\
=-\frac{1}{2} \left\langle JT(e_i, e_k), e_l \right\rangle
-\frac{1}{2} \left\langle Je_k,T(e_i,  e_l) \right\rangle\\
+ \left\langle J A(P^{TX_G}e_i)P^{N_G}e_k+ J A(P^{TX_G}e_k)P^{N_G}e_i,
P^{TX_G}e_l\right\rangle\\
+ \left\langle JP^{TX_G}e_k,  A(P^{TX_G}e_i)P^{N_G}e_l
+  A(P^{TX_G}e_l)P^{N_G}e_i\right\rangle.
\end{multline}

By  \eqref{3.11a}, (\ref{a3.28}), (\ref{3.30}) and (\ref{3.29}),
for $U\in T_{x_0}B$,
\begin{multline}\label{3.31}
 \frac{\sqrt{-1}}{2\pi} (\partial_U R^{L_{B}})_{x_0}(U,e_l)
= \frac{3}{2}\left\langle J T(U, e_l), U\right\rangle
- 2 \left\langle  A(P^{TX_G}U)P^{N_G}U, J P^{TX_G}e_l\right\rangle\\
+  \left\langle J P^{TX_G}U,  A(P^{TX_G}U)P^{N_G}e_l
+ A(P^{TX_G}e_l) P^{N_G}U\right\rangle\\
=  \frac{3}{2}\left\langle JT(U, e_l)- JT(P^{TX_G}U,P^{TX_G}e_l),
 U\right\rangle.
\end{multline}

Note that $(JTY)_G=N_G$ on $X_G$, by (\ref{3.31}), we get
(\ref{3.32}).

By \eqref{a3.18} and \eqref{a3.28}, one gets  at $x_0$,
\begin{align}\label{a3.31}
\frac{\sqrt{-1}}{\pi}B(Z, e_l)
= \frac{1}{2} \Big(\nabla \nabla \left\langle Je_k,e_l \right\rangle
+\nabla\nabla\left\langle \wi{\mu}, T(e_k,e_l)\right\rangle\Big)_{(\mR,\mR)}Z_k .
\end{align}

From \eqref{3.12a} we have
\begin{multline}\label{b3.30}
\Big(\nabla \nabla\left\langle J e^H_k,
e^{H}_l\right\rangle\Big)_{(\mR,\mR)} Z_k
= \langle J\mR, (\nabla^{TX}\nabla^{TX}e_l^{H} )_{(\mR,\mR)} \rangle\\
+\langle J (\nabla^{TX}\nabla^{TX}e_k^H )_{(\mR,\mR)},e_l^{H} \rangle Z_k
+2 \langle J\nabla^{TX}_\mR e_k^H,\nabla^{TX}_\mR e_l^{H}
\rangle Z_k.
\end{multline}

From \eqref{h01}, \eqref{b3.16},  one finds at $x_0$ that
\begin{align}\label{a3.311}
&J\mR^\bot\in TY,\quad J\mR^0\in TX_G, \nonumber\\
&\nabla^{TB}_{\mR} e^0_i = A(e^0_i)\mR, \quad
\nabla^{TB}_{\mR} e^\bot_i = A(\mR^0)e^\bot_i,\\
&(\nabla^{T^HX}_{e_j^H}e_i^H)Z_iZ_j=(\nabla^{TB}_{e_j}e_i)^H Z_i Z_j
 =2A(\mR^0)\mR^\perp+A(\mR^0)\mR^0. \nonumber
\end{align}

Now by  (\ref{g3.29}),
\begin{align}\label{3.33}
&(\nabla ^{TX}_{e^H_j} \nabla ^{TX}_{e^H_i} e^H_k)_{x_0}
= \nabla ^{T^HX}_{e^H_j} \nabla ^{T^HX}_{e^H_i} e^H_k
-\frac{1}{2} T(e^H_j, \nabla ^{T^HX}_{e^H_i} e^H_k)
-\frac{1}{2} \nabla ^{TX}_{e^H_j} (T(e^H_i,e^H_k)).
\end{align}

By (\ref{b3.18}), we get
\begin{align}\label{b3.31}
&(\nabla^{TB}\nabla^{TB} e_k)_{(\mR,\mR)}Z_k
= \nabla_{\mR^0}^{TB}(A(e^0_j)e^0_i)Z^0_j Z^0_i
+  3 A(\mR^0)A(\mR^0)\mR^\bot\\
&\hspace*{15mm}+ 3\nabla_{\mR^0}^{TX_G}(A(e^0_i)\mR^\bot)Z^0_i
+2R^{TB}(\mR^\bot,\mR^0)\mR^\bot
+ R^{TB}(\mR^\bot,\mR^0)\mR^0.\nonumber
\end{align}

From \eqref{b3.18}, \eqref{a3.311}, \eqref{3.33}, \eqref{b3.31},
 the anti-symmetric property of the torsion tensor $T$
and the fact that $A$ exchanges $TX_G$ and $N_G$, we get
\begin{align}\label{b3.32}
&\langle J\mR, (\nabla^{TX}\nabla^{TX}e_l^{0,H} )_{(\mR,\mR)} \rangle
= \left\langle \frac{1}{3}R^{TX_G}(\mR^0, e^0_l)\mR^0
+ \nabla^{TB}_{\mR^0}(A(e^0_j)e^0_l) Z^0_j, J \mR^0\right\rangle \\
&\hspace*{15mm} + \left\langle
2 \nabla ^{TX_G}_{\mR^0} (A(e^0_l)e^\bot_j)Z^\bot_j
+R^{TB}(\mR^\bot, e^0_l )\mR^\bot
+ R^{TB}(\mR^\bot,\mR^0) e^0_l, J \mR^0\right\rangle \nonumber\\
&\hspace*{15mm}
-\frac{1}{2}\left\langle J \mR^\bot, T(\mR,A(e^0_l)\mR)\right\rangle
-\frac{1}{2}\left\langle J \mR,
\nabla ^{TX}_{\mR} (T(e_i,e^0_l))Z_i\right\rangle, \nonumber\\
&\langle J (\nabla^{TX}\nabla^{TX}e_k^H )_{(\mR,\mR)},e_l^{0,H} \rangle Z_k
=\left\langle 2 J R^{TB}(\mR^\bot, \mR^0)\mR^\bot
+ J R^{TB}(\mR^\bot, \mR^0)\mR^0, e^0_l\right\rangle\nonumber\\
&\hspace*{10mm}+ \left\langle J \nabla^{TB}_{\mR^0}(A(e^0_j)e^0_i)
Z^0_jZ^0_i + 3 J \nabla^{TX_G}_{\mR^0} (A(e^0_i)e^\bot_j)
Z^0_iZ^\bot_j, e^0_l\right\rangle. \nonumber
\end{align}
Note that from \eqref{ah4}, \eqref{3.10}, \eqref{3.11a}, \eqref{a3.311}
 and $A$ exchanges $TX_G$ and $N_G$,
\begin{align} \label{b3.34}
&\left\langle J \mR,
\nabla ^{TX}_{\mR} (T(e_i,e^0_l)) Z_i\right\rangle
= \left\langle J \mR^\bot,
\nabla ^{TY}_{\mR} (T(e_i,e^0_l)) Z_i\right\rangle \\
&\hspace*{50mm}
+\frac{1}{2}\left\langle T(\mR, J\mR^0),T(\mR,e^0_l)\right\rangle,\nonumber\\
& \left\langle J \nabla^{TB}_{\mR^0}(A(e^0_j)e^0_i)Z^0_jZ^0_i,
e^0_l \right\rangle
=- \left\langle A(\mR^0)A(\mR^0)\mR^0, Je^0_l \right\rangle \nonumber\\
&\hspace*{30mm}
=-\frac{1}{4}\left\langle T(\mR^0, J\mR^0),T(\mR^0,e^0_l) \right\rangle,
\nonumber\\
&\left\langle \nabla^{TB}_{\mR^0}(A(e^0_j)e^0_l),  J\mR^0\right\rangle
= -\left\langle A(e^0_j)e^0_l, A(\mR^0)J\mR^0\right\rangle=0.\nonumber
\end{align}

\comment{
For $U\in TX_G, V\in N_G$, we have at $x_0$, by \eqref{3.11a},
\begin{align} \label{b3.36}
A(U)V= \left\langle A(U)V, e^0_j\right\rangle e^0_j
=-\left\langle V, A(U)e^0_j\right\rangle e^0_j
= \frac{1}{2} \left\langle T(U,Je^0_j), JV\right\rangle e^0_j.
\end{align}
}

By (\ref{g3.29}), (\ref{3.12a}), \eqref{33.36}, \eqref{a3.311} and
the fact that $A$ exchanges $TX_G$ and $N_G$, at $x_0$,
\begin{align}\label{b3.35}
&\left\langle J \nabla^{TX} _{\mR}e^H_k,
\nabla^{TX} _{\mR}e^{0,H}_l\right\rangle Z_k
=\Big\langle J \nabla^{TB} _{\mR}e_k,
A(e^{0}_l)\mR -\frac{1}{2} T(\mR,e^0_l)\Big\rangle Z_k  \\
&\hspace*{15mm} = \Big\langle J  A(\mR^0) \mR^0,  -\frac{1}{2}
T(\mR,e^0_l)\Big\rangle +2 \left\langle J  A(\mR^0) \mR^\bot,
A(e^0_l) \mR^\bot\right\rangle\nonumber\\
& = \frac{1}{4}\left\langle  T(\mR^0,J \mR^0), T(\mR,e^0_l)\right\rangle
+\frac{1}{2} \left\langle J\mR^\perp,T(\mR^0,e_j^0)\right\rangle\left\langle
J\mR^\perp,T(e_l^0,Je_j^0)\right\rangle.\nonumber
\end{align}

By \eqref{b3.30},
\eqref{b3.32}-\eqref{b3.35}, at $x_0$,
\begin{multline}\label{3.36}
\Big(\nabla \nabla\left\langle J e^H_k,
e^{0,H}_l\right\rangle\Big)_{(\mR,\mR)} Z_k
=\frac{1}{3}\left\langle R^{TX_G}(\mR^0, e^0_l)\mR^0, J \mR^0\right\rangle\\
+\left\langle 2\nabla^{TX_G}_{\mR^0} (A(e^0_l)e^\bot_j) Z^\bot_j
+ R^{TB}(\mR^\bot, e^0_l ) \mR^\bot
+ R^{TB}(\mR^\bot, \mR^0)e^0_l, J\mR^0\right\rangle\\
- \left\langle 2 R^{TB}(\mR^\bot, \mR^0)\mR^\bot+ R^{TB}(\mR^\bot, \mR^0)\mR^0
+ 3  \nabla^{TX_G}_{\mR^0} (A(e^0_i)e^\bot_j) Z^0_iZ^\bot_j
, J e^0_l\right\rangle\\
-\frac{1}{2} \left\langle J\mR^\bot, T(\mR, A(e^0_l)\mR)+
\nabla ^{TY}_{\mR} (T(e_i,e^0_l))Z_i\right\rangle
+\frac{1}{4} \left\langle  T(\mR^0,J\mR^0), T(\mR^\bot, e^0_l)\right\rangle\\
-\frac{1}{4} \left\langle T(\mR^\bot, J\mR^0), T(\mR,e^0_l)\right\rangle
+ \left\langle J\mR^\bot,  T(\mR^0, e^0_j)\right\rangle
\left\langle   J\mR^\bot, T(e^0_l, J e^0_j)\right\rangle .
\end{multline}

Observe that $A(e^0_i)\mR^0\in N_{G}$, $A(e^0_i)\mR^\bot\in TX_G$.
By \eqref{3.11a}, \eqref{3.11b}, \eqref{3.11d} and \eqref{33.36},
\begin{multline} \label{b3.37}
\left\langle J\mR^\bot, T(\mR, A(e^0_l)\mR)\right\rangle
=\langle J\mR^\perp,T(\mR,A(e_l^0)\mR^0)\rangle+\langle
J\mR^\perp,T(\mR,A(e_l^0)\mR^\perp)\rangle\\
= \frac{1}{2}\left\langle JT(e^0_l, J\mR^0), e^\bot_j\right\rangle
\left\langle J\mR^\bot, T(\mR, e^\bot_j)\right\rangle
+ \left\langle J \mR^\bot, T(\mR, A(e^0_l)\mR^\bot)\right\rangle\\
=-\frac{1}{2} \left\langle T(e^0_l, J\mR^0), T(\mR^0, \mR^\bot)\right\rangle
+ \frac{1}{2} \left\langle J\mR^\bot, T(\mR, e^0_j)\right\rangle
\left\langle   J\mR^\bot, T(e^0_l, J e^0_j)\right\rangle\\
+ \frac{1}{2} \left\langle JT(e^0_l, J\mR^0), e^\bot_j\right\rangle
\left\langle J\mR^\bot, T(\mR^\bot, e^\bot_j)\right\rangle.
\end{multline}

From (\ref{a3.43}), at $x_0$,
\begin{multline}\label{3.38}
(\nabla \nabla  \left\langle  \wi{\mu}, T(e_k,e_l)\right\rangle)_{(\mR,\mR)}\\
=  \left\langle (\nabla ^{TY}_\cdot\nabla^{TY}_\cdot \wi{\mu})_{(\mR,\mR)},
T(e_k,e_l)\right\rangle
+ 2 \left\langle \nabla^{TY}_{\mR} \wi{\mu},
\nabla ^{TY}_{\mR} (T(e_k,e_l))\right\rangle\\
=\left\langle T(\mR^\bot, J \mR^\bot) ,T(e_k,e_l)
\right\rangle-2\left\langle \nabla ^{TY}_{\mR} (T(e_k,e_l)), J
\mR^\bot\right\rangle .
\end{multline}

Finally, by (\ref{a3.15}), (\ref{a3.31}), (\ref{3.36}),
  (\ref{b3.37}) and (\ref{3.38}), we get \eqref{a3.32}.
\end{proof}

We now examine the coefficients in the expansion of terms
involving the moment map $\widetilde{\mu}$.

 Set
\begin{multline}\label{a3.48}
\mO''_2=
-\frac{1}{3}\left\langle (\nabla^{TY}_\cdot \dot{g}^{TY}_\cdot)_{(\mR,\mR)}
J\mR^\bot, J\mR^\bot\right\rangle
+ \frac{1}{6}\left\langle \nabla ^{TY}_\mR (T(e_j, J_{x_0}e^0_i)),
J\mR^\bot\right\rangle Z_j Z^0_i\\
+\frac{1}{3}\left\langle   \nabla ^{N_G}_{\mR^0} (A(e^0_j) e^0_i)Z^0_j Z^0_i
+ R^{TB}(\mR^\bot,\mR^0)\mR^0,\mR^\bot \right\rangle\\
-\frac{1}{12}\sum_l \left\langle T(\mR, e_l),J\mR^\bot\right\rangle^2
+\frac{1}{4}\left\langle J\mR^\bot,T(\mR^\bot, e^0_l)\right\rangle
\left\langle J\mR^\bot,T(\mR^0, e^0_l)\right\rangle\\
+ \frac{7}{12} | T(\mR^\bot, J\mR^\bot)|^2
+\frac{1}{3}\left\langle  T(\mR^0, J\mR^\bot),
T(\mR^\bot, J\mR^\bot)\right\rangle .
\end{multline}
\begin{lemma}\label{t4.5} For $|t|\leq 1$, we have
\begin{align}\label{3.48}
&|\frac{1}{t}\wi{\mu}|^2_{g^{TY}}(tZ)= |Z^\bot|^2
 - t \left\langle T(\mR^\bot, J\mR^\bot), J\mR^\bot \right\rangle
+ t^2\mO''_2 +\cO(t^3),\\
&\left\langle  \wi{\mu}, \wi{\mu}^E\right\rangle_{g^{TY}} (tZ)
=-t \left\langle J \mR^\bot, \wi{\mu}^E_{x_0}\right\rangle\nonumber\\
&\hspace*{10mm}
+t^2 \Big(\frac{1}{2}  \left\langle T(\mR^\bot, J\mR^\bot),
\wi{\mu}^E_{x_0}\right\rangle
- \left\langle J\mR^\bot, \nabla^{TY}_{\mR}\wi{\mu}^E\right\rangle_{x_0}\Big)
+\cO(t^3). \nonumber
\end{align}
 \end{lemma}
\begin{proof}\comment{ First, the following two formulas are clear,
\begin{align}\label{3.481}
\frac{1}{2}\left. \frac{\partial^2}{\partial t^2}
|\widetilde{\mu}(tZ)|^2\right|_{t=0}=\frac{1}{2}
\left.\left(\nabla\nabla|\widetilde{\mu}(tZ)|^2\right)_{(\mR,\mR)}\right|_{t=0}
=\langle\nabla\widetilde{\mu},\nabla\widetilde{\mu}\rangle_{(\mR,\mR)}
=|\nabla_\mR\widetilde{\mu}|^2,
\end{align}
\begin{align}\label{3.482}
\frac{1}{3!}\left. \frac{\partial^3}{\partial t^3}
|\widetilde{\mu}(tZ)|^2\right|_{t=0}=\frac{1}{6}
\left.\left(\nabla\nabla\nabla|\widetilde{\mu}(tZ)|^2
\right)_{(\mR,\mR,\mR)}\right|_{t=0}
=\langle\nabla\nabla\widetilde{\mu},\nabla\widetilde{\mu}
\rangle_{(\mR,\mR,\mR)}.
\end{align}
}
 By  \eqref{g3.26}, \eqref{g3.41},  \eqref{g3.42}, \eqref{3.12a},
 \eqref{a3.311}, $\bJ=J$ and
$\wi{\mu}=0$ on $P$, we get, at $x_0$,
\begin{multline}\label{a3.44}
(\nabla ^{TY}_{e^H_k}\nabla ^{TY}_{e^H_j}\nabla ^{TY}_{e^H_i} \wi{\mu})_{x_0}
=- P^{TY} J \nabla ^{TX}_{e^H_k} \nabla ^{TX}_{e^H_j} e^H_i
- \frac{1}{2} T(e^H_k,P^{T^HX} J \nabla ^{TX}_{e^H_j} e^H_i )\\
-\frac{1}{2} \nabla ^{TY}_{e^H_k} (T(e^H_j,P^{T^HX} J e^H_i ))
-\frac{1}{2} (\nabla ^{TY}_{e^H_j} \dot{g}^{TY}_{e^H_i})
(\nabla ^{TY}_{e^H_k}\wi{\mu})\\
- \frac{1}{2}(\nabla ^{TY}_{e^H_k} \dot{g}^{TY}_{e^H_i})
(\nabla ^{TY}_{e^H_j}\wi{\mu})
- \frac{1}{2} \dot{g}^{TY}_{e^H_i}
(\nabla ^{TY}_{e^H_k} \nabla ^{TY}_{e^H_j}\wi{\mu}).
\end{multline}
From \eqref{g3.29}, \eqref{a3.43}, \eqref{a3.311}, \eqref{3.33}, \eqref{b3.31}
 and \eqref{a3.44}, we have
 \begin{multline}\label{a3.45}
(\nabla ^{TY}_{\cdot}\nabla ^{TY}_{\cdot}\nabla ^{TY}_{\cdot} \wi{\mu})
_{(\mR,\mR,\mR)}:=
(\nabla ^{TY}_{e^H_k}\nabla ^{TY}_{e^H_j}\nabla ^{TY}_{e^H_i} \wi{\mu})_{x_0}
Z_k Z_j Z_i \\
=-J  \nabla ^{N_G}_{\mR^0}(A(e^0_j) e^0_i)Z^0_j Z^0_i
- 3 J A(\mR^0)A(\mR^0)\mR^\bot
- 2 P^{TY}JR^{TB}(\mR^\bot, \mR^0)\mR^\bot\\
-  P^{TY}JR^{TB}(\mR^\bot, \mR^0)\mR^0
- T(\mR, JA(\mR^0)\mR^\bot) \\
-\frac{1}{2}  \nabla ^{TY}_{\mR} (T(e^H_j,P^{T^HX} J e^H_i ))Z_j Z_i
+ (\nabla ^{TY}_{\cdot} \dot{g}^{TY}_{\cdot})_{(\mR,\mR)} J\mR^\bot
- \frac{1}{2} \dot{g}^{TY}_{\mR}
( T(\mR^\bot, J\mR^\bot)).
\end{multline}
Now by \eqref{g3.481},  \eqref{a3.43},
\eqref{a3.45}, and $\wi{\mu}=0$ on $P$,  we have
 \begin{multline}\label{a3.46}
|\frac{1}{t}\wi{\mu}|_{g^{TY}}^2 (tZ) =
\sum_{k=2}^4 \frac{1}{k!}
\frac{\partial ^k}{\partial t^k}\Big(|\wi{\mu}|_{g^{TY}}^2 (tZ)\Big)|_{t=0}\,
 t^{k-2}+\cO(t^3)\\
=|\nabla ^{TY}_{\mR}\wi{\mu}|^2_{x_0}
+ t \left\langle (\nabla ^{TY}_{\cdot}\nabla ^{TY}_{\cdot} \wi{\mu})
_{(\mR,\mR)}, \nabla ^{TY}_{\mR}\wi{\mu}\right\rangle_{x_0} \\
+  \frac{t^2}{4!}  \Big( 8 \left\langle (\nabla ^{TY}_{\cdot}
\nabla ^{TY}_{\cdot}\nabla ^{TY}_{\cdot} \wi{\mu})
_{(\mR,\mR,\mR)}, \nabla ^{TY}_{\mR}\wi{\mu}\right\rangle_{x_0}
+ 6 |(\nabla ^{TY}_{\cdot}\nabla ^{TY}_{\cdot} \wi{\mu})
_{(\mR,\mR)}|^2_{x_0} \Big) +\cO(t^3).
\end{multline}

By \eqref{3.11c},
\begin{align}\label{aa3.42}
T(\mR^0,J\mR^\bot)= \frac{1}{2} T(\mR^\bot, J\mR^0).
\end{align}

From \eqref{h3}, \eqref{33.36}, \eqref{a3.43},
\eqref{a3.45}, \eqref{a3.46} and \eqref{aa3.42},
we get the coefficients of $t^0, t^1$
in the expansion of $|\frac{1}{t}\wi{\mu}|_{g^{TY}}^2 (tZ)$
 in \eqref{3.48},  and the coefficients of $t^2$ is
 \begin{multline}\label{a3.47}
 \frac{1}{3}\left\langle J  \nabla ^{N_G}_{\mR^0} (A(e^0_j) e^0_i)Z^0_j Z^0_i
+3 J A(\mR^0)A(\mR^0)\mR^\bot +JR^{TB}(\mR^\bot, \mR^0)\mR^0,
J\mR^\bot \right\rangle \\
+ \frac{1}{3}\left\langle  2 JR^{TB}(\mR^\bot, \mR^0)\mR^\bot
+ T(\mR, JA(\mR^0)\mR^\bot) ,J\mR^\bot \right\rangle \\
-\frac{1}{3}\left\langle (\nabla^{TY}_\cdot \dot{g}^{TY}_\cdot)_{(\mR,\mR)}
J\mR^\bot, J\mR^\bot\right\rangle
+ \frac{1}{6} \left\langle\nabla ^{TY}_{\mR}
(T(e^H_j,P^{T^HX} J e^H_i ))Z_j Z_i, J\mR^\bot\right\rangle \\
+ \frac{1}{3} \left\langle T(\mR, J\mR^\bot),
 T(\mR^\bot, J\mR^\bot) \right\rangle
+ \frac{1}{4} \Big|T(\mR^\bot, J\mR^\bot)\Big|^2\\
=-\frac{1}{3}\left\langle (\nabla^{TY}_\cdot \dot{g}^{TY}_\cdot)_{(\mR,\mR)}
J\mR^\bot, J\mR^\bot\right\rangle
+ \frac{1}{6} \left\langle\nabla ^{TY}_{\mR}
(T(e^H_j,P^{T^HX} J e^H_i ))Z_j Z_i, J\mR^\bot\right\rangle \\
+\frac{1}{3}\left\langle \nabla ^{N_G}_{\mR^0} (A(e^0_j) e^0_i)Z^0_j Z^0_i
+  R^{TB}(\mR^\bot,\mR^0)\mR^0, \mR^\bot\right\rangle \\
-\frac{1}{4} \sum_j \left\langle T(\mR^0,
e^0_j), J\mR^\bot\right\rangle^2 +\frac{1}{6}\left\langle T(\mR,
e^0_j),J\mR^\bot\right\rangle
\left\langle T(\mR^0, e^0_j),J\mR^\bot\right\rangle\\
+  \frac{7}{12}\Big | T(\mR^\bot, J\mR^\bot)\Big|^2 +
\frac{1}{3}\left\langle  T(\mR^0, J\mR^\bot), T(\mR^\bot,
J\mR^\bot)\right\rangle.
\end{multline}

To get \eqref{3.48} from \eqref{a3.47}, we need to compute $\nabla
^{TY}_{e^H_k} (T(e^H_j,P^{T^HX} J e^H_i ))$.

For $W$ a section of $TX$, $U$ a section of $TB$, we have
\begin{multline}\label{aa3.48}
 \left\langle  \nabla ^{T^H X}_{e^H_k}P^{T^H X} W, U^H \right\rangle
=e^H_k \left\langle W, U^H \right\rangle
- \left\langle P^{T^H X} W, \nabla ^{TX}_{e^H_k}U^H\right\rangle \\
=\left\langle  P^{T^H X} \nabla ^{TX}_{e^H_k} W, U^H \right\rangle
+ \left\langle P^{TY} W, \nabla ^{TX}_{e^H_k}U^H \right\rangle .
\end{multline}

From \eqref{h4}, \eqref{aa3.48}, we get at $x_0$,
\begin{align}\label{aa3.49}
\nabla ^{T^H X}_{e^H_k}P^{T^H X} W
= P^{T^H X} \nabla ^{TX}_{e^H_k} W
 -\frac{1}{2}\left\langle T(e^H_k,e^H_l), P^{TY} W \right\rangle e^H_l.
\end{align}

Remark that $Je^{\bot,H}_i\in TY, Je_i^0\in T^HX$ only hold on
$P$.

 From \eqref{g3.29}, \eqref{3.11b}, \eqref{3.12a},
\eqref{33.36}, \eqref{b3.16} and \eqref{aa3.49},
\begin{align}\label{aa3.50}
\begin{split}
&(\nabla ^{T^H X}_{e^H_k}P^{T^H X} J e^{\bot,H}_i)_{x_0}
=JA(P^{TX_G}e_k) e^{\bot}_i -\frac{1}{2} JT(e_k,e^{\bot}_i)
- \frac{1}{2}\left\langle T(e_k,e_l),
J e^{\bot}_i \right\rangle e_l\\
&\hspace*{15mm}=  -\frac{1}{2} JT (e_k, e^{\bot}_i)
- \frac{1}{2}\left\langle T(e_k,e_l)- T(P^{TX_G}e_k,P^{TX_G}e_l),
J e^{\bot}_i \right\rangle e_l,\\
&(\nabla ^{T^H X}_{e_k}P^{T^H X} J e^{0}_i)_{x_0}
=  P^{T^HX} J \nabla ^{TX}_{e^H_k}e^{0,H}_i
= JA(e^0_i) P^{N_G}e_k -\frac{1}{2} JT (e_k, e^{0}_i)\\
&\hspace*{15mm}= -\frac{1}{2} JT (e_k, e^{0}_i)
+ \frac{1}{2} \left\langle JP^{N_G}e_k,T ( e^{0}_i, e^{0}_l) \right\rangle
e^{0}_l,\\
&(\nabla ^{TB}_{e_k} J_{x_0} e^0_i)_{x_0}
= A(J_{x_0} e^0_i) e_k
= -\frac{1}{2} JT (P^{TX_G}e_k, e^{0}_i)
+ \frac{1}{2} \left\langle JP^{N_G}e_k,T ( e^{0}_i, e^{0}_l)
\right\rangle  e^{0}_l.
\end{split}\end{align}

From \eqref{aa3.50}, we get at $x_0$ that
\begin{multline}\label{aa3.51}
 \left\langle\nabla ^{TY}_{\mR}
(T(e^H_j,P^{T^HX} J e^H_i ))Z_j Z_i, J\mR^\bot\right\rangle -
\left\langle\nabla ^{TY}_{\mR}
(T(e_j, J_{x_0} e^0_i ))Z_j Z^0_i, J\mR^\bot\right\rangle\\
=\left\langle
T(e_j,\nabla^{T^HX}_\mR P^{T^HX}Je_i^H-\nabla^{T^HX}_\mR
(J_{x_0}P^{T^HX}e_i)^H)Z_jZ_i,J\mR^\perp\right\rangle\\
=\left\langle T\left(\mR,-\frac{1}{2}JT(\mR,\mR^\perp)-\frac{1}{2}\left\langle
T(\mR,e_l)-T(\mR^0,P^{TX_G}e_l),J\mR^\perp\right\rangle
e_l\right),J\mR^\perp\right\rangle \\
+\left\langle T\left(e_j,- \frac{1}{2}JT(e_k,e_i^0)
+\frac{1}{2}JT(P^{TX_G}e_k,e_i^0)\right)Z_kZ_jZ_i^0,J\mR^\perp\right\rangle\\
= - \frac{1}{2}\sum_l \left\langle T(\mR,e_l),J\mR^\bot
\right\rangle ^2 +\frac{1}{2}\left\langle T(\mR,
e^0_l),J\mR^\bot\right\rangle \left\langle T(\mR^0,
e^0_l),J\mR^\bot\right\rangle.
\end{multline}

From  \eqref{a3.47} and \eqref{aa3.51},
$\mO''_2$ is the coefficient
of $t^2$ in the expansion of $|\frac{1}{t}\wi{\mu}|_{g^{TY}}^2 (tZ)$.

By \eqref{a3.43}, we get also the second equation of \eqref{3.48}.

The proof of Lemma \ref{t4.5} is completed.
\end{proof}

The following   is the main result of this Subsection.
\begin{thm} \label{t4.6} The following identities hold,
\begin{align}\label{a3.51}
\begin{split}
\mO_1=& 2\pi \sqrt{-1} \left\langle JT(\mR^\perp, e_i^0), \mR^\bot
\right\rangle\nabla_{0,e_i^0}+2\pi \sqrt{-1} \left\langle
JT(\mR, e_i^\perp),
\mR^\bot \right\rangle\nabla_{0,e_i^\perp}\\
& +\pi \sqrt{-1} \left\langle JT(\mR^0, e^\bot_i), e^\bot_i
\right\rangle - \left \langle JT(e^0_i, Je^0_j),\mR^\bot \right
\rangle
\nabla_{0,e^0_i} \nabla_{0,e^0_j}\\
&+ 4\pi^2\left\langle J T(\mR^\bot, J\mR^\bot), \mR^\bot \right\rangle
+4 \pi \sqrt{-1} \left\langle J \mR^\bot,
\wi{\mu}^E_{x_0}\right \rangle,\\
\mO_2 =& \mO^\prime_2+ 4\pi ^2 \mO''_2
-4\pi \sqrt{-1}\Big(\frac{1}{2}\left\langle T(\mR^\bot, J\mR^\bot),
\wi{\mu}^E_{x_0}\right \rangle
- \left\langle J  \mR^\bot,
\nabla^{TY}_{\mR}\wi{\mu}^E\right \rangle\Big)\\
&\hspace*{25mm}
- \left\langle \wi{\mu}^E_{x_0}, \wi{\mu}^E_{x_0}\right \rangle_{g^{TY}}.
\end{split}
\end{align}
\end{thm}
\begin{proof} By \eqref{3.11e},
 \begin{align}\label{a3.52}
\left\langle JT(\mR, e_i), e_i \right\rangle=
 \left\langle JT(\mR^0, e^\bot_i), e^\bot_i \right\rangle.
\end{align}
By \eqref{3.32}, \eqref{3.31} and \eqref{a3.52},
\begin{align}\label{aa3.52}
\begin{split}
&-\frac{2}{3}(\partial_\mR R^{L_B})_{x_0} (\mR,e_i) \nabla_{0,e_i}\\
&\hspace*{15mm}=2\pi \sqrt{-1}\Big( \left\langle JT(\mR^\perp, e_i^0), \mR^\bot
\right\rangle\nabla_{0,e_i^0}+ \left\langle
JT(\mR, e_i^\perp),
\mR^\bot \right\rangle\nabla_{0,e_i^\perp} \Big),\\
&-\frac{1}{3} (\partial_i R^{L_B})_{x_0} (\mR, e_i)
=\pi \sqrt{-1} \left\langle JT(\mR^0, e^\bot_i), e^\bot_i
\right\rangle.
\end{split}\end{align}
From \eqref{a3.17}, \eqref{a3.20}, \eqref{3.48} and \eqref{aa3.52},
we get \eqref{a3.51}.
\end{proof}


\subsection{Computation of the coefficient $\Phi_1$}\label{s4.5}
Recall that the operator $\cL^0_2$ is defined in \eqref{a3.17},
 $P_{\cL^\bot}$ is the orthogonal projection from $L^2(\bR^{n_0})$
onto $\Ker \cL^\bot$ and $P_{\cL}$ is the orthogonal projection from
$L^2(\bR^{2n-2n_0})$ onto $\Ker \cL$ as in \eqref{g16}.

 For $Z^\bot\in \bR^{n_0}$, set
 \begin{align}\label{b3.52}
\begin{split}&\Psi_{1,1}(Z^\bot) = \left((\cL^0_2)^{-1}P^{N^\bot}
\mO_1 (\cL^0_2)^{-1}P^{N^\bot}\mO_1 P^N \right)
\left((0,Z^\bot),(0,Z^\bot)\right),\\
&\Psi_{1,2}(Z^\bot) = -\left((\cL^0_2)^{-1}P^{N^\bot} \mO_2
P^N\right)\left((0,Z^\bot),(0,Z^\bot)\right),\\
&\Psi_{1,3}(Z^\bot) = \left((\cL^0_2)^{-1}P^{N^\bot}\mO_1 P^N
\mO_1 (\cL^0_2)^{-1}P^{N^\bot}\right)
\left((0,Z^\bot),(0,Z^\bot)\right),\\
&\Psi_{1,4}(Z^\bot) =\left(P^N \mO_1 (\cL^0_2)^{-2} P^{N^\bot}
\mO_1 P^N\right)\left((0,Z^\bot),(0,Z^\bot)\right),\\
&\wi{\Psi}_{1,1}(Z^\bot) = \left((\cL^0_2)^{-1}P^{N^\bot}P_{\cL^\bot}
\mO_1 (\cL^0_2)^{-1}P^{N^\bot}\mO_1 P^N \right)
\left((0,Z^\bot),(0,Z^\bot)\right), \\
&\wi{\Psi}_{1,2}(Z^\bot) = -\left((\cL^0_2)^{-1}P^{N^\bot}P_{\cL^\bot} \mO_2
P^N\right)\left((0,Z^\bot),(0,Z^\bot)\right),\\
&\Phi_{1,i} = \int_{\bR^{n_0}}\Psi_{1,i}(Z^\bot) dv_{N_G}(Z^\bot),
\quad {\rm for} \,  \,  \,  i=1,2, 3,4.
\end{split}\end{align}

\begin{prop}\label{t4.7}
The following two identities hold for $i=1,2$,
\begin{align}\label{3.532}
&\int_{\bR^{n_0}}\wi{\Psi}_{1,i}(Z^\bot)
 dv_{N_G}(Z^\bot)=\Phi_{1,i}.
\end{align}
\end{prop}
\begin{proof}  In fact, in our case, by \eqref{g17},
$P^N= P_{\cL}\otimes  P_{\cL^\bot}\otimes {\rm Id}_{E}$.

By \eqref{g15}, \eqref{g16},
 \begin{align}\label{b3.53}
\left((\cL^0_2)^{-1}P^{N^\bot}
\mO_2 P^N\right)(Z,(0,Z^{\prime \bot}))
= \left((\cL^0_2)^{-1}P^{N^\bot}
\mO_2 P_{\cL}(\cdot, 0) G^\bot\right) (Z) G^\bot(Z^{\prime \bot}).
\end{align}

From Theorem \ref{t3.4}, \eqref{b3.53},
 \begin{multline}\label{b3.54}
\Phi_{1,2}
=\left\langle \left(-(\cL^0_2)^{-1}P^{N^\bot}
\mO_2 P_{\cL}(\cdot, 0) G^\bot\right) (0, Z^\bot),
 G^\bot(Z^{\bot})\right \rangle _{L^2(\bR^{n_0})}\\
= \left\langle \left(-(\cL^0_2)^{-1}P^{N^\bot}P_{\cL^\bot}
\mO_2 P_{\cL}(\cdot, 0) G^\bot\right) (0, Z^\bot),
 G^\bot(Z^{\bot})\right \rangle _{L^2(\bR^{n_0})}\\
=\int_{\bR^{n_0}}\wi{\Psi}_{1,2}(Z^\bot) dv_{N_G}(Z^\bot).
\end{multline}

In the same way, we get \eqref{3.532} for $i=1$.
\end{proof}

 Note that the restriction of $\norm{\,\cdot\,}_{t,0}$ in \eqref{u4}
on $\cC^\infty(\bR^{2n-n_0}, E_{G,x_0})$
does not depend on $t$ and we denote it by $\norm{\,\cdot\,}_{0}$.

Since $\cL^t_2$ in \eqref{a3.17} is a self-adjoint  elliptic
operator with respect to $\norm{\,\cdot\,}_{0}$ as we conjugated
the operator with $\kappa^{1/2}$, $\cL^0_2$ and $\mO_r$  are also
formally  self-adjoint with respect to $\norm{\,\cdot\,}_{0}$.
Thus in the right hand side of \eqref{g22.7},
the third and fourth terms are the adjoints of the first two terms.

From \eqref{g22.7}, \eqref{6.560} and \eqref{b3.52},  we get
 \begin{align}\label{a3.53}
\Phi_1= \Phi_{1,1}+\Phi_{1,2} + (\Phi_{1,1}+\Phi_{1,2})^*
+ \Phi_{1,3}-\Phi_{1,4}.
\end{align}

From \eqref{b3.52}, \eqref{3.532}, \eqref{a3.53}, we learn that in
order to compute $\Phi_1$, we only need to evaluate $\Psi_{1,i}$
and $\wi{\Psi}_{1,i}$ $(i\in\{1,2,3,4\})$.

\begin{lemma} \label{t4.8} The following identity holds,
\begin{align} \label{3.65}
\wi{\Psi}_{1,1}(Z^{\bot})=
- \frac{1}{8\pi} \Big|T(\tfrac{\partial}{\partial
\ov{z}^0_j}, e^\bot_k)\Big|^2 P_{\cL^\bot}(Z^\bot,Z^\bot) .
\end{align}
\end{lemma}
\begin{proof}

 Recall that the operators $b_i$, $b_i^+$, $b_j^\perp$ and
 $b_j^{\perp +}$ have been defined in \eqref{g6}. In particular,
 by \eqref{b3.11}, one has
 \begin{align}\label{3.533}
 4\pi Z_j^\perp=b_j^\perp+b_j^{\perp +}, \quad
\nabla_{0, e^\bot_j}
= \tfrac{\partial}{\partial Z^\bot_j}
= \tfrac{1}{2} ( b^{\bot +}_j-  b^{\bot}_j).
 \end{align}

By \eqref{g6}, \eqref{g7} and \eqref{3.533}, set
\begin{align} \label{3.53}
\begin{split}
&B^\bot_{jk} = (4\pi)^2 Z^\bot_j Z^\bot_k = b^{\bot+}_jb^{\bot+}_k
+ b^{\bot}_k b^{\bot+}_j + b^{\bot}_j b^{\bot+}_k
+ b^{\bot}_j b^{\bot}_k + 4\pi \delta_{jk},\\
&B^\bot_{ijk}= b^{\bot}_i b^{\bot}_j b^{\bot}_k
+ 3b^{\bot}_i b^{\bot}_j b^{\bot+}_k+ 3b^{\bot}_i b^{\bot+}_j b^{\bot+}_k
+b^{\bot+}_i b^{\bot+}_j b^{\bot+}_k.
\end{split}\end{align}

If $a_{ijk}$ is symmetric on $i,j,k$, then by \eqref{g6},
\eqref{g7}, \eqref{3.533} and \eqref{3.53}, one verifies
\begin{align} \label{3.54}
a_{ijk}(4\pi)^3 Z^\bot_i Z^\bot_j Z^\bot_k
=a_{ijk}B^\bot_{ijk} +12\pi a_{ijj} (b^{\bot}_i+ b^{\bot+}_i).
\end{align}

By \eqref{g7}, \eqref{3.11e},  \eqref{34.21}, \eqref{3.533},
\eqref{3.53} and the fact that $T(\, ,\,)$ is anti-symmetric, we
get
\begin{multline}\label{b3.55}
2\pi \left\langle
JT(\mR^\bot, e_i^\perp),
\mR^\bot \right\rangle\nabla_{0,e_i^\perp}
= \frac{1}{16\pi} \wi{\mT}_{jik} B^\bot_{jk} (b^{\bot+}_i- b^{\bot}_i)\\
=  \frac{1}{16\pi} \wi{\mT}_{jik}
\left[ (b^{\bot}_j b^{\bot+}_k+ b^{\bot}_j b^{\bot}_k)b^{\bot+}_i
- (b^{\bot+}_jb^{\bot+}_k+ b^{\bot}_k b^{\bot+}_j + b^{\bot}_j b^{\bot+}_k)
b^{\bot}_i\right]\\
=-\frac{1}{8\pi} \wi{\mT}_{ijk}(b^{\bot}_j b^{\bot+}_k
+ b^{\bot}_j b^{\bot}_k)b^{\bot+}_i.
\end{multline}

By Theorem \ref{t4.2}, Remark \ref{r4.3},
 \eqref{g7}, \eqref{0g6}, \eqref{34.21}, \eqref{a3.51},
 \eqref{3.53}-\eqref{b3.55},
we  can reformulate $\mO_1$ as follows by using the creation and
annihilation operators introduced in \eqref{g6},
\begin{multline} \label{3.56}
\mO_1= -\frac{\sqrt{-1}}{8\pi} \left\langle
JT(\tfrac{\partial}{\partial z^0_i}, e^\bot_j),
e^\bot_k\right\rangle B^\bot_{jk} b^{+}_i
+ b_i\frac{\sqrt{-1}}{8\pi}\left\langle
JT(\tfrac{\partial}{\partial \ov{z}^0_i}, e^\bot_j),
e^\bot_k \right\rangle B^\bot_{jk}\\
+ \frac{\sqrt{-1}}{4} \left\langle JT(\mR^0, e^\bot_i),
e^\bot_j \right\rangle (b^{\bot +}_ib^{\bot +}_j- b^{\bot}_ib^{\bot}_j)
-\frac{\sqrt{-1} }{8\pi} \wi{\mT}_{ijk}(b^{\bot}_j b^{\bot+}_k
+ b^{\bot}_j b^{\bot}_k)b^{\bot+}_i\\
-\frac{\sqrt{-1}}{4\pi}  \left\langle JT(\tfrac{\partial}{\partial z^0_i},
\tfrac{\partial}{\partial \ov{z}^0_j}), e^\bot_k \right\rangle
(b^{\bot +}_k+ b^{\bot }_k) (2 b_j b^+_i +4\pi \delta_{ij})
+\sqrt{-1}  \left\langle J e^\bot_j,
\wi{\mu}^E_{x_0}\right \rangle (b^{\bot +}_j+ b^{\bot }_j)\\
+\frac{1}{16\pi} \left\langle JT(e^\bot_i, Je^\bot_j), e^\bot_k \right\rangle
[B^\bot_{ijk} +12\pi \delta_{ik} (b^{\bot +}_j+ b^{\bot }_j)]\\
= -\frac{\sqrt{-1}}{8\pi}\mT_{jk}(\tfrac{\partial}{\partial z^0_i})
B^\bot_{jk} b^{+}_i
+\frac{\sqrt{-1}}{8\pi} \mT_{jk}(\tfrac{\partial}{\partial \ov{z}^0_i})
b_i B^\bot_{jk}
+\frac{\sqrt{-1}}{4}  \mT_{ij}(\mR^0)
(b^{\bot +}_ib^{\bot +}_j- b^{\bot}_ib^{\bot}_j)\\
+\sqrt{-1}  \left\langle J e^\bot_j,
\wi{\mu}^E_{x_0}\right \rangle (b^{\bot +}_j+ b^{\bot }_j)
-\frac{\sqrt{-1}}{4\pi}  \left\langle JT(\tfrac{\partial}{\partial z^0_i},
\tfrac{\partial}{\partial \ov{z}^0_j}), e^\bot_k \right\rangle
(b^{\bot +}_k+ b^{\bot }_k) (2 b_j b^+_i +4\pi \delta_{ij})\\
-\frac{\sqrt{-1} }{8\pi} \wi{\mT}_{ijk}(b^{\bot}_j b^{\bot+}_k
+ b^{\bot}_j b^{\bot}_k)b^{\bot+}_i+\frac{1}{16\pi}\mT_{ijk}
[B^\bot_{ijk} +12\pi \delta_{ik} (b^{\bot +}_j+ b^{\bot }_j)] .
\end{multline}

\comment{
By \eqref{g6} and \eqref{g16},
\begin{align} \label{3.57}
&b^+_iP^N =b^{\bot+}_iP^N =0\,,\quad \,
b^{\bot}_i P^N(Z,Z^{\prime})=2a^{\bot}_i Z^\bot_i P^N(Z,Z^{\prime}),\\
&b_iP^N(Z,Z^{\prime})=a_i(\ov{z}^0_i-\ov{z}^{\prime 0}_i)P^N(Z,Z^{\prime}).
\nonumber
\end{align}
We learn from \eqref{3.57} that for any polynomial $g(Z^\bot)$ in $Z^\bot$
we can write $g(Z^\bot)P^N(Z,Z^{\prime})$ as
sums of $g_{\beta^\bot}(b^\bot)^{\beta^\bot} P^N(Z,Z^{\prime})$
with constants $g_{\beta^\bot}$.
By Theorem \ref{t3.4},
\begin{equation}\label{3.59}
 P_{\cL^\bot} (b^\bot)^{\alpha^\bot}   g(Z^\bot)P^N=0\,,
\quad\text{for $|\alpha^\bot|>0$}.
\end{equation}
}

From Theorem \ref{t3.4}, \eqref{g3.57}, \eqref{3.53}, \eqref{3.56}
and $a_i=a^+_i=2\pi$, we get
\begin{multline} \label{3.58}
\left((\cL^0_2)^{-1}\mO_1 P^N\right)(Z, Z^{\prime})=
\sqrt{-1} \Big \{ \frac{b_l}{8\pi}
 \mT_{kk}(\tfrac{\partial}{\partial \ov{z}^0_l})
+ \left\langle J e^\bot_k,
\wi{\mu}^E_{x_0}\right \rangle \frac{b^{\bot}_k}{4\pi}\\
-  \left\langle JT(\tfrac{\partial}{\partial z^0_l},
\tfrac{\partial}{\partial \ov{z}^0_l}), e^\bot_k \right\rangle
\frac{b^{\bot}_k}{4\pi}
- \frac{b^{\bot}_lb^{\bot}_k}{32\pi}
\mT_{kl}(z^0+\ov{z}^{\prime 0}) \\
-\frac{\sqrt{-1}}{16\pi} \mT_{klm}
\Big[\frac{b^{\bot}_mb^{\bot}_lb^{\bot}_k}{12\pi}
+3 b^{\bot}_k \delta_{lm}\Big]\Big\}P^N(Z,Z^{\prime}).
\end{multline}

By  Theorem \ref{t3.4}, \eqref{g3.59}, \eqref{3.53} and
\eqref{3.56},
\begin{multline} \label{3.60}
P^{N^\bot} P_{\cL^\bot}\mO_1=\sqrt{-1} P^{N^\bot} P_{\cL^\bot}
\Big\{  -\frac{1}{2}\mT_{jj}(\tfrac{\partial}{\partial z^0_i}) b_i^+
+\frac{1}{2} \mT_{jj}(\tfrac{\partial}{\partial \ov{z}^0_i}) b_i
+\left\langle J e^\bot_j, \wi{\mu}^E_{x_0}\right \rangle  b^{\bot +}_j\\
-\frac{1}{4\pi}  \left\langle JT(\tfrac{\partial}{\partial z^0_i},
\tfrac{\partial}{\partial \ov{z}^0_j}), e^\bot_{j'} \right\rangle
b^{\bot +}_{j'} (2 b_j b^+_i +4\pi \delta_{ij})\\
+ \frac{1}{4} \Big(\mT_{jj'}(\mR^0)
-  \mT_{jj'}(\tfrac{\partial}{\partial z^0_i})  \frac{b^{+}_i}{2\pi}
+  \mT_{jj'}(\tfrac{\partial}{\partial \ov{z}^0_i})\frac{b_i}{2\pi}\Big)
 b^{\bot +}_jb^{\bot +}_{j'}\\
- \frac{\sqrt{-1}}{16\pi}\mT_{ijj'}
[b^{\bot +}_ib^{\bot +}_jb^{\bot +}_{j'} +12\pi \delta_{i{j'}} b^{\bot +}_j].
\end{multline}

In the following equation, by \eqref{g7}, \eqref{g3.57},
\eqref{g3.59}, we only need to pair the terms in \eqref{3.58} and
\eqref{3.60} which have the same length on $b^{\bot +}_j$ and
$b^\bot_j$, and the total degree on $b_i,b_i^+, z^0, \ov{z}^0$
should not be zero. Thus by  \eqref{g7}, \eqref{g3.57},
\eqref{3.58} and \eqref{3.60},
\begin{multline} \label{3.61}
\left(P^{N^\bot} P_{\cL^\bot}\mO_1(\cL^0_2)^{-1}\mO_1 P^N\right)(Z, (0,Z^{\prime\bot}))=
 \left\{P^{N^\bot}\Big[- \frac{1}{16\pi} \Big(\sum_{ij} b_i \mT_{jj}(\tfrac{\partial}{\partial \ov{z}^0_i}) \Big)^2\right.\\
+  \frac{1}{128\pi} \Big( \mT_{jj'}(\mR^0)
+ \frac{b_i}{2\pi} \mT_{jj'}(\tfrac{\partial}{\partial \ov{z}^0_i})\Big)
 b^{\bot +}_jb^{\bot +}_{j'} \left.
\cdot b^{\bot}_lb^{\bot}_k \mT_{kl}(z^0)
 \Big]P^N\right\}(Z, (0,Z^{\prime\bot})).
\end{multline}

From \eqref{g7}, \eqref{g3.57}, \eqref{3.11d}, \eqref{34.21},
\eqref{3.61} and $a_i=a_i^+=2\pi$, one gets
\begin{multline} \label{3.62}
\left(P^{N^\bot} P_{\cL^\bot}\mO_1(\cL^0_2)^{-1}\mO_1
P^N\right)(Z, (0,Z^{\prime\bot}))=
 \left\{P^{N^\bot}\Big[- \frac{1}{16\pi} \Big(\sum_{ij} b_i \mT_{jj}(\tfrac{\partial}{\partial \ov{z}^0_i}) \Big)^2\right.\\
+ \left. \frac{1}{8} \left\langle 2\pi J T(\mR^0, e^\bot_l)
+b_i J T( \tfrac{\partial}{\partial \ov{z}^0_i}, e^\bot_l),
 JT(z^0, e^\bot_l) \right\rangle
\Big]P^N\right\}(Z, (0,Z^{\prime\bot})).
\end{multline}

Set $P^\bot_{\cL}= {\rm Id}_{L^2(\bR^{2n-2n_0})}-P_{\cL}$.

Let $h_i(Z^0)$ (resp. $F(Z^0)$) be polynomials in $Z^0$ with
degree $1$ (resp. $2$) and $a_{ij}\in \bC$.

By Theorem \ref{t3.4}, \eqref{g7} and \eqref{g3.57},
\begin{align}\label{3.63}
&\left(F(Z^0)P_{\cL}\right)(Z^0,0) = \Big(\frac{1}{2}
\frac{\partial ^2 F}{\partial z^0_i\partial z^0_j} z^0_iz^0_j +
\frac{\partial ^2 F}{\partial z^0_i\partial \ov{z}^0_j} z^0_i
\frac{b_j}{a_j} + \frac{1}{2}\frac{\partial ^2 F}{\partial
\ov{z}^0_i\partial \ov{z}^0_j}
 \frac{b_ib_j}{a_ia_j} \Big) P_{\cL}(Z^0,0).
\end{align}
By Theorem \ref{t3.4}, (\ref{g6}),  (\ref{g7}), \eqref{g16},
(\ref{g3.57}), (\ref{3.63}) and $a_j=2\pi$, we have
\begin{equation}\label{3.64}
\begin{split}
&(P^{\bot}_{\cL} FP_{\cL})(0,0)
=-\frac{1}{\pi} \frac{\partial ^2 F}{\partial z^0_i\partial \ov{z}^0_i},\\
&\left(\cL^{-1} P^{\bot}_{\cL} a_{ij}b_ib_j P_{\cL}\right)(0,0) =
\left(\cL^{-1} P^{\bot}_{\cL} h_iP_{\cL}\right)(0,0) = 0,\\
&\left(\cL^{-1} P^{\bot}_{\cL} h_ib_i P_{\cL}\right)(0,0)=
\left(\cL^{-1} P^{\bot}_{\cL}b_ih_i P_{\cL}\right)(0,0)
=-\frac{1}{2\pi}\frac{\partial h_i}{\partial z^0_i},\\
&\left(\cL^{-1} P^{\bot}_{\cL} F P_{\cL}\right)(0,0)
= -\frac{1}{4\pi ^2}\frac{\partial ^2 F}{\partial z^0_i\partial \ov{z}^0_i},\\
&\left(\cL^{-1} P^{\bot}_{\cL}b_i F b_j P_{\cL}\right)(0,0)=
-\left(\cL^{-1} P^{\bot}_{\cL} b_i b_j F  P_{\cL}\right)(0,0)
= -\frac{1}{2\pi} \frac{\partial ^2 F}{\partial z^0_i\partial z^0_j},\\
&\left(\cL^{-1} P^{\bot}_{\cL}F b_i b_j P_{\cL}\right)(0,0)
= -\frac{3}{2\pi} \frac{\partial ^2 F}{\partial z^0_i\partial z^0_j},\\
&\left(\cL^{-1} P^{\bot}_{\cL} \Big(\sum_i b_i h_i \Big)^2
P_{\cL}\right)(0,0)= -\frac{1}{2\pi} \Big(\frac{\partial
h_i}{\partial z^0_j} \frac{\partial h_j}{\partial z^0_i}
-\Big(\sum_i \frac{\partial h_i}{\partial z^0_i}\Big)^2\Big).
\end{split}
\end{equation}

Finally by \eqref{3.532}, \eqref{3.62}, \eqref{3.64}
and $\cL^0_2= \cL +\cL^\bot$ , we get \eqref{3.65}.
\end{proof}

\begin{lemma}\label{at4.9} The following identity holds,
\begin{align}\label{a3.54}
\Phi_{1,3}= \Phi_{1,4}.
\end{align}
\end{lemma}
\begin{proof} Let $ \mF_2\in T^*_{x_0}X_G$ with values
in real polynomials on $Z^\bot$ with even degree,
$\mF_1\in N^*_{G,x_0}$, $\mF_3(Z^\bot)$ a polynomial on $Z^\bot$
 with odd degree, be defined by
\begin{align}\label{a3.55}
\begin{split}
&\mF_1 (e^\bot_k)=\sqrt{-1} \left\langle J e^\bot_k,
\wi{\mu}^E_{x_0}\right \rangle
- \sqrt{-1} \left\langle JT(\tfrac{\partial}{\partial z^0_l},
\tfrac{\partial}{\partial \ov{z}^0_l}), e^\bot_k \right\rangle
+ \frac{3}{4} \mT_{llk},\\
& \mF_2 (\cdot, Z^\bot)P^N(Z,Z^\prime)=
\Big( \mT_{kl}(\cdot)
\frac{b^{\bot}_lb^{\bot}_k}{32\pi}P^N\Big)(Z,Z^\prime),\\
&\mF_3(Z^\bot)P^N(Z,Z^\prime)= \frac{1}{16\pi}
\Big(\mT_{klm}
\frac{b^{\bot}_mb^{\bot}_lb^{\bot}_k}{12\pi}P^N\Big)(Z,Z^\prime).
\end{split}\end{align}
Then from \eqref{g3.57}, \eqref{3.58} and \eqref{a3.55},
\begin{multline}\label{a3.56}
\left((\cL^0_2)^{-1}\mO_1 P^N\right)(Z, Z^{\prime})=
\Big(\frac{\sqrt{-1}}{4} \mT_{kk}(\ov{z}^0-\ov{z}^{\prime 0})
- \sqrt{-1} \mF_2(z^0+\ov{z}^{\prime 0}, Z^\bot) \\
+ \Big(\mF_1 + \mF_3\Big)(Z^\bot)\Big)P^N(Z, Z^{\prime}).
\end{multline}

Observe that $\mF_i(Z^\bot)^*= \mF_i(Z^\bot)$ for $i=1,3$, thus
from \eqref{3.58} and \eqref{a3.56},
\begin{multline}\label{a3.57}
\left(P^N\mO_1(\cL^0_2)^{-1}\right)(Z^{\prime}, Z)=
\left(\left((\cL^0_2)^{-1}\mO_1 P^N\right)(Z, Z^{\prime})\right)^*\\
= \Big(-\frac{\sqrt{-1}}{4} \mT_{kk}(z^0- z^{\prime 0})
+\sqrt{-1} \mF_2(\ov{z}^0+z^{\prime 0}, Z^\bot)
+ \Big(\mF_1 + \mF_3\Big)(Z^\bot)\Big)P^N( Z^{\prime},Z).
\end{multline}

For $h_1(z^0), h_2(\ov{z}^0)$ two linear functions on $z^0,\ov{z}^0$,
by Theorem \ref{t3.4}, \eqref{g3.57},
\begin{align}\label{a3.58}
(P_{\cL} h_1(z^0) h_2(\ov{z}^0)P_{\cL})(0,0)
= \Big( P_{\cL} h_1(z^0) \frac{\partial h_2 }{\partial \ov{z}^0_i}
\frac{b_i}{2\pi} P_{\cL}  \Big)(0,0)
= \frac{1}{\pi} \frac{\partial h_1}{\partial z^0_i}
\frac{\partial h_2 }{\partial \ov{z}^0_i} .
\end{align}

From \eqref{b3.52}, \eqref{a3.56}-\eqref{a3.58},
\begin{align}\label{a3.59}
\Psi_{1,3}(Z^\bot)= \Big[ \Big((\mF_1+\mF_3)(Z^\bot)\Big)^2
+ \frac{1}{\pi} \Big|\frac{1}{4}\sum_k \mT_{kk}(\tfrac{\partial}{\partial \ov{z}^0_i})
+ \mF_2( \tfrac{\partial}{\partial \ov{z}^0_i}, Z^\bot)\Big|^2 \Big]
G^\bot(Z^\bot)^2.
\end{align}

By Theorem \ref{t3.4}, \eqref{g15}, \eqref{a3.54}, $\mF_j G^\bot$,
$(j=1,3)$, $\mF_2(\tfrac{\partial}{\partial \ov{z}^0_i},\cdot)
G^\bot$ are eigenfunctions of $\cL^\bot$ with  eigenvalues $4\pi
j$, $8\pi$,
 thus they are orthogonal to each other.

From \eqref{b3.52}, \eqref{a3.56}-\eqref{a3.58}, we have
\begin{multline}\label{a3.60}
\Psi_{1,4}(Z^\bot)=G^\bot(Z^\bot)^2 \int_{\bR^{n_0}}
\Big\{ \Big((\mF_1 G^\bot)(Z^{\prime \bot}) \Big)^2
+ \left((\mF_3 G^\bot)(Z^{\prime \bot})\right)^2 \\
 + \frac{1}{16\pi}
\Big|\sum_k \mT_{kk}(\tfrac{\partial}{\partial \ov{z}^0_i})
G^\bot\Big|^2 (Z^{\prime \bot})
+ \frac{1}{\pi}  \left|\mF_2(\tfrac{\partial}{\partial \ov{z}^0_i},
\cdot)G^\bot\right|^2 (Z^{\prime \bot})\Big\} dv_{N_G}(Z^{\prime \bot}).
\end{multline}

From \eqref{g15}, \eqref{b3.52}, \eqref{a3.59}, \eqref{a3.60} and
the above discussion, we get \eqref{a3.54}.
\end{proof}


Now we need to compute the contribution from
$-(\cL^0_2)^{-1} P^{N^\bot}\mO_2 P^N$.
\begin{lemma}\label{t4.9} The following identity holds,
\begin{multline}\label{a3.64}
\wi{\Psi}_{1,2}(Z^\bot)= \Big\{ \frac{1}{2\pi}
 \left\langle  R^{TX_G}(\tfrac{\partial}{\partial z^0_j},
\tfrac{\partial}{\partial \ov{z}^0_i})\tfrac{\partial}{\partial z^0_i},
\tfrac{\partial}{\partial \ov{z}^0_j} \right\rangle
+ \frac{1}{48\pi} \left\langle
R^{TB}(e^\bot_k, \tfrac{\partial}{\partial z^0_j})e^\bot_k,
\tfrac{\partial}{\partial \ov{z}^0_j}\right\rangle\\
 + \frac{1}{96\pi}\Big| T(\tfrac{\partial}{\partial z^0_i},
\tfrac{\partial}{\partial \ov{z}^0_j})\Big|^2
- \frac{\sqrt{-1}}{16\pi} \left\langle T(e^\bot_k,Je^\bot_k),
T(\tfrac{\partial}{\partial z^0_j},
\tfrac{\partial}{\partial \ov{z}^0_j}) \right\rangle
+ \frac{13}{192\pi}
\Big|T(e^\bot_k,\tfrac{\partial}{\partial \ov{z}^0_j})\Big|^2\\
+\frac{\sqrt{-1}}{96\pi} \left\langle
11\,   \nabla ^{TY}_{\tfrac{\partial}{\partial z^0_j}}
(T(e^\bot_k,\tfrac{\partial}{\partial \ov{z}^0_j}))
+4\,  \nabla ^{TY}_{\tfrac{\partial}{\partial \ov{z}^0_j}}
(T(e^\bot_k, \tfrac{\partial}{\partial z^0_j}))
+7\, \nabla ^{TY}_{e^\bot_k} (T(\tfrac{\partial}{\partial z^0_j},
\tfrac{\partial}{\partial \ov{z}^0_j})), Je^\bot_k \right\rangle\\
-\frac{2}{3\pi}\nabla_{\tfrac{\partial}{\partial z^0_j}}
\nabla_{\tfrac{\partial}{\partial \ov{z}^0_j}} \log h
+\frac{1}{2\pi}R^{E_B}(\tfrac{\partial}{\partial z^0_j},
\tfrac{\partial}{\partial \ov{z}^0_j})
\Big \} P_{\cL^\bot}(Z^\bot,Z^\bot).
\end{multline}
\end{lemma}
\begin{proof}
By \eqref{g7}, \eqref{0g6}, \eqref{g3.57}, \eqref{a3.18} and \eqref{3.533},
\begin{align}\label{a3.67}
 I_1 P^N= \Big\{
\frac{1}{2} b^{\bot}_i B(Z,\tfrac{\partial}{\partial Z^\bot_i})
+ b_j B(Z,\tfrac{\partial}{\partial \ov{z}^0_j})
+\tfrac{\partial}{\partial z^0_j}
\Big(B(Z,\tfrac{\partial}{\partial \ov{z}^0_j})\Big)
-\tfrac{\partial}{\partial \ov{z}^0_j}
\Big(B(Z,\tfrac{\partial}{\partial z^0_j}) \Big)
\Big\}P^N.
\end{align}

By \eqref{g3.59} and \eqref{a3.67},
\begin{align}\label{3.67}
P_{\cL^\bot} I_1 P^N=
P_{\cL^\bot} \Big\{ b_j B(Z,\tfrac{\partial}{\partial \ov{z}^0_j})
+ \tfrac{\partial}{\partial z^0_j}
\Big(B(Z,\tfrac{\partial}{\partial \ov{z}^0_j})\Big)
-\tfrac{\partial}{\partial \ov{z}^0_j}
\Big(B(Z,\tfrac{\partial}{\partial z^0_j}) \Big)
\Big\}P^N.
\end{align}

By \eqref{a3.32}, and observe that from Theorem \ref{t3.4},
 only the monomials which have even degree on $Z^\bot$
and $\nabla_{e ^\bot_j}$, and which have also strictly positive
 degree on $Z^0$ and $\nabla_{0,e ^0_j}$,
have contributions in $P^{N^\bot}P_{\cL^\bot} I_1 P^N$.

By Remark \ref{r4.3}, \eqref{g3.59} and \eqref{a3.32},
\begin{multline}\label{3.68}
P^{N^\bot}P_{\cL^\bot}\Big( \tfrac{\partial}{\partial z^0_j}
\Big(B(Z,\tfrac{\partial}{\partial \ov{z}^0_j})\Big)
-\tfrac{\partial}{\partial \ov{z}^0_j}
\Big(B(Z,\tfrac{\partial}{\partial z^0_j}) \Big) \Big)P^N
=-\pi \sqrt{-1} P^{N^\bot}P_{\cL^\bot} \\
\frac{1}{6} \Big\{\tfrac{\partial}{\partial z^0_j}
\left\langle R^{TX_G}(\mR^0,J\mR^0)\mR^0,
\tfrac{\partial}{\partial \ov{z}^0_j}\right\rangle
-\tfrac{\partial}{\partial \ov{z}^0_j}
\left\langle R^{TX_G}(\mR^0,J\mR^0)\mR^0,
\tfrac{\partial}{\partial z^0_j}\right\rangle \Big\}P^N\\
= - \frac{\pi}{3} P^{N^\bot} \left\langle 2 R^{TX_G}(z^0,\ov{z}^0)
\tfrac{\partial}{\partial z^0_j}
+ R^{TX_G}(\tfrac{\partial}{\partial z^0_j},\mR^0)z^0
+   R^{TX_G}(\tfrac{\partial}{\partial z^0_j},\ov{z}^0)\mR^0,
\tfrac{\partial}{\partial \ov{z}^0_j}\right\rangle  P^N.
\end{multline}

By \eqref{a3.17}, \eqref{3.64} and \eqref{3.68},
\begin{multline}\label{3.69}
-\left((\cL^0_2)^{-1}P^{N^\bot}P_{\cL^\bot}
\Big(\tfrac{\partial}{\partial z^0_j}
\Big(B(Z,\tfrac{\partial}{\partial \ov{z}^0_j}) \Big)
-\tfrac{\partial}{\partial \ov{z}^0_j}
\Big(B(Z,\tfrac{\partial}{\partial z^0_j}) \Big)  \Big)
P^N\right)((0,Z^\bot),(0,Z^\bot))\\
=-\frac{1}{6\pi}
 \left\langle  R^{TX_G}(\tfrac{\partial}{\partial z^0_i},
\tfrac{\partial}{\partial \ov{z}^0_i})
\tfrac{\partial}{\partial z^0_j}
+ R^{TX_G}(\tfrac{\partial}{\partial z^0_j},
\tfrac{\partial}{\partial \ov{z}^0_i})
\tfrac{\partial}{\partial z^0_i},
\tfrac{\partial}{\partial \ov{z}^0_j}\right\rangle
P_{\cL^\bot}(Z^\bot,Z^\bot).
\end{multline}

Observe that if $Q$ is an odd degree monomial on
$b_j,b^+_j, z^0_j,\ov{z}^0_j$, then
\begin{align}\label{a3.68}
\left(QP^N\right)\left(\left(0,Z^\bot\right), \left(0, Z^{\prime \bot}\right)\right)=0.
\end{align}
By using this observation and \eqref{a3.32}, we get
\begin{multline}\label{3.70}
-\left((\cL^0_2)^{-1}P^{N^\bot}
b_jB(Z,\tfrac{\partial}{\partial \ov{z}^0_j}) P^N\right)((0,Z^\bot), (0, Z^{\prime \bot}))\\
= \pi \sqrt{-1} \left\{(\cL^0_2)^{-1}P^{N^\bot}b_j
\Big[ \frac{1}{6} \left\langle R^{TX_G}(\mR^0,J\mR^0)\mR^0 ,
\tfrac{\partial}{\partial \ov{z}^0_j}\right\rangle \right.\\
-\frac{5}{4} \left\langle  \nabla ^{TY}_{\mR^0}
(T(e^\bot_k,\tfrac{\partial}{\partial \ov{z}^0_j}))Z^\bot_k
+\nabla ^{TY}_{\mR^\bot} (T(e^0_k,\tfrac{\partial}{\partial
\ov{z}^0_j}))Z^0_k,
 J\mR^\bot \right\rangle\\
+ \left\langle \frac{1}{2}R^{TB}(\mR^\bot,J\mR^0)\mR^\bot
+\sqrt{-1} R^{TB}(\mR^\bot, \mR^0)\mR^\bot,
\tfrac{\partial}{\partial \ov{z}^0_j}\right\rangle\\
-\frac{3}{8}\sqrt{-1} \left\langle J\mR^\bot,T(\mR^0, e^0_i)\right\rangle
\left\langle J\mR^\bot,
T( e^0_i,\tfrac{\partial}{\partial \ov{z}^0_j} )\right\rangle \\
- \frac{1}{8}\left\langle T(\mR^\bot,J\mR^0),
 T(\mR^\bot,\tfrac{\partial}{\partial \ov{z}^0_j}) \right\rangle
+ \frac{1}{2} \left\langle  T(\mR^\bot,J\mR^\bot),
T(\mR^0,\tfrac{\partial}{\partial \ov{z}^0_j}) \right\rangle \\
\left.
-\frac{1}{8} \left\langle JT(\tfrac{\partial}{\partial \ov{z}^0_j},J\mR^0),
e^\bot_j \right\rangle
\left\langle J \mR^\bot,T(\mR^\bot,e^\bot_j) \right\rangle
\Big] P^N\right\}((0,Z^\bot), (0, Z^{\prime \bot})).
\end{multline}

From \eqref{0g4}, \eqref{g3.57}, \eqref{3.11b} and \eqref{3.53}, we have
\begin{subequations}
\begin{align}  \label{a3.70}
& \left\langle  T(\tfrac{\partial}{\partial z^0_j}, e^0_i),
T(e^0_i, \tfrac{\partial}{\partial \ov{z}^0_j})\right\rangle
= -2 \Big|T(\tfrac{\partial}{\partial z^0_i},
\tfrac{\partial}{\partial \ov{z}^0_j})\Big|^2,\\
&P_{\cL^\bot} Z^\bot_k Z^\bot_l P_{\cL^\bot}
=\frac{\delta_{kl}}{4\pi} P_{\cL^\bot}. \label{a3.70b}
\end{align}
\end{subequations}

By  \eqref{g3.57}, \eqref{3.11e}, \eqref{3.64}, \eqref{3.70} and \eqref{a3.70},
\begin{multline}\label{3.71}
-\left((\cL^0_2)^{-1}P^{N^\bot}P_{\cL^\bot}b_j
B(Z,\tfrac{\partial}{\partial \ov{z}^0_j}) P^N\right)((0,Z^\bot),(0,Z^\bot))\\
=\left\{(\cL^0_2)^{-1} P^{N^\bot}b_j
\Big[\frac{\pi}{3}  \left\langle R^{TX_G}(z^0,\ov{z}^0)\mR^0,
\tfrac{\partial}{\partial \ov{z}^0_j}\right\rangle  \right.\\
-\frac{5\sqrt{-1}}{16} \left\langle
\nabla ^{TY}_{\mR^0} (T(e^\bot_k,\tfrac{\partial}{\partial \ov{z}^0_j}))
+\nabla ^{TY}_{e^\bot_k} (T(e^0_i,\tfrac{\partial}{\partial \ov{z}^0_j}))Z^0_i,
 Je^\bot_k \right\rangle\\
+\frac{1}{8}  \left\langle \sqrt{-1} R^{TB}(e^\bot_k, J \mR^0)e^\bot_k
-2 R^{TB}(e^\bot_k, \mR^0)e^\bot_k,
\tfrac{\partial}{\partial \ov{z}^0_j}\right\rangle\\
+\frac{3}{32} \left\langle  T(\mR^0, e^0_i),
T(e^0_i, \tfrac{\partial}{\partial \ov{z}^0_j})\right\rangle
- \frac{\sqrt{-1}}{32}\left\langle T(e^\bot_k,J\mR^0),
T(e^\bot_k,\tfrac{\partial}{\partial \ov{z}^0_j}) \right\rangle\\
\left.
 + \frac{\sqrt{-1}}{8}\left\langle   T(e^\bot_k,Je^\bot_k),
T(\mR^0,\tfrac{\partial}{\partial \ov{z}^0_j}) \right\rangle
 \Big]P^N \right\}((0,Z^\bot), (0, Z^{\bot}))\\
=\Big\{ -\frac{1}{12\pi}
  \left\langle R^{TX_G}(\tfrac{\partial}{\partial z^0_j},
\tfrac{\partial}{\partial \ov{z}^0_i})\tfrac{\partial}{\partial z^0_i}
+ R^{TX_G}(\tfrac{\partial}{\partial z^0_i},
\tfrac{\partial}{\partial \ov{z}^0_i})\tfrac{\partial}{\partial z^0_j},
\tfrac{\partial}{\partial \ov{z}^0_j}\right\rangle\\
+\frac{5\sqrt{-1}}{32\pi} \Big\langle
\nabla ^{TY}_{\tfrac{\partial}{\partial z^0_j}}
(T(e^\bot_k,\tfrac{\partial}{\partial \ov{z}^0_j}))
+\nabla ^{TY}_{e^\bot_k} (T(\tfrac{\partial}{\partial z^0_j},
\tfrac{\partial}{\partial \ov{z}^0_j})),
 Je^\bot_k \Big\rangle\\
+ \frac{3}{16\pi} \left\langle
R^{TB}(e^\bot_k, \tfrac{\partial}{\partial z^0_j})e^\bot_k,
\tfrac{\partial}{\partial \ov{z}^0_j}\right\rangle
+\frac{3}{32\pi}\Big|T(\tfrac{\partial}{\partial z^0_i},
\tfrac{\partial}{\partial \ov{z}^0_j})\Big|^2\\
- \frac{1}{64\pi}
\Big|T(e^\bot_k,\tfrac{\partial}{\partial \ov{z}^0_j})\Big|^2
- \frac{\sqrt{-1}}{16\pi} \left\langle T(e^\bot_k,Je^\bot_k),
T(\tfrac{\partial}{\partial z^0_j},
\tfrac{\partial}{\partial \ov{z}^0_j}) \right\rangle
\Big\} P_{\cL^\bot}(Z^\bot,Z^\bot).
\end{multline}

For $G_1(Z)$ (resp. $G_2(Z)$) polynomials on $Z$ with degree $1$
(resp. $2$) and $F\in T^*_{x_0}X_G\otimes T^*_{x_0}X_G$, by
Theorem \ref{t3.4}, \eqref{g7}, \eqref{0g6}, \eqref{g16},
\eqref{g3.57} and \eqref{g3.59}, for any $k,l$, $k',l'$,
\begin{align}\label{3.72}
\begin{split}
&\nabla_{0,e^\bot_j}P^N =-2\pi Z^\bot_j P^N,\\
&P^{N^\bot} P_{\cL^\bot} (G_1(Z) b^\bot_k +G_2(Z)b^\bot_k b^\bot_l
+Z^\bot_{k'}b_{l'}) P^N=0,\\
&\frac{1}{3}\left\langle  R^{TB} (\mR^\bot, e^\bot_i)\mR^\bot,
e^\bot_j\right\rangle  \nabla_{0,e^\bot_i} \nabla_{0,e^\bot_j} P^N
= -\frac{2\pi}{3}\left\langle  R^{TB} (\mR^\bot, e^\bot_j)\mR^\bot,
 e^\bot_j\right \rangle P^N,\\
&F(e^0_i,e^0_j) \nabla_{0,e^0_i} \nabla_{0,e^0_j}P^N
= \left[F(\tfrac{\partial}{\partial \ov{z}^0_i},
 \tfrac{\partial}{\partial \ov{z}^0_j})b_ib_j
-4\pi F(\tfrac{\partial}{\partial z^0_j},
\tfrac{\partial}{\partial \ov{z}^0_j})\right]P^N.
\end{split}\end{align}

By \eqref{a3.18}, \eqref{3.72}, we get
\begin{multline}\label{a3.71}
I_2 P^N= \Big \{ \Big(\Big\langle\frac{1}{3}
R^{TX_G} (\mR^0,\tfrac{\partial}{\partial \ov{z}^0_i})\mR^0
+ R^{TB} ( \mR^\bot, \tfrac{\partial}{\partial \ov{z}^0_i})\mR^\bot
+\nabla_{\mR^0}^{TX_G}(A(\tfrac{\partial}{\partial \ov{z}^0_i})\mR^\bot),
\tfrac{\partial}{\partial \ov{z}^0_j}\Big \rangle\\
-3  \left \langle A(\tfrac{\partial}{\partial \ov{z}^0_i}) \mR^\bot,
A(\tfrac{\partial}{\partial \ov{z}^0_j}) \mR^\bot\right \rangle
+\Big\langle \tfrac{\partial}{\partial \ov{z}^0_i},
\nabla_{\mR^0}^{TX_G}(A(\tfrac{\partial}{\partial \ov{z}^0_j})\mR^\bot)
\Big \rangle \Big)  b_i b_j \\
- 4\pi  \Big\langle \frac{1}{3}R^{TX_G} (\mR^0,
\tfrac{\partial}{\partial z^0_j})\mR^0
+ R^{TB} ( \mR^\bot, \tfrac{\partial}{\partial z^0_j})\mR^\bot
+\nabla_{\mR^0}^{TX_G}(A(\tfrac{\partial}{\partial z^0_j})\mR^\bot),
\tfrac{\partial}{\partial \ov{z}^0_j}\Big\rangle\\
+ 12\pi |A(\tfrac{\partial}{\partial \ov{z}^0_j})\mR^\bot|^2
-4\pi \Big\langle \tfrac{\partial}{\partial z^0_j},
\nabla_{\mR^0}^{TX_G}(A(\tfrac{\partial}{\partial \ov{z}^0_j})\mR^\bot)
\Big \rangle
 -\frac{2\pi}{3}\left\langle  R^{TB} (\mR^\bot, e^\bot_j)\mR^\bot,
 e^\bot_j\right \rangle\Big \} P^N.
\end{multline}

Observe that as $A(e^0_i)e^0_i\in N_{G}$, we have at $x_0$,
\begin{align}\label{a3.72}
\left \langle \nabla ^{TB}_{\mR^0}(A(e^0_i)e^0_i), e^0_j\right \rangle
=\left \langle A(\mR^0)A(e^0_i)e^0_i, e^0_j\right \rangle.
\end{align}

Thus by \eqref{0g6}, \eqref{g3.57},
\eqref{g3.59},  \eqref{b3.15}, \eqref{a3.70b},
 \eqref{3.72}-\eqref{a3.72}, $a_j=a^\bot_j=2\pi$,
and the arguments above \eqref{3.68},
\begin{subequations}
\begin{align}\label{3.73a}
&P^{N^\bot} P_{\cL^\bot}\left \langle \Gamma_{i i}(\mR),
 e_l \right \rangle\nabla_{0,e_l} P^N
= -\frac{2}{3}P^{N^\bot} \left\langle R^{TX_G} (\mR^0,e^0_i)e^0_i,
\tfrac{\partial}{\partial \ov{z}^0_j}\right \rangle b_j P^N,\\
&\label{3.73}
P^{N^\bot} P_{\cL^\bot}I_2 P^N
=  P^{N^\bot}\Big \{ \Big(
 \Big\langle \frac{1}{3}
R^{TX_G} (\mR^0,\tfrac{\partial}{\partial \ov{z}^0_i})\mR^0
+\frac{1}{4\pi}
R^{TB} (e^\bot_k,\tfrac{\partial}{\partial \ov{z}^0_i})e^\bot_k,
\tfrac{\partial}{\partial \ov{z}^0_j}\Big\rangle\\
&\hspace*{10mm}-\frac{3}{4\pi}\left \langle
A(\tfrac{\partial}{\partial \ov{z}^0_i}) e^\bot_k,
A(\tfrac{\partial}{\partial \ov{z}^0_j}) e^\bot_k\right\rangle \Big) b_i b_j
- \frac{4\pi}{3}\left\langle R^{TX_G}(\mR^0,
\tfrac{\partial}{\partial z^0_j})\mR^0,
\tfrac{\partial}{\partial \ov{z}^0_j}\right \rangle  \Big \} P^N.\nonumber
\end{align}
\end{subequations}

By \eqref{0g4}, \eqref{a3.15}, \eqref{3.64}, \eqref{3.73a},
\eqref{3.73} and the fact that $R^{TX_G}(\, , \,)$ is
$(1,1)$-form, we get
\begin{multline} \label{3.74}
-\left((\cL^0_2)^{-1}P^{N^\bot} P_{\cL^\bot}(I_2+ \left \langle
\Gamma_{i i}(\mR),
 e_l \right \rangle\nabla_{0,e_l}) P^N \right)((0,Z^{\bot}), (0,Z^{\bot}))\\
= \frac{1}{6\pi}\Big \{
3 \left\langle R^{TX_G} (\tfrac{\partial}{\partial z^0_i},
\tfrac{\partial}{\partial \ov{z}^0_i})\tfrac{\partial}{\partial z^0_j}
+R^{TX_G} (\tfrac{\partial}{\partial z^0_j},
\tfrac{\partial}{\partial \ov{z}^0_i})\tfrac{\partial}{\partial z^0_i},
\tfrac{\partial}{\partial \ov{z}^0_j}\right\rangle\\
-2 \left \langle R^{TX_G} (\tfrac{\partial}{\partial z^0_j},e^0_i)e^0_i
+R^{TX_G}(\tfrac{\partial}{\partial \ov{z}^0_i},
\tfrac{\partial}{\partial z^0_j})\tfrac{\partial}{\partial z^0_i},
\tfrac{\partial}{\partial \ov{z}^0_j}\right \rangle
 \Big \} P_{\cL^\bot}(Z^\bot,Z^\bot) \\
= \frac{2}{3\pi}\left\langle R^{TX_G} (\tfrac{\partial}{\partial z^0_j},
\tfrac{\partial}{\partial \ov{z}^0_i})\tfrac{\partial}{\partial z^0_i},
\tfrac{\partial}{\partial \ov{z}^0_j}\right\rangle
 P_{\cL^\bot}(Z^\bot,Z^\bot).
\end{multline}

Now by \eqref{3.32}, \eqref{3.53}, \eqref{a3.70b} and \eqref{3.72},
\begin{align} \label{3.75}
\begin{split}
&-P^{N^\bot} P_{\cL^\bot}\frac{1}{9} \sum_i
\Big[\sum_j(\partial_j R^{L_B})_{x_0} (\mR,e_i)Z_j\Big]^2 P^N
= \frac{\pi}{4} P^{N^\bot}  \Big | T (\mR^0, e^\bot_i)\Big |^2P^N,\\
&-\frac{1}{4}(\cL^0_2)^{-1}  P^{N^\bot} P_{\cL^\bot} [K_2(\mR),\cL^0_2]P^N
=\frac{1}{4} P^{N^\bot} P_{\cL^\bot}K_2(\mR) P^N\\
&\hspace*{25mm}
= \frac{1}{12}  P^{N^\bot} \left \langle R^{TX_G} (\mR^0,e^0_i)
\mR^0, e^0_i\right \rangle P^N. 
\end{split}\end{align}

By \eqref{33.36}, \eqref{a3.28}, \eqref{3.30} and \eqref{3.29}, we get
\begin{align} \label{a3.75}
 \frac{\sqrt{-1}}{2\pi}(\partial^0_j R^{L_B})_{x_0} (\mR, e^0_i)
= -\frac{1}{2}  \left \langle J\mR^\bot, T(e^0_j,e^0_i)\right \rangle
+ \left \langle JA(e^0_j)\mR^\bot,e^0_i \right \rangle =0.
\end{align}

Thus by \eqref{g7}, \eqref{a3.20}, \eqref{3.32},
\eqref{3.75} and \eqref{a3.75}, we get
\begin{multline} \label{3.76}
- P^{N^\bot} P_{\cL^\bot} \mO^\prime_2 P^N
=P^{N^\bot} P_{\cL^\bot}\Big\{
-I_1-(I_2+ \left \langle \Gamma_{i i}(\mR),
 e_l \right \rangle\nabla_{0,e_l}) \\
-\frac{1}{4} [K_2(\mR),\cL^0_2]
-  R^{E_B}(\mR,\tfrac{\partial}{\partial \ov{z}^0_j})b_j
-\frac{\pi}{4} \Big | T (\mR^0, e^\bot_i)\Big |^2\Big\} P^N.
\end{multline}

Note that $R^{TX_G} (\cdot,\cdot)$ is a $(1,1)$-form,
 by \eqref{g3.57}, \eqref{a3.15}, \eqref{3.64}, \eqref{3.67},
\eqref{3.69}, \eqref{3.71}, \eqref{3.74} and \eqref{3.76},
\begin{multline} \label{3.77}
-\left((\cL^0_2)^{-1} P^{N^\bot} P_{\cL^\bot} \mO^\prime_2 P^N\right)
((0,Z^{\bot}), (0,Z^{\bot}))
\\
 = -\left((\cL^0_2)^{-1} P^{N^\bot} P_{\cL^\bot}(I_1+I_2
+ \left \langle \Gamma_{i i}(\mR), e_l \right \rangle\nabla_{0,e_l})
P^N \right) ((0,Z^{\bot}), (0,Z^{\bot}))\\
+  \frac{1}{2\pi}\Big\{ R^{E_B}(\tfrac{\partial}{\partial z^0_j},
\tfrac{\partial}{\partial \ov{z}^0_j})
+\frac{1}{3}  \left \langle
R^{TX_G} (\tfrac{\partial}{\partial z^0_j},e^0_i)
e^0_i, \tfrac{\partial}{\partial \ov{z}^0_j}\right \rangle
+ \frac{1}{4}
\Big |T (\tfrac{\partial}{\partial \ov{z}^0_j}, e^\bot_i)\Big |^2\Big \} P_{\cL^\bot}(Z^\bot,Z^\bot)\\
= \Big\{ \frac{1}{2\pi}
 \left\langle  R^{TX_G}(\tfrac{\partial}{\partial z^0_j},
\tfrac{\partial}{\partial \ov{z}^0_i})\tfrac{\partial}{\partial z^0_i},
\tfrac{\partial}{\partial \ov{z}^0_j} \right\rangle
+ \frac{3}{16\pi} \left\langle
R^{TB}(e^\bot_k, \tfrac{\partial}{\partial z^0_j})e^\bot_k,
\tfrac{\partial}{\partial \ov{z}^0_j}\right\rangle\\
 + \frac{3}{32\pi}\Big| T(\tfrac{\partial}{\partial z^0_i},
\tfrac{\partial}{\partial \ov{z}^0_j})\Big|^2
- \frac{\sqrt{-1}}{16\pi} \left\langle T(e^\bot_k,Je^\bot_k),
T(\tfrac{\partial}{\partial z^0_j},
\tfrac{\partial}{\partial \ov{z}^0_j}) \right\rangle
+ \frac{7}{64\pi}
\Big|T(e^\bot_k,\tfrac{\partial}{\partial \ov{z}^0_j})\Big|^2\\
+\frac{5\sqrt{-1}}{32\pi} \left\langle
\nabla ^{TY}_{\tfrac{\partial}{\partial z^0_j}}
(T(e^\bot_k,\tfrac{\partial}{\partial \ov{z}^0_j}))
+\nabla ^{TY}_{e^\bot_k} (T(\tfrac{\partial}{\partial z^0_j},
\tfrac{\partial}{\partial \ov{z}^0_j})), Je^\bot_k \right\rangle
+\frac{1}{2\pi}R^{E_B}(\tfrac{\partial}{\partial z^0_j},
\tfrac{\partial}{\partial \ov{z}^0_j})
\Big \} P_{\cL^\bot}(Z^\bot,Z^\bot).
\end{multline}

By  \eqref{g3.57}, \eqref{a3.48}, \eqref{3.53}, \eqref{a3.70b},
\eqref{3.72} and the arguments above \eqref{3.68},
\begin{multline} \label{3.78}
4\pi ^2 P^{N^\bot} P_{\cL^\bot} \mO''_2 P^N
=4\pi ^2 P^{N^\bot} P_{\cL^\bot} \Big\{
-\frac{1}{3}\left\langle (\nabla^{TY}_\cdot
\dot{g}^{TY}_\cdot)_{(\mR^0,\mR^0)}
J\mR^\bot, J\mR^\bot\right\rangle \\
+\frac{1}{6}\left\langle \nabla ^{TY}_{\mR^0} (T(e^\bot_j, J_{x_0}e^0_i))
Z^\bot_j Z^0_i
 +\nabla ^{TY}_{\mR^\bot} (T(e^0_j, J_{x_0}e^0_i))Z^0_j Z^0_i,
 J\mR^\bot\right\rangle \\
+\frac{1}{3} \left\langle R^{TB}(\mR^\bot,\mR^0)\mR^0,\mR^\bot\right\rangle
- \frac{1}{12} \sum_l \left\langle T(\mR^0, e_l),J\mR^\bot\right\rangle^2
 \Big\} P^N \\
=\frac{\pi}{3}  P^{N^\bot} \Big\{ \frac{1}{2}\left\langle \nabla ^{TY}_{\mR^0}
(T(e^\bot_k, J_{x_0}e^0_i)) Z^0_i
+\nabla ^{TY}_{e^\bot_k} (T(e^0_j, J_{x_0}e^0_i))Z^0_j Z^0_i,
 Je^\bot_k\right\rangle  \\
-\left\langle (\nabla^{TY}_\cdot \dot{g}^{TY}_\cdot)
_{(\mR^0,\mR^0)} Je^\bot_k, Je^\bot_k \right\rangle
+ \left\langle R^{TB}(e^\bot_k,\mR^0)\mR^0,e^\bot_k\right\rangle
- \frac{1}{4} |T(\mR^0, e_l)|^2
\Big\} P^N .
\end{multline}

Let $\{f_l\}$ be an orthonormal frame of $TY$ on $X$.

As $\nabla^{TY}$ preserves the metric $g^{TY}$,
by  \eqref{h1},  \eqref{h14},
\begin{align} \label{3.79}
 \left\langle (\nabla^{TY}_{e^0_i}\dot{g}^{TY}_{e^0_j}) f_l,
f_l \right\rangle = \nabla_{e^0_i}\left\langle
\dot{g}^{TY}_{e^0_j} f_l, f_l \right\rangle = 4
\nabla_{e^0_i}\nabla_{e^0_j} \log h.
\end{align}

\comment{
By \eqref{3.57}, \eqref{3.64} and  \eqref{a3.68},
\begin{multline}\label{3.80}
-\pi ^2 \cL^{-1} P^\bot_{\cL}\Big|T(z^0, \ov{z}^0)\Big|^2 P_{\cL}(0,0)
= \frac{1}{4} \cL^{-1} P^\bot_{\cL}
\left\langle T(z^0, \tfrac{\partial}{\partial \ov{z}^0_i}),
T(z^0,\tfrac{\partial}{\partial \ov{z}^0_j})\right\rangle b_i b_j
P_{\cL}(0,0)\\
= \frac{3}{8\pi}  \Big|\sum_{i} T(\tfrac{\partial}{\partial z^0_i},
\tfrac{\partial}{\partial \ov{z}^0_i})\Big|^2
+\frac{3}{8\pi} \sum_{i j}\Big|T(\tfrac{\partial}{\partial z^0_i},
\tfrac{\partial}{\partial \ov{z}^0_j})\Big|^2.
\end{multline}
}

Now $\{J e^\bot_k\}$ is an orthonormal basis of $TY$ along the fiber $Y_{x_0}$
and $\{e_l\}= \{e^0_i\}\cup \{e^\bot_k\}$.

By \eqref{g3.57}, \eqref{3.64}, \eqref{a3.70}, \eqref{3.78} and \eqref{3.79},
\begin{multline}\label{3.81}
-4\pi ^2 \left((\cL^0_2)^{-1}P^{N^\bot} P_{\cL^\bot} \mO''_2
P^N\right)
 ((0,Z^{\bot}), (0,Z^{\bot}))\\
=\frac{1}{4\pi}\Big\{\frac{\sqrt{-1}}{6}\left\langle
-\nabla ^{TY}_{\tfrac{\partial}{\partial z^0_j}}
(T(e^\bot_k, \tfrac{\partial}{\partial \ov{z}^0_j}))
+ \nabla ^{TY}_{\tfrac{\partial}{\partial \ov{z}^0_j}}
(T(e^\bot_k, \tfrac{\partial}{\partial z^0_j}))
-2\nabla ^{TY}_{e^\bot_k} (T(\tfrac{\partial}{\partial z^0_j},
\tfrac{\partial}{\partial \ov{z}^0_j})),Je^\bot_k\right\rangle \\
-\frac{8}{3}\nabla_{\tfrac{\partial}{\partial z^0_j}}
\nabla_{\tfrac{\partial}{\partial \ov{z}^0_j}} \log h
-\frac{1}{3} \Big|T(\tfrac{\partial}{\partial z^0_i},
\tfrac{\partial}{\partial \ov{z}^0_j})\Big|^2
-\frac{1}{6} \Big|T(e ^\bot_k,
\tfrac{\partial}{\partial \ov{z}^0_j})\Big|^2 \\
- \frac{2}{3}\left\langle R^{TB}(e^\bot_k,
\tfrac{\partial}{\partial z^0_j})e^\bot_k,
\tfrac{\partial}{\partial \ov{z}^0_j}\right\rangle
\Big \} P_{\cL^\bot}(Z^\bot,Z^\bot).
\end{multline}

By \eqref{a3.51}, \eqref{b3.52}, \eqref{3.77} and \eqref{3.81},
we get \eqref{a3.64}.
The proof of Lemma \ref{t4.9} is complete.
\end{proof}

\subsection{Final computations: the proof of Theorem \ref{t0.6}}\label{s4.6}

By \eqref{g3.29}, \eqref{3.10}, \eqref{3.11a}, \eqref{3.12a} and
\eqref{b3.16},
 as $Je^\bot_k\in TY$ on $P$, we get at $x_0$,
\begin{align}\label{4.1}
\begin{split}
&\nabla ^{TY}_{e^{0}_i}Je^\bot_k = P^{TY} \nabla ^{TX}_{e^{0}_i}J e^\bot_k
= P^{TY}J \nabla ^{TX}_{e^{0}_i}e^\bot_k =0,\\
&\nabla ^{TB}_{e^{0}_i} J  e^0_j
= \nabla ^{TX_G}_{e^{0}_i}J e^0_j +A(e^{0}_i)J e^0_j
= -\frac{1}{2} JT(e^{0}_i,e^0_j) = \nabla ^{TB}_{e^{0}_i} (J_{x_0}  e^0_j).
\end{split}\end{align}
By \eqref{h3}, \eqref{h14}, \eqref{3.11c} and \eqref{4.1},
as in \eqref{3.79}, at $x_0$,
\begin{multline}\label{4.6}
\left\langle \nabla ^{TY}_{e^0_i} (T(e^\bot_k,e^{0}_j) ),
J e^\bot_k \right\rangle_{x_0}
=-2 \left\langle \nabla ^{TY}_{e^0_i}
 (T(J e^{0}_j, J e^\bot_k)) , J e^\bot_k \right\rangle\\
= - \left\langle (\nabla ^{TY}_{e^0_i}\dot{g}^{TY}_{J e^0_j}) J e^\bot_k ,
J e^\bot_k \right\rangle
=- 4  \nabla_{e^0_i}\nabla_{J_{x_0} e^0_j} \log h.
\end{multline}
By  \eqref{r6} and \eqref{4.6}, we get
\begin{align}\label{4.7}
\begin{split}
&\sqrt{-1}\Big\langle \nabla ^{TY}_{\tfrac{\partial}{\partial z^0_j}}
(T(e^\bot_k,\tfrac{\partial}{\partial \ov{z}^0_j}) ),
J e^\bot_k \Big\rangle
= - 4  \nabla_{\tfrac{\partial}{\partial z^0_j}}
\nabla_{\tfrac{\partial}{\partial \ov{z}^0_j}} \log h
=\Delta_{X_G}\log h ,\\
&\sqrt{-1}\Big\langle \nabla ^{TY}_{\tfrac{\partial}{\partial \ov{z}^0_j}}
(T(e^\bot_k,\tfrac{\partial}{\partial z^0_j})),J e^\bot_k \Big\rangle
=-\Delta_{X_G}\log h
.
\end{split}\end{align}
Note that $T(e_i,e_j)= -[e^{H}_i,e^{H}_j]$,
as $[e_i,e_j]=0$. By \eqref{h1} and \eqref{h3},
\begin{multline}\label{4.8}
 \nabla ^{TY}_{e^{\bot,H}_k} (T(e^{0,H}_i,e^{0,H}_j))
= - \Big[e^{\bot,H}_k, [e^{0,H}_i,e^{0,H}_j]\Big]
 + T(e^{\bot,H}_k, T(e^{0,H}_i,e^{0,H}_j))\\
= L_{e^{0,H}_i} (T(e^{\bot,H}_k, e^{0,H}_j))
- L_{e^{0,H}_j} (T(e^{\bot,H}_k, e^{0,H}_i))
+ T(e^{\bot,H}_k, T(e^{0,H}_i,e^{0,H}_j))\\
= \nabla ^{TY}_{e^{0,H}_i} (T(e^{\bot,H}_k, e^{0,H}_j))
-\nabla ^{TY}_{e^{0,H}_j} (T(e^{\bot,H}_k, e^{0,H}_i))
- T(e^{0,H}_i, T(e^{\bot,H}_k, e^{0,H}_j))  \\
+T(e^{0,H}_j, T(e^{\bot,H}_k, e^{0,H}_i))
+ T(e^{\bot,H}_k, T(e^{0,H}_i,e^{0,H}_j)).
\end{multline}
Thus by Theorem \ref{t4.2}, \eqref{4.7} and \eqref{4.8},
\begin{multline}\label{4.9}
\sqrt{-1}\left\langle \nabla ^{TY}_{e^\bot_k}
(T(\tfrac{\partial}{\partial z^0_j},\tfrac{\partial}{\partial \ov{z}^0_j}) ),
J e^\bot_k \right\rangle
=\sqrt{-1}\Big \{ 2 \Big\langle \nabla ^{TY}_{\tfrac{\partial}{\partial z^0_j}}
(T(e^\bot_k,\tfrac{\partial}{\partial \ov{z}^0_j}) ),
J e^\bot_k \Big\rangle \\
- \left\langle T(\tfrac{\partial}{\partial z^0_j}, J e^\bot_k ),
T(e^\bot_k, \tfrac{\partial}{\partial \ov{z}^0_j} )\right\rangle
+  \left\langle T(\tfrac{\partial}{\partial \ov{z}^0_j}, J e^\bot_k ),
T(e^\bot_k, \tfrac{\partial}{\partial z^0_j} )\right\rangle
+  \left\langle T(e^\bot_k,J e^\bot_k ),
T(\tfrac{\partial}{\partial z^0_j},
\tfrac{\partial}{\partial \ov{z}^0_j})\right\rangle \Big \}\\
= 2  \Delta_{X_G} \log h
+  | T(e^\bot_k,\tfrac{\partial}{\partial \ov{z}^0_j})|^2
+ \sqrt{-1}
 \left\langle T(e^\bot_k,J e^\bot_k ), T(\tfrac{\partial}{\partial z^0_j},
\tfrac{\partial}{\partial \ov{z}^0_j})\right\rangle.
\end{multline}

By $T(e_i,e_j)=-[e^H_i,e^H_j]$, \eqref{g3.29} and \eqref{3.33}, we have
\begin{align}\label{4.10}
&R^{TX}(e^H_k,e^H_j) e^H_i
= \nabla ^{TX}_{e^H_k}\nabla ^{TX}_{e^H_j} e^H_i
-\nabla ^{TX}_{e^H_j}\nabla ^{TX}_{e^H_k} e^H_i
- \nabla ^{TX}_{[e^H_k,e^H_j]}e^H_i \\
&\hspace*{5mm}
= R^{TB}(e_k,e_j) e_i -\frac{1}{2} T(e_k, \nabla ^{TB}_{e_j} e_i)
+\frac{1}{2} T(e_j, \nabla ^{TB}_{e_k} e_i)\nonumber\\
&\hspace*{25mm}-\frac{1}{2} \nabla ^{TX}_{e^H_k}(T(e_j, e_i))
+ \frac{1}{2} \nabla ^{TX}_{e^H_j}(T(e_k, e_i))
+  \nabla ^{TX}_{T(e^H_k,e^H_j)}e^H_i,\nonumber\\
&\left\langle  R^{TX}(e^{\bot,H}_k,e^{0,H}_j) (J_{x_0} e^{0}_j)^H,
J_{x_0} e^{\bot,H}_k \right\rangle
=\left\langle  R^{TX}(e^{\bot,H}_k,e^{0,H}_j) e^{0,H}_j,
e^{\bot,H}_k \right\rangle. \nonumber
\end{align}
By \eqref{3.11a}, \eqref{3.12a}, \eqref{33.36}, \eqref{4.1} and $T(e^\bot_k, e^0_i)\in TY$,
 at $x_0$, we get
\begin{align}  \label{a4.11}
\begin{split}
&\nabla^{TB} _{e^0_j}(J_{x_0}e^0_j)=0,\quad
  \nabla^{TB} _{e^{\bot}_k}(J_{x_0}e^0_j)
= \frac{1}{2} \left\langle T( e^0_j, e^0_l), J e^{\bot}_k\right\rangle e^0_l,\\
&\left\langle \nabla ^{TX}_{T(e^\bot_k, e^0_i)}
(J_{x_0} e^{0}_i)^H,J_{x_0} e^\bot_k\right\rangle
=\left\langle \nabla ^{TX}_{T(e^\bot_k, e^0_i)}
e^{0,H}_i, e^\bot_k\right\rangle .
\end{split}\end{align}
We apply now the first equation of \eqref{4.10} into the second equation
of \eqref{4.10}, by using \eqref{ah4} and \eqref{4.17}
and $T(\, , \,)$ is $(1,1)$-form, we get at $x_0$,
\begin{multline}\label{4.11}
\frac{1}{4}  |T( e^0_j, e^0_l)|^2
-  \frac{1}{2}\left\langle \nabla ^{TY}_{e^\bot_k}
(T(e^0_j, J_{x_0}e^{0}_j)),J e^\bot_k \right\rangle
+\frac{1}{2} \left\langle \nabla ^{TY}_{e^0_j} (T(e^\bot_k, J_{x_0}e^{0}_j)),
J e^\bot_k \right\rangle\\
=\left\langle  R^{TB}(e^\bot_k,e^{0}_j) e^{0}_j, e^\bot_k\right\rangle
+ \frac{1}{2} \left\langle \nabla ^{TX}_{e^{0}_j} (T(e^\bot_k,e^{0}_j)),
 e^\bot_k \right\rangle \\
=\left\langle  R^{TB}(e^\bot_k,e^{0}_j) e^{0}_j, e^\bot_k\right\rangle
- \frac{1}{4}  |T(e^\bot_k, e^0_l)|^2.
\end{multline}
Finally from \eqref{0g4}, \eqref{4.7}, \eqref{4.9} and \eqref{4.11}
and $T(\, ,\,)$ is $(1,1)$-form, we get
\begin{multline}\label{4.12}
2\left\langle  R^{TB}(e^\bot_k,\tfrac{\partial}{\partial z^0_j})
 \tfrac{\partial}{\partial \ov{z}^0_j}, e^\bot_k\right\rangle
=\sqrt{-1} \left\langle \nabla ^{TY}_{e^\bot_k}
(T(\tfrac{\partial}{\partial z^0_j},\tfrac{\partial}{\partial \ov{z}^0_j}) ),
J e^\bot_k \right\rangle\\
-\sqrt{-1} \Big\langle \nabla ^{TY}_{\tfrac{\partial}{\partial z^0_j}}
(T(e^\bot_k,\tfrac{\partial}{\partial \ov{z}^0_j}) ),
J e^\bot_k \Big\rangle
+\frac{1}{2} | T(e^\bot_k,\tfrac{\partial}{\partial \ov{z}^0_j})|^2
+  \Big| T(\tfrac{\partial}{\partial z^0_i},
\tfrac{\partial}{\partial \ov{z}^0_j})\Big|^2\\
=\Delta_{X_G} \log h
+ \frac{3}{2}| T(e^\bot_k,\tfrac{\partial}{\partial \ov{z}^0_j})|^2
+  \Big| T(\tfrac{\partial}{\partial z^0_i},
\tfrac{\partial}{\partial \ov{z}^0_j})\Big|^2
+\sqrt{-1} \left\langle T(e^\bot_k,J e^\bot_k ),
T(\tfrac{\partial}{\partial z^0_j},
\tfrac{\partial}{\partial \ov{z}^0_j})\right\rangle .
\end{multline}
From \eqref{4.7}-\eqref{4.12},
\begin{multline}\label{4.16}
\frac{\sqrt{-1}}{96 \pi}\left\langle
11\,   \nabla ^{TY}_{\tfrac{\partial}{\partial z^0_j}}
(T(e^\bot_k,\tfrac{\partial}{\partial \ov{z}^0_j}))
+4\,  \nabla ^{TY}_{\tfrac{\partial}{\partial \ov{z}^0_j}}
(T(e^\bot_k, \tfrac{\partial}{\partial z^0_j}))
+7\, \nabla ^{TY}_{e^\bot_k} (T(\tfrac{\partial}{\partial z^0_j},
\tfrac{\partial}{\partial \ov{z}^0_j})), Je^\bot_k \right\rangle\\
+ \frac{1}{48 \pi}\left\langle
R^{TB}(e^\bot_k,\tfrac{\partial}{\partial z^0_j})e^\bot_k,
 \tfrac{\partial}{\partial \ov{z}^0_j}\right\rangle
-\frac{\sqrt{-1}}{16 \pi}  \left\langle T(e^\bot_k,J e^\bot_k ),
 T(\tfrac{\partial}{\partial z^0_j},
\tfrac{\partial}{\partial \ov{z}^0_j})\right\rangle\\
= \frac{5}{24\pi} \Delta_{X_G}\log h
+ \frac{11}{192\pi}
\Big|T(e^\bot_k,\tfrac{\partial}{\partial \ov{z}^0_j})\Big|^2
-\frac{1}{96 \pi}\Big| T(\tfrac{\partial}{\partial z^0_i},
\tfrac{\partial}{\partial \ov{z}^0_j})\Big|^2
.
\end{multline}
By \eqref{g16}, \eqref{b3.52}, \eqref{3.65}, \eqref{a3.64} and \eqref{4.16},
\begin{align}\label{4.15}
\Phi_{1,1}+\Phi_{1,2}
&=\frac{1}{2\pi}
 \left\langle  R^{TX_G}(\tfrac{\partial}{\partial z^0_j},
\tfrac{\partial}{\partial \ov{z}^0_i})\tfrac{\partial}{\partial z^0_i},
\tfrac{\partial}{\partial \ov{z}^0_j} \right\rangle
+ \frac{3}{8\pi} \Delta_{X_G}\log h
+\frac{1}{2\pi}R^{E_B}(\tfrac{\partial}{\partial z^0_j},
\tfrac{\partial}{\partial \ov{z}^0_j})\\
&=\frac{1}{16\pi} r^{X_G}_{x_0} + \frac{3}{8\pi} \Delta_{X_G}\log h
+\frac{1}{4\pi}R^{E_G}_{x_0}(w^0_j, \ov{w}^0_j). \nonumber
\end{align}

From Lemma \ref{at4.9}, \eqref{a3.53} and \eqref{4.15}, we get \eqref{0.15}.

Recall that we compute everything on $\cC^\infty(X, L^p\otimes
E)$.

From \eqref{3.19}, \eqref{a3.16}, \eqref{a3.17}, comparing to
\eqref{1ue1}, we know that in \eqref{0.14},  $\Phi_r(x_0)\in  \End
(E_G)_{x_0}$, and the term $r^X$, $R^{\det}$ will not appear here,
and $\tau = 2\pi n$,
 thus we get the remainder part of Theorem \ref{t0.6}
from Corollary \ref{t0.5}.

The proof of Theorem \ref{t0.6}  is complete.

\newpage

\subsection{Coefficient  $P^{(2)}(0,0)$}\label{s4.7}

As in \eqref{a3.53}, we have
 \begin{align}\label{4.20}
P^{(2)}(0,0)= (\Psi_{1,1}+\Psi_{1,2})(0) + (\Psi_{1,1}+\Psi_{1,2})^*(0)
+ (\Psi_{1,3}-\Psi_{1,4})(0).
\end{align}

For $k\in \bN$, let $H_k(x)$ be the Hermite polynomial,
\begin{align}\label{4.17}
H_k(x) =\sum_{j=0}^{[k/2]} (-1)^j \frac{k!\,  (2x)^{k-2j}}{j!\, (k-2j)!}.
\end{align}
By \cite[\S 8.6]{Taylor96}, \eqref{g6} and $a^\bot_l=2\pi$, we have
\begin{align}\label{4.18}
(b^\bot_l)^k e ^{-\pi |Z^\bot_l|^2}
= (2\pi )^{k/2} H_k(\sqrt{2\pi} Z^\bot_l)e ^{-\pi |Z^\bot_l|^2}.
\end{align}
Especially, for $l$ fixed, $i\in \bN$,
\begin{align}\label{4.19}
&((b^\bot_l)^{2i+1} e ^{-\pi |Z^\bot_l|^2})(0)=0, \nonumber\\\
&((b^\bot_l)^2 e ^{-\pi |Z^\bot_l|^2})(0)= -4\pi,\quad
((b^\bot_l)^4 e ^{-\pi |Z^\bot_l|^2})(0)=3\cdot (4\pi)^2,\\
& ((b^\bot_l)^6 e ^{-\pi |Z^\bot_l|^2})(0)=15\cdot (-4\pi)^3.\nonumber
\end{align}

\comment{
For $U\in T_{x_0}X_G$, set
 \begin{align}\label{4.21}
\mT_{ijk}=  \left\langle JT(e^\bot_i, J e^\bot_j), e^\bot_k \right\rangle,
\quad \mT_{jk}(U)= \left\langle JT(U,e^\bot_j), e^\bot_k \right\rangle
\end{align}
Then by Theorem \ref{t3.4}, $\mT_{ijk}$ is symmetric on $i,j,k$ and
 $\mT_{jk} \in T_{x_0}^*X_G$ is  symmetric on $j,k$.
}

Recall that when we meet the operation $|\quad|^2$,
 we will first do this operation, then take the sum of the indices.
Thus $|\mT_{ijk}|^2$ means
$\sum_{ijk}|\mT_{ijk}|^2$, etc.

By \eqref{g18},  \eqref{a3.55} and \eqref{4.19},
 \begin{align}\label{4.22}
\mF_2 (\cdot, 0)
= - \frac{1}{8}  \mT_{kk};
\quad P^N(0,0)=2^{n_0/2}.
\end{align}

By \eqref{a3.59}, \eqref{4.19} and \eqref{4.22}, we know
 \begin{align}\label{4.23}
\Psi_{1,3}(0)= \frac{2^{n_0/2}}{\pi}
\Big|\frac{1}{4}  \sum_k \mT_{kk}(\tfrac{\partial}{\partial \ov{z}^0_i})
+ \mF_2 (\tfrac{\partial}{\partial \ov{z}^0_i}, 0)\Big|^2
=  \frac{2^{n_0/2}}{64 \pi}
\Big| \sum_k \mT_{kk}(\tfrac{\partial}{\partial \ov{z}^0_i})\Big|^2 .
\end{align}
and from \eqref{g14}, \eqref{g15}, \eqref{g3.57},
 \eqref{a3.60}, \eqref{4.22} and $a^\bot_i=2\pi$,
 \begin{multline}\label{4.24}
\Psi_{1,4}(0)=G^\bot(0)^2
\Big\{ \frac{1}{4 \pi} \sum_k \mF_1(e^\bot_k)^2
+ \frac{6\cdot (4\pi)^3}{(192 \pi^2)^2} | \mT_{klm}|^2 \\
+ \frac{1}{16\pi}  \Big| \sum_k \mT_{kk}(\tfrac{\partial}{\partial \ov{z}^0_i}) \Big|^2
+ \frac{2\cdot (4\pi)^2}{\pi\cdot (32\pi)^2}
\Big|\mT_{kl}(\tfrac{\partial}{\partial \ov{z}^0_i}) \Big|^2\Big\}\\
=\frac{2^{n_0/2}}{4\pi} \Big\{ \sum_k\mF_1(e^\bot_k)^2
+  \frac{1}{24}| \mT_{klm}|^2
+ \frac{1}{4} \Big|\sum_k
\mT_{kk}(\tfrac{\partial}{\partial \ov{z}^0_i}) \Big|^2
+ \frac{1}{8}
\Big|\mT_{kl}(\tfrac{\partial}{\partial \ov{z}^0_i}) \Big|^2\Big\}.
\end{multline}
\begin{lemma}\label{t4.11} The following identity holds,
 \begin{multline}\label{4.25}
\Psi_{1,1}(0)= \Big\{
-\frac{19}{2^6\cdot 3\pi}|\mT_{jj'}(\tfrac{\partial}{\partial \ov{z}^0_i})|^2
 -\frac{11}{2^7\cdot 3\pi}\mT_{klm}^2 +\frac{1}{2^8 \pi} \mT_{kkm}\mT_{llm}\\
-\frac{5}{2 ^7 \pi}\mT_{jj}(\tfrac{\partial}{\partial z^0_l})
\mT_{kk}(\tfrac{\partial}{\partial \ov{z}^0_l})
- \frac{1}{8\pi} \sum_k  \mF_1(e^\bot_k)^2
- \frac{1}{16\pi}  \mF_1(e^\bot_k)\mT_{kll} \Big\}P^N(0,0).
\end{multline}
\end{lemma}
\begin{proof} Set
\begin{align}\label{4.26}
\begin{split}
&\mI_1=-\sqrt{-1}\Big(\mT_{jj'}(\tfrac{\partial}{\partial z^0_i})\frac{b^{+}_i}{8\pi}
+\frac{1}{4} \mT_{jj'}(z^0) \Big) b^{\bot}_j b^{\bot}_{j'}
\frac{\sqrt{-1}}{8\pi} b_l \mT_{kk}(\tfrac{\partial}{\partial \ov{z}^0_l})\\
&\mI_2=\sqrt{-1}\Big(
 \mT_{jj'}(\tfrac{\partial}{\partial \ov{z}^0_i}) b_i
\frac{B^\bot_{jj'}}{8\pi}
+ \frac{1}{4}  \mT_{jj'}(\ov{z}^0)
(b^{\bot +}_j b^{\bot +}_{j'} -b^{\bot}_j b^{\bot}_{j'})\Big)
\frac{-\sqrt{-1}}{32\pi} \mT_{kl}(z^0)b^{\bot}_k b^{\bot}_{l},\\
&\mI_3=
-\frac{\sqrt{-1}}{8\pi} \wi{\mT}_{ijj'}(b^{\bot}_j b^{\bot+}_{j'}
+ b^{\bot}_j b^{\bot}_{j'})b^{\bot+}_i
\Big(\mF_1(e^\bot_k)\frac{b^{\bot }_k}{4\pi} + \mT_{klm}
\frac{b^{\bot}_k b^{\bot}_l b^{\bot}_m}{192\pi ^2}\Big).
\end{split}\end{align}

Observe that by \eqref{3.64}, when we evaluate $\Psi_{1,1}$ in
\eqref{b3.52}, in each monomial, if the total degree of  $b_l$,
$\ov{z}^0$ is not as same as the total degree of $b^+_l$, $z^0$,
then the contribution of this term is $0$. Thus by \eqref{g7},
\eqref{g3.57},  \eqref{b3.52}, \eqref{3.53}, \eqref{3.56},
\eqref{3.58}, \eqref{a3.55} and \eqref{4.26},
 \begin{multline}\label{4.27}
\Psi_{1,1}(Z^\bot)=\Big\{(\cL^0_2)^{-1} P^{N^\bot} \Big[\mI_1 +\mI_2+\mI_3\\
+ \Big(\mF_1(e^\bot_j) (b^{\bot +}_j+ b^{\bot }_j)
+  \mT_{ijj'}\frac{B^\bot_{ijj'}}{16\pi}\Big)
\Big(\mF_1(e^\bot_k)\frac{b^{\bot }_k}{4\pi} + \mT_{klm}
\frac{b^{\bot}_k b^{\bot}_l b^{\bot}_m}{192\pi ^2}\Big)
\Big]P^N\Big\}(Z^\bot,Z^\bot).
\end{multline}

By  \eqref{g6}, \eqref{g16} and \eqref{4.19},
\begin{align}\label{a4.34}
&(b_jz^0_i P^N)(0,0)= -2 \delta_{ij}P^N(0,0),\quad
(b^\bot_kb^\bot_l b_jz^0_i P^N)(0,0)
= 8\pi  \delta_{ij}\delta_{kl}P^N(0,0).
\end{align}

From Theorem \ref{t3.4}, \eqref{g7},  \eqref{g3.57},
\eqref{4.19}, \eqref{4.26} and  \eqref{a4.34},
\begin{multline}\label{4.35}
((\cL^0_2)^{-1} P^{N^\bot} \mI_1P^N)(0,0)\\
= \frac{1}{32\pi} \left\{(\cL^0_2)^{-1} P^{N^\bot}
\mT_{kk}(\tfrac{\partial}{\partial \ov{z}^0_l})
\Big ( 4 \mT_{jj'}(\tfrac{\partial}{\partial z^0_l})b^{\bot}_j b^{\bot}_{j'}
+ b_l b^{\bot}_j b^{\bot}_{j'} \mT_{jj'}(z^0)\Big)P^N\right\}(0,0)\\
= \frac{1}{32\pi} \mT_{kk}(\tfrac{\partial}{\partial \ov{z}^0_l})
\Big\{\Big(\mT_{jj'}(\tfrac{\partial}{\partial z^0_l})
\frac{b^{\bot}_j b^{\bot}_{j'}}{2\pi}
+ \frac{b_l b^{\bot}_j b^{\bot}_{j'}}{12\pi}\mT_{jj'}(z^0)\Big)P^N\Big\}(0,0)\\
=-\frac{1}{24\pi} \mT_{jj}(\tfrac{\partial}{\partial z^0_l})
\mT_{kk}(\tfrac{\partial}{\partial \ov{z}^0_l})P^N(0,0).
\end{multline}

By  \eqref{g7}, \eqref{g3.57}, \eqref{3.53} and \eqref{4.26},
\begin{multline}\label{4.36}
 (P^{N^\bot} \mI_2 P^N)(Z,(0,Z^{\prime\bot})) =\frac{1}{2^8 \pi^2}
\Big\{P^{N^\bot}
\mT_{jj'}(\tfrac{\partial}{\partial \ov{z}^0_i})\\
\Big[b_i \mT_{kl}(z^0) B^\bot_{jj'}
+ (b^{\bot +}_j b^{\bot +}_{j'} -b^{\bot}_j b^{\bot}_{j'})
\mT_{kl}(z^0) b_i\Big] b^{\bot}_k b^{\bot}_{l}P^N\Big\}(Z,(0,Z^{\prime\bot}))\\
=\frac{1}{2^8\pi^2}  \mT_{jj'}(\tfrac{\partial}{\partial \ov{z}^0_i})
\Big\{P^{N^\bot} \Big[ b_i \mT_{kl}(z^0)(2 b^{\bot +}_j b^{\bot +}_{j'}
+ 2 b^{\bot}_j b^{\bot +}_{j'}+4\pi\delta_{jj'}) \\
+ 2 \mT_{kl}(\tfrac{\partial}{\partial z^0_i})(b^{\bot +}_j b^{\bot +}_{j'}
-b^{\bot}_j b^{\bot}_{j'})  \Big]b^{\bot}_k b^{\bot}_{l}P^N
\Big\}(Z,(0,Z^{\prime\bot}))\\
=\frac{1}{2^8\pi^2}  \mT_{jj'}(\tfrac{\partial}{\partial \ov{z}^0_i})
 \Big\{ b_i \Big(64\pi^2 \mT_{jj'}(z^0)
+16 \pi\mT_{kj'}(z^0)b^{\bot}_j b^{\bot}_k
+ 4\pi \delta_{jj'}\mT_{kl}(z^0)b^{\bot}_k b^{\bot}_{l}\Big)\\
-2 \mT_{kl}(\tfrac{\partial}{\partial z^0_i})
 b^{\bot}_j b^{\bot}_{j'}b^{\bot}_k b^{\bot}_{l}P^N
\Big\}(Z,(0,Z^{\prime\bot})).
\end{multline}
Thus by Theorem \ref{t3.4},  \eqref{g6},
 \eqref{4.19}, \eqref{a4.34}, \eqref{4.36}
and use a similar equation as \eqref{4.32} for
$\mT_{jj'}(\tfrac{\partial}{\partial \ov{z}^0_i})
\mT_{kl}(\tfrac{\partial}{\partial z^0_i})
 b^{\bot}_j b^{\bot}_{j'}b^{\bot}_k b^{\bot}_{l}$, we get
\begin{multline}\label{4.37}
 ((\cL^0_2)^{-1} P^{N^\bot} \mI_2 P^N)(0,0)
=\frac{1}{2^8\pi^2} \mT_{jj'}(\tfrac{\partial}{\partial \ov{z}^0_i})
\left[\Big(16\pi b_i\mT_{jj'}(z^0) \right. \\
\left.+\frac{4}{3}b_i\mT_{kj'}(z^0)b^{\bot}_j b^{\bot}_k
+\frac{1}{3}\delta_{jj'} b_i\mT_{kl}(z^0)b^{\bot}_k b^{\bot}_l
- \frac{1}{8\pi}\mT_{kl}(\tfrac{\partial}{\partial z^0_i})
b^{\bot}_j b^{\bot}_{j'}b^{\bot}_k b^{\bot}_{l}\Big)P^N\right](0,0)\\
= \frac{1}{2^8\pi^2}\left[
-\frac{64\pi}{3}|\mT_{jj'}(\tfrac{\partial}{\partial \ov{z}^0_i})|^2
+ \frac{8\pi}{3} \mT_{jj}(\tfrac{\partial}{\partial \ov{z}^0_i})
\mT_{kk}(\tfrac{\partial}{\partial z^0_i})\right. \\
\left.-2\pi\Big(2 |\mT_{jj'}(\tfrac{\partial}{\partial \ov{z}^0_i})|^2
+ \mT_{jj}(\tfrac{\partial}{\partial \ov{z}^0_i})
\mT_{kk}(\tfrac{\partial}{\partial z^0_i}) \Big)\right] P^N(0,0)\\
= \frac{1}{2^8\cdot 3\pi}
\left[-76|\mT_{jj'}(\tfrac{\partial}{\partial \ov{z}^0_i})|^2
+2 \mT_{jj}(\tfrac{\partial}{\partial \ov{z}^0_i})
\mT_{kk}(\tfrac{\partial}{\partial z^0_i})\right] P^N(0,0).
\end{multline}

By \eqref{g7}, \eqref{g3.57},  \eqref{4.26}, we get
\begin{align}\label{a4.36}
\mI_3 P^N=- \frac{\sqrt{-1}}{8\pi} \wi{\mT}_{ijj'}
\Big[  b^{\bot}_j b^{\bot}_{j'} \mF_1(e^\bot_i)
+ \mT_{ilm}  b^{\bot}_j b^{\bot}_{j'} \frac{b^{\bot}_l b^{\bot}_m}{16\pi}
+ \frac{1}{2}  \mT_{ilj'}b^{\bot}_j b^{\bot}_l\Big]P^N.
\end{align}
By  \eqref{3.11e},  \eqref{34.21},
 \eqref{4.19},  \eqref{a4.36} and a similar equation as \eqref{4.32} for
$\wi{\mT}_{jij'}\mT_{kli}  b^{\bot}_j b^{\bot}_{j'} b^{\bot}_k b^{\bot}_l$,
 we get
\begin{align}\label{a4.37}
\left((\cL^0_2)^{-1} P^{N^\bot} \mI_3 P^N\right)(0,0)
=\frac{\sqrt{-1}}{64\pi} \wi{\mT}_{ijj'}\mT_{ijj'} P^N(0,0) =0,
\end{align}
as $ \wi{\mT}_{ijj'}$ is anti-symmetric on $i,j$
and $\mT_{ijj'}$ is symmetric on $i,j$.

By Theorem \ref{t3.4}, \eqref{g7}, \eqref{g3.57} and \eqref{4.19},
\begin{multline}\label{4.28}
\Big((\cL^0_2)^{-1} P^{N^\bot} \mF_1(e^\bot_j) (b^{\bot +}_j+ b^{\bot }_j)
\mF_1(e^\bot_k)\frac{b^{\bot }_k}{4\pi} P^N\Big)(0,0)\\
= \frac{1}{32\pi^2} \Big(\mF_1(e^\bot_j)^2 (b^{\bot }_j)^2  P^N\Big)(0,0)
=-\frac{1}{8\pi} \sum_j \mF_1(e^\bot_j)^2 P^N(0,0).
\end{multline}

Recall that $\mT_{klm}$ is symmetric on $k,l,m$.

 By Theorem \ref{t3.4}, \eqref{g7}, \eqref{g3.57},
\eqref{3.53} and \eqref{4.19},
\begin{multline}\label{4.29}
\Big \{(\cL^0_2)^{-1} P^{N^\bot}
\Big ( \mF_1(e^\bot_j) (b^{\bot +}_j+ b^{\bot }_j)
\mT_{klm} \frac{b^{\bot}_mb^{\bot}_lb^{\bot}_k}{192\pi ^2}
+ \mT_{ijj'}\frac{B^\bot_{ijj'}}{64\pi ^2}\mF_1(e^\bot_k)b^{\bot }_k
\Big )P^N\Big \}(0,0)\\
= \Big\{(\cL^0_2)^{-1} P^{N^\bot}  \mF_1(e^\bot_j)\Big (b^{\bot }_j
\mT_{klm} \frac{b^{\bot}_k b^{\bot}_l b^{\bot}_m}{48\pi ^2}
+ \mT_{jlm}
\frac{b^{\bot}_l b^{\bot}_m}{4\pi}\Big )P^N\Big \}(0,0)\\
=\frac{1}{32\pi ^2}\Big\{ \mF_1(e^\bot_j)\Big (\mT_{klm}
\frac{b^{\bot }_j b^{\bot}_k b^{\bot}_l b^{\bot}_m}{24\pi}
+ \mT_{jlm} b^{\bot}_l b^{\bot}_m
\Big )P^N\Big \}(0,0)\\
= \frac{1}{32\pi ^2}\Big\{ \mF_1(e^\bot_j)\Big(
 \sum_{l\neq j}\mT_{jll}
\frac{(b^{\bot }_j)^2 (b^{\bot}_l)^2}{8\pi}
+ \mT_{jjj}\frac{(b^{\bot }_j)^4}{24 \pi}
+ \mT_{jll} (b^{\bot}_l)^2
\Big )P^N\Big \}(0,0)\\
=-\frac{1}{16\pi}  \mF_1(e^\bot_j)\mT_{jll} P^N(0,0).
\end{multline}

As $\mT_{klm}$ is symmetric on $k,l,m$, we know that
\begin{align}\label{4.30}
&\mT_{klm}^2=6\sum_{k<l< m}\mT_{klm}^2+ 3\sum_{k\neq m}\mT_{kkm}^2
+\mT_{mmm}^2,\\
&\mT_{kkm}\mT_{llm}= \sum_{k\neq l\neq m\neq k}\mT_{kkm}\mT_{llm}
+  \sum_{k\neq m}(2 \mT_{kkm}\mT_{mmm} +  \mT_{kkm}^2) + \mT_{mmm}^2.
\nonumber
\end{align}
 From \eqref{4.19},  \eqref{4.30},
\begin{multline}\label{4.31}
\Big(\mT_{ijj'}\mT_{klm}
b^{\bot}_i b^{\bot}_j b^{\bot}_{j'}
b^{\bot}_k b^{\bot}_l b^{\bot}_m P^N\Big)(0,0)
= \left\{\Big( 36\sum_{k<l< m}\mT_{klm}^2 (b^{\bot}_k)^2
 (b^{\bot}_l)^2 (b^{\bot}_m)^2\right.\\
 +9 \sum_{k\neq l\neq m\neq k}\mT_{kkm}\mT_{llm}
(b^{\bot}_k)^2 (b^{\bot}_l)^2 (b^{\bot}_m)^2
+ 6 \sum_{k\neq m} \mT_{kkm}\mT_{mmm}  (b^{\bot}_k)^2(b^{\bot}_m)^4\\
\left.+ 9  \sum_{k\neq m}\mT_{mmk}\mT_{mm k} (b^{\bot}_k)^2(b^{\bot}_m)^4
+  \mT_{mmm}^2 (b^{\bot}_m)^6\Big)P^N\right\}(0,0)\\
= (-4\pi)^3 \Big(36 \sum_{k<l< m}\mT_{klm}^2
 +9  \sum_{k\neq l\neq m\neq k}\mT_{kkm}\mT_{llm}\\
+  3 \sum_{k\neq m}
(6 \mT_{kkm}\mT_{mmm} +9 \mT_{mmk}\mT_{mm k})
+  15 \mT_{mmm}^2\Big)P^N(0,0)\\
=  (-4\pi)^3   \cdot 3  \Big( 2 \mT_{klm}^2
+ 3  \mT_{kkm}\mT_{llm}\Big)P^N(0,0).
\end{multline}

By  \eqref{4.19},
\begin{multline}\label{4.32}
\Big(\mT_{ijm}\mT_{klm}b^{\bot}_i b^{\bot}_jb^{\bot}_k b^{\bot}_l
P^N\Big)(0,0)\\
= \left\{\left[\sum_{k\neq l}\Big(2\mT_{klm}^2
+ \mT_{kkm}\mT_{llm}\Big)(b^{\bot}_k)^2 (b^{\bot}_l)^2
+ \mT_{llm}^2 (b^\bot_l)^4 \right]P^N\right\}(0,0)\\
= (4\pi)^2 \Big(\sum_{k\neq l}(2\mT_{klm}^2
+ \mT_{kkm}\mT_{llm})
+3 \mT_{llm}^2\Big) P^N(0,0)\\
=  (4\pi)^2 (2\mT_{klm}^2 + \mT_{kkm}\mT_{llm})P^N(0,0).
\end{multline}

By \eqref{g7}, \eqref{g3.57} and \eqref{3.53},
we have also
\begin{multline}\label{4.33}
P^{N^\bot} \mT_{ijj'}B^\bot_{ijj'}\mT_{klm}
b^{\bot}_k b^{\bot}_l b^{\bot}_mP^N
=\Big(\mT_{ijj'}\mT_{klm}
b^{\bot}_i b^{\bot}_j b^{\bot}_{j'}
b^{\bot}_k b^{\bot}_l b^{\bot}_m\\
+ 36\pi \mT_{ijm}\mT_{klm}  b^{\bot}_i b^{\bot}_j b^{\bot}_k b^{\bot}_l
+ 36\pi \cdot 8\pi \mT_{ilm}\mT_{klm}  b^{\bot}_i b^{\bot}_k\Big)P^N.
\end{multline}
Thus from Theorem \ref{t3.4}, \eqref{4.31}-\eqref{4.33},
\begin{multline}\label{4.34}
\left\{\Big((\cL^0_2)^{-1} P^{N^\bot}
\frac{1}{16\pi} \mT_{ijj'}B^\bot_{ijj'}\mT_{klm}
\frac{b^{\bot}_k b^{\bot}_l b^{\bot}_m}{192\pi ^2}\Big)P^N\right\}(0,0)\\
= \frac{1}{2^{10}\cdot 3\pi^3}
\left\{\Big( \frac{1}{24\pi}\mT_{ijj'}\mT_{klm}
b^{\bot}_i b^{\bot}_j b^{\bot}_{j'}
b^{\bot}_k b^{\bot}_l b^{\bot}_m
+ \frac{9}{4} \mT_{ijm}\mT_{klm}
b^{\bot}_i b^{\bot}_j b^{\bot}_k b^{\bot}_l\right.\\
\left.+  36\pi \mT_{ilm}\mT_{klm}  b^{\bot}_i b^{\bot}_k\Big)P^N\right\}(0,0)\\
= \frac{1}{2^{10}\cdot 3\pi}
\left\{-8 (2\mT_{klm}^2 +3 \mT_{kkm}\mT_{llm})
+  36 (2\mT_{klm}^2 + \mT_{kkm}\mT_{llm})
-144 \mT_{klm}^2\right\}P^N(0,0)\\
=\frac{1}{2^{8}\cdot 3\pi} \left(-22\mT_{klm}^2
+3 \mT_{kkm}\mT_{llm} \right)P^N(0,0).
\end{multline}

From \eqref{4.27}, \eqref{4.35}, \eqref{4.37}, \eqref{a4.37},
\eqref{4.28}, \eqref{4.29} and \eqref{4.34},  we get
\begin{multline}\label{4.38}
\Psi_{1,1}(0)= \Big\{
\frac{1}{2^8\cdot 3\pi}
\left[-76|\mT_{jj'}(\tfrac{\partial}{\partial \ov{z}^0_i})|^2
+2 \mT_{jj}(\tfrac{\partial}{\partial \ov{z}^0_i})
\mT_{kk}(\tfrac{\partial}{\partial z^0_i})
-22\mT_{klm}^2 +3 \mT_{kkm}\mT_{llm} \right]\\
-\frac{1}{24\pi} \mT_{jj}(\tfrac{\partial}{\partial z^0_l})
\mT_{kk}(\tfrac{\partial}{\partial \ov{z}^0_l})
- \frac{1}{8\pi}\sum_j \mF_1(e^\bot_j)^2
- \frac{1}{16\pi} \mF_1(e^\bot_j)\mT_{jll}
\Big\}P^N(0,0).
\end{multline}

From \eqref{4.38} we get \eqref{4.25}.
\end{proof}

Recall that $B(Z, e^\bot_l)$ was defined in \eqref{a3.18}.
\begin{lemma}\label{t4.12} The following identity holds,
\begin{multline}\label{4.39}
\frac{\sqrt{-1}}{\pi}B(Z, e^\bot_l)
=\frac{1}{2} \left\langle R^{TB}(\mR^\bot,\mR^0)e^\bot_l, J\mR^0\right\rangle
-\frac{5}{4}\left\langle \nabla ^{TY}_{\mR} (T(e_k,e^\bot_l)),
J\mR^\bot\right\rangle Z_k \\
+ \frac{1}{2} \left\langle \frac{1}{3} R^{TB}(\mR^\bot,e^\bot_l)\mR^\bot
 + \nabla_{\mR^0}^{TX_G}(A(e^0_k)e^\bot_l)Z^0_k,
J\mR^0\right\rangle\\
+\frac{1}{8}\left\langle T(\mR^0, e^0_j), Je^\bot_l\right\rangle
\left\langle T(\mR^\bot-\mR^0,J e^0_j), J\mR^\bot\right\rangle\\
+\frac{1}{4}\left\langle T(\mR^\bot, e^0_j), Je^\bot_l\right\rangle
\left\langle T(\mR^0,J e_j^0),J\mR^\perp\right\rangle\\
+\frac{1}{8}\left\langle T(\mR^0,J \mR^0), T(\mR^\perp, e^\bot_l)\right\rangle
-\frac{1}{8}\left\langle T(\mR, e^\bot_l), T(\mR^\bot,J \mR^0)\right\rangle\\
+\frac{1}{8}\left\langle  T(\mR^\perp, JT(\mR^0,J \mR^0)),
J e^\bot_l\right\rangle
+ \frac{1}{2}\left\langle T(\mR^\bot,J \mR^\bot),T(\mR,e^\bot_l)\right\rangle.
\end{multline}
\end{lemma}
\begin{proof}
By \eqref{b3.18}, \eqref{3.33} and  $A(\mR^0)A(\mR^0)e^\bot_l\in N_G$,
as $A$ exchanges $TX_G$ and $N_G$, we get
\begin{multline}\label{4.41}
\langle J\mR, (\nabla^{TX}\nabla^{TX}e_l^{\bot,H} )_{(\mR,\mR)} \rangle
=-\frac{1}{2}\left\langle J\mR, T(\mR, \nabla ^{TB}_{\mR}e^\bot_l)
+\nabla ^{TX}_{\mR} (T(e^H_i,e^\bot_l))Z_i\right\rangle\\
+  \left\langle J\mR^0,
\frac{1}{3} R^{TB}(\mR^\bot,e^\bot_l)\mR^\bot
+ R^{TB}(\mR^\bot,\mR^0)e^\bot_l
+ \nabla_{\mR^0}^{TX_G}(A(e^0_k)e^\bot_l)Z^0_k\right\rangle.
\end{multline}
By \eqref{ah4},  \eqref{33.36}, \eqref{a3.311}, we have at $x_0$,
\begin{align}\label{4.42}
& -\frac{1}{2}\left\langle J\mR^\bot,
T(\mR, \nabla ^{TB}_{\mR}e^\bot_l)\right\rangle
= \frac{1}{4}\left\langle J e^\bot_l, T(\mR^0, e^0_j)  \right\rangle
\left\langle J\mR^\bot, T(\mR, J e^0_j) \right\rangle,\\
&-\frac{1}{2}\left\langle J\mR^0,  \nabla ^{TX}_{\mR}
(T(e^H_i,e^{\bot,H}_l))Z_i\right\rangle
= -\frac{1}{4} \left\langle T(\mR, e^\bot_l), T(\mR, J\mR^0) \right\rangle
.\nonumber
\end{align}

By \eqref{3.11a},  \eqref{3.11d},  \eqref{33.36}, \eqref{a3.311}, \eqref{3.33}
and $ \nabla ^{TX}_{\mR} (T(e^H_i,e^H_k))Z_iZ_k=0$, we have
\begin{multline}\label{4.43}
\langle J(\nabla^{TX}\nabla^{TX}e_k^H )_{(\mR,\mR)},e_l^{\bot} \rangle Z_k
=\frac{1}{2}\left\langle T(\mR, \nabla ^{TB}_{\mR}e_k),
 Je^\bot_l\right\rangle Z_k\\
=  \frac{1}{2}\left\langle T(\mR, 2A(\mR^0)\mR^\bot +A(\mR^0)\mR^0),
Je^\bot_l\right\rangle\\
=\frac{1}{2}\left\langle T(\mR, e^0_j), Je^\bot_l\right\rangle
\left\langle T(\mR^0, J e^0_j), J\mR^\bot\right\rangle\\
-\frac{1}{4}\left\langle  T(\mR^0,e^\bot_l),
 T(\mR^0,J \mR^0)\right\rangle
+\frac{1}{4}\left\langle  T(\mR^\perp, JT(\mR^0,J \mR^0)),
J e^\bot_l\right\rangle .
\end{multline}

From \eqref{g3.29}, \eqref{3.11a}, \eqref{33.36}, \eqref{a3.311}
 and the fact that $A$ exchanges $TX_G$ and $N_G$, we get
\begin{align}\label{4.44}
&\left\langle J \nabla^{TX} _{\mR}e^H_k,
\nabla^{TX} _{\mR}e^{\bot,H}_l\right\rangle Z_k
=\left\langle J \nabla^{TB} _{\mR}e_k, A(\mR^0)  e^{\bot}_l
 -\frac{1}{2} T(\mR,e^\bot_l)\right\rangle Z_k  \\
&\hspace*{15mm} = \left\langle J  A(\mR^0) \mR^0,  -\frac{1}{2}
T(\mR,e^\bot_l)\right\rangle +2 \left\langle J  A(\mR^0) \mR^\bot,
A(\mR^0) e^\bot_l\right\rangle,\nonumber\\
& = \frac{1}{4}\left\langle T(\mR^0,J \mR^0),   T(\mR,e^\bot_l)\right\rangle
-\frac{1}{2} \left\langle J e^\bot_l,T(\mR^0,e_j^0)\right\rangle
\left\langle J\mR^\perp,T(\mR^0, J e_j^0)\right\rangle.\nonumber
\end{align}

From \eqref{a3.31}, \eqref{b3.30}, \eqref{3.38},
\eqref{4.41}-\eqref{4.44}, we get
\begin{multline}\label{4.45}
\frac{\sqrt{-1}}{\pi}B(Z, e^\bot_l)
=\frac{1}{8}\left\langle J e^\bot_l, T(\mR^0,  e^0_j) \right\rangle
\left\langle J\mR^\bot, T(\mR,J e^0_j) \right\rangle \\
-\frac{1}{4}\left\langle J\mR^\bot,
\nabla ^{TY}_{\mR} (T(e_i,e^\bot_l))Z_i\right\rangle
-\frac{1}{8} \left\langle T(\mR, e^\bot_l), T(\mR, J\mR^0) \right\rangle\\
+ \frac{1}{2} \left\langle J\mR^0,
\frac{1}{3} R^{TB}(\mR^\bot,e^\bot_l)\mR^\bot
+ R^{TB}(\mR^\bot,\mR^0)e^\bot_l
+ \nabla_{\mR^0}^{TX_G}(A(e^0_k)e^\bot_l)Z^0_k\right\rangle\\
+\frac{1}{4}\left\langle T(\mR, e^0_j), Je^\bot_l\right\rangle
\left\langle T(\mR^0, J e^0_j), J\mR^\bot\right\rangle
-\frac{1}{8}\left\langle T(\mR^0,e^\bot_l),
 T(\mR^0,J \mR^0)\right\rangle\\
+\frac{1}{8}\left\langle  T(\mR^\perp, JT(\mR^0,J \mR^0)),
J e^\bot_l\right\rangle
+\frac{1}{4}\left\langle T(\mR^0,J \mR^0),   T(\mR,e^\bot_l)\right\rangle \\
- \frac{1}{2}\left\langle J e^\bot_l,T(\mR^0,e_j^0)\right\rangle
\left\langle J\mR^\perp,T(\mR^0,J e_j^0)\right\rangle\\
+\frac{1}{2}\left\langle T(\mR^\bot, J \mR^\bot) ,T(\mR,e^\bot_l)\right\rangle
- \left\langle \nabla ^{TY}_{\mR} (T(e_k,e^\bot_l)),
J\mR^\bot\right\rangle Z_k  .
\end{multline}

From \eqref{4.45} we get \eqref{4.39}.
\end{proof}

Now we need to compute the contribution from
$-(\cL^0_2)^{-1} P^{N^\bot}\mO_2 P^N$.
\begin{lemma}\label{t4.13} The following identity holds,
\begin{multline}\label{4.46}
\Psi_{1,2}(0)=\left\{\frac{1}{16\pi} r^{X_G}_{x_0}
 +\frac{1}{2\pi} R^{E_G}(\tfrac{\partial}{\partial z^0_j},
\tfrac{\partial}{\partial \ov{z}^0_j}) + \frac{1}{2\pi}\Delta_{X_G} \log h
+\frac{29}{2^5\cdot 3\pi}|T(e^\bot_k,\tfrac{\partial}{\partial \ov{z}^0_j})|^2
\right.\\
+\frac{\sqrt{-1}}{16\pi}\left\langle T(e^\bot_k,J e^\bot_k ),
T(\tfrac{\partial}{\partial z^0_j},
\tfrac{\partial}{\partial \ov{z}^0_j})\right\rangle
+ \frac{1}{4\pi}\Big|T(\tfrac{\partial}{\partial z^0_i},
\tfrac{\partial}{\partial \ov{z}^0_j})\Big|^2
+\frac{1}{32\pi} \Big|\sum_j
\mT_{jj}(\tfrac{\partial}{\partial \ov{z}^0_i})\Big|^2\\
-\frac{1}{2^6 \pi} \left\langle (\nabla^{TY}_\cdot
\dot{g}^{TY}_\cdot)_{(e^\bot_j,e^\bot_j)}
Je^\bot_k, Je^\bot_k\right\rangle
-\frac{1}{2^5 \pi} \left\langle (\nabla^{TY}_\cdot
\dot{g}^{TY}_\cdot)_{(e^\bot_j,e^\bot_k)}
Je^\bot_j, Je^\bot_k\right\rangle\\
+ \frac{1}{2^7 \pi}\wi{\mT}_{ijk} (\wi{\mT}_{kji} + \wi{\mT}_{ijk})
+ \frac{7}{2^8 \pi} \Big(2\mT_{jkm}^2 + \mT_{jjm}\mT_{kkm} \Big)\\
\left. - \frac{\sqrt{-1}}{16\pi} \Big(\left\langle T(e^\bot_j, Je^\bot_j),
\wi{\mu}^E\right\rangle
-2\left\langle Je^\bot_j,
\nabla^{TY}_{e^\bot_j}\wi{\mu}^E\right\rangle\Big)\right\} P^N (0,0).
\end{multline}
\end{lemma}
\begin{proof} From Theorem \ref{t3.4}, \eqref{g7}, \eqref{g3.57},
\eqref{b3.11} and \eqref{4.19},
\begin{multline}\label{4.48}
4\pi \left((\cL^0_2)^{-1} P^{N^\bot} Z^\bot_kZ^\bot_l P^N\right)(0,0)
=  \left((\cL^0_2)^{-1} P^{N^\bot}Z^\bot_l b^\bot_k P^N\right)(0,0)\\
=  \left((\cL^0_2)^{-1} b^\bot_kZ^\bot_l P^N\right)(0,0)
= \Big(\frac{b^\bot_k b^\bot_l}{32\pi^2} P^N\Big)(0,0)
 =-\frac{\delta_{kl}}{8\pi} P^N(0,0).
\end{multline}

Set
\begin{align}\label{a4.47}
\mI_4=- \left\{(\cL^0_2)^{-1}P^{N^\bot}
\Big(\tfrac{\partial}{\partial z^0_j}
\Big(B(Z,\tfrac{\partial}{\partial \ov{z}^0_j}) \Big)
-\tfrac{\partial}{\partial \ov{z}^0_j}
\Big(B(Z,\tfrac{\partial}{\partial z^0_j}) \Big)  \Big)P^N\right\}(0,0).
\end{align}
At first, if $Q$ is a monomial on $b_i,b_i^+$,
$b^\bot_j, b^{\bot+}_j, Z_i$ and the total degree of $b_i,b_i^+, Z^0_i$
 or $b^\bot_j, b^{\bot+}_j, Z^\bot_j$ is odd, then by Theorem \ref{t3.4},
\begin{align}\label{a4.45}
\left((\cL^0_2)^{-1}P^{N^\bot}Q P^N\right)(0,0)=0.
\end{align}

By \eqref{a4.45}, only the monomials of $B(Z,e^0_l)$
with odd degree on $Z^0$ have contributions for $\mI_4$.

If we denote by $\wi{B}_Z(e^0_l)$ the odd degree component on $Z^0$ of
the difference of $B(Z,e^0_l)$
and of the sum the the first two and the last terms of $B(Z,e^0_l)$
 in \eqref{a3.32},
 then by \eqref{a3.32} we know that $\wi{B}_Z(e^0_l)$ is a linear
function on $Z^0$ and
$\tfrac{\partial}{\partial z^0_j}
\Big(\wi{B}_Z(\tfrac{\partial}{\partial \ov{z}^0_j}) \Big)$
and $-\tfrac{\partial}{\partial \ov{z}^0_j}
\Big(\wi{B}_Z(\tfrac{\partial}{\partial z^0_j}) \Big)$  are equal.

Moreover, by \eqref{3.11e}, \eqref{4.48}, we know the contribution of
the last term of $B(Z,e^0_l)$ in \eqref{a3.32} is zero in $\mI_4$.

Thus by Remark \ref{r4.3},  \eqref{a3.32} and \eqref{a4.47},
\begin{multline}\label{4.47}
\mI_4=\pi \sqrt{-1}\left\{(\cL^0_2)^{-1} P^{N^\bot}
\left[\frac{1}{6}\tfrac{\partial}{\partial z^0_j}
\left\langle R^{TX_G}(\mR^0,J\mR^0)\mR^0,
\tfrac{\partial}{\partial \ov{z}^0_j}\right\rangle \right.  \right.\\
- \frac{1}{6}\tfrac{\partial}{\partial \ov{z}^0_j}
\left\langle R^{TX_G}(\mR^0,J\mR^0)\mR^0,
\tfrac{\partial}{\partial z^0_j}\right\rangle \\
-\frac{5}{4}\Big\langle  J\mR^\bot,
2 \nabla ^{TY}_{\mR^\bot} (T(\tfrac{\partial}{\partial z^0_j},
\tfrac{\partial}{\partial \ov{z}^0_j}))
+ \nabla ^{TY}_{\tfrac{\partial}{\partial z^0_j}} (T(e^\bot_i,
\tfrac{\partial}{\partial \ov{z}^0_j}))Z^\bot_i
-\nabla ^{TY}_{\tfrac{\partial}{\partial \ov{z}^0_j}} (T(e^\bot_i,
\tfrac{\partial}{\partial z^0_j}))Z^\bot_i \Big\rangle\\
+ 3 \sqrt{-1} \left\langle  R^{TB}(\mR^\bot,
\tfrac{\partial}{\partial z^0_j}) \mR^\bot,
\tfrac{\partial}{\partial \ov{z}^0_j}\right\rangle
-\frac{3\sqrt{-1}}{4} \left\langle J\mR^\bot,
T(\tfrac{\partial}{\partial z^0_j}, e^0_i)
\right\rangle
 \left\langle J\mR^\bot,
T(e^0_i, \tfrac{\partial}{\partial \ov{z}^0_j})\right\rangle\\
-\frac{\sqrt{-1}}{4}\left\langle  T(\mR^\bot, \tfrac{\partial}{\partial \ov{z}^0_j}),
 T(\mR^\bot,\tfrac{\partial}{\partial z^0_j})\right\rangle
\left.\left.+ \left\langle T(\mR^\bot, J \mR^\bot),
 T(\tfrac{\partial}{\partial z^0_j},
\tfrac{\partial}{\partial \ov{z}^0_j})\right\rangle
\right]P^N\right\}(0,0).
\end{multline}

By \eqref{3.64}, \eqref{a3.70}, \eqref{4.48} and \eqref{4.47},
comparing with \eqref{3.68} and \eqref{3.69}, we get
\begin{multline}\label{4.49}
\mI_4 = \left\{-\frac{1}{6\pi}
 \left\langle  R^{TX_G}(\tfrac{\partial}{\partial z^0_i},
\tfrac{\partial}{\partial \ov{z}^0_i})
\tfrac{\partial}{\partial z^0_j}
+ R^{TX_G}(\tfrac{\partial}{\partial z^0_j},
\tfrac{\partial}{\partial \ov{z}^0_i})
\tfrac{\partial}{\partial z^0_i},
\tfrac{\partial}{\partial \ov{z}^0_j}\right\rangle\right.\\
+\frac{5\sqrt{-1}}{2^7\pi}\left\langle  Je^\bot_k,
2 \nabla ^{TY}_{e^\bot_k} (T(\tfrac{\partial}{\partial z^0_j},
\tfrac{\partial}{\partial \ov{z}^0_j}))
+\nabla ^{TY}_{\tfrac{\partial}{\partial z^0_j}} (T(e^\bot_k,
\tfrac{\partial}{\partial \ov{z}^0_j}))
-\nabla ^{TY}_{\tfrac{\partial}{\partial \ov{z}^0_j}} (T(e^\bot_k,
\tfrac{\partial}{\partial z^0_j})) \right\rangle\\
+ \frac{3}{32\pi}  \left\langle  R^{TB}(e^\bot_k,
\tfrac{\partial}{\partial z^0_j}) e^\bot_k,
\tfrac{\partial}{\partial \ov{z}^0_j}\right\rangle
+\frac{3}{64\pi} |T(\tfrac{\partial}{\partial z^0_i},
\tfrac{\partial}{\partial \ov{z}^0_j})|^2\\
\left. - \frac{1}{2^7\pi}|T(e^\bot_k, \tfrac{\partial}{\partial \ov{z}^0_j})|^2
- \frac{\sqrt{-1}}{32\pi}\left\langle T(e^\bot_k, J e^\bot_k),
  T(\tfrac{\partial}{\partial z^0_j},
\tfrac{\partial}{\partial \ov{z}^0_j})\right\rangle \right\}
P^N(0,0).
\end{multline}
By  \eqref{g7},
 \eqref{g3.57} and \eqref{3.53},
\begin{align}\label{a4.50}
\begin{split}
&(z^0_i\ov{z}^0_j P^N)(Z,0)= (z^0_i\frac{b_j}{2\pi} P^N)(Z,0)
= \frac{1}{2\pi}((b_j z^0_i +2\delta_{ij})P^N)(Z,0),\\
&Z^\bot_kZ^\bot_l P^N= \frac{1}{16\pi^2}
(b^\bot_kb^\bot_l+4\pi \delta_{kl}) P^N,\\
&(4\pi)^4 (Z^\bot_k)^4 P^N
= \left((b^\bot_k)^4 +24\pi (b^\bot_k)^2 + 3\cdot (4\pi)^2\right)P^N.
\end{split}\end{align}
From Theorem \ref{t3.4}, \eqref{g7}, \eqref{g3.57}, \eqref{3.64},
 \eqref{4.19}, \eqref{a4.34} and \eqref{a4.50},
\begin{align}\label{4.50}
& (P^{N^\bot}Z^\bot_kZ^\bot_l P^N)(0,0)
= \frac{1}{16\pi^2} (b^\bot_kb^\bot_l P^N)(0,0)=-\frac{\delta_{kl}}{4\pi},\nonumber\\
&\left((\cL^0_2)^{-1}P^{N^\bot}b_j z^0_i Z^\bot_kZ^\bot_l P^N\right)(0,0)\\
&\hspace*{5mm}
= \frac{1}{16\pi^2}\Big\{\Big(\frac{1}{12\pi}b^\bot_kb^\bot_l b_j z^0_i
+ \delta_{kl}b_j z^0_i\Big)  P^N\Big\}(0,0)
=-\frac{1}{12\pi^2}\delta_{ij}\delta_{kl}P^N(0,0),\nonumber\\
&\left((\cL^0_2)^{-1} b^\bot_l Z^\bot_k z^0_i\ov{z}^0_j P^N\right)(0,0)
= \frac{1}{8\pi^2} \left\{b^\bot_l b^\bot_k \Big(\frac{b_j}{12\pi} z^0_i
+\frac{2}{8\pi}\delta_{ij}\Big)P^N\right\}(0,0)\nonumber\\
&\hspace*{5mm}
= - \frac{1}{24\pi^2}\delta_{ij}\delta_{kl}P^N(0,0),\nonumber\\
&\left((\cL^0_2)^{-1}P^{N^\bot} Z^\bot_l Z^\bot_k z^0_i\ov{z}^0_j
P^N\right)(0,0)\nonumber\\
&\hspace*{5mm}= \frac{1}{4\pi} \left\{(\cL^0_2)^{-1}P^{N^\bot}
(b^\bot_l Z^\bot_k+ \delta_{kl}) z^0_i\ov{z}^0_j P^N\right\}(0,0)
=\frac{-7}{96\pi^3} \delta_{ij}\delta_{kl}P^N(0,0).\nonumber
\end{align}

By \eqref{3.11e}, \eqref{3.70}, \eqref{a3.70}, \eqref{4.50}
and comparing with \eqref{3.71}, we get
\begin{multline}\label{4.51}
-\left((\cL^0_2)^{-1}P^{N^\bot}
b_jB(Z,\tfrac{\partial}{\partial \ov{z}^0_j}) P^N\right)(0,0)\\
= \Big\{- \frac{1}{12\pi}
  \left\langle R^{TX_G}(\tfrac{\partial}{\partial z^0_j},
\tfrac{\partial}{\partial \ov{z}^0_i})\tfrac{\partial}{\partial z^0_i}
+ R^{TX_G}(\tfrac{\partial}{\partial z^0_i},
\tfrac{\partial}{\partial \ov{z}^0_i})\tfrac{\partial}{\partial z^0_j},
\tfrac{\partial}{\partial \ov{z}^0_j}\right\rangle\\
+\frac{5\sqrt{-1}}{48\pi} \Big\langle
\nabla ^{TY}_{\tfrac{\partial}{\partial z^0_j}}
(T(e^\bot_k,\tfrac{\partial}{\partial \ov{z}^0_j}))
+\nabla ^{TY}_{e^\bot_k} (T(\tfrac{\partial}{\partial z^0_j},
\tfrac{\partial}{\partial \ov{z}^0_j})),
 Je^\bot_k \Big\rangle\\
+ \frac{1}{8\pi} \left\langle
R^{TB}(e^\bot_k, \tfrac{\partial}{\partial z^0_j})e^\bot_k,
\tfrac{\partial}{\partial \ov{z}^0_j}\right\rangle
+ \frac{1}{16\pi}\Big|T(\tfrac{\partial}{\partial z^0_i},
\tfrac{\partial}{\partial \ov{z}^0_j})\Big|^2\\
- \frac{1}{96\pi}
\Big|T(e^\bot_k,\tfrac{\partial}{\partial \ov{z}^0_j})\Big|^2
- \frac{\sqrt{-1}}{24\pi} \left\langle T(e^\bot_k,Je^\bot_k),
T(\tfrac{\partial}{\partial z^0_j},
\tfrac{\partial}{\partial \ov{z}^0_j}) \right\rangle
\Big\} P(0,0).
\end{multline}

From \eqref{4.39} and \eqref{a4.45},
\begin{multline}\label{4.52}
\left((\cL^0_2)^{-1} b^\bot_lB(Z,e^\bot_l) P^N\right)(0,0)
=-\pi\sqrt{-1}\left\{(\cL^0_2)^{-1} b^\bot_l \right.\\
\Big[
\frac{1}{2} \left\langle R^{TB}(\mR^\bot,\mR^0)e^\bot_l, J\mR^0\right\rangle
-\frac{5}{4}\left\langle \nabla ^{TY}_{\mR^\bot} (T(e^\bot_k,e^\bot_l)),
J\mR^\bot\right\rangle Z_k^\bot  \\
-\frac{5}{4}\left\langle \nabla ^{TY}_{\mR^0} (T(e^0_k,e^\bot_l)),
J\mR^\bot\right\rangle Z_k^0
-\frac{1}{8}\left\langle T(\mR^0, e^0_j), Je^\bot_l\right\rangle
\left\langle T(\mR^0,J e^0_j), J\mR^\bot\right\rangle\\
+\frac{1}{8}\left\langle T(\mR^0,J \mR^0), T(\mR^\perp, e^\bot_l)\right\rangle
-\frac{1}{8}\left\langle T(\mR^0, e^\bot_l), T(\mR^\bot,J \mR^0)\right\rangle\\
+\frac{1}{8}\left\langle  T(\mR^\perp, JT(\mR^0,J \mR^0)),
J e^\bot_l\right\rangle
\left.+ \frac{1}{2}\left\langle T(\mR^\bot,J \mR^\bot),T(\mR^\bot,e^\bot_l)\right\rangle\Big] P^N\right\}(0,0).
\end{multline}
As $T$ is anti-symmetric, from \eqref{g7}, \eqref{g3.57},  we get
\begin{align}\label{4.53}
b^\bot_l\left\langle \nabla ^{TY}_{\mR^\bot} (T(e^\bot_k,e^\bot_l)),
J\mR^\bot\right\rangle Z_k^\bot P^N
= -\Big(\tfrac{\partial}{\partial Z^\bot_l}
\left\langle \nabla ^{TY}_{\mR^\bot} (T(e^\bot_k,e^\bot_l)),
J\mR^\bot\right\rangle\Big) Z_k^\bot P^N,\\
b^\bot_l\left\langle T(\mR^\bot,J \mR^\bot),
T(\mR^\bot,e^\bot_l)\right\rangle P^N
= -\left\langle T(\mR^\bot,Je^\bot_l)+ T(e^\bot_l,J \mR^\bot),
T(\mR^\bot,e^\bot_l)\right\rangle P^N.\nonumber
\end{align}

From \eqref{3.11e}, \eqref{4.7}, \eqref{4.48}, \eqref{4.50},
 \eqref{4.52}, \eqref{4.53} and the anti-symmetric property of $T$, we get
\begin{multline}\label{4.54}
-\frac{1}{2}\left((\cL^0_2)^{-1} b^\bot_l B(Z,e^\bot_l) P^N\right)(0,0)\\
= \frac{\sqrt{-1}}{2\pi} \left\{
-\frac{5}{2^7}\Big(\left\langle \nabla ^{TY}_{e^\bot_k}(T(e^\bot_k,e^\bot_l)),
Je^\bot_l\right\rangle
+ \left\langle \nabla ^{TY}_{e^\bot_l} (T(e^\bot_k,e^\bot_l)),
Je^\bot_k\right\rangle\Big) \right.\\
+ \frac{5}{96}\Big\langle \nabla ^{TY}_{\tfrac{\partial}{\partial z^0_j}}
 (T(\tfrac{\partial}{\partial \ov{z}^0_j},e^\bot_l))
+\nabla ^{TY}_{\tfrac{\partial}{\partial \ov{z}^0_j}}
 (T(\tfrac{\partial}{\partial z^0_j},e^\bot_l)),J e^\bot_l\Big\rangle\\
\left.+\frac{1}{2^6} \left\langle T(e^\bot_k, J e^\bot_l)
+ T(e^\bot_l, J e^\bot_k),
T(e^\bot_k,  e^\bot_l)\right\rangle\right\}P^N(0,0)=0.
\end{multline}

By \eqref{a3.67}, \eqref{4.7}, \eqref{a4.47}, \eqref{4.49},
\eqref{4.51}, \eqref{4.54} and since $R^{TX_G}(\cdot,\cdot)$ is a
$(1,1)$-form, comparing with \eqref{3.69} and \eqref{3.71}, we get
\begin{multline} \label{a4.54}
-\left((\cL^0_2)^{-1}P^{N^\bot}I_1 P^N \right)(0,0)
= \left\{-\frac{1}{2\pi}
 \left\langle  R^{TX_G}(\tfrac{\partial}{\partial z^0_i},
\tfrac{\partial}{\partial \ov{z}^0_i})
\tfrac{\partial}{\partial z^0_j},
\tfrac{\partial}{\partial \ov{z}^0_j}\right\rangle   \right.\\
+\frac{7}{6} \left[
\frac{5\sqrt{-1}}{2^5 \pi}\left\langle  Je^\bot_k,
 \nabla ^{TY}_{e^\bot_k} (T(\tfrac{\partial}{\partial z^0_j},
\tfrac{\partial}{\partial \ov{z}^0_j}))
+ \nabla ^{TY}_{\tfrac{\partial}{\partial z^0_j}} (T(e^\bot_k,
\tfrac{\partial}{\partial \ov{z}^0_j}))\right\rangle  \right.\\
+ \frac{3}{16\pi}  \left\langle  R^{TB}(e^\bot_k,
\tfrac{\partial}{\partial z^0_j}) e^\bot_k,
\tfrac{\partial}{\partial \ov{z}^0_j}\right\rangle
+\frac{3}{32\pi} |T(\tfrac{\partial}{\partial z^0_i},
\tfrac{\partial}{\partial \ov{z}^0_j})|^2\\
 \left. \left.
-\frac{1}{2^6 \pi}|T(e^\bot_k, \tfrac{\partial}{\partial \ov{z}^0_j})|^2
-\frac{\sqrt{-1}}{16\pi}\left\langle T(e^\bot_k, J e^\bot_k),
 T(\tfrac{\partial}{\partial z^0_j},
\tfrac{\partial}{\partial \ov{z}^0_j}) \right\rangle
 \right]\right\}P^N(0,0). \\
\end{multline}

Recall that from \eqref{0g4}, \eqref{3.11a}, \eqref{3.11b} and \eqref{33.36},
\begin{align} \label{4.55}
\begin{split}&|A(e^0_i)e^\bot_k|^2
=4|A(\tfrac{\partial}{\partial z^0_i})e^\bot_k|^2
= |T(\tfrac{\partial}{\partial z^0_i},J e^0_j)|^2
=2 |T(\tfrac{\partial}{\partial z^0_i},
\tfrac{\partial}{\partial \ov{z}^0_j})|^2,\\
&\left\langle A(e^0_i)e^0_i, A(e^0_j)e^0_j\right\rangle
=4 \Big|\sum_i T(\tfrac{\partial}{\partial z^0_i},
\tfrac{\partial}{\partial \ov{z}^0_i})\Big|^2,\\
&|A(e^0_i)e^0_j|^2 = \frac{1}{4} |T(e^0_i, Je^0_j)|^2
= |T(\tfrac{\partial}{\partial z^0_i},  Je^0_j)|^2
= 2|T(\tfrac{\partial}{\partial z^0_i},
\tfrac{\partial}{\partial \ov{z}^0_j})|^2.
\end{split}\end{align}

From \eqref{3.64}, \eqref{a3.71}, \eqref{4.48}, \eqref{4.55} and
since $R^{TX_G}(\cdot,\cdot)$ is a $(1,1)$-form (comparing with
\eqref{3.73}, \eqref{3.74}), we get
\begin{multline} \label{4.56}
-\left((\cL^0_2)^{-1}P^{N^\bot}I_2 P^N \right)(0,0)
= \left\{
 \frac{4}{3\pi}\left\langle R^{TX_G} (\tfrac{\partial}{\partial z^0_j},
\tfrac{\partial}{\partial \ov{z}^0_i})\tfrac{\partial}{\partial z^0_i},
\tfrac{\partial}{\partial \ov{z}^0_j}\right\rangle\right.\\
-\frac{1}{8\pi} \left\langle
R^{TB} (e^\bot_k, \tfrac{\partial}{\partial z^0_j})e^\bot_k,
\tfrac{\partial}{\partial \ov{z}^0_j}\right\rangle
 \left. -\frac{1}{48\pi}\left\langle R^{TB} (e^\bot_k, e^\bot_j)e^\bot_k,
e^\bot_j\right\rangle
 + \frac{3}{16\pi} |T(\tfrac{\partial}{\partial z^0_i},
\tfrac{\partial}{\partial \ov{z}^0_j})|^2 \right\} P^N(0,0).
\end{multline}

By \eqref{0g4}, \eqref{g3.57}, \eqref{b3.15}, \eqref{3.533},
\eqref{3.64}, \eqref{a3.72},  \eqref{4.48}, \eqref{4.55} and since
$R^{TX_G}(\cdot,\cdot)$ is a $(1,1)$-form (comparing with
\eqref{3.73a}), we get
\begin{multline} \label{4.57}
-\left((\cL^0_2)^{-1}P^{N^\bot}\left\langle \Gamma_{i i}(\mR),
 e_l \right\rangle\nabla_{0,e_l} P^N \right)(0,0)\\
= \left\{(\cL^0_2)^{-1}P^{N^\bot} \right.
 \Big(\frac{2}{3} \left\langle R^{TX_G} (\mR^0,e^0_i)e^0_i,
\tfrac{\partial}{\partial \ov{z}^0_j}\right\rangle b_j\\
\left.+ \frac{1}{2}  \left\langle R^{TB} (\mR^\bot,e^0_i)e^0_i
+ A(e^0_i) A(e^0_i)\mR^\bot, e^\bot_k\right \rangle b^\bot_k\Big)
P^N \right\}(0,0)\\
= \left\{
-\frac{1}{3\pi}\left\langle R^{TX_G} (\tfrac{\partial}{\partial z^0_j},
e^0_i)e^0_i,\tfrac{\partial}{\partial \ov{z}^0_j}\right\rangle
- \frac{1}{16\pi} \left\langle R^{TB} (e^\bot_k, e^0_i)e^0_i,
e^\bot_k\right\rangle
+ \frac{1}{16\pi} |A(e^0_i)e^\bot_k|^2  \right\}P^N(0,0)\\
=\left\{
\left\langle - \frac{2}{3\pi}R^{TX_G} (\tfrac{\partial}{\partial z^0_j},
\tfrac{\partial}{\partial \ov{z}^0_i})\tfrac{\partial}{\partial z^0_i}
+ \frac{1}{4\pi} R^{TB} (e^\bot_k, \tfrac{\partial}{\partial z^0_j})e^\bot_k,
\tfrac{\partial}{\partial \ov{z}^0_j}\right\rangle  \right.
\left. + \frac{1}{8\pi} |T(\tfrac{\partial}{\partial z^0_i},
\tfrac{\partial}{\partial \ov{z}^0_j})|^2
 \right\}P^N(0,0).
\end{multline}

\comment{
Thus from \eqref{0g4}, \eqref{g3.57}, \eqref{4.48}, \eqref{4.55}-\eqref{4.57}
and comparing with \eqref{3.74}, we get
\begin{multline} \label{4.58}
-\left((\cL^0_2)^{-1}P^{N^\bot} (I_2+ \left \langle
\Gamma_{i i}(\mR), e_l \right \rangle\nabla_{0,e_l}) P^N \right)(0,0)\\
= \left\{
 \frac{2}{3\pi}\left\langle R^{TX_G} (\tfrac{\partial}{\partial z^0_j},
\tfrac{\partial}{\partial \ov{z}^0_i})\tfrac{\partial}{\partial z^0_i},
\tfrac{\partial}{\partial \ov{z}^0_j}\right\rangle
-\frac{5}{6\pi}\left\langle R^{TX_G} (\tfrac{\partial}{\partial z^0_j},
\tfrac{\partial}{\partial z^0_i})\tfrac{\partial}{\partial \ov{z}^0_i} ,
\tfrac{\partial}{\partial \ov{z}^0_j}\right\rangle\right.\\
 + \frac{5}{16\pi}|T(\tfrac{\partial}{\partial z^0_i},
\tfrac{\partial}{\partial \ov{z}^0_j})|^2
+\frac{1}{8\pi} \left\langle
R^{TB} (e^\bot_k, \tfrac{\partial}{\partial z^0_j})e^\bot_k,
\tfrac{\partial}{\partial \ov{z}^0_j}\right\rangle\\
\left.-\frac{1}{48\pi}\left\langle R^{TB} (e^\bot_k, e^\bot_j)e^\bot_k,
e^\bot_j\right\rangle \right\}P^N(0,0).
\end{multline}
}

By $\cL^0_2 P^N=0$, \eqref{b3.15}, \eqref{3.64}, \eqref{4.50},
\eqref{4.55} and since $R^{TX_G}(\cdot,\cdot)$ is a $(1,1)$-form
(comparing with \eqref{3.75}), we get
\begin{multline}\label{4.59}
- \Big\{(\cL^0_2)^{-1}P^{N^\bot}
\Big[\frac{1}{4} K_2(\mR)-\frac{3}{8}\Big(\sum_l \left \langle A(e^0_l)e^0_l,
\mR^\bot\right \rangle\Big)^2, \cL_2^0\Big]P^N\Big\}(0,0)\\
= \Big\{P^{N^\bot}\Big[\frac{1}{4} K_2(\mR)
-\frac{3}{8}\Big(\sum_l \left \langle A(e^0_l)e^0_l,
\mR^\bot\right \rangle\Big)^2 \Big]P^N\Big\}(0,0)\\
=  \frac{1}{4}\left\{P^{N^\bot}  \left[
 \left\langle\frac{1}{3} R^{TX_G}(\mR^0, e^0_i)\mR^0
+  R^{TB}(\mR^\bot, e^0_i)\mR^\bot, e^0_i\right\rangle \right.\right.\\
+\frac{1}{3} \left \langle R^{TB} (\mR^\bot, e^\bot_i)\mR^\bot,
e^\bot_i\right \rangle
\left.\left. + \frac{1}{2}\Big(\sum_i \left \langle A(e^0_i)e^0_i,
 \mR^\bot\right \rangle\Big)^2
- | A(e^0_i)\mR^\bot|^2\right] P^N\right\}(0,0)\\
=  \Big(\frac{1}{6\pi}  \left \langle
R^{TX_G} (\tfrac{\partial}{\partial z^0_j},e^0_i)
e^0_i, \tfrac{\partial}{\partial \ov{z}^0_j}\right \rangle
-\frac{1}{16\pi}\left\langle R^{TB}(e^\bot_k, e^0_i)e^\bot_k,
e^0_i\right\rangle\\
-\frac{1}{48\pi}\left \langle R^{TB} (e^\bot_k, e^\bot_i)e^\bot_k,
e^\bot_i\right \rangle
-\frac{1}{32\pi}\Big|\sum_i A(e^0_i)e^0_i\Big|^2
+ \frac{1}{16\pi}| A(e^0_i)e^\bot_k|^2\Big)P^N(0,0),\\
=  \Big(\frac{1}{3\pi}  \left \langle
R^{TX_G} (\tfrac{\partial}{\partial z^0_j},
\tfrac{\partial}{\partial \ov{z}^0_i})\tfrac{\partial}{\partial z^0_i},
\tfrac{\partial}{\partial \ov{z}^0_j}\right \rangle
-\frac{1}{4\pi}\left\langle R^{TB}(e^\bot_k,
\tfrac{\partial}{\partial z^0_j})e^\bot_k,
\tfrac{\partial}{\partial \ov{z}^0_j}\right\rangle\\
-\frac{1}{8\pi}\Big|\sum_i T(\tfrac{\partial}{\partial z^0_i},
\tfrac{\partial}{\partial \ov{z}^0_i})\Big|^2
+ \frac{1}{8\pi}| T(\tfrac{\partial}{\partial z^0_i},
\tfrac{\partial}{\partial \ov{z}^0_j})|^2
-\frac{1}{48\pi}\left \langle R^{TB} (e^\bot_k, e^\bot_j)e^\bot_k,
e^\bot_j\right \rangle\Big)P^N(0,0).
\end{multline}
By \eqref{0g6}, \eqref{g3.57}, \eqref{3.533},
\eqref{3.64}, \eqref{4.48} and \eqref{4.55},
\begin{align} \label{4.60}
&-\Big\{(\cL^0_2)^{-1}P^{N^\bot}\Big(- \frac{1}{2}\left\langle A(e^0_l)e^0_l,
\mR^\bot\right \rangle \nabla_{A(e^0_k)e^0_k}
+2 \left \langle A(e^0_i)e^0_j, \mR^\bot\right \rangle \nabla_{A(e^0_i)e^0_j}\\
&\hspace*{30mm}+ \frac{2}{3} \left \langle R^{TB}(\mR^\bot, e^\bot_i)e^\bot_i,
 e_j\right \rangle  \nabla_{0,e_j}\Big)P^N\Big\}(0,0)\nonumber\\
&=- \frac{1}{16\pi}\Big( - \frac{1}{2} \Big|\sum_l A(e^0_l)e^0_l\Big|^2
+2 | A(e^0_i)e^0_j|^2
+ \frac{2}{3} \left \langle R^{TB}(e^\bot_j, e^\bot_i)e^\bot_i,
 e^\bot_j\right \rangle \Big)P^N(0,0),\nonumber\\
&=\Big( \frac{1}{8\pi}\Big|\sum_i T(\tfrac{\partial}{\partial z^0_i},
\tfrac{\partial}{\partial \ov{z}^0_i})\Big|^2
- \frac{1}{4\pi} \Big|T(\tfrac{\partial}{\partial z^0_i},
\tfrac{\partial}{\partial \ov{z}^0_j})\Big|^2
+ \frac{1}{24\pi} \left \langle R^{TB}(e^\bot_k, e^\bot_j)e^\bot_k,
 e^\bot_j\right \rangle \Big)P^N(0,0),\nonumber\\
&-\Big\{(\cL^0_2)^{-1}P^{N^\bot}(-R^{E_B}(\mR, e_i))
\nabla_{0,e_i}P^N\Big\}(0,0)
= \frac{1}{2\pi}R^{E_B}(\tfrac{\partial}{\partial z^0_j},
\tfrac{\partial}{\partial \ov{z}^0_j})P^N(0,0).\nonumber
\end{align}
For $F_{ij;kl}\in \bC$,
 from Theorem \ref{t3.4}, \eqref{b3.11}, \eqref{4.19}, \eqref{a4.50}
and comparing with \eqref{4.32}, we get
\begin{multline} \label{a4.60}
\left\{(\cL^0_2)^{-1}P^{N^\bot}F_{ij;kl} Z^\bot_i Z^\bot_j Z^\bot_k Z^\bot_l
 P^N\right\}(0,0)\\
= \Big\{(\cL^0_2)^{-1}P^{N^\bot}
 \Big[\sum_{j\neq k} (F_{jj;kk}+ F_{kj;kj}+ F_{kj;jk}) (Z^\bot_j)^2(Z^\bot_k)^2
+ F_{kk;kk} (Z^\bot_k)^4\Big]P^N\Big\}(0,0)\\
= \frac{1}{2^8 \pi^4}\left\{P^{N^\bot}
 \Big[\sum_{j\neq k} (F_{jj;kk}+ F_{kj;kj}+ F_{kj;jk})
 \Big(\frac{(b^\bot_j)^2(b^\bot_k)^2}{16\pi}+\frac{1}{2} ((b^\bot_j)^2+(b^\bot_k)^2) \Big)\right.\\
\left. + F_{kk;kk} \Big(\frac{(b^\bot_k)^4}{16\pi}+3(b^\bot_k)^2\Big) \Big]P^N\right\}(0,0)\\
=  \frac{-3}{2^8 \pi^3}(F_{jj;kk}+ F_{kj;kj}+ F_{kj;jk}) P^N(0,0).
\end{multline}
By \eqref{3.32}, \eqref{a3.52},
 \begin{align} \label{4.61}
\frac{1}{9} \sum_i
\Big[(\partial_\mR R^{L_B})_{x_0} (\mR,e_i)\Big]^2
=& -\pi^2 \sum_i\left\langle JT(\mR^\bot, e^0_i), \mR^\bot\right\rangle ^2 \\
&-\pi^2\sum_j \left\langle JT(\mR, e^\bot_j), \mR^\bot\right\rangle ^2 .
\nonumber
\end{align}
By \eqref{0g4}, \eqref{34.21}, \eqref{a4.60}
and $\mT_{kl}(e^0_i)$ is symmetric on $k,l$, we get
\begin{multline} \label{4.62}
-\pi^2 \sum_i\left((\cL^0_2)^{-1} P^{N^\bot}
\left\langle JT(\mR^\bot, e^0_i), \mR^\bot\right\rangle ^2 P^N\right)(0,0) \\
=-\pi^2 \left((\cL^0_2)^{-1} P^{N^\bot} \mT_{jj'}(e^0_i)\mT_{kl}(e^0_i)
Z^\bot_j Z^\bot_{j'}Z^\bot_k Z^\bot_l P^N\right)(0,0) \\
= \frac{3}{2^8 \pi} \Big(2 \mT_{jk}(e^0_i)^2
+ \mT_{jj}(e^0_i)\mT_{kk}(e^0_i)\Big)  P^N(0,0)\\
= \frac{3}{2^6 \pi}\Big(2 |T(e^\bot_k,
\tfrac{\partial}{\partial \ov{z}^0_j})|^2
+ \Big|\sum_j \mT_{jj}(\tfrac{\partial}{\partial \ov{z}^0_i})
\Big|^2\Big)P^N(0,0).
\end{multline}
In the same way, by \eqref{3.11e}, \eqref{34.21}, \eqref{a4.60}, we get
 \begin{multline} \label{a4.62}
-\pi^2 \sum_j \left((\cL^0_2)^{-1} P^{N^\bot}
\left\langle JT(\mR^\bot, e^\bot_j), \mR^\bot\right\rangle ^2
P^N\right)(0,0)\\
=\frac{3}{2^8 \pi}\wi{\mT}_{ijk} (\wi{\mT}_{ijk}+ \wi{\mT}_{kji})P^N(0,0).
\end{multline}

By \eqref{34.21}, \eqref{4.50},
\begin{multline} \label{4.64}
- \pi^2 \sum_j \left((\cL^0_2)^{-1} P^{N^\bot}
\left\langle JT(\mR^0, e^\bot_j), \mR^\bot\right\rangle ^2 P^N\right)(0,0)\\
=\frac{7}{48\pi}|\mT_{jk}(\tfrac{\partial}{\partial \ov{z}^0_i})|^2 P^N(0,0)
=\frac{7}{48\pi}|T(e^\bot_k,\tfrac{\partial}{\partial \ov{z}^0_j})|^2 P^N(0,0).
\end{multline}

By \eqref{3.32} and \eqref{a3.75}, the total degree of $Z^0$,
$\nabla_{0,e^0_i}$ in the fourth term of $\mO^\prime_2$ in
\eqref{a3.20} is $1$, thus the contribution of the fourth term of
$\mO^\prime_2$ in \eqref{a3.20} for
 $-((\cL^0_2)^{-1}P^{N^\bot} \mO^\prime_2 P^N)(0,0)$ is zero.
By \eqref{a3.20}, \eqref{a4.54}, \eqref{4.56}-\eqref{4.60}
and \eqref{4.61}-\eqref{4.64}, comparing with \eqref{3.77}, we get
\begin{multline} \label{4.66}
-\left((\cL^0_2)^{-1} P^{N^\bot}  \mO^\prime_2 P^N\right)(0,0)
=\left\{\frac{1}{2\pi} \left\langle R^{TX_G}(\tfrac{\partial}{\partial z^0_j},
\tfrac{\partial}{\partial \ov{z}^0_i})
\tfrac{\partial}{\partial z^0_i},
\tfrac{\partial}{\partial \ov{z}^0_j}\right\rangle\right.\\
+\frac{7}{6}\left[
\frac{5\sqrt{-1}}{2^5 \pi}\left\langle  Je^\bot_k,
 \nabla ^{TY}_{e^\bot_k} (T(\tfrac{\partial}{\partial z^0_j},
\tfrac{\partial}{\partial \ov{z}^0_j}))
+ \nabla ^{TY}_{\tfrac{\partial}{\partial z^0_j}} (T(e^\bot_k,
\tfrac{\partial}{\partial \ov{z}^0_j}))\right\rangle
+\frac{3}{32\pi} |T(\tfrac{\partial}{\partial z^0_i},
\tfrac{\partial}{\partial \ov{z}^0_j})|^2 \right.\\
+ \frac{3}{16\pi}  \left\langle  R^{TB}(e^\bot_k,
\tfrac{\partial}{\partial z^0_j}) e^\bot_k,
\tfrac{\partial}{\partial \ov{z}^0_j}\right\rangle
+\frac{7}{2^6 \pi}|T(e^\bot_k, \tfrac{\partial}{\partial \ov{z}^0_j})|^2
 \left. -\frac{\sqrt{-1}}{16\pi}\left\langle T(e^\bot_k, J e^\bot_k),
 T(\tfrac{\partial}{\partial z^0_j},
\tfrac{\partial}{\partial \ov{z}^0_j}) \right\rangle
 \right]\\
-\frac{1}{8\pi} \left\langle  R^{TB}(e^\bot_k,
\tfrac{\partial}{\partial z^0_j}) e^\bot_k,
\tfrac{\partial}{\partial \ov{z}^0_j}\right\rangle
+\frac{3}{16\pi} |T(\tfrac{\partial}{\partial z^0_i},
\tfrac{\partial}{\partial \ov{z}^0_j})|^2
+ \frac{3}{32\pi} |T(e^\bot_k, \tfrac{\partial}{\partial \ov{z}^0_j})|^2 \\
+\frac{3}{64\pi}\Big|\sum_j \mT_{jj}(\tfrac{\partial}{\partial \ov{z}^0_i})
\Big|^2
+\frac{3}{2^8 \pi}\wi{\mT}_{ijk} (\wi{\mT}_{ijk}+ \wi{\mT}_{kji})
\left. +\frac{1}{2\pi} R^{E_G}(\tfrac{\partial}{\partial z^0_j},
\tfrac{\partial}{\partial \ov{z}^0_j})
\right\}P^N(0,0).
\end{multline}

By \eqref{a3.48} and \eqref{a4.45},
\begin{multline} \label{4.67}
-4\pi ^2 \left((\cL^0_2)^{-1} P^{N^\bot} \mO''_2 P^N\right)(0,0)
= -4\pi ^2 \left\{(\cL^0_2)^{-1} P^{N^\bot} \right.\\
 \Big[ -\frac{1}{3}\left\langle (\nabla^{TY}_\cdot
\dot{g}^{TY}_\cdot)_{(\mR^0,\mR^0)}
J\mR^\bot + (\nabla^{TY}_\cdot
\dot{g}^{TY}_\cdot)_{(\mR^\bot,\mR^\bot)}
J\mR^\bot, J\mR^\bot\right\rangle \\
+\frac{1}{6}\left\langle \nabla ^{TY}_{\mR^0} (T(e^\bot_j, J_{x_0}e^0_i))
Z^\bot_j Z^0_i
 +\nabla ^{TY}_{\mR^\bot} (T(e^0_j, J_{x_0}e^0_i))Z^0_j Z^0_i,
 J\mR^\bot\right\rangle \\
+\frac{1}{3} \left\langle R^{TB}(\mR^\bot,\mR^0)\mR^0,\mR^\bot\right\rangle
- \frac{1}{12} \sum_l \left\langle T(\mR^0, e_l),J\mR^\bot\right\rangle^2\\
- \frac{1}{12} \sum_l \left\langle T(\mR^\bot, e_l),J\mR^\bot\right\rangle^2
\left.+ \frac{7}{12} |T(\mR^\bot, J\mR^\bot)|^2
 \Big] P^N  \right\}(0,0).
\end{multline}

Now $\{e_l\}= \{e^0_i\}\cup \{e^\bot_k\}$.

By \eqref{a3.70}, \eqref{3.79}, \eqref{4.7}, \eqref{4.50}, \eqref{a4.60},
\eqref{4.62}, \eqref{a4.62}, \eqref{4.67}
and comparing with \eqref{3.81},
\begin{multline}\label{4.70}
-4\pi ^2 \left((\cL^0_2)^{-1}P^{N^\bot} \mO''_2 P^N\right)(0,0)
=\left\{\frac{7}{24\pi}\left[
-\frac{8}{3}\nabla_{\tfrac{\partial}{\partial z^0_j}}
\nabla_{\tfrac{\partial}{\partial \ov{z}^0_j}} \log h\right.\right.\\
+\frac{\sqrt{-1}}{3}\Big\langle
-\nabla ^{TY}_{\tfrac{\partial}{\partial z^0_j}}
(T(e^\bot_k, \tfrac{\partial}{\partial \ov{z}^0_j}))
-\nabla ^{TY}_{e^\bot_k} (T(\tfrac{\partial}{\partial z^0_j},
\tfrac{\partial}{\partial \ov{z}^0_j})),Je^\bot_k\Big\rangle \\
\left.  -\frac{1}{3} \Big|T(\tfrac{\partial}{\partial z^0_i},
\tfrac{\partial}{\partial \ov{z}^0_j})\Big|^2
-\frac{1}{6} \Big|T(e ^\bot_k,
\tfrac{\partial}{\partial \ov{z}^0_j})\Big|^2
 - \frac{2}{3}\left\langle R^{TB}(e^\bot_k,
\tfrac{\partial}{\partial z^0_j})e^\bot_k,
\tfrac{\partial}{\partial \ov{z}^0_j}\right\rangle \right ]\\
-\frac{1}{2^6 \pi} \left\langle (\nabla^{TY}_\cdot
\dot{g}^{TY}_\cdot)_{(e^\bot_j,e^\bot_j)}
Je^\bot_k, Je^\bot_k\right\rangle
- \frac{1}{2^5 \pi}\left\langle (\nabla^{TY}_\cdot
\dot{g}^{TY}_\cdot)_{(e^\bot_j,e^\bot_k)}
Je^\bot_j, Je^\bot_k\right\rangle \\
 - \frac{1}{2^8 \pi}
  \Big(8 |T(e^\bot_k,\tfrac{\partial}{\partial \ov{z}^0_j})|^2
+ 4 \Big|\sum_j \mT_{jj}(\tfrac{\partial}{\partial \ov{z}^0_i})\Big|^2
+ \wi{\mT}_{ijk} (\wi{\mT}_{ijk} + \wi{\mT}_{kji})\Big)\\
\left.+ \frac{7}{2^8 \pi} \Big(2\mT_{jkm}^2 + \mT_{jjm}\mT_{kkm} \Big)
\right \} P^N(0,0).
\end{multline}

By \eqref{a3.51}, \eqref{b3.52}, \eqref{4.48}, \eqref{4.66} and \eqref{4.70},
comparing with \eqref{a3.64},
we have
\begin{multline}\label{4.71}
\Psi_{1,2}(0)= - \left((\cL^0_2)^{-1}P^{N^\bot} (\mO^\prime_2+4\pi ^2\mO''_2)
 P^N\right)(0,0) \\
- \frac{\sqrt{-1}}{16\pi} \Big(\left\langle T(e^\bot_j, Je^\bot_j),
\wi{\mu}^E\right\rangle
-2\left\langle Je^\bot_j,
\nabla^{TY}_{e^\bot_j}\wi{\mu}^E\right\rangle\Big) P^N (0,0) \\
=\left\{\frac{1}{2\pi} \left\langle R^{TX_G}(\tfrac{\partial}{\partial z^0_j},
\tfrac{\partial}{\partial \ov{z}^0_i})
\tfrac{\partial}{\partial z^0_i},
\tfrac{\partial}{\partial \ov{z}^0_j}\right\rangle\right.
+\frac{1}{2\pi} R^{E_G}(\tfrac{\partial}{\partial z^0_j},
\tfrac{\partial}{\partial \ov{z}^0_j})\\
+\frac{7}{6}\left[
\frac{1}{6\pi}\Delta_{X_G} \log h
+  \frac{1}{48\pi} \left\langle
R^{TB}(e^\bot_k, \tfrac{\partial}{\partial z^0_j})e^\bot_k,
\tfrac{\partial}{\partial \ov{z}^0_j}\right\rangle  \right.
 + \frac{1}{96\pi}\Big| T(\tfrac{\partial}{\partial z^0_i},
\tfrac{\partial}{\partial \ov{z}^0_j})\Big|^2\\
- \frac{\sqrt{-1}}{16\pi} \left\langle T(e^\bot_k,Je^\bot_k),
T(\tfrac{\partial}{\partial z^0_j},
\tfrac{\partial}{\partial \ov{z}^0_j}) \right\rangle
+ \frac{13}{192\pi}
\Big|T(e^\bot_k,\tfrac{\partial}{\partial \ov{z}^0_j})\Big|^2\\
\left. +\frac{7\sqrt{-1}}{96\pi} \left\langle
\nabla ^{TY}_{\tfrac{\partial}{\partial z^0_j}}
(T(e^\bot_k,\tfrac{\partial}{\partial \ov{z}^0_j}))
+ \nabla ^{TY}_{e^\bot_k} (T(\tfrac{\partial}{\partial z^0_j},
\tfrac{\partial}{\partial \ov{z}^0_j})), Je^\bot_k \right\rangle
\right]\\
- \frac{1}{8\pi}  \left\langle  R^{TB}(e^\bot_k,
\tfrac{\partial}{\partial z^0_j}) e^\bot_k,
\tfrac{\partial}{\partial \ov{z}^0_j}\right\rangle
+\frac{3}{16 \pi} |T(\tfrac{\partial}{\partial z^0_i},
\tfrac{\partial}{\partial \ov{z}^0_j})|^2
+\frac{1}{16\pi}
|T(e^\bot_k, \tfrac{\partial}{\partial \ov{z}^0_j})|^2 \\
+\frac{1}{32\pi} \Big|\sum_j
\mT_{jj}(\tfrac{\partial}{\partial \ov{z}^0_i})\Big|^2
+\frac{1}{2^7 \pi}\wi{\mT}_{ijk} (\wi{\mT}_{ijk} + \wi{\mT}_{kji})
+ \frac{7}{2^8 \pi} \Big(2\mT_{jkm}^2 + \mT_{jjm}\mT_{kkm} \Big)\\
-\frac{1}{2^6 \pi} \left\langle (\nabla^{TY}_\cdot
\dot{g}^{TY}_\cdot)_{(e^\bot_j,e^\bot_j)}
Je^\bot_k, Je^\bot_k\right\rangle
-\frac{1}{2^5 \pi}\left\langle (\nabla^{TY}_\cdot
\dot{g}^{TY}_\cdot)_{(e^\bot_j,e^\bot_k)}
Je^\bot_j, Je^\bot_k\right\rangle\\
\left.- \frac{\sqrt{-1}}{16\pi} \Big(\left\langle T(e^\bot_j, Je^\bot_j),
\wi{\mu}^E\right\rangle
-2\left\langle Je^\bot_j,
\nabla^{TY}_{e^\bot_j}\wi{\mu}^E\right\rangle\Big)\right \}  P^N (0,0) .
\end{multline}

By \eqref{4.7}, \eqref{4.16}, the term $\frac{7}{6}[\cdots]$ in \eqref{4.71}
is $\frac{7}{6} \left(\frac{3}{8\pi} \Delta_{X_G}\log h
+\frac{1}{8\pi}
|T(e^\bot_k,\tfrac{\partial}{\partial \ov{z}^0_j})|^2\right)$.

By  \eqref{4.12} and \eqref{4.71}, we get \eqref{4.46}.

 The proof of Lemma \ref{t4.13} is complete.
\end{proof}

\begin{lemma}\label{t4.14} The following identity holds,
\begin{align}\label{f4.1}
\left\langle (\nabla^{TY}_{e^\bot_k} \dot{g}^{TY}_{e^\bot_k})
Je^\bot_l, Je^\bot_l\right\rangle=
& 4 \nabla_{e^\bot_k}\nabla_{e^\bot_k}\log h, \nonumber\\
 \left\langle (\nabla^{TY}_{e^\bot_k} \dot{g}^{TY}_{e^\bot_l})
Je^\bot_l, Je^\bot_k\right\rangle=
& 4 \nabla_{e^\bot_k}\nabla_{e^\bot_k}\log h
+ 2 \Big|\sum_l \mT_{ll}(\tfrac{\partial}{\partial \ov{z}^0_j})\Big|^2\\
&-2 |T(e^\bot_k, \tfrac{\partial}{\partial \ov{z}^0_j}) |^2
+\frac{1}{2} (\wi{\mT}_{jki} + \wi{\mT}_{ijk})\wi{\mT}_{ijk}. \nonumber
\end{align}
\end{lemma}
\begin{proof} By using the same argument as in \eqref{3.79},
we get the first equation of \eqref{f4.1}.

Recall that $P^{T^HX}, P^{TY}$ are the projections from $TX=T^HX\oplus TY$
 onto $T^HX$,$TY$.

 By \eqref{h0}, \eqref{h4}, \eqref{g1},
\eqref{g3.29} and \eqref{g3.31} (cf. also \eqref{b3.16}),
\begin{subequations}
 \begin{align}\label{f4.2a}
& (P^{T^HX}J e^{\bot,H}_l)|_{\mu^{-1}(0)}=0,
\quad  (J e^{\bot,H}_l)_{x_0} \in TY,  \\
& (\nabla ^{TX}_{e^{\bot,H}_k} e^{\bot,H}_l)_{x_0}
= -\frac{1}{2} T(e^{\bot}_k, e^{\bot}_l),\label{f4.2}\\
& (\nabla ^{TX}_{e^{0}_j} e^{\bot,H}_l)_{x_0}
= (A(e^0_j)e^\bot_l)^H -\frac{1}{2} T(e^{0}_j, e^{\bot}_l),\nonumber\\
&(\nabla ^{TX}_{J e^{\bot,H}_l} e^H_k)_{x_0}
=\frac{1}{2} \left\langle T(e_k, e_j),
Je^{\bot}_l\right\rangle  e^H_j
+ \left\langle T(e_k, Je^{\bot}_l),
Je^{\bot}_j\right\rangle Je^\bot_j.\label{f4.2b}
\end{align}
\end{subequations}
From \eqref{f4.2a}, we get
\begin{align}\label{0f4.1}
 \nabla ^{TX}_{e^0_i} P^{T^HX}J e^{\bot,H}_l
& = \nabla ^{TX}_{J e^{\bot,H}_k} P^{T^HX}J e^{\bot,H}_l=0.
\end{align}

By \eqref{g3.29}, \eqref{34.21}, \eqref{aa3.50} and \eqref{f4.2},
we get at $x_0$,
\begin{align}\label{f4.3}
  \nabla ^{TX}_{e^{\bot,H}_k} P^{T^HX}J e^{\bot,H}_l
& =  \nabla ^{T^HX}_{e^{\bot,H}_k} P^{T^HX}J e^{\bot,H}_l\\
& = -\frac{1}{2} JT(e^\bot_k, e^\bot_l)
+ \frac{1}{2} \left\langle JT(e^{\bot}_k, e_j),
e^{\bot}_l \right\rangle  e_j\nonumber\\
&= -\frac{1}{2}(\wi{\mT}_{klj}- \wi{\mT}_{kjl})e^\bot_j
+ \frac{1}{2} \left\langle JT(e^{\bot}_k, e^0_j),
e^{\bot}_l \right\rangle  e^0_j . \nonumber
 \end{align}
By \eqref{3.12a}, \eqref{f4.2} and \eqref{f4.2b}, at $x_0$,
\begin{multline}\label{0f4.3}
[e^{\bot,H}_k, J e^{\bot,H}_l]= \nabla ^{TX}_{e^{\bot,H}_k}J e^{\bot,H}_l
-  \nabla ^{TX}_{J e^{\bot,H}_l}e^{\bot,H}_k
=J \nabla ^{TX}_{e^{\bot,H}_k} e^{\bot,H}_l
-  \nabla ^{TX}_{J e^{\bot,H}_l}e^{\bot,H}_k\\
= -\frac{1}{2} JT(e^\bot_k, e^\bot_l)
+\frac{1}{2} \left\langle JT(e^{\bot}_k, e_j),
e^{\bot}_l \right\rangle  e_j
- \left\langle T(e^{\bot}_k, J e^{\bot}_j),
Je^{\bot}_l \right\rangle J e^{\bot,H}_j.
\end{multline}
By  \eqref{34.21}, \eqref{0f4.1},
\eqref{f4.3} and \eqref{0f4.3}, we get at $x_0$,
\begin{multline}\label{0f4.4}
\left\langle \nabla ^{TX}_{[e^{\bot,H}_k, J e^{\bot,H}_l]}
P^{T^HX}J e^{\bot,H}_l, e^{\bot,H}_k \right\rangle\\
= -\frac{1}{2} \left\langle \nabla ^{TX}_{JT(e^\bot_k, e^\bot_l)
- \left\langle JT(e^{\bot}_k, e^{\bot}_j),
e^{\bot}_l \right\rangle e^{\bot}_j }
P^{T^HX}J e^{\bot,H}_l, e^{\bot,H}_k \right\rangle \\
=\frac{1}{4}(\wi{\mT}_{klj} - \wi{\mT}_{kjl})
(\wi{\mT}_{jlk} - \wi{\mT}_{jkl}) .
\end{multline}
By  \eqref{3.12a}, \eqref{f4.2}, \eqref{0f4.1} and \eqref{f4.3}, at $x_0$,
\begin{align}\label{f4.4}
 &\nabla ^{TX}_{e^{0,H}_j}P^{TY}J e^{\bot,H}_l
=  J \nabla ^{TX}_{e^{0,H}_j}  e^{\bot,H}_l
= JA(e^{0}_j)e^{\bot}_l -\frac{1}{2} JT(e^{0}_j, e^{\bot}_l) , \nonumber\\
& \nabla ^{TX}_{e^{\bot,H}_k}P^{TY}J e^{\bot,H}_l
= J \nabla ^{TX}_{e^{\bot,H}_k}e^{\bot,H}_l
- \nabla ^{TX}_{e^{\bot,H}_k}P^{T^HX}J e^{\bot,H}_l\\
&\hspace*{20mm}=- \frac{1}{2} \left\langle J T(e^{\bot}_k, e_j),
e^{\bot}_l \right\rangle  e_j
=- \frac{1}{2}  \wi{\mT}_{kjl} e^{\bot}_j
-\frac{1}{2} \left\langle J T(e^{\bot}_k, e^0_j),
e^{\bot}_l \right\rangle  e^0_j. \nonumber
 \end{align}
Thus by  \eqref{f4.4}, at $x_0$,
\begin{align}\label{f4.5}
 \nabla ^{TY}_{e^{\bot,H}_k}P^{TY}J e^{\bot,H}_l
= P^{TY}\nabla ^{TX}_{e^{\bot,H}_k}P^{TY}J e^{\bot,H}_l=0.
 \end{align}

By \eqref{h0}, \eqref{h3}, \eqref{h4} and \eqref{f4.5}, at $x_0$,
\begin{multline}\label{f4.6}
 \left\langle (\nabla^{TY}_{e^\bot_k}
\dot{g}^{TY}_{e^\bot_l})
Je^\bot_l, Je^\bot_k\right\rangle
= e^\bot_k  \left\langle \dot{g}^{TY}_{e^\bot_l} P^{TY} Je^\bot_l,
 P^{TY} Je^\bot_k\right\rangle
= 2 e^\bot_k  \left\langle \nabla ^{TX}_{P^{TY} Je^\bot_l}e^\bot_l,
 P^{TY} Je^\bot_k\right\rangle\\
 = 2 e^\bot_k  \left\langle \nabla ^{TX}_{ Je^\bot_l}e^\bot_l,
 Je^\bot_k\right\rangle
-2  e^\bot_k \left\langle \nabla ^{TX}_{P^{T^HX} Je^\bot_l}e^\bot_l,
  Je^\bot_k\right\rangle
-2  e^\bot_k \left\langle \nabla ^{TX}_{P^{TY} Je^\bot_l}e^\bot_l,
 P^{T^HX} Je^\bot_k\right\rangle.
 \end{multline}

By \eqref{3.11d}, \eqref{34.21},
\eqref{f4.2a}, \eqref{f4.2} and \eqref{f4.3}, at $x_0$
\begin{multline}\label{f4.7}
-2  e^\bot_k \left\langle \nabla ^{TX}_{P^{T^HX} Je^\bot_l}e^\bot_l,
  Je^\bot_k\right\rangle
= -2\Big\langle \nabla ^{TX}_{\nabla ^{TX}_{ e^\bot_k}P^{T^HX} Je^\bot_l}
e^\bot_l, Je^\bot_k\Big\rangle\\
= \left\langle T\Big(-\frac{1}{2}(\wi{\mT}_{klj}- \wi{\mT}_{kjl})e^\bot_j
+ \frac{1}{2}\left\langle JT(e^{\bot}_k, e^0_j),
e^{\bot}_l \right\rangle  e^0_j, e^\bot_l\Big), Je^\bot_k\right\rangle\\
=\frac{1}{2}(\wi{\mT}_{klj} - \wi{\mT}_{kjl})\wi{\mT}_{jlk}
+\frac{1}{2} \left\langle T (e^0_j,e^\bot_l), Je^\bot_k\right\rangle
\left\langle JT (e^\bot_k, e^0_j), e^\bot_l\right\rangle\\
= \frac{1}{2}(\wi{\mT}_{klj} - \wi{\mT}_{kjl})\wi{\mT}_{jlk}
+ \frac{1}{2} |T (e^\bot_k, e^0_j)|^2.
 \end{multline}

By \eqref{3.11e}, \eqref{34.21},  \eqref{f4.2a},
 \eqref{f4.2} and \eqref{f4.3},
 at $x_0$, we have
\begin{multline}\label{f4.8}
-2  e^\bot_k \left\langle \nabla ^{TX}_{P^{TY} Je^\bot_l}e^\bot_l,
 P^{T^HX} Je^\bot_k\right\rangle
= -2 \left\langle \nabla ^{TX}_{P^{TY} Je^\bot_l}e^\bot_l,
\nabla ^{TX}_{e^\bot_k} P^{T^HX} Je^\bot_k\right\rangle\\
= - \frac{1}{2}
 \left\langle JT(e^\bot_k, e_j), e^\bot_k\right\rangle
 \left\langle  T(e^\bot_l, e_j), Je^\bot_l\right\rangle
= \frac{1}{2} \mT_{ll}(e^0_j)\mT_{kk}(e^0_j).
\end{multline}

Now by \eqref{3.12a},
\begin{multline}\label{f4.9}
 e^\bot_k  \left\langle \nabla ^{TX}_{ Je^\bot_l}e^\bot_l,
 Je^\bot_k\right\rangle
= - e^\bot_k  \left\langle \nabla ^{TX}_{ Je^\bot_l}Je^\bot_l,
 e^\bot_k\right\rangle\\
=-  e^\bot_k \left\langle \nabla ^{TX}_{ P^{TY}Je^\bot_l}P^{TY}J e^\bot_l
+ \nabla ^{TX}_{ P^{T^HX}Je^\bot_l}P^{TY}J e^\bot_l
+ \nabla ^{TX}_{ Je^\bot_l} P^{T^HX}Je^\bot_l,e^\bot_k\right\rangle .
 \end{multline}
By Theorem \ref{t4.2}, \eqref{h4}, \eqref{f4.2a} and \eqref{f4.3},
at $x_0$,
\begin{multline}\label{f4.10}
-2 e^\bot_k   \left\langle\nabla ^{TX}_{ P^{T^HX}Je^\bot_l} P^{TY}Je^\bot_l,
e^\bot_k\right\rangle
=-2   \Big\langle\nabla ^{TX}_{\nabla ^{TX}_{e^\bot_k}P^{T^HX}Je^\bot_l}
 P^{TY}Je^\bot_l,e^\bot_k\Big\rangle\\
= -\left\langle T\Big(-\frac{1}{2}(\wi{\mT}_{klj}- \wi{\mT}_{kjl})e^\bot_j
+ \frac{1}{2} \left\langle JT(e^{\bot}_k, e^0_j),
e^{\bot}_l \right\rangle  e^0_j, e^{\bot}_k\Big), J e^{\bot}_l\right\rangle\\
= -\frac{1}{2}(\wi{\mT}_{klj} - \wi{\mT}_{kjl})\wi{\mT}_{jkl}
- \frac{1}{2} |T (e^\bot_k, e^0_j)|^2.
\end{multline}
And by \eqref{3.11d}, \eqref{f4.2a}-\eqref{f4.2b}, \eqref{0f4.1},
 \eqref{f4.3}, \eqref{0f4.4} and  \eqref{f4.10}, at $x_0$,
\begin{multline}\label{f4.11}
-2 e^\bot_k \left\langle\nabla ^{TX}_{Je^\bot_l} P^{T^HX}Je^\bot_l,
e^\bot_k\right\rangle
=-2 \left\langle  \Big(\nabla ^{TX}_{Je^\bot_l}
\nabla ^{TX}_{e^\bot_k}
+\nabla ^{TX}_{[e^{\bot,H}_k, J e^{\bot,H}_l]} \Big)
P^{T^HX}J e^{\bot,H}_l , e^\bot_k\right\rangle   \\
= -\left\langle T\Big(-\frac{1}{2}(\wi{\mT}_{klj}- \wi{\mT}_{kjl})e^\bot_j
+ \frac{1}{2} \left\langle JT(e^{\bot}_k, e_j),
e^{\bot}_l \right\rangle  e_j, e^{\bot}_k\Big), J e^{\bot}_l\right\rangle
- \frac{1}{2}(\wi{\mT}_{klj} - \wi{\mT}_{kjl})
(\wi{\mT}_{jlk} - \wi{\mT}_{jkl}) \\
= -\frac{1}{2}(\wi{\mT}_{klj} - \wi{\mT}_{kjl}) \wi{\mT}_{jlk}
- \frac{1}{2} |T (e^\bot_k, e^0_j)|^2.
\end{multline}

Finally, by \eqref{h1}, \eqref{h4}, \eqref{h14} and \eqref{f4.5},
as in \eqref{3.79},
\begin{multline}\label{f4.12}
-2 e^\bot_k \left\langle \nabla ^{TX}_{ P^{TY}Je^\bot_l}P^{TY}J e^\bot_l,
e^\bot_k\right\rangle
=2 e^\bot_k \left\langle T(e^\bot_k, P^{TY}Je^\bot_l),
P^{TY}Je^\bot_l\right\rangle\\
=  \left\langle (\nabla ^{TY}_{e^\bot_k} \dot{g}^{TY}_{e^\bot_k})
Je^\bot_l, Je^\bot_l\right\rangle
=4 \nabla_{e^\bot_k}\nabla_{e^\bot_k}\log h.
\end{multline}

Thus by \eqref{f4.6}-\eqref{f4.12},
\begin{multline}\label{f4.13}
 \left\langle (\nabla^{TY}_{e^\bot_k}
\dot{g}^{TY}_{e^\bot_l})
Je^\bot_l, Je^\bot_k\right\rangle
= 4 \nabla_{e^\bot_k}\nabla_{e^\bot_k}\log h
- \frac{1}{2}|T(e^\bot_k, e^0_j) |^2\\
+ \frac{1}{2} \mT_{ll}(e^0_j)\mT_{kk}(e^0_j)
-\frac{1}{2}(\wi{\mT}_{klj} - \wi{\mT}_{kjl}) \wi{\mT}_{jkl}.
\end{multline}
From \eqref{0g4} and \eqref{f4.13}, we get \eqref{f4.1}.
\end{proof}

\begin{proof}[Proof of Theorem \ref{t4.15}] By \eqref{34.21}, \eqref{a3.55},
\begin{align}\label{4.75}
\begin{split}
&\sum_k\mF_1 (e^\bot_k)^2 = - \left\langle \wi{\mu}^E_{x_0},
\wi{\mu}^E_{x_0}\right\rangle_{g^{TY}}
- \left\langle \wi{\mu}^E, \frac{3}{2} \sqrt{-1}  T(e^\bot_l, Je^\bot_l)
+2T(\tfrac{\partial}{\partial z^0_j},
\tfrac{\partial}{\partial \ov{z}^0_j})\right\rangle\\
&\hspace*{5mm} + \Big |\sum_j T(\tfrac{\partial}{\partial z^0_j},
\tfrac{\partial}{\partial \ov{z}^0_j})\Big |^2
+ \frac{9}{16} \mT_{llm}\mT_{kkm}
-\frac{3\sqrt{-1}}{2}  \left\langle T(e^\bot_l, Je^\bot_l),
T(\tfrac{\partial}{\partial z^0_j},
\tfrac{\partial}{\partial \ov{z}^0_j})\right\rangle,\\
&\mF_1 (e^\bot_k) \mT_{kll} = -\sqrt{-1}\left\langle T(e^\bot_l, Je^\bot_l),
\wi{\mu}^E
+ T(\tfrac{\partial}{\partial z^0_j},\tfrac{\partial}{\partial \ov{z}^0_j})
\right\rangle
+ \frac{3}{4} \mT_{llm}\mT_{kkm} . 
\end{split}\end{align}

By \eqref{4.23}, \eqref{4.24}, \eqref{4.25} and \eqref{4.75}, we have
\begin{multline}\label{4.76}
(\Psi_{1,1}+ \Psi_{1,1}^* + \Psi_{1,3}-\Psi_{1,4})(0)
= \left\{- \frac{1}{2\pi} \sum_k \mF_1 (e^\bot_k)^2
-\frac{1}{8\pi} \mF_1 (e^\bot_k)\mT_{kll}   \right.\\
- \frac{13}{2^6 \cdot 3\pi}\mT_{klm}^2
+ \frac{1}{2^7 \pi} \mT_{kkm}\mT_{llm}
\left. -\frac{11}{48 \pi}
\Big|\mT_{kl}(\tfrac{\partial}{\partial \ov{z}^0_i})\Big|^2
-\frac{1}{8 \pi}\Big|\sum_{k}\mT_{kk}
(\tfrac{\partial}{\partial \ov{z}^0_i})\Big|^2\right\} P^N(0,0)\\
= \left\{ \frac{1}{2\pi}\left\langle \wi{\mu}^E_{x_0},
\wi{\mu}^E_{x_0}\right\rangle_{g^{TY}}
+ \frac{1}{\pi} \left\langle \wi{\mu}^E,
\frac{7}{8} \sqrt{-1} T(e^\bot_l, Je^\bot_l)
+T(\tfrac{\partial}{\partial z^0_l},
\tfrac{\partial}{\partial \ov{z}^0_l})\right\rangle  \right.\\
-\frac{1}{2\pi} \Big| \sum_j T(\tfrac{\partial}{\partial z^0_j},
\tfrac{\partial}{\partial \ov{z}^0_j})\Big|^2
+ \frac{7\sqrt{-1}}{8\pi} \left\langle  T(e^\bot_l, Je^\bot_l),
T(\tfrac{\partial}{\partial z^0_j},
\tfrac{\partial}{\partial \ov{z}^0_j})\right\rangle
-  \frac{47}{2^7 \pi}\mT_{kkm}\mT_{llm}\\
- \frac{13}{2^6 \cdot 3\pi}\mT_{klm}^2
\left.    - \frac{11}{48 \pi}
\Big|T(e^\bot_k,\tfrac{\partial}{\partial \ov{z}^0_j})\Big|^2
-\frac{1}{8 \pi}\Big|\sum_{k}\mT_{kk}
(\tfrac{\partial}{\partial \ov{z}^0_i})\Big|^2\right\} P^N(0,0).
\end{multline}
By \eqref{4.46} and \eqref{f4.1}, we get
\begin{multline}\label{a4.76}
\Psi_{1,2}(0)+ \Psi_{1,2}(0)^*
=\left\{\frac{1}{8\pi} r^{X_G}_{x_0}
 +\frac{1}{\pi} R^{E_G}(\tfrac{\partial}{\partial z^0_j},
\tfrac{\partial}{\partial \ov{z}^0_j}) + \frac{1}{\pi}\Delta_{X_G} \log h
\right. \\
-\frac{3}{8\pi} \nabla_{e^\bot_k}\nabla_{e^\bot_k}\log h
+ \frac{35}{48 \pi}
|T(e^\bot_k,\tfrac{\partial}{\partial \ov{z}^0_j})|^2
+ \frac{\sqrt{-1}}{8 \pi} \left\langle T(e^\bot_l, Je^\bot_l),
T(\tfrac{\partial}{\partial z^0_j},
\tfrac{\partial}{\partial \ov{z}^0_j}) \right\rangle \\
+\frac{1}{2\pi} \Big| T(\tfrac{\partial}{\partial z^0_i},
\tfrac{\partial}{\partial \ov{z}^0_j})\Big|^2
- \frac{1}{16 \pi}\Big|\sum_{k}\mT_{kk}
(\tfrac{\partial}{\partial \ov{z}^0_i})\Big|^2
+\frac{1}{2^6 \pi}\Big[\wi{\mT}_{ijk}(\wi{\mT}_{kji} + \wi{\mT}_{ijk})
-2 (\wi{\mT}_{jki} + \wi{\mT}_{ijk})\wi{\mT}_{ijk} \Big]\\
+ \frac{7}{2^7 \pi} \Big(2\mT_{jkm}^2 + \mT_{jjm}\mT_{kkm} \Big)
\left.- \frac{\sqrt{-1}}{8 \pi} \Big(\left\langle T(e^\bot_l, Je^\bot_l),
 \wi{\mu}^E\right\rangle
-2  \left\langle J e^\bot_k,
\nabla^{TY}_{e^\bot_k} \wi{\mu}^E\right\rangle\Big)
 \right\} P^N (0,0).
\end{multline}

Thus by \eqref{4.20}, \eqref{4.76} and \eqref{a4.76}, as
$\wi{\mT}_{ijk}$ is anti-symmetric on $i,j$, we get
\begin{multline}\label{4.77}
P^{(2)}(0,0)=\left\{\frac{1}{8\pi} r^{X_G}_{x_0}
 +\frac{1}{\pi} R^{E_G}(\tfrac{\partial}{\partial z^0_j},
\tfrac{\partial}{\partial \ov{z}^0_j}) + \frac{1}{\pi}\Delta_{X_G} \log h
-\frac{3}{8\pi} \nabla_{e^\bot_k}\nabla_{e^\bot_k}\log h \right. \\
+ \frac{1}{2 \pi}
|T(e^\bot_k,\tfrac{\partial}{\partial \ov{z}^0_j})|^2
+\frac{1}{2\pi} \Big| T(\tfrac{\partial}{\partial z^0_i},
\tfrac{\partial}{\partial \ov{z}^0_j})\Big|^2
-\frac{1}{2\pi} \Big| \sum_j T(\tfrac{\partial}{\partial z^0_j},
\tfrac{\partial}{\partial \ov{z}^0_j})\Big|^2  \\
+ \frac{\sqrt{-1}}{\pi} \left\langle T(e^\bot_l, Je^\bot_l),
T(\tfrac{\partial}{\partial z^0_j},
\tfrac{\partial}{\partial \ov{z}^0_j}) \right\rangle
- \frac{3}{16 \pi}\Big|\sum_{k}\mT_{kk}
(\tfrac{\partial}{\partial \ov{z}^0_i})\Big|^2
+ \frac{1}{24 \pi}\mT_{klm}^2\\
-\frac{5}{16 \pi}\mT_{kkm}\mT_{llm}
+ \frac{1}{2^6 \pi} \wi{\mT}_{ijk} (3 \wi{\mT}_{kji} - \wi{\mT}_{ijk})
+ \frac{1}{2\pi}\left\langle \wi{\mu}^E_{x_0},
\wi{\mu}^E_{x_0}\right\rangle_{g^{TY}} \\
+ \frac{1}{\pi} \left\langle \wi{\mu}^E,
\frac{3}{4} \sqrt{-1} T(e^\bot_l, Je^\bot_l)
+T(\tfrac{\partial}{\partial z^0_l},
\tfrac{\partial}{\partial \ov{z}^0_l})\right\rangle
\left.+ \frac{\sqrt{-1}}{4 \pi} \left\langle J e^\bot_k,
\nabla^{TY}_{e^\bot_k} \wi{\mu}^E\right\rangle
 \right\} P^N (0,0).
\end{multline}
By Theorem \ref{t4.2}, \eqref{h1}, \eqref{h14},
\eqref{3.11c} and \eqref{34.21}, as same as in \eqref{3.79},
we get for $U\in T_{x_0}X_G$,
\begin{align}\label{4.78}
\begin{split}
&\mT_{llm}=\left\langle T(e^\bot_m, J e^\bot_l), J e^\bot_l\right\rangle
= 2 \nabla_{e^\bot_m} \log h,\quad
T(e^\bot_l, J e^\bot_l) = 2 (\nabla_{e^\bot_k}\log h) Je^\bot_k,\\
&\mT_{kk}(U)= -2\left\langle T(JU, J e^\bot_k), J e^\bot_k\right\rangle
= -\left\langle \dot{g}^{TY}_{JU}
Je^\bot_k, Je^\bot_k\right\rangle =-4 \nabla_{JU}\log h.
\end{split}\end{align}
By \eqref{4.22}, \eqref{4.77} and \eqref{4.78},
we get Theorem \ref{t4.15}.
\end{proof}

\subsection{Coefficient $\Phi_1$: general case}\label{s8.1}
We use the general assumption at the beginning of this Section,
but we do not suppose that $\bJ=J$ in \eqref{0.1}.


Let $\overline{\partial} ^{L^p\otimes E,*}$  be the formal adjoint of
the Dolbeault operator $\overline{\partial} ^{L^p\otimes E}$
on the Dolbeault complex
 $\Omega ^{0,\bullet}(X, L^p\otimes E)$ with the scalar product
$\left\langle \quad \right \rangle$
induced by $g^{TX}$, $h^L$, $h^E$ as in Section \ref{s3.01}.

Set
\begin{align}\label{bf12}
D_p = \sqrt{2}\left( \overline{\partial} ^{L^p\otimes E} +
\overline{\partial} ^{L^p\otimes E,*}\right)
\end{align}
 Then
\begin{align}\label{bf14}
D_p^2= 2\left( \overline{\partial} ^{L^p\otimes E}\overline{\partial}
^{L^p\otimes E,*} +\overline{\partial} ^{L^p\otimes
E,*}\overline{\partial} ^{L^p\otimes E}\right)
\end{align}
 preserves the $\bZ$-grading of $\Omega ^{0,\bullet}(X, L^p\otimes E)$.

For $p$ large enough,
\begin{equation} \label{f12}
\Ker D_p  =\Ker D_p^2  = H^0 (X,L^p\otimes E).
\end{equation}
Here $D_{p}$ need not be a spin$^c$ Dirac operator on
 $\Omega ^{0,\bullet}(X, L^p\otimes E)$.

Let $P_p^G(x,x')$ $(x,x'\in X)$ be the smooth kernel of the orthogonal
projection $P_p^G$ from 
$(\cC ^\infty(X, L^p\otimes E), \left\langle \quad \right \rangle)$ 
onto $(\Ker D_p^2)^G$
with respect to the Riemannian volume form $dv_X(x')$ for $p$ large enough.

We explain now how to reduce the study of the asymptotic expansion of
$P_p^G(x,x')$ to the $\bJ=J$ case.

Let $g^{TX}_{\om}(\cdot,\cdot):= \om(\cdot,J\cdot)$ be the metric on $TX$
induced by $\om, J$.
We will use a subscript $\om$ to indicate the objects corresponding to
 $g^{TX}_{\om}$, especially  $r^{X}_{\om}$ is the scalar curvature of
$(TX,g^{TX}_{\om})$, and  $\Delta_{X_G,\om}$ is the Bochner-Laplace
 operator on $X_G$ as in \eqref{r6} associated  to $g^{TX_G}_{\om}$.

Let $\det_{\bC} $ denote the determinant function on the complex
 bundle $T^{(1,0)}X$, and $|{\bf J}|= (-{\bf J}^2)^{-1/2}$.

Let $h^E_\om:= ({\det}_{\bC}|{\bf J}|)^{-1} h^E$ define a metric
on $E$. Let $R^E_\om$ be the curvature associated to the
holomorphic Hermitian connection on $(E, h^E_\om)$.

Let $\left\langle \quad \right \rangle_\om$ be the Hermitian product
on $\cC^\infty (X, L^p\otimes E)$
induced by $g^{TX}_\om, h^L, h^E_\om$ as in \eqref{h10}, then
\begin{equation}\label{af19}
(\cC^\infty (X, L^p\otimes E), \left\langle \quad \right \rangle_\om)
= (\cC^\infty (X, L^p\otimes E), \left\langle \quad \right \rangle),
 \quad dv_{X,\om} =({\det}_{\bC}|{\bf J}|) dv_{X}.
\end{equation}

Observe that $H^0(X, L^p\otimes E)$ does not depend on $g^{TX},  h^L, h^E$.

Let $P_{\om,p}^G(x,x')$ $(x,x'\in X)$ be the smooth kernel of the
orthogonal projection $P^G_{\om,p}$ from $(\cC^\infty (X,
L^p\otimes E), \left\langle \quad \right \rangle_\om)$ onto
$H^0(X, L^p\otimes E)^G$ with respect to $dv_{X,\om}(x)$.

By \eqref{af19},
\begin{equation}\label{af20}
P_{p}^G(x,x')= ({\det}_{\bC}|{\bf J}|)(x') P_{\om,p}^G(x,x').
\end{equation}

We will use the trivialization in Introduction corresponding
 to $g^{TX}_{\om}$.

 Since
$g^{TX}_\om(\cdot, \cdot)= \om(\cdot, J\cdot)$ is a K\"ahler metric on $TX$,
$D_{\om,p}$ is a Dirac operator (cf. Def. \ref{Dirac}).
Thus Theorems \ref{t0.0}, \ref{t0.1} hold for $P^G_{\om, p}(x,x')$.


Let $dv_B$ be the volume form on $B$ induced by $g^{TX}$
as in Introduction.

As in \eqref{0.9}, let $\wi{\kappa}\in \cC^\infty (TB|_{X_G},\bR)$
be defined by for $Z\in T_{x_0}B$, $x_0\in X_G$,
\begin{align}\label{af21}
dv_B (x_0, Z)
=\wi{\kappa}(x_0, Z) dv_{X_G,\om}(x_0)dv_{N_{G,\om, x_0}}.
\end{align}
As in \eqref{0.10}, we introduce $\mathscr{I}_p(x_0)$ a section of $\End(E_G)$
on $X_G$,
\begin{align}\label{af22}
 \mathscr{I}_p(x_0)=  \int_{Z\in N_{G}, |Z|\leq \var_0}
h^2 (x_0,Z) P^G_{p}\circ \Psi((x_0,Z),(x_0,Z))
 \wi{\kappa}(x_0, Z)dv_{N_{G,\om, x_0}}.
\end{align}
Then \eqref{0.11} still holds.

Summarizes, we have the following result,

\begin{thm}\label{nonkahler}
The smooth kernel $P_p^G(x,x')$
has a full off--diagonal asymptotic expansion analogous to \eqref{0.8}
with $\mQ_0=(\det_\bC |{\bf J}|) \Id_{E_G}$ as $p\to\infty$\,.
There exist $\Phi_r(x_0)\in  \End (E_G)_{x_0}$
 polynomials in $A_\om$, $R^{TB}_\om$, $R^{E_B}$,
$\mu^E$, $R^E$ {\rm (}resp. $h_\om$, $R^{L_B}$; resp. $\mu${\rm )}
and their derivatives at $x_0$ to order $2r-1$
 {\rm (}resp. $2r$, resp. $2r+1${\rm )},
 and $\Phi_0= {\rm Id}_{E_G}$
such that \eqref{0.15} holds for $\mathscr{I}_p$. Moreover
\begin{align}\label{af12}
\Phi_1(x_0)= \frac{1}{8\pi}\Big[r^{X_G}_\om
+ 6 \Delta_{X_G,\om}\log \left(h_\om|_{X_G}\right)
-2 \Delta_{X_G,\om} \Big(\log({\det}_\bC |{\bf J}|)\Big)
+ 4 R^{E_G} (w_{\om,j}^0,\ov{w}_{\om,j}^0)\Big].
\end{align}
Here
$\{w_{\om,j}\}$ is an orthogonal basis of $(T^{(1,0)}X_G, g^{TX_G}_\om)$.
\end{thm}
\begin{proof} By \eqref{af19}-\eqref{af22},
\begin{align}\label{af23}
 \mathscr{I}_p(x_0)=  \int_{Z\in N_{G}, |Z|\leq \var_0}
h^2_\om (x_0,Z) P^G_{\om,p}\circ \Psi((x_0,Z),(x_0,Z))
 \kappa_\om (x_0, Z) dv_{N_G,\om}(Z).
\end{align}
From the above discussion, only \eqref{af12} reminds to be proved. But
\begin{equation}\label{af17}
R^{E_G}_\om =R^{E_G} -\ov{\partial}\partial \log ({\det}_{\bC}|{\bf J}|),
\end{equation}
Thus
\begin{equation}\label{af18}
2R^{E_G}_\om(w_{\om,j}^0,\ov{w}_{\om,j}^0)
=2 R^{E_G}(w_{\om,j}^0,\ov{w}_{\om,j}^0)
-\Delta_{X_G,\om} \log ({\det}_{\bC}|{\bf J}|),
\end{equation}
and \eqref{af12} is  from \eqref{0.5} and \eqref{af23}.
\end{proof}

\section{Bergman kernel and geometric quantization}\label{s8}

In this Section, we prove Theorems \ref{t0.9}, \ref{t0.7}.

\begin{proof}[Proof of Theorem \ref{t0.9}]
We use the notation in Section \ref{s6.4}.

By Theorem \ref{t6.7} and Lemma \ref{at6.8}, we know that
$p^{-\frac{n_0}{4}}(\sigma_p\circ \sigma_p^*)^{\frac{1}{2}}$ is a
Toeplitz operator with principal symbol
$(2^{\frac{n_0}{4}}/\wi{h}(x_0)) \Id_{E_G}$ in the sense of Def.
\ref{d6.6}, and its kernel has an expansion analogous  to
\eqref{a6.32} and $Q_{0,0}$ therein is
$2^{\frac{n_0}{4}}/\wi{h}(x_0)$.

We claim that
\begin{align}\label{b0.20}
\bT_p=p^{-\frac{n_0}{2}}(\sigma_p\circ \sigma_p^*)^{\frac{1}{2}}
\wi{h}^2(\sigma_p\circ \sigma_p^*)^{\frac{1}{2}}
\end{align}
is a Toeplitz operator with principal symbol $2^{\frac{n_0}{2}} \Id_{E_G}$.

Indeed, when $E=\bC$, this is a consequence of \cite{BMS94} on the
composition of the Toeplitz operators.

To get the above claim for general $E$, we need just keep in mind
 that the kernel $\bT_p(x_0,x'_0)$ of $\bT_p$ with respect $dv_{X_G}(x'_0)$
has the expansion analogous to \eqref{a6.32} and $Q_{0,0}$ therein
is $2^{\frac{n_0}{2}} \Id_{E_G}$.

 Our claim then  follows from
the composition of the expansion of the kernel of
$p^{-\frac{n_0}{4}}(\sigma_p\circ \sigma_p^*)^{\frac{1}{2}}$, as
well as the Taylor expansion of $\wi{h}^2$ (cf. also the recent
book \cite{MM05b} for a more detailed proof).

Now we still denote by $\langle \, , \, \rangle$ the $L^2$-scalar
product on $\cC^\infty(X_G,L^p_G\otimes E_G)$ induced by
$h^{L^p_G}$, $h^{E_G}$, $g^{TX_G}$ as in \eqref{h10}.

Let $\{s^p_i\}$ be an orthonormal basis of
 $(H^0(X,L^p\otimes E)^G,\langle \, , \, \rangle)$,
then $\varphi^p_i= (\sigma_p\circ \sigma_p^*)^{-\frac{1}{2}}\sigma_p s^p_i$
is an orthonormal basis of
$(H^0(X_G,L^p_G\otimes E_G),\langle \, , \, \rangle)$.

From Def. \ref{d6.6}, \eqref{0.26}, \eqref{h10} and \eqref{b0.20},
we get
\begin{multline}\label{b0.21}
(2p)^{-\frac{n_0}{2}} \left\langle \sigma_p s^p_i,\sigma_p
s^p_j\right\rangle_{\widetilde{h}} =(2p)^{-\frac{n_0}{2}}
\left\langle (\sigma_p\circ \sigma_p^*)^{\frac{1}{2}} \varphi^p_i,
(\sigma_p\circ \sigma_p^*)^{\frac{1}{2}}
\varphi^p_j\right\rangle_{\widetilde{h}}\\
= 2^{-\frac{n_0}{2}}\left\langle \bT_p
\varphi^p_i,\varphi^p_j\right\rangle = \delta_{ij} +
\cO\left(\frac{1}{p}\right).
\end{multline}

The proof of Theorem \ref{t0.9} is complete.
\end{proof}



\begin{proof}[Proof of Theorem \ref{t0.7}]
Set  $$\widetilde{h}^{E_G}=\widetilde{h}^2 h^{E_G}.$$ Then
$\wi{P}^{X_G}_{p}$ is the orthogonal projection from $\cC^\infty
(X_G,L_G^p\otimes E_G)$ onto $H^0(X,L_G^p\otimes E_G)$,
 associated to the Hermitian product on  $\cC^\infty
(X_G,L_G^p\otimes E_G)$ induced by
 the metrics $h^{L_G}$, $\widetilde{h}^{E_G}$, $g^{TX_G}$ as in
\eqref{h10}.

Let $\wi{P}^{X_G}_{p,\om}(x_0,x_0')$ be the smooth kernel of
$\wi{P}^{X_G}_{p}$ with respect to $dv_{X_G}(x_0')$. Then
\begin{align}\label{b0.17}
\wi{P}^{X_G}_{p,\om}(x_0,x_0')
= \widetilde{h}^2(x'_0)\wi{P}^{X_G}_{p}(x_0,x_0').
\end{align}


Let $\wi{\nabla} ^{E_G}$ be the Hermitian holomorphic connection on
$(E_G, \widetilde{h}^{E_G})$ with curvature $\wi{R}^{E_G}$.

Then
\begin{align}\label{b0.16}
\wi{\nabla} ^{E_G}= {\nabla} ^{E_G} + \partial \log (\widetilde{h}^2),\quad
\wi{R}^{E_G} = {R}^{E_G} + 2\ov{\partial}\partial \log \widetilde{h}.
\end{align}
Thus from \eqref{b0.16},
\begin{align}\label{b0.29}
&\wi{R}^{E_G}(w^0_j, \ov{w}^0_j)
= 2 \wi{R}^{E_G}(\tfrac{\partial}{\partial z^0_j},
\tfrac{\partial}{\partial \ov{z}^0_j})
= R^{E_G}(w^0_j, \ov{w}^0_j) +  \Delta_{X_G} \log \widetilde{h}.
\end{align}

By \eqref{3.19}, \eqref{b0.17} and \eqref{b0.29},
Theorem \ref{t0.7} is a direct
consequence of \cite[Theorem 1.3]{DLM04a} (or Theorem \ref{t0.6}
with $G=\{1\}$) for $\wi{P}^{X_G}_{p,\om}(x_0,x_0)$.
\end{proof}

\newpage

\subsection*{Acknowledgments}
We thank Professor Jean-Michel Bismut for useful conversations.
The work of the second author was partially supported
by MOEC, MOSTC and NNSFC. Part of work was done while the second author was
visiting IHES during February and March, 2005. He would like to thank
Professor Jean-Pierre Bourguignon for the hospitality.


\def\cprime{$'$} \def\cprime{$'$}
\providecommand{\bysame}{\leavevmode\hbox to3em{\hrulefill}\thinspace}
\providecommand{\MR}{\relax\ifhmode\unskip\space\fi MR }
\providecommand{\MRhref}[2]{%
  \href{http://www.ams.org/mathscinet-getitem?mr=#1}{#2}
}
\providecommand{\href}[2]{#2}

\end{document}